\documentclass{book}

\usepackage{amssymb}
\usepackage{amsmath}
\usepackage{latexsym}
\usepackage{eufrak}
\usepackage{euscript}
\usepackage{xypic}
\usepackage{makeidx}
\usepackage[backref=true, pagebackref=true]{hyperref}

\makeindex

\begin{document}

\frontmatter

\title{Exterior Differential Systems and \\ Euler-Lagrange Partial
  Differential Equations}
\author{Robert Bryant \and Phillip Griffiths \and Daniel
  Grossman}
\maketitle

\tableofcontents

\newtheorem{Definition}{Definition}[chapter]
\newtheorem{Fact}{Fact}[chapter]
\newtheorem{Theorem}{Theorem}[chapter]
\newtheorem{Example}{Example}[chapter]
\newtheorem{Proposition}{Proposition}[chapter]
\newtheorem{Lemma}{Lemma}[chapter]
\newtheorem{Corollary}{Corollary}[chapter]
\newtheorem{Conjecture}{Conjecture}[chapter]

\newenvironment{Proof}{\noindent\textbf{Proof.}}{\hfill$\square$}

\newcommand{\B}{{\mathbf B}}
\newcommand{\C}{{\mathbf C}}
\newcommand{\I}{{\mathcal I}}
\newcommand{\nn}[1]{#1n}
\newcommand{\p}{{\partial}}
\newcommand{\R}{{\mathbf R}}
\newcommand{\F}{{\mathcal F}}

\def\lefthook{\mathbin{
\hbox spread 0pt{\vrule height1.4pt depth-1pt width 4pt
                 \vrule height6pt depth-1pt}}}

\def\innerprod{\,\,\hbox to 4pt{\hrulefill}\kern-3.25pt\hbox{
     \vrule height7pt\,\,}}

\def\bw#1{\textstyle\bigwedge^{#1}\displaystyle}
\def\lie#1{{\EuFrak{#1}}}
\def\sf#1#2{\textstyle\frac{#1}{#2}\displaystyle}
\def\ss{\textstyle \sum\displaystyle}

\chapter{Preface}

During the 1996-97 academic year, Phillip Griffiths and Robert Bryant
conducted a seminar at the Institute for Advanced Study in Princeton,
NJ, outlining
their recent work (with Lucas Hsu) on a geometric
approach to the calculus of variations in several variables. 
The present work is an outgrowth of that project;
it includes all of the material
presented in the seminar, with numerous additional details and
a few extra topics of interest.

The material can be viewed as a chapter in the ongoing development of
a theory of the geometry of differential equations.  The relative
importance among PDEs of second-order
Euler-Lagrange equations suggests that their geometry should be
particularly rich, as does the geometric character of their
conservation laws, which we discuss at length.  

A second purpose for the present work is to give an exposition of
certain aspects of the theory of exterior differential systems, which
provides the language and the techniques for the entire study.  Special
emphasis is placed on the method of equivalence, which plays a central
role in uncovering geometric properties of differential equations.
The Euler-Lagrange PDEs of the calculus of variations have turned out
to provide excellent illustrations of the general theory.


\chapter{Introduction}

In the classical calculus of variations, one studies
functionals\index{functional|(} 
of the form
\begin{equation}
  {\mathcal F}_L(z) = \int_\Omega L(x,z,\nabla z)\,dx,
  \qquad \Omega\subset\R^n,
\label{ClassicalFunctional}
\end{equation}
where $x=(x^1,\ldots,x^n),\ dx=dx^1\wedge\cdots\wedge dx^n,\ 
  z=z(x)\in C^1(\bar\Omega)$ (for example), and the 
{\em Lagrangian}\index{Lagrangian|(}
$L=L(x,z,p)$ is a smooth function of $x$, $z$, and $p=(p_1,\ldots,p_n)$.
Examples frequently encountered in physical field theories are
Lagrangians of the form
$$
  L=\sf12||p||^2+F(z),
$$
usually interpreted as a kind of 
energy\index{energy}.  
The 
{\em Euler-Lagrange equation}\index{Euler-Lagrange!equation} 
describing functions $z(x)$ that are 
stationary\index{stationary|(} 
for such a
functional\index{functional} 
is the second-order partial differential equation
$$
  \Delta z(x) = F^\prime(z(x)).
$$
For another example, we may identify a function $z(x)$ with its graph
$N\subset\R^{n+1}$, and take the Lagrangian
$$
  L = \sqrt{1+||p||^2},
$$
whose associated functional ${\mathcal F}_L(z)$ equals the 
area\index{area functional} 
of the graph, regarded as a hypersurface in 
Euclidean space.\index{Euclidean!space}
The Euler-Lagrange equation\index{Euler-Lagrange!equation} 
describing functions $z(x)$ stationary for this functional is $H=0$,
where $H$ is the 
mean curvature\index{mean curvature} 
of the graph $N$.

To study these Lagrangians and Euler-Lagrange equations geometrically,
one has to choose a class of admissible coordinate changes, and there
are four natural candidates.  In increasing order of generality, they
are: 
\begin{itemize}
\item Classical transformations\index{transformation!classical}, 
of the form $x^\prime=x^\prime(x)$, $z^\prime=z^\prime(z)$; in this
situation, we think of $(x,z,p)$ as coordinates on the space
$J^1(\R^n,\R)$ of $1$-jets of maps $\R^n\to\R$.\footnote{A 
{\em $1$-jet}\index{jets|nn} 
is an equivalence class of functions having
the same value and the same first derivatives at some designated point
of the domain.}
\item Gauge transformations\index{transformation!gauge}, 
of the form $x^\prime=x^\prime(x)$, $z^\prime=z^\prime(x,z)$; here, we
think of $(x,z,p)$ as coordinates on the space of $1$-jets of sections
of a bundle $\R^{n+1}\to\R^n$, where $x=(x^1,\ldots,x^n)$ are
coordinates on the base $\R^n$ and $z\in\R$ is a fiber coordinate.
\item Point transformations\index{transformation!point}, 
of the form $x^\prime=x^\prime(x,z)$, $z^\prime=z^\prime(x,z)$; here,
we think of $(x,z,p)$ as coordinates on the space of tangent
hyperplanes
$$
  \{dz-p_idx^i\}^\perp\subset T_{(x^i,z)}(\R^{n+1})
$$ 
of the manifold $\R^{n+1}$ with coordinates $(x^1,\ldots,x^n,z)$.
\item Contact transformations\index{transformation!contact}, 
of the form $x^\prime=x^\prime(x,z,p)$, $z^\prime=z^\prime(x,z,p)$,
$p^\prime=p^\prime(x,z,p)$, satisfying the equation of differential
$1$-forms 
$$
  dz^\prime-\ss p_i^\prime dx^{i\prime} = f\cdot(dz-\ss p_idx^i)
$$
for some function $f(x,z,p)\neq 0$.
\end{itemize}
We will be studying the geometry of 
functionals\index{functional}
${\mathcal F}_L(z)$ subject to the class of 
contact transformations\index{transformation!contact}, 
which is strictly larger than the other three classes.  The effects of
this choice will become clear as we proceed.  Although contact
transformations were recognized classically, appearing most notably in
studies of surface geometry, they do not seem to have been extensively
utilized in the calculus of variations.

\

Classical calculus of variations primarily concerns the following
features of a functional ${\mathcal F}_L$.

The {\em first variation}\index{first variation} 
$\delta{\mathcal F}_L(z)$ is analogous to
the derivative of a function, where $z=z(x)$ is thought of as an
independent variable in an infinite-dimensional space of functions.
The analog of the condition that a point be critical is the condition
that $z(x)$
be {\em stationary}\index{stationary|)} 
for all fixed-boundary variations.  Formally, one writes
$$
  \delta{\mathcal F}_L(z) = 0,
$$
and as we shall explain, this gives a second-order scalar partial
differential equation for the unknown function $z(x)$ of the form
$$
  \frac{\partial L}{\partial z} - \sum
    \frac{d}{dx^i}\left(\frac{\partial L}{\partial p_i}\right)=0.
$$
This is the 
{\em Euler-Lagrange equation}\index{Euler-Lagrange!equation} 
of the Lagrangian $L(x,z,p)$, and we will study it in an invariant,
geometric setting.  This seems 
especially promising in light of the fact that, although it is not
obvious, the process by which we associate an Euler-Lagrange equation
to a Lagrangian is invariant under the large class of 
contact transformations\index{transformation!contact}.  
Also, note that the
Lagrangian $L$ determines the functional ${\mathcal F}_L$, but not
vice versa.  To see this, observe that if we add to $L(x,z,p)$ a
``divergence term'' and consider 
$$
  L^\prime(x,z,p) = L(x,z,p) + \sum
    \left(\frac{\partial K^i(x,z)}{\partial x^i} + 
      \frac{\partial K^i(x,z)}{\partial z}p^i\right)
$$
for functions $K^i(x,z)$, then by Green's 
theorem\index{Green's theorem}, 
the functionals ${\mathcal F}_L$ and ${\mathcal F}_{L^\prime}$
differ by a constant depending only on values of $z$ on
$\partial\Omega$.  For many purposes, such functionals should be
considered equivalent; in particular, $L$ and $L^\prime$ have the same
Euler-Lagrange equations.

Second, there is a relationship between 
{\em symmetries}\index{symmetry} 
of a Lagrangian $L$ and 
{\em conservation laws}\index{conservation law} 
for the corresponding Euler-Lagrange equations, described by a
classical theorem of 
Noether\index{Noether's theorem}.  
A subtlety here is that the group of symmetries of an equivalence
class of Lagrangians may be strictly larger than the group of
symmetries of any particular representative.  We will investigate how
this discrepancy is reflected in the space of conservation laws, in a
manner that involves global topological issues.

Third, one considers the 
{\em second variation}\index{second variation} 
$\delta^2{\mathcal F}_L$, analogous to the Hessian of a smooth
function, usually with the goal of identifying local minima of the
functional.  There has been a great deal of analytic work done in this
area for classical variational problems, reducing the problem of local
minimization to understanding the behavior of certain 
{\em Jacobi operators}\index{Jacobi operator}, 
but the geometric theory is not as well-developed as that of the first
variation and the Euler-Lagrange equations. 
\index{functional|)}
\index{Lagrangian|)}

\

We will consider these issues and several others in a geometric
setting as suggested above, using various methods from the subject of 
exterior differential systems, to be explained along the way.
Chapter 1 begins with an introduction to
{\em contact manifolds}\index{contact!manifold}, 
which provide the geometric setting for the study of first-order
functionals (\ref{ClassicalFunctional}) subject to contact
transformations.  We then construct an object that is central to
the entire theory: the 
{\em Poincar\'e-Cartan form}\index{Poincar\'e-Cartan form}, 
an explicitly computable differential form that is associated to the
equivalence class of any Lagrangian, where the notion of
equivalence includes that alluded to above for classical Lagrangians.
We then carry out a calculation using the
Poincar\'e-Cartan form to associate to any Lagrangian on a
contact manifold an exterior differential system---the 
{\em Euler-Lagrange system}\index{Euler-Lagrange!system}---whose 
integral manifolds are 
{\em stationary}\index{stationary} 
for the associated functional; in the classical case, these
correspond to solutions of the Euler-Lagrange equation.
The Poincar\'e-Cartan form also makes it quite easy to state and
prove 
{\em Noether's theorem}\index{Noether's theorem}, 
which gives an
isomorphism between a space of symmetries of a Lagrangian and
a space of conservation laws for the Euler-Lagrange equation; exterior
differential systems provides a particularly natural setting for
studying the latter objects.  We illustrate all of this theory in the
case of 
minimal hypersurfaces\index{minimal surface} 
in Euclidean space ${\mathbf E}^n$, and in the case of more general
linear Weingarten surfaces\index{Weingarten equation} 
in ${\mathbf E}^3$, providing intuitive and computationally simple
proofs of known results.

In Chapter 2, we consider the geometry of Poincar\'e-Cartan forms more
closely.  The main tool for this is 
\'E.~Cartan's\index{Cartan, \'E.} 
{\em method of equivalence}\index{equivalence method}, 
by which one develops an algorithm for associating to certain
geometric structures their differential invariants under a specified
class of equivalences.  We explain the various steps of this method
while illustrating them in several major cases.  First, we apply the
method to 
{\em hyperbolic Monge-Ampere 
systems}\index{Monge-Ampere system!hyperbolic} 
in two independent variables; these exterior differential systems
include many important
Euler-Lagrange systems that arise from classical problems, and
among other results, we find a characterization of those PDEs that are
contact-equivalent to the homogeneous linear wave equation.  We then
turn to the case of $n\geq 3$ independent variables, and carry out
several steps of the equivalence method for Poincar\'e-Cartan forms,
after isolating those of the algebraic type arising from classical
problems.  Associated to such a 
{\em neo-classical form}\index{Poincar\'e-Cartan form!neo-classical} 
is a
field of hypersurfaces in the fibers of a vector bundle, well-defined
up to 
affine transformations\index{affine!transformation}.  
This motivates a digression on the
affine geometry of hypersurfaces, conducted using Cartan's {\em method
of moving frames}, which we will illustrate but not discuss in any
generality.  After identifying a number of differential invariants for
Poincar\'e-Cartan forms in this manner, we show that they are
sufficient for characterizing those
Poincar\'e-Cartan forms associated to the PDE for hypersurfaces having
prescribed mean curvature\index{mean curvature!prescribed}.

A particularly interesting branch of the equivalence problem for
neo-classical Poincar\'e-Cartan forms includes some highly symmetric
Poincar\'e-Cartan forms corresponding to 
Poisson equations\index{Poisson equation}, 
discussed in Chapter 3.  Some of
these equations have good invariance properties under the group of
conformal transformations of the $n$-sphere, 
and we find that the corresponding branch of the
equivalence problem reproduces a construction
that is familiar in conformal
geometry.  We will discuss the relevant aspects of conformal geometry
in some detail; these include another application of the equivalence
method, in which the important conceptual step of 
{\em prolongation}\index{prolongation!of a $G$-structure}
of $G$-structures appears for the first time.  This
point of view allows us to apply
Noether's theorem in a particularly simple way to the most symmetric
of non-linear Poisson equations, the one with the critical exponent:
$$
  \Delta u = Cu^{\frac{n+2}{n-2}}.
$$
Having calculated the conservation laws for this equation, we also
consider the case of wave equations, and in particular the very
symmetric example:
$$
  \square z = Cz^{\frac{n+3}{n-1}}.
$$
Here, conformal geometry with Lorentz signature is the appropriate
background, and we present the conservation laws corresponding to the
associated symmetry group, along with a few elementary applications.

The final chapter addresses certain matters which
are thus far not so well-developed.  First, we consider the 
second variation\index{second variation} 
of a functional, with the goal of understanding
which integral manifolds of an Euler-Lagrange system are local
minima.  We give an interesting geometric formula for the second
variation, in which conformal geometry makes another appearance
(unrelated to that in the preceding chapter). 
Specifially, we find that the critical
submanifolds for certain variational problems inherit a canonical
conformal structure\index{conformal!structure}, 
and the 
second variation\index{second variation} 
can be expressed in
terms of this structure and an additional scalar curvature invariant.
This interpretation does not seem to appear in the classical literature. 
Circumstances under
which one can carry out in an invariant manner the usual ``integration
by parts'' in the second-variation formula, which is crucial for the
study of local minimization, turn out to be somewhat limited.  We
discuss the reason for this, and illustrate the optimal situation by
revisiting the example of 
prescribed mean curvature\index{mean curvature!prescribed}
systems.

We also consider the problem of finding an analog of the
Poincar\'e-Cartan form in the case of functionals on vector-valued
functions and their Euler-Lagrange PDE systems.  Although there is no
analog of proper contact transformations in this case, we will present
and describe the merits of 
D. Betounes'\index{Betounes form} 
construction of such an
analog, based on some rather involved multi-linear algebra.  An
illuminating special case is that of 
harmonic maps\index{harmonic map} 
between Riemannian manifolds, for which we find the associated forms
and conservation laws.

Finally, we consider the appearance of higher-order conservation laws
for first-order variational problems.  The geometric setting for these
is the {\em infinite prolongation} of an Euler-Lagrange system, which
has come to play a major role in classifying conservation laws.  We
will propose a generalized version of Noether's theorem appropriate to
our setting, but we do not have a proof of our statement.  In any
case, there are other ways to illustrate two of the
most well-known but intriguing examples: the system describing
Euclidean surfaces of Gauss curvature $K=-1$, and that corresponding
to the sine-Gordon equation\index{sine-Gordon equation}, 
$\square z = \sin z$.  We will generate examples of higher-order
conservation laws by relating these two systems, first in the
classical manner, and then more systematically using the notions of 
{\em prolongation} and 
{\em integrable extension}\index{integrable extension}, 
which come from the subject of exterior differential systems.
Finally, having explored these systems this far, it is convenient to
exhibit and relate the 
{\em B\"acklund transformations}\index{B\"acklund transformation} 
that act on each.

One particularly appealing aspect of this study is that one sees in
action so many aspects of the subject of exterior differential
systems.  There are particularly beautiful instances of the method of
equivalence, a good illustration of the method of moving frames
(for affine hypersurfaces), essential use of prolongation both of
$G$-structures and of differential systems, and a use of the notion of
integrable extension to clarify a confusing issue. 

\

Of course, the study of Euler-Lagrange equations
by means of exterior differential forms and the method of
equivalence is not new.  In fact, much of the 19th century
material in this area is so naturally formulated in terms
of differential forms (cf. the Hilbert form in the one-variable
calculus of variations) that it is difficult to say exactly
when this approach was initiated.

However, there is no doubt that 
\'Elie Cartan's\index{Cartan, \'E.}
1922 work {\em Le\c cons sur les invariants 
int\'egraux}~\cite{Cartan:Lecons} serves both as an elegant summary of
the known material at the time and as a remarkably forward-looking
formulation of the use of differential forms in the
calculus of variations.  At that time, Cartan did
not bring his method of equivalence (which he had developed
beginning around 1904 as a tool to study the geometry of
pseudo-groups) to bear on the subject.
It was not until his 1933 work {\em Les espaces m\'etriques
fond\'es sur la notion d'aire}~\cite{Cartan:Surface} and his
1934 monograph {\em Les espaces de Finsler}~\cite{Cartan:Finsler}
that Cartan began to explore the geometries that one could
attach to a Lagrangian for surfaces or for curves.  Even in these
works, any explicit discussion of the full method of equivalence
is supressed and Cartan contents himself with deriving the
needed geometric structures by seemingly ad hoc methods.

After the modern formulation of jet spaces and their
contact systems was put into place, Cartan's approach was extended
and further developed by several people. One might particularly
note the 1935 work of Th.~de Donder~\cite{Donder} and its development.
Beginning in the early 1940s,
Th.~Lepage~\cite{Lepage:Classe,Lepage:Equation} undertook a study
of first order Lagrangians that made extensive use of the
algebra of differential forms on a contact manifold.
Beginning in the early 1950s, this point of view was developed
further by P. Dedecker~\cite{MR56:16680}, who undertook a serious
study of the calculus of variations via tools of homological algebra.
All of these authors are concerned in one way or another with
the canonical construction of differential geometric (and other) 
structures
associated to a Lagrangian, but the method of equivalence is
not utilized in any extensive way.  Consequently, they deal
primarily with first-order linear-algebraic invariants of variational
problems.  Only with the method of equivalence can one 
uncover the full set of higher-order geometric
invariants.   This is one of the central themes
of the present work; without the equivalence method, 
for example, one could not give our unique
characterizations of certain classical, ``natural'' systems
(cf.~\S\ref{Section:SmallEquiv},
\S\ref{Section:PrescribedH}, and \S\ref{Section:ConfEquiv}).

In more modern times, numerous works of
I. Anderson, D. Betounes,
R. Hermann, N. Kamran, V. Lychagin, P. Olver, H. Rund, A. Vinogradov,
and their coworkers, just to name a few,
all concern themselves with geometric aspects and invariance properties
of the calculus of variations.  Many of the results expounded
in this monograph can be found in one form or another in works by
these or earlier authors.  We certainly make no pretext of giving a
complete historical account of the work in this area in the 20th
century.  Our
bibliography lists those works of which we were aware that
seemed most relevant to our approach, if not necessarily to the
results themselves, and it identifies
only a small portion of the work done in these areas.
The most substantially developed alternative theory in this area is that of the
{\em variational bicomplex} associated to the algebra of differential
forms on a fiber bundle.  The reader can learn this material
from Anderson's works~\cite{Anderson:Introduction}
and~\cite{Anderson:Variational}, and references therein, which contain
results heavily overlapping those of our Chapter~\ref{Chapter:Additional}.

\

Some terminology and notation that we will use follows, with more
introduced in the text.  An {\em exterior differential 
system}\index{exterior differential system}
(EDS) is a pair
$(M,{\mathcal E})$ consisting of a smooth manifold $M$ and a
homogeneous, differentially closed ideal ${\mathcal E}\subseteq
\Omega^*(M)$ in the algebra of smooth differential forms on $M$.  Some
of the EDSs that we study are differentially generated by the sections
of a smooth subbundle $I\subseteq T^*M$ of the cotangent bundle of
$M$; this subbundle, and sometimes its space of sections, is called a
{\em Pfaffian system}\index{Pfaffian system} on $M$.  It will be
useful to use the notation $\{\alpha,\beta,\ldots\}$ for the
(two-sided) {\em algebraic} ideal generated by forms $\alpha$,
$\beta$,\ldots, and to use the notation $\{I\}$ for the algebraic
ideal generated by the sections of a Pfaffian system $I\subseteq T^*M$.
An 
{\em integral manifold}\index{integral manifold} 
of an EDS
$(M,{\mathcal E})$ is a submanifold immersion $\iota:N\hookrightarrow
M$ for which $\varphi_N\stackrel{\mathit{def}}{=} \iota^*\varphi = 0$
for all $\varphi\in{\mathcal E}$.  Integral manifolds of Pfaffian
systems are defined similarly.

A differential form $\varphi$ on the total space of a fiber bundle
$\pi:E\to B$ is said to be
{\em semibasic}\index{semibasic form} 
if its contraction with any vector field tangent to the fibers of
$\pi$ vanishes, or equivalently, if its value at each point $e\in E$
is the pullback via $\pi^*_e$ of some form at $\pi(e)\in B$.  Some
authors call such a form {\em horizontal}.
A stronger condition is that $\varphi$ be 
{\em basic}\index{basic form}, 
meaning that it is locally (in open subsets of $E$) the pullback via
$\pi^*$ of a form on the base $B$. 

Our computations will frequently require the following multi-index
notation.  If $(\omega^1,\ldots,\omega^n)$ is an ordered basis for a
vector space $V$, then corresponding to a multi-index
$I = (i_1,\ldots,i_k)$ is the $k$-vector
$$
  \omega^I = \omega^{i_1}\wedge\cdots\wedge\omega^{i_k}
    \in \bw{k}(V),
$$
and for the complete multi-index we simply define
$$
  \omega = \omega^1\wedge\cdots\wedge\omega^n.
$$
Letting $(e_1,\ldots,e_n)$ be a dual basis for $V^*$,
we also define the $(n-k)$-vector
$$
  \omega_{(I)} = e_I\innerprod \omega = 
    e_{i_k}\innerprod(e_{i_{k-1}}\innerprod
    \cdots(e_{i_1}\innerprod\omega)\cdots).
$$
This $\omega_{(I)}$ is, up to sign, just $\omega^{I_c}$, where $I_c$
is a multi-index complementary to $I$.
For the most frequently occurring cases $k=1,2$ we have the formulae 
(with ``hats'' $\,\hat{}\,$ indicating omission of a
factor)
\begin{eqnarray*}
  \omega_{(i)} & = & (-1)^{i-1}\omega^1\wedge\cdots\wedge
    \hat\omega^i\wedge\cdots\wedge\omega^n,\\
  \omega_{(ij)} & = & (-1)^{i+j-1}\omega^1\wedge\cdots\wedge
    \hat\omega^i\wedge\cdots\wedge\hat\omega^j\wedge\cdots
    \wedge\omega^n \\
  & = & -\omega_{(ji)},\quad\mbox{for $i<j$},
\end{eqnarray*}
and the identities
\begin{eqnarray*}
  \omega^i\wedge\omega_{(j)} & = & \delta^i_j\omega, \\
  \omega^i\wedge\omega_{(jk)} & = & \delta^i_k\omega_{(j)} -
    \delta^i_j\omega_{(k)}.
\end{eqnarray*}

We will often, but not always, use without comment the convention of
summing over repeated indices.  Always, $n\geq 2$.\footnote{For the
  case $n=1$, an analogous geometric approach to the calculus of
  variations for curves may be found in \cite{Griffiths:Exterior}.}

\mainmatter

\chapter{Lagrangians and Poincar\'e-Cartan Forms}
\label{Chapter:Lagrangians}

In this chapter, we will construct and illustrate our basic objects of
study.  The geometric setting that one uses for studying Lagrangian
functionals subject to contact transformations is a {\em contact
manifold}, and we will begin with its definition and relevant
cohomological properties.  These properties allow us to formalize an
intuitive notion of equivalence for functionals, and more importantly,
to replace such an equivalence class by a more concrete differential
form, the 
{\em Poincar\'e-Cartan form}\index{Poincar\'e-Cartan form}, 
on which all of our later calculations depend.  In particular, we will
first use it to derive the 
Euler-Lagrange differential system\index{Euler-Lagrange!system}, 
whose integral manifolds correspond to stationary points of a given
functional.  We then use it to give an elegant version of the solution
to the 
inverse problem\index{inverse problem}, 
which asks when a differential system of the appropriate algebraic
type is the Euler-Lagrange system of some functional.  Next, we use it
to define the isomorphism between a certain Lie algebra of
infinitesimal symmetries of a variational problem and the space
conservation laws for the Euler-Lagrange system, as described in
Noether's theorem\index{Noether's theorem}.  
All of this will be illustrated at an elementary level using examples
from Euclidean hypersurface geometry.

\section{Lagrangians and Contact Geometry}

We begin by introducing the geometric setting in which we will study
Lagrangian functionals and their Euler-Lagrange systems.

\begin{Definition}
A {\em contact manifold}\index{contact!manifold|(} 
$(M, I)$ 
is a smooth manifold $M$ of dimension $2n+1\ (n\in\mathbf{Z}^+)$, 
with a distinguished line
sub-bundle $I\subset T^*M$ of the cotangent bundle which is
non-degenerate\index{non-degenerate!$1$-form} 
in the sense that for any local $1$-form $\theta$ generating $I$,
$$
  \theta\wedge(d\theta)^n \neq 0.
$$
\end{Definition}
Note that the non-degeneracy criterion is independent of the choice of 
$\theta$; this is because if $\bar\theta=f\theta$ for some function
$f\neq 0$, then we find 
$$
  \bar\theta\wedge(d\bar\theta)^n = 
    f^{n+1}\theta\wedge(d\theta)^n.
$$

For example, on the space $J^1(\R^n,\R)$ of $1$-jets\index{jets} of
functions, we can take coordinates $(x^i,z,p_i)$ corresponding to the
jet at $(x^i)\in\R^n$ of the linear function $f(\bar x)=z+\sum
p_i(\bar x^i-x^i)$.  Then we define the 
{\em contact form}\index{contact!form}
$$
  \theta = dz-\sum p_idx^i,
$$
for which
$$
  d\theta=-\sum dp_i\wedge dx^i,
$$
so the non-degeneracy condition $\theta\wedge(d\theta)^n\neq 0$ is
apparent.  In fact, the 
Pfaff theorem\index{Pfaff theorem} 
(cf.~Ch.~I, \S3 of~\cite{Bryant:Exterior})
implies that every contact manifold is locally isomorphic to this
example; that is, every contact manifold\index{contact!manifold}
$(M,I)$ has local coordinates $(x^i,z,p_i)$ for which the form
$\theta=dz-\sum p_idx^i$ generates $I$.

More relevant for differential geometry is the example
$G_n(TX^{n+1})$, the
Grassmannian\index{Grassmannian} 
bundle parameterizing $n$-dimensional oriented subspaces of the
tangent spaces of an $(n+1)$-dimensional manifold $X$.  It is 
naturally a contact manifold, and will be considered in more detail
later.

Let $(M,I)$ be a contact manifold of dimension $2n+1$, and assume
that $I$ is generated by a global, non-vanishing section
$\theta\in\Gamma(I)$; this assumption only simplifies our notation,
and would in any case hold on a double-cover of $M$.
Sections of $I$ generate the {\em contact
differential ideal}\index{ideal!differential}
$$ 
  \I=\{\theta,d\theta\}\subset\Omega^*(M)
$$ 
in the exterior algebra of
differential forms on $M$.\footnote{Recall our convention that braces
  $\{\cdot\}$ denote the 
  algebraic ideal\index{ideal!algebraic|nn} 
  generated by an object; for instance, $\{\theta\}$ consists of exterior
  multiples of any contact form $\theta$, and is smaller than $\I$.
  We sometimes use $\{I\}$ as alternate notation for $\{\theta\}$.}    
A {\em Legendre submanifold}\index{Legendre submanifold|(} 
of $M$ is an immersion $\iota:N\hookrightarrow M$ of an
$n$-dimensional submanifold $N$ such that $\iota^*\theta = 0$ for any
contact form $\theta\in\Gamma(I)$; in this case $\iota^*d\theta=0$ as
well, so a Legendre submanifold is the same thing as
an integral manifold\index{integral manifold} 
of the 
differential ideal\index{ideal!differential} $\I$.
In Pfaff coordinates with $\theta=dz-\sum p_idx^i$,
one such integral manifold is given by
$$
  N_0=\{z=p_i=0\}.
$$
To see other Legendre submanifolds ``near'' this one, note than any
submanifold $C^1$-close to $N_0$ satisfies the 
independence condition\index{independence condition}
$$
  dx^1\wedge\cdots\wedge dx^n \neq 0,
$$
and can therefore be described locally as a graph
$$
  N=\{(x^i,z(x),p_i(x))\}.
$$
In this case, we have
$$
  \theta|_N=0\quad\mbox{if and only if}\quad 
    p_i(x)=\frac{\partial z}{\partial x^i}(x).
$$
Therefore, $N$ is determined by the function
$z(x)$, and conversely, every function $z(x)$ determines such an $N$;
 we informally say that ``the generic Legendre submanifold depends locally 
on one arbitrary function of $n$ variables.''  Legendre submanifolds
of this form, with $dx|_N\neq 0$, will often be described 
as {\em transverse}.\index{Legendre submanifold!transverse}

Motivated by (\ref{ClassicalFunctional}) in the Introduction, we are
primarily interested in 
functionals\index{functional|(} 
given by triples $(M,I,\Lambda)$, where $(M,I)$ is a $(2n+1)$-dimensional
contact manifold,
and $\Lambda\in\Omega^n(M)$ is a differential form of degree $n$ on
$M$; such a $\Lambda$ will be referred to as a {\em Lagrangian} on $(M,I)$.\footnote{In the Introduction, we used the
term {\em Lagrangian} for a function, 
rather than for a differential form, but we will not do so again.}
We then define a functional on the set of smooth, compact Legendre
submanifolds $N\subset M$, possibly with boundary $\partial N$, 
by
$$
  {\mathcal F}_\Lambda(N) = \int_N\Lambda.
$$
The classical variational problems described above may be recovered
from this notion by taking
$M=J^1(\R^n,\R)\cong\R^{2n+1}$ with coordinates $(x^i,z,p_i)$, $I$
generated by $\theta = dz-\sum
p_idx^i$, and $\Lambda = L(x^i,z,p_i)dx$.  This formulation
also admits certain functionals depending on second derivatives of
$z(x)$, because there may be $dp_i$-terms in $\Lambda$.  Later, we
will restrict attention to a class of functionals which, possibly after a
contact transformation, can be expressed without second derivatives.

There are two standard notions of 
equivalence\index{equivalence!of Lagrangians} for 
Lagrangians\index{Lagrangian|(}
$\Lambda$.  First, note that if the difference
$\Lambda-\Lambda^\prime$ of two Lagrangians lies in the contact
ideal\index{contact!ideal} $\I$ then the
functionals ${\mathcal F}_\Lambda$ and ${\mathcal F}_{\Lambda^\prime}$ 
are equal, because they are defined only for Legendre
submanifolds, on which all forms in $\I$ vanish.  Second, suppose that 
the difference of two Lagrangians is an exact $n$-form,
$\Lambda-\Lambda^\prime=d\varphi$ for some
$\varphi\in\Omega^{n-1}(M)$.  Then we find
$$
  {\mathcal F}_\Lambda(N)={\mathcal F}_{\Lambda^\prime}(N)
    + \int_{\partial N}\varphi
$$
for all 
Legendre submanifolds\index{Legendre submanifold|)} 
$N$.  One typically studies the variation
of ${\mathcal F}_\Lambda$ along $1$-parameter families $N_t$ with
fixed boundary, and the 
preceding equation shows that ${\mathcal F}_\Lambda$ and ${\mathcal
F}_{\Lambda^\prime}$ differ only by a constant on such a family.  Such
$\Lambda$ and $\Lambda^\prime$ are sometimes said to be {\em
  divergence-equivalent}.

These two notions of equivalence suggest that we consider the class
$$
  [\Lambda] \in \Omega^n(M)/(\I^n+d\Omega^{n-1}(M)),
$$
where $\I^n=\I\,\cap\,\Omega^n(M)$.
The natural setting for this space is the quotient $(\bar\Omega^*,
\bar d)$ of the de Rham
complex $(\Omega^*(M),d)$,
where $\bar\Omega^n = \Omega^n(M)/\I^n$, and $\bar d$ is induced by the
usual exterior derivative $d$ on this quotient.
We then have {\em characteristic cohomology
  groups}\index{characteristic cohomology} $\bar
H^n=H^n(\bar\Omega^*,\bar
d)$.  We will show in a moment that (recalling $\mbox{dim}(M)=2n+1$):
\begin{equation}
  \mbox{for }k>n,\ \ \I^k=\Omega^k(M).
\label{HighCohomology}
\end{equation}
In other words, all forms on $M$ of degree greater than $n$ lie in the 
contact ideal\index{contact!ideal}; one consequence is that $\I$ can have no integral
manifolds\index{integral manifold} of dimension greater than $n$.  The importance of
(\ref{HighCohomology}) is that it implies that $d\Lambda\in\I^{n+1}$,
and we can therefore regard our equivalence class of functionals as a
characteristic cohomology\index{characteristic cohomology} class
$$
  [\Lambda]\in \bar H^n.
$$
This class is almost, but not quite, our fundamental object 
of study.

To prove both (\ref{HighCohomology}) and several later results, we need to
describe some of the pointwise linear algebra associated with the
contact ideal\index{contact!ideal} $\I=\{\theta,d\theta\}\subset\Omega^*(M)$.  Consider the
tangent distribution of rank $2n$
$$
  I^\perp\subset TM
$$
given by the annihilator of the contact line 
bundle\index{contact!line bundle}.  Then the
non-degeneracy condition on $I$ implies that the $2$-form
$$
  \Theta \stackrel{\mathit{def}}{=} d\theta
$$ 
restricts fiberwise to $I^\perp$ as a non-degenerate,
alternating bilinear form, determined by $I$ up to scaling.  This
allows one to use tools from symplectic
linear algebra\index{symplectic!linear algebra|(}; the main fact is the
following.
 
\begin{Proposition} 
Let $(V^{2n},\Theta)$ be a symplectic vector space, where
$\Theta\in\bigwedge^2V^*$ is a non-degenerate alternating bilinear
form.  Then 

\noindent
(a) for $0\leq k\leq n$, the map
\begin{equation}
  (\Theta\wedge)^k:\bw{n-k}V^*\to \bw{n+k}V^*
\label{SymplecticMap}
\end{equation}
is an isomorphism, and

\noindent
(b) if we
define the space of {\em primitive forms}\index{primitive} to be
$$
  P^{n-k}(V^*)=\mbox{{\em Ker} }((\Theta\wedge)^{k+1}:\bw{n-k}V^*\to
     \bw{n+k+2}V^*),
$$
then
we have a decomposition of $Sp(n,\R)$-modules
$$
  \bw{n-k}V^*=P^{n-k}(V^*)\oplus\left(\Theta\wedge\bw{n-k-2}V^*\right).
\footnote{Of course, we can
extend this decomposition inductively to obtain the {\em Hodge-Lepage
decomposition}\index{Hodge-Lepage decomposition}
$$
  \bw{n-k}V^* \cong P^{n-k}(V^*)\oplus
    P^{n-k-2}(V^*)\oplus\cdots\oplus \left\{\begin{array}{l}
       P^1(V^*)=V^* \mbox{ for }n-k\mbox{ odd,} \\
       P^0(V^*)=\R \mbox{ for }n-k\mbox{ even,}
  \end{array}\right.
$$
under which any element $\xi\in\bigwedge^{n-k}V^*$ can be
written uniquely as
$$
  \xi = \xi_0 + (\Theta\wedge\xi_1) + (\Theta^2\wedge\xi_2) +
   \cdots + \left(\Theta^{\lfloor\frac{n-k}{2}\rfloor}\wedge
      \xi_{\lfloor\frac{n-k}{2}\rfloor}\right),
$$
with each $\xi_i\in P^{n-k-2i}(V^*)$.
What we will not prove here is that the representation of
$Sp(n,\R)$ on $P^{n-k}(V^*)$ is irreducible for each $k$, 
so this gives
the complete irreducible decomposition of $\bigwedge^{n-k}(V^*)$.}
$$
\label{Prop:Symplectic}
\end{Proposition}

Proposition~\ref{Prop:Symplectic}
implies in particular (\ref{HighCohomology}), for it says that
modulo $\{\theta\}$ (equivalently, restricted to $I^\perp$), every
form $\varphi$ of degree greater than $n$ is a multiple
of $d\theta$, which is exactly to say that $\varphi$ is in the
algebraic ideal generated by $\theta$ and $d\theta$.

\

\begin{Proof}
(a) Because $\bigwedge^{n-k}V^*$ and $\bigwedge^{n+k}V^*$ have the same
dimension, it suffices to show that the map (\ref{SymplecticMap}) is
injective.  We proceed by induction on $k$, downward from $k=n$ to $k=0$.
In case $k=n$, the (\ref{SymplecticMap}) is just multiplication
$$
  (\Theta^n)\cdot:\R\to\bw{2n}V^*,
$$
which is obviously injective, because $\Theta$ is non-degenerate.  

Now suppose that the statement is proved for some $k$, and suppose that
$\xi\in\bigwedge^{n-(k-1)}$ satisfies
$$
  \Theta^{k-1}\wedge\xi = 0.
$$
This implies that
$$
  \Theta^k\wedge\xi = 0,
$$
so that for every vector $X\in V$, we have
$$
  0 = X \innerprod (\Theta^k\wedge\xi) 
    = k(X\innerprod\Theta)\wedge\Theta^{k-1}\wedge\xi
            + \Theta^k\wedge(X\innerprod\xi).
$$
Now, the first term on the right-hand side vanishes by our assumption
on $\xi$ (our second use of this assumption), so we must have
$$
  0 = \Theta^k\wedge(X\innerprod\xi),
$$
and the induction hypothesis then gives
$$
  X \innerprod \xi = 0.
$$
This is true for every $X\in V$, so we conclude that $\xi = 0$.

\

\noindent
(b) We will show that any $\xi\in\bigwedge^{n-k}V^*$ has a unique
decomposition as the sum of a primitive form and a multiple of $\Theta$.
For the existence of such a decomposition, we apply the surjectivity in
part (a) to the element
$\Theta^{k+1}\wedge\xi\in\bigwedge^{n+k+2}V^*$, and find
$\eta\in\bigwedge^{n-k-2}V^*$ for which
$$
  \Theta^{k+2}\wedge\eta = \Theta^{k+1}\wedge\xi.
$$
Then we can decompose
$$
  \xi = (\xi - \Theta\wedge\eta) + (\Theta\wedge\eta),
$$
where the first summand is primitive by construction.

To prove uniqueness, we need to show that if $\Theta\wedge\eta$ is
primitive for some $\eta\in\bigwedge^{n-k-2}(V^*)$, then
$\Theta\wedge\eta = 0$.  In fact, primitivity means
$$
  0 = \Theta^k\wedge\Theta\wedge\eta,
$$
which implies that $\eta=0$ by the injectivity in part (a).
\end{Proof}
\index{symplectic!linear algebra|)}

\

Returning to our discussion of Lagrangian functionals,
observe that there is a short exact sequence of complexes
$$
  0\to\I^*\to\Omega^*(M)\to\bar\Omega^*\to 0
$$
giving a long exact cohomology sequence
$$
  \cdots\to H^n_{dR}(M)\to \bar H^n
        \stackrel{\delta}{\to}H^{n+1}(\I)\to H^{n+1}_{dR}(M)
    \to\cdots,
$$
where $\delta$ is essentially exterior differentiation.  Although an
equivalence class $[\Lambda]\in\bar H^n$ generally has no canonical
representative differential form, we can now show that 
its image $\delta([\Lambda])\in H^{n+1}(\I)$ does.

\begin{Theorem}
Any class $[\Pi]\in H^{n+1}(\I)$ has a unique global representative closed
form $\Pi\in\I^{n+1}$ satisfying $\theta\wedge\Pi=0$ for any contact
form $\theta\in\Gamma(I)$, or equivalently, $\Pi\equiv 0\mbox{ (mod
$\{I\}$)}$.  
\end{Theorem}

\begin{Proof}
Any $\Pi\in\I^{n+1}$ may be written locally as
$$
  \Pi = \theta\wedge\alpha+d\theta\wedge\beta
$$
for some $\alpha\in\Omega^n(M)$, $\beta\in\Omega^{n-1}(M)$.  But this
is the same as
$$
  \Pi = \theta\wedge(\alpha+d\beta)+d(\theta\wedge\beta),
$$
so replacing $\Pi$ with the equivalent (in $H^{n+1}(\I)$) form
$\Pi-d(\theta\wedge\beta)$, we have the {\em local} existence of a 
representative as claimed.

For uniqueness, suppose that $\Pi_1-\Pi_2 = d(\theta\wedge\gamma)$ for
some $(n-1)$-form $\gamma$
(this is exactly equivalence in $H^{n+1}(\I)$), and that
$\theta\wedge\Pi_1=\theta\wedge\Pi_2=0$.
Then $\theta\wedge
d\theta\wedge\gamma=0$, so $d\theta\wedge\gamma\equiv
0\mbox{ (mod $\{I\}$)}$.  By symplectic linear algebra, this
implies that $\gamma\equiv 0\mbox{ (mod $\{I\}$)}$, so $\Pi_1-\Pi_2=0$.

Finally, global existence follows from local existence and uniqueness.
\end{Proof}

We can now define our main object of study.
\begin{Definition}
For a contact manifold $(M,I)$ with Lagrangian $\Lambda$,
the unique representative $\Pi\in\I^{n+1}$ of $\delta([\Lambda])$
satisfying $\Pi\equiv0\mbox{ (mod $\{I\}$)}$ is called
the {\em Poincar\'e-Cartan form}\index{Poincar\'e-Cartan form|(} of $\Lambda$.
\end{Definition}

Poincar\'e-Cartan forms of Lagrangians will be the main object of
study in these lectures, and there are two computationally useful ways to
think of them.  The first
is as above:  given a representative Lagrangian $\Lambda$, express
$d\Lambda$ locally as
$\theta\wedge(\alpha+d\beta)+d(\theta\wedge\beta)$, and
then 
$$
\boxed{
  \Pi = \theta\wedge(\alpha+d\beta).
}
$$  
The second, which will be
important for computing the first variation\index{first variation} and
the Euler-Lagrange system\index{Euler-Lagrange!system} of $[\Lambda]$,
is as an exact form:
$$
\boxed{
   \Pi = d(\Lambda-\theta\wedge\beta).
}
$$
In fact, $\beta$ is the unique $(n-1)$-form modulo $\{I\}$ such
that 
$$
  d(\Lambda-\theta\wedge\beta)\equiv 0 \pmod{\{I\}}. 
$$
This observation will be used later, in the proof of Noether's
theorem.\index{Noether's theorem}
\index{Poincar\'e-Cartan form|)}
\index{contact!manifold|)}
\index{functional|)}
\index{Lagrangian|)}

\section{The Euler-Lagrange System}
\label{Section:EulerLagrange}

In the preceding section, we showed how one can associate to an
equivalence class $[\Lambda]$ of Lagrangians on a contact
manifold\index{contact!manifold} $(M,I)$ a canonical
$(n+1)$-form $\Pi$.  In this section, we use this Poincar\'e-Cartan
form\index{Poincar\'e-Cartan form|(} to find an exterior differential system whose integral manifolds
are precisely the
stationary\index{stationary|(} Legendre submanifolds
for the functional $\mathcal F_{\Lambda}$.  This requires
us to calculate the {\em first variation}\index{first variation} of $\mathcal F_{\Lambda}$,
which gives the derivative of $\mathcal F_{\Lambda}(N_t)$ for any
$1$-parameter family $N_t$ of Legendre submanifolds of $(M,I)$.  The
Poincar\'e-Cartan form enables us to carry out the usual integration
by parts for this calculation in an invariant manner.

We also consider the relevant version of the {\em inverse
  problem}\index{inverse problem} of
the calculus of
variations, which asks whether a given PDE of the appropriate type is
equivalent to the Euler-Lagrange
equation\index{Euler-Lagrange!equation} for some functional.  We
answer this by giving a necessary and sufficient condition for an EDS
of the appropriate type to be locally equivalent to the Euler-Lagrange
system\index{Euler-Lagrange!system} of some $[\Lambda]$.  We find
these conditions by reducing the
problem to a search for a Poincar\'e-Cartan form.

\subsection{Variation of a Legendre Submanifold}
\index{Legendre submanifold|(}
\index{Legendre variation|(}

Suppose that we have a 1-parameter family $\{N_t\}$ of Legendre
submanifolds of a 
contact manifold\index{contact!manifold}
$(M,I)$; more precisely, this is given by a compact manifold with
boundary $(N,\partial N)$ and a 
smooth map
$$
  F:N\times [0,1]\to M
$$
which is a Legendre submanifold $F_t$ for each fixed $t\in[0,1]$ and is
independent of $t\in[0,1]$ on $\partial N\times[0,1]$.  Because
$F_t^*\theta=0$ for any contact form\index{contact!form} $\theta\in\Gamma(I)$, we must have locally
\begin{equation}
   F^*\theta=G\,dt
\label{FlowGenerator}
\end{equation}
for some function $G$ on $N\times[0,1]$.  We let $g=G|_{N\times\{0\}}$
be the restriction to the initial submanifold.

It will be useful to know that given a Legendre submanifold
$f:N\hookrightarrow M$, every function $g$ 
may be realized as in (\ref{FlowGenerator}) for some fixed-boundary
variation and some contact form $\theta$, locally in the interior
$N^o$.  This may be seen in Pfaff coordinates
$(x^i,z,p_i)$ on $M$, for which $\theta=dz-\sum
p_idx^i$ generates $I$ and such that our
given $N$ is a $1$-jet graph $\{(x^i,z(x),p_i(x)=z_{x^i}(x))\}$.
Then $(x^i)$ give coordinates on $N$, and a variation of 
$N$ is of the form
$$
  F(x,t) = (x^i, z(x,t), z_{x^i}(x,t)).
$$
Now $F^*(dz-\sum p_idx^i) = z_tdt$; and given $z(x,0)$, we can
always extend to $z(x,t)$ with $g(x)=z_t(x,0)$ prescribed arbitrarily,
which is what we claimed.

\subsection{Calculation of the Euler-Lagrange System}
\label{Subsection:EulerLagrange}
\index{Euler-Lagrange!system|(}

We can now carry out a calculation that is fundamental for the whole
theory.  Suppose given a Lagrangian $\Lambda\in\Omega^n(M)$ on a
contact manifold\index{contact!manifold} $(M,I)$, and a fixed-boundary variation of Legendre
submanifold $F:N\times[0,1]\to M$; we wish to compute
$\frac{d}{dt}(\int_{N_t}\Lambda)$.

To do this, first recall the calculation of the Poincar\'e-Cartan form 
for the equivalence class $[\Lambda]\in \bar H^n$.  Because 
$\I^{n+1}=\Omega^{n+1}(M)$,
we can always write
\begin{eqnarray*}
  d\Lambda & = & \theta\wedge\alpha+d\theta\wedge\beta \\
  & = & \theta\wedge(\alpha+d\beta) + d(\theta\wedge\beta),
\end{eqnarray*}
and then
\begin{equation}
  \Pi = \theta\wedge(\alpha+d\beta) =
      d(\Lambda-\theta\wedge\beta).
\label{DefPC2}
\end{equation}
We are looking for conditions on a Legendre submanifold
$f:N\hookrightarrow M$ to be
{\em stationary} for $[\Lambda]$ under all fixed-boundary variations, in
the sense that
$\left.\frac{d}{dt}\right|_{t=0}(\int_{N_t}\Lambda)=0$ whenever $F|_{t=0}=f$.
We compute (without writing the $F^*$s)
\begin{eqnarray*}
   \frac{d}{dt}\int_{N_t}\Lambda  & = &
      \frac{d}{dt}\int_{N_t}(\Lambda-\theta\wedge\beta) \\
   & = & \int_{N_t}{\mathcal L}_{\frac{\partial}{\partial t}}
          (\Lambda-\theta\wedge\beta) \\
   & = & \int_{N_t}\left(\sf{\partial}{\partial t}\innerprod
       d(\Lambda-\theta\wedge\beta)\right)
       + \int_{N_t}d\left(\sf{\partial}{\partial t}\innerprod
         (\Lambda-\theta\wedge\beta)\right) \\
   & = & \int_{N_t}\sf{\partial}{\partial t}\innerprod\Pi
       \quad \mbox{(because $\partial N$ is fixed)}.
\end{eqnarray*}
One might express this result as
$$
  \delta({\mathcal F}_\Lambda)_N(v) = \int_N v\innerprod f^*\Pi,
$$
where the variational vector field $v$, lying in the space
$\Gamma_0(f^*TM)$ of sections of $f^*TM$ vanishing along $\partial
N$, plays the role of $\frac{\partial}{\partial t}$.
The condition $\Pi\equiv 0\mbox{ (mod $\{I\}$)}$ allows us to write
$\Pi = \theta\wedge\Psi$ for some $n$-form $\Psi$, not
uniquely determined, and we have
$$
  \left.\frac{d}{dt}\right|_{t=0}\int_{N_t}\Lambda = \int_{N}g\,f^*\Psi,
$$
where $g = (\frac{\partial}{\partial t}\innerprod F^*\theta)|_{t=0}$.
It was shown previously that this $g$
could locally be chosen arbitrarily in the interior $N^o$, so the necessary and
sufficient condition for a Legendre submanifold
$f:N\hookrightarrow M$ to be stationary for
${\mathcal F}_\Lambda$ is that $f^*\Psi=0$.
\index{Legendre variation|)}

\begin{Definition}
  The {\em Euler-Lagrange system} of the Lagrangian $\Lambda$ is the
differential ideal\index{ideal!differential} generated algebraically as
$$
   {\mathcal E}_\Lambda = \{\theta, d\theta, \Psi\}
      \subset \Omega^*(M).
$$
A {\em stationary Legendre submanifold} of $\Lambda$ is an integral
manifold of ${\mathcal E}_\Lambda$.  The functional is said to be {\em 
non-degenerate}\index{non-degenerate!functional} if its
Poincar\'e-Cartan form $\Pi=\theta\wedge\Psi$
has no degree-1 divisors (in the exterior algebra of $T^*M$)
other than multiples of $\theta$.
\label{ELDef2}
\end{Definition}

Note first that ${\mathcal E}_\Lambda$ is uniquely determined by $\Pi$,
even though $\theta$ and $\Psi$ may not be.\footnote{Actually, given $\Pi$ we
  have not only a
  well-defined ${\mathcal E}_\Lambda$, but a well-defined $\Psi$
  modulo $\{I\}$ which is primitive\index{primitive|nn} on $I^\perp$.
  There is a canonical map ${\mathcal
    E}:H^n(\bar\Omega^*)\to P^n(T^*M/I)$ to the space of primitive
  forms, taking a Lagrangian class $[\Lambda]$ to the corresponding
  $\Psi$ in its Euler-Lagrange system; and this map fits into a full
  resolution of the constant sheaf
  $$
    0\to\R\to\bar\Omega^0\to\cdots\to\bar\Omega^{n-1}\to
    H^n(\bar\Omega^*)\stackrel{{\mathcal E}}{\to}P^n(T^*M/I)\to
    \cdots\to P^0(T^*M/I)\to 0.
  $$
  This has been developed and applied in the context of CR geometry in
  \index{Rumin, M.|nn}\cite{Rumin:Complexe}.}
Note also that the ideal in $\Omega^*(M)$ algebraically generated by
$\{\theta,d\theta,\Psi\}$ is already differentially closed, because
$d\Psi\in\Omega^{n+1}(M)=\I^{n+1}$.  

We can examine this for the classical situation where 
$M=\{(x^i,z,p_i)\}$, $\theta=dz-\sum p_idx^i$, and $\Lambda=L(x,z,p)dx$.  We find
\begin{eqnarray*}
  d\Lambda & = & L_z\theta\wedge dx + \ss L_{p_i}dp_i\wedge dx \\
    & = & \theta\wedge L_zdx-d\theta\wedge \ss L_{p_i}dx_{(i)},
\end{eqnarray*}
so referring to (\ref{DefPC2}),
$$
  \Pi =  \theta\wedge(L_zdx-\ss d(L_{p_i}dx_{(i)})) = \theta\wedge\Psi.
$$
Now, for a transverse Legendre submanifold\index{Legendre
submanifold!transverse}
$N=\{(x^i,z(x),z_{x^i}(x))\}$, we have $\Psi|_N=0$ if and only if
along $N$
$$
  \frac{\partial L}{\partial z} - \sum
    \frac{d}{dx^i}\left(\frac{\partial L}{\partial p_i}\right)=0,
$$
where 
$$
  \frac{d}{dx^i} = \frac{\partial}{\partial x^i} +
    z_{x^i}\frac{\partial}{\partial z} + \sum_j z_{x^ix^j}
    \frac{\partial}{\partial p_j}
$$
is the {\em total derivative}\index{total derivative}.  This is
the usual 
Euler-Lagrange equation\index{Euler-Lagrange!equation}, a
second-order, quasi-linear PDE for $z(x^1,\ldots,x^n)$ having
symbol\index{symbol} $L_{p_ip_j}$.  It is an exercise to show that this
symbol matrix is invertible at $(x^i,z,p_i)$ if and only if $\Lambda$ is
non-degenerate\index{non-degenerate!functional} in the sense of
Definition~\ref{ELDef2}.
\index{Legendre submanifold|)}
\index{stationary|)}

\subsection{The Inverse Problem}
\label{Subsection:Inverse}
\index{inverse problem|(}

There is a reasonable model for exterior differential systems of
``Euler-Lagrange type''.
\begin{Definition}
\index{Monge-Ampere system|(}
  A {\em Monge-Ampere differential system} $(M,{\mathcal E})$ consists
  of a contact manifold\index{contact!manifold}
$(M,I)$ of dimension $2n+1$, together with a differential ideal 
${\mathcal E}\subset\Omega^*(M)$\index{ideal!differential}, 
generated locally by the
contact ideal\index{contact!ideal} 
$\I$ and an $n$-form $\Psi\in\Omega^n(M)$.
\end{Definition}
Note that in this definition, the contact line
bundle\index{contact!line bundle} $I$ can be
recovered from ${\mathcal E}$ as its degree-$1$ part.
We can now pose a famous question.

\

\noindent
\textbf{Inverse Problem:}  {\em When is a given Monge-Ampere system
${\mathcal E}$ on $M$ equal to the
Euler-Lagrange system ${\mathcal E}_\Lambda$
of some Lagrangian $\Lambda\in\Omega^n(M)$?}

\

Note that if a given ${\mathcal E}$ does equal ${\mathcal E}_\Lambda$
for some $\Lambda$, then for some local generators $\theta,\Psi$ of
${\mathcal E}$ we must have $\theta\wedge\Psi=\Pi$, the
Poincar\'e-Cartan form of $\Lambda$.
Indeed, we can say that $(M,{\mathcal E})$ is
Euler-Lagrange if and only if there is an exact form
$\Pi\in\Omega^{n+1}(M)$, locally of the form $\theta\wedge\Psi$ for some
generators $\theta,\Psi$ of ${\mathcal E}$.
However, we face the difficulty that
$(M,{\mathcal E})$ does not determine either
$\Psi\in\Omega^n(M)$ or $\theta\in\Gamma(I)$ uniquely.

This can be partially overcome by normalizing $\Psi$ as follows.
Given only $(M,{\mathcal E}=\{\theta,d\theta,\Psi\})$, $\Psi$ is
determined as an element of $\bar\Omega^n=\Omega^n(M)/\I^n$.  We
can obtain a representative $\Psi$ that is unique modulo $\{I\}$
by adding the unique multiple of
$d\theta$ that yields a {\em primitive}\index{primitive}
form on $I^\perp$, referring to the symplectic decomposition of
$\bw{n}(T^*M/I)$ (see 
Proposition~\ref{Prop:Symplectic})\index{symplectic!linear algebra}.  
With this choice, we have a form
$\theta\wedge\Psi$ which is uniquely determined up to scaling; the
various multiples $f\theta\wedge\Psi$, where $f$ is a locally defined
function on $M$, are the candidates to be Poincar\'e-Cartan form.
Note that using a primitive normalization is reasonable, because our
actual Poincar\'e-Cartan forms $\Pi = \theta\wedge\Psi$ satisfy
$d\Pi=0$, which in particular implies that $\Psi$ is primitive on
$I^\perp$.  The proof of Noether's theorem\index{Noether's theorem} in
the next section will use a
more refined normalization of $\Psi$.

The condition for a Monge-Ampere system to be Euler-Lagrange is
therefore that there should be a globally defined exact $n$-form 
$\Pi$, locally of the form $f\theta\wedge\Psi$ with $\Psi$ normalized as
above.  This suggests the more accessible {\em local inverse problem},
which asks whether there is a {\em closed} $n$-form that is 
locally expressible as
$f\theta\wedge\Psi$.  It is for this local version that we give
a criterion.

We start with any candidate
Poincar\'e-Cartan form $\Xi = \theta\wedge\Psi$, and consider the
following criterion on $\Xi$:
\begin{equation}
\boxed{
  d\Xi = \varphi\wedge\Xi \quad \mbox{for some }\varphi
    \mbox{ with }d\varphi\equiv 0\ (\mbox{mod }\I).
}
\label{InverseCriterion}
\end{equation}

We first note that if this holds for some choice of
$\Xi=\theta\wedge\Psi$, then it holds for all other choices $f\Xi$;
this is easily verified.

Second, we claim that if (\ref{InverseCriterion}) holds,
then we can find $\tilde\varphi$ also
satisfying $d\Xi=\tilde\varphi\wedge\Xi$, and in addition,
$d\tilde\varphi = 0$.  To see this, write
$$
  d\varphi = \theta\wedge\alpha + \beta\,d\theta
$$
(here $\alpha$ is a $1$-form and $\beta$ is a function),
and differentiate using $d^2=0$, modulo the algebraic ideal $\{I\}$, to obtain
$$
  0 \equiv d\theta\wedge(\alpha+d\beta)\pmod{\{I\}}.
$$
But with the standing assumption $n\geq 2$, symplectic linear
algebra\index{symplectic!linear algebra}
implies that the $1$-form $\alpha+d\beta$ must vanish modulo
$\{I\}$.  As a result,
$$
  d(\varphi-\beta \, \theta)=\theta\wedge(\alpha+d\beta) = 0,
$$
so we can take $\tilde\varphi = \varphi-\beta\,\theta$, verifying
the claim.

Third, once we know that $d\Xi=\varphi\wedge\Xi$ with $d\varphi=0$,
then on a possibly smaller neighborhood, we use the Poincar\'e
lemma\index{Poincar\'e lemma} to
write $\varphi = du$ for a function $u$, and then
$$
  d(e^{-u}\Xi) = e^{-u}(\varphi\wedge\Xi-du\wedge\Xi) = 0.
$$  This proves the following.

\begin{Theorem}
A Monge-Ampere system $(M,{\mathcal E}=\{\theta,d\theta,\Psi\})$ on a
$(2n+1)$-dimensional 
contact manifold $M$ with $n\geq 2$, where $\Psi$ is assumed to 
be primitive\index{primitive} modulo $\{I\}$, is locally equal to an Euler-Lagrange system
$\mathcal E_\Lambda$ if and
only if it satisfies (\ref{InverseCriterion}).
\end{Theorem}
\index{Poincar\'e-Cartan form|)}

\noindent
\textbf{Example 1.}
\index{Poisson equation|(}
Consider a scalar PDE of the form 
\begin{equation}
  \Delta z=f(x,z,\nabla z),
\label{LinPoisson}
\end{equation}
where $\Delta = \sum\frac{\partial^2}{\partial x^{i2}}$;
we ask which functions $f:\R^{2n+1}\to\R$ are such that
(\ref{LinPoisson}) is contact-equivalent to an Euler-Lagrange
equation.  To apply our framework, we let
$M=J^1(\R^n,\R)$, $\theta=dz-\sum p_idx^i$ so $d\theta=-\sum
dp_i\wedge dx^i$, and set
$$
  \Psi  =  \ss dp_i\wedge dx_{(i)} - f(x,z,p)dx.
$$
Restricted to a Legendre submanifold\index{Legendre submanifold} of the form
$N=\{(x^i,z(x),\frac{\partial z}{\partial x^i}(x)\}$, we find
$$
  \Psi|_N = (\Delta z-f(x,z,\nabla z))dx.
$$
Evidently $\Psi$ is primitive\index{primitive} modulo $\{I\}$, and $\mathcal
E=\{\theta,d\theta,\Psi\}$ is a
Monge-Ampere system whose transverse integral
manifolds\index{integral manifold!transverse} 
(i.e., those on which $dx^1\wedge\cdots\wedge dx^n\neq 0$) correspond to
solutions of the equation (\ref{LinPoisson}).  To apply our
test, we start with the candidate $\Xi=\theta\wedge\Psi$, for which
$$
  d\Xi = -\theta\wedge d\Psi = \theta\wedge df\wedge dx.
$$
Therefore, we consider $\varphi$ satisfying
$$
  \theta\wedge df\wedge dx = \varphi\wedge\Xi,
$$
or equivalently
$$
  df\wedge dx\equiv -\varphi\wedge\Psi \pmod{\{I\}},
$$
and find that they are exactly those $1$-forms of the form
$$
  \varphi = \ss f_{p_i}dx^i + c\,\theta
$$
for an arbitrary function $c$.  The problem is reduced to describing
those $f(x,z,p)$ for which there exists some $c(x,z,p)$ so that $\varphi
= \sum f_{p_i}dx^i+c\,\theta$ is closed.  We can
determine all such forms explicitly, as follows.  
The condition that $\varphi$ be closed expands to
\begin{eqnarray*}
  0 & = & c_{p_i}dp_i\wedge dz \\
      & & \ + (f_{p_ip_j}-c\delta_i^j-c_{p_j}p_i)dp_j\wedge dx^i \\
      & & \ + \sf12(f_{p_ix^j}-f_{p_jx^i}-c_{x^j}p_i+c_{x^i}p_j)
                   dx^j\wedge dx^i \\
      & & \ + (f_{p_iz}-c_{x^i}-c_zp_i)dz\wedge dx^i.
\end{eqnarray*}
These four terms must vanish separately.  The vanishing of the first
term
implies that $c=c(x^i,z)$ does not depend on any $p_i$.  Given this,
the vanishing of the second term implies that $f(x^i,z,p_i)$ is
quadratic in the $p_i$, with diagonal leading term:
$$
  f(x^i,z,p_i) = \textstyle\frac12\displaystyle c(x^i,z)\ss p_j^2
       +\ss e^j(x^i,z)p_j + a(x^i,z)
$$
for some functions $e^j(x^i,z)$ and $a(x^i,z)$.  Now the vanishing of
the third term reduces to
$$
  0 = e^i_{x^j}-e^j_{x^i},
$$
implying that for some function $b(x^j,z)$, 
$$
  e^j(x^i,z) = \frac{\partial b(x^i,z)}{\partial x^j};
$$
this $b(x^j,z)$ is uniquely determined only up to addition of a
function of $z$. Finally, the vanishing of the fourth term reduces to
$$
  (b_z - c)_{x^i} = 0,
$$
so that $c(x^i,z)$ differs from $b_z(x^i,z)$ by a function of $z$
alone.  By adding an antiderivative of this difference to $b(x^i,z)$
and relabelling the result as $b(x^i,z)$, we see that our criterion for the
Monge-Ampere system to be Euler-Lagrange is that $f(x^i,z,p_i)$ be of
the form
$$
  f(x^i,z,p_i) = \sf12b_z(x,z)\ss p_i^2 + \ss b_{x^i}(x,z)p_i+a(x,z)
$$
for some functions $b(x,z)$, $a(x,z)$.  These describe exactly those
Poisson equations that are locally contact-equivalent to
Euler-Lagrange equations.
\index{Poisson equation|)}

\

\noindent
\textbf{Example 2.}
An example that is not quasi-linear is given by
$$
  \mbox{det}(\nabla^2z) - g(x,z,\nabla z) = 0.
$$
The $n$-form $\Psi = dp - g(x,z,p)dx$ and the standard contact system
generate a Monge-Ampere system whose transverse integral manifolds
correspond to solutions of this equation.  A calculation similar to
that in the preceding example shows that this Monge-Ampere system is
Euler-Lagrange if and only if $g(x,z,p)$ is of the form
$$
  g(x,z,p) = g_0(x,z)\, g_1(p, z-\ss p_ix^i).
$$

\

\noindent
\textbf{Example 3.}
The linear Weingarten equation\index{Weingarten equation} $aK+bH+c=0$
for a surface in
Euclidean space having Gauss curvature\index{Gauss curvature} $K$ and
mean curvature\index{mean curvature} $H$ is
Euler-Lagrange for all choices of constants $a,b,c$, as we shall see
in \S\ref{Subsection:Euclidean-invariant}.  
In this case, the appropriate contact manifold for the problem is
$M=G_2(T{\mathbf E}^3)$, the Grassmannian\index{Grassmannian} of
oriented tangent planes of Euclidean space.

\

\noindent
\textbf{Example 4.}
Here is an example of a Monge-Ampere system which is locally, but not
globally, Euler-Lagrange, suitable for those readers familiar with
some complex algebraic geometry.  Let $X$ be a K3 surface\index{K3
  surface}; that is, $X$ is a 
simply connected, compact, complex manifold of complex
dimension $2$ with trivial
canonical bundle, necessarily of K\"ahler type.
Suppose also that there is a positive holomorphic
line bundle $L\to X$ with a Hermitian metric having positive first
Chern form
$\omega\in\Omega^{1,1}(X)$.  Our contact manifold $M$ is the unit
circle subbundle of $L\to X$, a smooth manifold of real dimension
$5$; the contact form is
$$
  \theta = \sf{i}{2\pi}\alpha,\qquad d\theta = \omega,
$$
where $\alpha$ is the $\lie{u}(1)$-valued Hermitian connection form on $M$.
Note $\theta\wedge(d\theta)^2\neq 0$, because the $4$-form
$(d\theta)^2=\omega^2$ is actually a volume form on $M$ (by
positivity) and $\theta$ is non-vanishing on fibers of $M\to X$,
unlike $(d\theta)^2$.

Now we trivialize the canonical bundle of $X$ with a holomorphic
$2$-form $\Phi = \Psi + i\Sigma$, and take for our Monge-Ampere system
$$
  {\mathcal E}=\{\theta,d\theta=\omega,\Psi=\mbox{Re}(\Phi)\}.
$$
We can see that $\mathcal E$ is locally Euler-Lagrange as follows.  First, by
reasons of type, $\omega\wedge\Phi=0$; and $\omega$ is
real, so $0=\mbox{Re}(\omega\wedge\Phi)=\omega\wedge\Psi$.  In
particular, $\Psi$ is primitive.  With
$\Xi=\theta\wedge\Psi$, we compute
$$
  d\Xi = \omega\wedge\Psi-\theta\wedge d\Psi = -\theta\wedge d\Psi,
$$
but $d\Psi = \mbox{Re}(d\Phi) = 0$, because $\Phi$ is holomorphic and
therefore closed.

On the other hand, $(M,{\mathcal E})$ cannot be globally
Euler-Lagrange; that is,
$\Xi=\theta\wedge\Psi$ cannot be exact, for if $\Xi=d\xi$, then
$$
  \int_M\Xi\wedge\Sigma = \int_Md(\xi\wedge\Sigma) = 0,
$$
but also
$$
  \int_M\Xi\wedge\Sigma = \int_M\theta\wedge\Psi\wedge\Sigma
     = c\int_X\Phi\wedge\bar\Phi,
$$
for some number $c\neq 0$.
\index{inverse problem|)}
\index{Euler-Lagrange!system|)}
\index{Monge-Ampere system|)}

\section{Noether's Theorem}
\label{Section:Noether}

\index{Noether's theorem|(}
\index{symmetry|(}
\index{conservation law|(}

The classical theorem of Noether describes an isomorphism between a Lie
algebra of infinitesimal symmetries
associated to a variational problem, and a space of
conservation laws for its Euler-Lagrange equations.  We will often
assume without comment that our Lagrangian is non-degenerate in the
sense discussed earlier.

There are four reasonable Lie algebras of symmetries that we might
consider in our setup.  Letting $\mathcal V(M)$ denote the Lie algebra
of all vector fields on $M$, they are the following.
\begin{itemize}
\item
  Symmetries of $(M, I, \Lambda)$:
  $$
    \lie{g}_\Lambda = \{v\in{\mathcal V}(M):{\mathcal L}_vI\subseteq I,\ 
       {\mathcal L}_v\Lambda = 0\}.
  $$
\item
  Symmetries of $(M, I, [\Lambda])$:
  $$
    \lie{g}_{[\Lambda]} = \{v\in{\mathcal V}(M):{\mathcal L}_vI\subseteq I,\
      {\mathcal L}_v[\Lambda] = 0\}.
  $$
\item
  Symmetries of $(M, \Pi)$:
  $$
    \lie{g}_\Pi = \{v\in{\mathcal V}(M):{\mathcal L}_v\Pi = 0\}.
   $$
      (Note that ${\mathcal L}_v\Pi = 0$ implies
           ${\mathcal L}_vI\subseteq I$ for non-degenerate $\Lambda$.)
\item 
  Symmetries of $(M, {\mathcal E}_\Lambda)$:
  $$
    \lie{g}_{{\mathcal E}_\Lambda} = 
       \{v\in{\mathcal V}(M):{\mathcal L}_v{\mathcal E}_\Lambda \subseteq 
                {\mathcal E}_\Lambda\}.
    $$
     (Note that ${\mathcal L}_v{\mathcal E}_\Lambda \subseteq
                          {\mathcal E}_\Lambda$ implies
                        ${\mathcal L}_vI\subseteq I$.)
\end{itemize}

We comment on the relationship between these spaces.  Clearly, there
are inclusions
$$
  \lie{g}_\Lambda \subseteq \lie{g}_{[\Lambda]} \subseteq
    \lie{g}_\Pi \subseteq \lie{g}_{{\mathcal E}_\Lambda}.
$$
Any of the three inclusions may be strict.  For example, we
{\em locally} have $\lie{g}_{[\Lambda]}=\lie{g}_\Pi$ because $\Pi$ is
the image of $[\Lambda]$ under the coboundary
$\delta\!:\!H^{n}(\Omega^*/\I)\to H^{n+1}(\I)$, which is invariant under
diffeomorphisms of $(M,I)$ and is an isomorphism on contractible
open sets.  However, we shall see later that globally there is an
inclusion
$$
  \lie{g}_\Pi/\lie{g}_{[\Lambda]} \hookrightarrow H^n_{dR}(M),
$$
and this discrepancy between the two symmetry algebras introduces some 
subtlety into Noether's theorem.

Also, there is a bound
\begin{equation}
  \mbox{dim}\left(\lie{g}_{{\mathcal E}_\Lambda}/\lie{g}_\Pi\right)
    \leq 1.
\label{SymmBound}
\end{equation}
This follows from noting that if a vector field $v$ preserves
${\mathcal E}_\Lambda$, then it preserves $\Pi$ up to multiplication
by a function; that is, ${\mathcal L}_v\Pi = f\Pi$.  Because $\Pi$ is
a closed form, we find that $df\wedge\Pi = 0$; in the non-degenerate
case, this implies $df=u\theta$ for some function $u$.  The definition 
of a contact form\index{contact!form} prohibits any $u\theta$ from
being closed unless
$u=0$, meaning that $f$ is a constant.  This constant gives a linear
functional on $\lie{g}_{{\mathcal E}_\Lambda}$ whose kernel is $\lie{g}_\Pi$,
proving (\ref{SymmBound}).
The area\index{area functional} functional and minimal 
surface\index{minimal surface} equation for Euclidean
hypersurfaces provide an example where the two spaces are
different.  In that case, the induced Monge-Ampere system is invariant not only
under Euclidean motions, but under dilations of Euclidean space as 
well; this is not true of the Poincar\'e-Cartan form.

The next step in introducing Noether's theorem is to describe the
relevant spaces of conservation laws.  In general, suppose that $(M,{\mathcal
J})$ is an exterior differential system with integral
manifolds\index{integral manifold} of
dimension $n$.  A {\em conservation law} for $(M,{\mathcal
J})$ is an $(n-1)$-form $\varphi\in\Omega^{n-1}(M)$ such that
$d(f^*\varphi) = 0$ for every integral manifold $f:N^n\hookrightarrow M$
of ${\mathcal J}$.  Actually, we will only consider as conservation
laws those $\varphi$
on $M$ such that $d\varphi\in{\mathcal J}$, which may be a strictly
smaller set.  This will not present any liability, as one can
always ``saturate'' ${\mathcal J}$ to remove this discrepancy.
The two apparent
ways in which a conservation law may be {\em
trivial}\index{conservation law!trivial} are when either $\varphi\in{\mathcal J}^{n-1}$ already or
$\varphi$ is exact on $M$.  Factoring out these cases
leads us to the following.
\begin{Definition}
The space of {\em conservation laws} for $(M,{\mathcal J})$ is
$$
  {\mathcal C} = H^{n-1}(\Omega^*(M)/{\mathcal J}).
$$
\end{Definition}

It also makes sense to factor out those conservation laws represented by
$\varphi\in\Omega^{n-1}(M)$ which are already closed on $M$, and not
merely on integral manifolds of ${\mathcal J}$.  This can be
understood using the long exact sequence:
$$
  \cdots\to H^{n-1}_{dR}(M) \stackrel{\pi}{\to} {\mathcal C} \to
     H^n({\mathcal J}) \to H^n_{dR}(M) \to\cdots.
$$
\begin{Definition}
The space of {\em proper conservation 
laws}\index{conservation law!proper} is $\bar{\mathcal C} =
{\mathcal C}/\pi(H^{n-1}_{dR}(M))$.
\end{Definition}
Note that there is an inclusion $\bar{\mathcal C}\hookrightarrow
H^n({\mathcal J})$.  In case $\mathcal J=\mathcal E_\Lambda$ is the
Euler-Lagrange system\index{Euler-Lagrange!system} of a non-degenerate
functional\index{non-degenerate!functional} $\Lambda$ on a
contact manifold $(M,I)$, we have the following.
\begin{Theorem}[Noether]
Let $(M,{\mathcal E}_\Lambda)$ be the Euler-Lagrange system of a
non-degenerate functional $\Lambda$.  There is a 
linear isomorphism
$$
  \eta:\lie{g}_\Pi \to H^n({\mathcal E}_\Lambda),
$$
taking the subalgebra $\lie{g}_{[\Lambda]}\subset\lie{g}_\Pi$ to the
subspace $\eta(\lie{g}_{[\Lambda]})=\bar{\mathcal C} \subset
H^n({\mathcal E}_\Lambda).$
\label{NoetherTheorem}
\end{Theorem}

Before proceeding to the proof, which will furnish an explicit formula
for $\eta$, we need to make a digression on the
algebra of infinitesimal contact
transformations\index{transformation!contact|(}
$$
  \lie{g}_I = \{v\in{\mathcal V}(M):{\mathcal L}_vI\subseteq I\}.
$$
The key facts are that on any neighborhood where $I$ has a non-zero
generator $\theta$, a contact symmetry $v$ is uniquely determined by
its so-called {\em generating function}\index{generating function|(}
$g=v\innerprod\theta$, and
that given such $\theta$, any function $g$ is the generating function of some
$v\in\lie{g}_I$.  This
can be seen on a possibly smaller neighborhood by taking Pfaff
coordinates with $\theta = dz-\sum p_idx^i$.  Working in a basis
$\partial_\theta,\partial^i,
\partial_i$ dual to the basis $\theta,dp_i,dx^i$ of $T^*M$,
we write
$$
  v=g\,\partial_\theta + \ss v^i\partial_i
      + \ss v_i\partial^i.
$$
Now the condition
$$
  {\mathcal L}_v\theta\equiv 0\pmod{\{I\}}
$$
can be made explicit, and it turns out to be
$$
  v_i = \partial_ig = \left(\frac{\partial}{\partial x^i} + p_i
     \frac{\partial}{\partial z}\right)g,\quad
  v^i = -\partial^ig = -\frac{\partial g}{\partial p_i}.
$$
This establishes our claim, because the correspondence between $v$ and 
$g$ is now given by
\begin{equation}
  v = g\partial_\theta - \ss g_{p_i}\partial_i +
    \ss (g_{x^i}+p_ig_z)\partial^i.
\label{CanonicalOperator}
\end{equation}

As we have presented it, the correspondence between infinitesimal
contact symmetries and their generating functions is local.  But a simple
patching argument shows that globally, as one moves between different
local generators $\theta$ for $I$, the different generating functions
$g$ glue together to give a global section $g\in\Gamma(M,I^*)$ of the
dual line bundle\index{contact!line bundle}.  In fact, the formula (\ref{CanonicalOperator})
describes a canonical splitting of the surjection
$$
  \Gamma(TM)\to \Gamma(I^*)\to 0.
$$
Note that this splitting is not a bundle map, but a differential operator.

Returning to Noether's theorem,
the proof that we present is slightly incomplete in that we assume
given a global non-zero contact form\index{contact!form}
$\theta\in\Gamma(I)$, or 
equivalently, that the contact line bundle\index{contact!line bundle}
is trivial.  This
allows us to treat generating functions of contact symmetries as
functions rather than as sections of $I^*$. 
It is an enlightening exercise to develop the patching
arguments needed to 
overcome this using sheaf cohomology.
Alternatively, one can
simply pull everything up to a double cover of $M$ on which $I$ has
a global generator, and little will be lost.
\index{transformation!contact|)}

\

\noindent
\textbf{Proof of Theorem \ref{NoetherTheorem}.}
{\em Step 1: Definition of the map $\eta$.}
The map in question is given by
$$
  \eta(v) = v\innerprod\Pi\qquad \mbox{for }v\in\lie{g}_\Pi\subset
     {\mathcal V}(M).
$$
Note that locally $v\innerprod\Pi = (v\innerprod\theta)\Psi -
\theta\wedge(v\innerprod\Psi)$, so that $v\innerprod\Pi$ lies in ${\mathcal
E}_\Lambda$.  Furthermore, the condition ${\mathcal L}_v\Pi = 0$ gives
$$
  0 = v\innerprod d\Pi + d(v\innerprod\Pi) = d(v\innerprod\Pi),
$$
so that $v\innerprod\Pi$ is closed, and gives a well-defined class
$\eta(v)\in H^n({\mathcal E}_\Lambda)$.

\

\noindent
{\em Step 2: $\eta$ is injective.}
Write
$$
  \eta(v) = (v\innerprod\theta)\Psi - \theta\wedge(v\innerprod\Psi).
$$
Suppose that this $n$-form is cohomologous to zero in $H^n({\mathcal 
E}_\Lambda)$; that is,
\begin{eqnarray*}
  (v\innerprod\theta)\Psi - \theta\wedge(v\innerprod\Psi) & = &
    d(\theta\wedge\alpha+d\theta\wedge\beta) \\
  & = & -\theta\wedge d\alpha + d\theta\wedge(\alpha+d\beta).
\end{eqnarray*}
Regarding this equation modulo $\{I\}$ and using the primitivity of
$\Psi$, we conclude that
$$
  v\innerprod\theta = 0.
$$
An infinitesimal symmetry $v\in\lie{g}_I$ of the contact system
is locally determined by its generating 
function\index{generating function|)} $v\innerprod\theta$
as in (\ref{CanonicalOperator}), so we conclude that $v=0$, proving
injectivity.

\

\noindent
{\em Step 3: $\eta$ is locally surjective.}
We start by representing a class in
$H^n({\mathcal E}_\Lambda)$ by a closed $n$-form
\begin{equation}
  \Phi = g\Psi + \theta\wedge\alpha.
\label{BadPresentation}
\end{equation}
We can choose the unique contact vector field $v$
such that
$v\innerprod\theta = g$, and our goals are to show that
$v\in\lie{g}_\Pi$ and that $\eta(v) =
[\Phi]\in H^n({\mathcal E}_\Lambda)$.

For this, we need a special choice of $\Psi$, which so far is
determined only modulo $\{I\}$; this is reasonable because the
presentation (\ref{BadPresentation}) is not unique.  In fact, we 
can further normalize $\Psi$ by the condition 
$$
  d\theta\wedge\Psi = 0.
$$
To see why this is so, first note that by symplectic linear
algebra\index{symplectic!linear algebra|(}
(Proposition~\ref{Prop:Symplectic}), 
\begin{equation}
  d\Psi\equiv d\theta\wedge\Gamma \pmod{\{I\}},
\label{Hypo1}
\end{equation}
for some $\Gamma$, because $d\Psi$ is of degree $n+1$.  Now suppose we
replace $\Psi$ by $\bar\Psi = \Psi - \theta\wedge\Gamma$, which
certainly preserves the essential condition $\Pi =
\theta\wedge\bar\Psi$.  Then we have
\begin{eqnarray*}
  d\theta\wedge\bar\Psi & = & d\theta\wedge(\Psi-\theta\wedge\Gamma)\\ 
    & = & (d\Pi +\theta\wedge d\Psi) - d\theta\wedge\theta\wedge\Gamma 
             \\
    & = & \theta\wedge(d\Psi - d\theta\wedge\Gamma) \\
    & = & 0,\mbox{ by (\ref{Hypo1})},
\end{eqnarray*}
and we have obtained our refined normalization.

Now we combine the following three equations modulo 
$\{I\}$:
\begin{itemize}
\item $0 \equiv {\mathcal L}_v\theta \equiv dg + v\innerprod d\theta$, when
multiplied by $\Psi$, gives
$$
  dg\wedge\Psi + (v\innerprod d\theta)\wedge\Psi \equiv 0 \pmod{\{I\}};
$$
\item $0 = d\Phi = d(g\Psi + \theta\wedge\alpha),$ so using our
normalization condition $d\theta\wedge\Psi=0$  (which implies $d\Psi
\equiv 0\ (\mbox{mod }\{I\})$),
$$
  dg\wedge\Psi + d\theta\wedge\alpha \equiv 0 \pmod{\{I\}};
$$
\item Ordinary primitivity gives $d\theta\wedge\Psi \equiv 0
  \ (\mbox{mod }\{I\})$, and contracting with $v$,
$$
  (v\innerprod d\theta)\wedge\Psi + d\theta\wedge(v\innerprod\Psi)
    \equiv 0 \pmod{\{I\}}.
$$
\end{itemize}
These three equations combine to give
$$
  d\theta\wedge(\alpha + v\innerprod\Psi) \equiv 0 \pmod{\{I\}},
$$
and from symplectic linear algebra\index{symplectic!linear algebra|)}, we have
$$
  \alpha + v\innerprod\Psi \equiv 0 \pmod{\{I\}}.
$$
This allows us to conclude
\begin{equation}
  v\innerprod\Pi = g\Psi + \theta\wedge\alpha
\label{Helpful3}
\end{equation}
which would complete the proof of surjectivity, except that we have
not yet shown that ${\mathcal L}_v\Pi = 0$.  However, by hypothesis
$g\Psi+\theta\wedge\alpha$ is closed; with (\ref{Helpful3}), 
this is enough to compute ${\mathcal L}_v\Pi=0$.

The global isomorphism asserted in the theorem follows easily from
these local conclusions, so long as we maintain the assumption that
there exists a global contact form. 

\

\noindent
{\em Step 4: $\eta$ maps symmetries of $[\Lambda]$ to proper
conservation laws.}\index{conservation law!proper|(}
For this, first note that there
is an exact sequence
$$
  0\to \bar{\mathcal C} \to H^n({\mathcal E}_\Lambda) \stackrel{i}{\to} H^n_{dR}(M)
    \to \cdots
$$
so it suffices to show that for $v\in\lie{g}_\Pi$,
\begin{equation}
  {\mathcal L}_v[\Lambda]=0\mbox{  if and only if }
    \eta(v)\in\mbox{Ker }i.
\label{StrongCondn}
\end{equation}
Recall that $\Pi = d(\Lambda-\theta\wedge\beta)$ for some
$\beta$, and we can therefore calculate
\begin{eqnarray*}
  \eta(v) & = & v\innerprod\Pi\\ & = &
    {\mathcal L}_v(\Lambda-\theta\wedge\beta) - 
    d(v\innerprod(\Lambda-\theta\wedge\beta))\\
    & \equiv & {\mathcal L}_v\Lambda \pmod {d\Omega^{n-1}(M) + \I^n}.
\end{eqnarray*}
This proves that $j\circ i(\eta(v)) = {\mathcal L}_v[\Lambda]$ in the
composition 
$$
  H^n({\mathcal E}_\Lambda)\stackrel{i}{\to} H^n_{dR}(M)
    \stackrel{j}{\to}
    H^n(\Omega^*(M)/\I).
$$
The conclusion (\ref{StrongCondn}) will follow if we can prove
that $j$ is injective.

To see that $j$ is injective, note that it occurs in the long exact 
cohomology sequence of
$$
  0 \to \I\to \Omega^*(M) \to \Omega^*(M)/\I \to 0;
$$
namely, we have
$$
  \cdots\to H^n(\I)\to H^n_{dR}(M)\stackrel{j}{\to}
      H^n(\Omega^*(M)/\I)\to\cdots.
$$
So it suffices to show that $H^n(\I)=0$, which we will do under the standing
assumption that there is a global contact form $\theta$.
Suppose that the $n$-form
$$
  \varphi = \theta\wedge\alpha+d\theta\wedge\beta
     =\theta\wedge(\alpha+d\beta)+d(\theta\wedge\beta)
$$
is closed.
Then regarding $0=d\varphi$ modulo $\{I\}$, we have by symplectic linear
algebra that
$$
  \alpha+d\beta\equiv 0\pmod{\{I\}}
$$
so that actually
$$
  \varphi = d(\theta\wedge\beta).
$$
This says that $\varphi\sim 0$ in $H^n(\I)$, and our proof is complete.
\hfill$\square$

\

It is important in practice to have a local formula for a
representative in $\Omega^{n-1}(M)$, closed modulo ${\mathcal
E}_\Lambda$, for the proper conservation 
law\index{conservation law!proper|)} $\eta(v)$.
This is obtained by first writing as usual
\begin{equation}
  \Pi = d(\Lambda - \theta\wedge\beta),
\label{NPrescription1}
\end{equation}
and also, for a given $v\in\lie{g}_{[\Lambda]}$,
\begin{equation}
  {\mathcal L}_v\Lambda \equiv d\gamma \pmod{\I}.
\label{NPrescription2}
\end{equation}
We will show that the $(n-1)$-form
\begin{equation}
\boxed{
  \varphi = -v\innerprod\Lambda + (v\innerprod\theta)\beta+\gamma
}
\label{NPrescription3}
\end{equation}
is satisfactory.  First, compute
\begin{eqnarray*}
  d\varphi & = & (-{\mathcal L}_v\Lambda + v\innerprod d\Lambda)
                 +d((v\innerprod\theta)\beta)
                 +d\gamma \\
           & \equiv & v\innerprod(\Pi + d(\theta\wedge\beta))
                  +d((v\innerprod\theta)\beta) \pmod{\I} \\
           & \equiv & \eta(v) + {\mathcal L}_v(\theta\wedge\beta)
                        \pmod{\I} \\
           & \equiv & \eta(v) \pmod{\I}.
\end{eqnarray*}
Now we have $d\varphi = \eta(v)+\Xi$ for some closed
$\Xi\in\I^n$.  We proved in the last part of the proof of Noether's
theorem that $H^n(\I)=0$, which implies that $\Xi=d\xi$ for some
$\xi\in\I^{n-1}$.  Now we have
$d(\varphi-\xi) = \eta(v)$, and $\varphi\sim\varphi-\xi$ in ${\mathcal 
C}=H^{n-1}(\Omega^*(M)/{\mathcal E}_\Lambda)$.  This justifies our
prescription (\ref{NPrescription3}).

Note that the prescription is especially simple when
$v\in\lie{g}_\Lambda\subseteq\lie{g}_{[\Lambda]}$, for then we can
take $\gamma = 0$.

\

\noindent
\textbf{Example.}
Let ${\mathbf L}^{n+1} = \{(t,y^1,\ldots,y^n)\}\cong\R^{n+1}$ be
Minkowski\index{Minkowski space} space, and let $M^{2n+3}=J^1({\mathbf
L}^{n+1},\R)$ be the standard
contact manifold, with coordinates $(t,y^i,z,p_a)$ (where
$0\leq a\leq n$), $\theta = dz-p_0dt-\sum p_idy^i$.  For a Lagrangian, 
take
$$
  \Lambda = \left(\sf12 ||p||^2+F(z)\right)dt\wedge dy
$$
for some ``potential'' function $F(z)$, where $dy =
  dy^1\wedge\cdots\wedge dy^n$ and  $||p||^2 =
-p_0^2+\sum p_i^2$ is the Lorentz-signature norm.  The local
symmetry group of this
functional is generated by two subgroups, the translations in
${\mathbf L}^{n+1}$ and the linear isometries $SO^o(1,n)$; as we shall
see in Chapter~\ref{Chapter:Conformal}, for certain $F(z)$ the
symmetry group of the associated Poincar\'e-Cartan form is strictly
larger.  For now, we 
calculate the conservation law corresponding to translation
in $t$, and begin by finding the Poincar\'e-Cartan form $\Pi$.
Letting $f(z) = F^\prime(z)$, we differentiate
\begin{eqnarray*}
  d\Lambda & = & (-p_0dp_0+\ss p_idp_i + f(z)\theta)\wedge
     dt\wedge dy \\
  & = & \theta\wedge(f(z)dt\wedge dy) + d\theta\wedge
        (p_0dy+\ss p_j dt\wedge dy_{(j)});
\end{eqnarray*}
with the usual recipe $\Pi = \theta\wedge(\alpha+d\beta)$ whenever
$d\Lambda = \theta\wedge\alpha + d\theta\wedge\beta$, we obtain
$$
  \Pi = \theta\wedge\left(f(z)dt\wedge dy + dp_0\wedge dy + 
         \ss dp_j\wedge dt\wedge dy_{(j)}\right).
$$
We see the Euler-Lagrange equation using
$$
  \Psi = f(z)dt\wedge dy + dp_0\wedge dy + \ss dp_j\wedge 
     dt\wedge dy_{(j)};
$$
an integral manifold of ${\mathcal
E}_\Lambda=\{\theta,d\theta,\Psi\}$ of the form
$$
  \{(t,y,z(t,y),z_t(t,y),z_{y^i}(t,y))\}
$$
 must satisfy
$$
  0 = \Psi|_N = \left(\frac{\partial^2z}{\partial t^2} -
    \sum\frac{\partial^2z}{\partial y^{i2}} + f(z)\right)
    dt\wedge dy.
$$
With the independence condition\index{independence condition}
$dt\wedge dy\neq 0$, we have the
familiar wave equation\index{wave equation}
$$
  \square z(t,y) = f(z).
$$
Now considering the time-translation symmetry $v=\frac{\partial}{\partial
t}\in\lie{g}_\Lambda$, the Noether prescription (\ref{NPrescription3}) gives
\begin{eqnarray*}
  \varphi & = & -v\innerprod\Lambda + (v\innerprod\theta)\beta \\
  & = & -\left(\frac{||p||^2}{2}+F(z)\right)dy - 
     p_0\left(p_0dy+\ss p_jdt\wedge dy_{(j)}\right) \\
  & = & -\left(\sf{1}{2}\ss p_a^2+F(z)\right)dy - 
      p_0dt\wedge(\ss p_jdy_{(j)}).
\end{eqnarray*}
One can verify that $\varphi$ is closed when restricted to a
solution of $\square z = f(z)$.
The question of how one might use this conservation law will be taken
up later.

\index{Noether's theorem|)}
\index{symmetry|)}
\index{conservation law|)}

\section{Hypersurfaces in Euclidean Space}
\label{Section:Hypersurface}
\index{Euclidean!space|(}

We will apply the the theory developed so far to the study of
hypersurfaces in Euclidean space
$$
  N^n\hookrightarrow{\mathbf E}^{n+1}.
$$
We are particularly interested in the study of those functionals on
such hypersurfaces which are invariant under the group $E(n+1)$ of
orientation-preserving Euclidean motions.\index{Euclidean!motion|(}

\subsection{The Contact Manifold over ${\mathbf E}^{n+1}$}

Points of ${\mathbf E}^{n+1}$ will be
denoted $x=(x^0,\ldots,x^n)$,
and each tangent space $T_x{\mathbf E}^{n+1}$ will be canonically
identified with ${\mathbf E}^{n+1}$ itself via translation.  A {\em
frame}\index{Euclidean!frame bundle} for ${\mathbf E}^{n+1}$ is a pair
$$
  f = (x,e)
$$
consisting of a point $x\in{\mathbf E}^{n+1}$ and a
positively-oriented orthonormal
basis $e=(e_0,\ldots,e_n)$ for $T_x{\mathbf E}^{n+1}$.
The set $\mathcal F$ of all such frames is a manifold, and the
right $SO(n+1,\R)$-action
$$
   (x,(e_0,\ldots,e_n))\cdot (g^a_b) = (x,(\ss e_ag^a_0,\ldots,
       \ss e_ag^a_n))
$$
gives the basepoint map
$$
  x:{\mathcal F}\to{\mathbf E}^{n+1}
$$
the structure of a principal bundle.\footnote{Throughout this section,
  we use index ranges $1\leq i,j\leq n$ and $0\leq a,b\leq n$.}  There
is also an obvious left-action of $E(n+1,\R)$ on 
${\mathcal F}$, and a choice of reference frame gives a left-equivariant
identification ${\mathcal F}\cong E(n+1)$ of the bundle of frames with 
the group of Euclidean motions.

The relevant contact manifold for studying hypersurfaces in ${\mathbf
E}^{n+1}$ is the manifold of {\em contact elements}
$$
  M^{2n+1} = \{(x,H):x\in{\mathbf E}^{n+1},\ H^n\subset T_x{\mathbf
    E}^{n+1} \mbox{ an oriented hyperplane}\}.
$$
This $M$ will be given the structure of a contact manifold in such a
way that transverse Legendre submanifolds correspond to arbitrary
immersed hypersurfaces in ${\mathbf E}^{n+1}$.  Note that $M$ may be
identified with the unit sphere bundle of ${\mathbf E}^{n+1}$ by
associating to a contact element $(x,H)$ its oriented orthogonal
complement $(x,e_0)$.  We will use this identification without further
comment.

The projection ${\mathcal F}\to M$ taking $(x,(e_a))\mapsto
(x,e_0)$ is $E(n+1,\R)$-equivariant (for the left-action).  To describe the
contact structure on $M$ and to carry out calculations, we will
actually work on ${\mathcal F}$ using the following structure equations.
First, we define canonical $1$-forms on ${\mathcal F}$ by differentiating the
vector-valued coordinate functions $x(f),e_a(f)$ on ${\mathcal F}$,
and decomposing the
resulting vector-valued $1$-forms at each $f\in{\mathcal F}$ with respect to
the frame $e_a(f)$:
\begin{equation}
   dx = \sum e_b\cdot\omega^b,\qquad de_a = \sum e_b\cdot\omega^b_a.
\label{EucForms}
\end{equation}
Differentiating the relations $\langle e_a(f),e_b(f)\rangle =
\delta_{ab}$ yields
$$
  \omega^a_b+\omega^b_a = 0.
$$
The forms $\omega^a$, $\omega^a_b$ satisfy no other linear algebraic relations,
giving a total of $(n+1)+\frac12n(n+1)=\mbox{dim}({\mathcal F})$ independent
$1$-forms.  By taking the derivatives of the defining relations
(\ref{EucForms}), we obtain the structure equations
\begin{equation}
  d\omega^a +\sum \omega^a_c\wedge\omega^c = 0,\qquad
  d\omega^a_b+\sum \omega^a_c\wedge\omega^c_b = 0.
\label{EucStructure}
\end{equation}
The forms $\omega^a$ are identified with the usual tautological
$1$-forms on the orthonormal frame 
bundle\index{Riemannian!frame bundle} of a Riemannian
manifold\index{Riemannian!manifold} (in
this case, of ${\mathbf E}^{n+1}$); and then the first equation
indicates that $\omega^a_b$ are components of the Levi-Civita
connection\index{connection!Levi-Civita} of ${\mathbf E}^{n+1}$, while
the second indicates that it
has vanishing Riemann curvature tensor\index{Riemann curvature tensor}.

In terms of these forms, the fibers of $x:{\mathcal F}\to{\mathbf
  E}^{n+1}$ are exactly the maximal connected integral manifolds of
the Pfaffian system $\{\omega^a\}$.  Note that $\{\omega^a\}$ and
$\{dx^a\}$ are alternative bases for the space of forms on ${\mathcal
  F}$ that are semibasic over ${\mathbf E}^{n+1}$,
but the former is $E(n+1)$-invariant, while the latter is not.

We return to an explanation of our contact manifold $M$, by first
distinguishing the $1$-form on ${\mathcal F}$
$$
  \theta \stackrel{\mathit{def}}{=} \omega^0.
$$
Note that its defining formula
$$
  \theta_f(v) = \langle dx(v),e_0(f)\rangle, \qquad v\in T_f{\mathcal F},
$$
shows that it is the pullback of a unique, globally defined $1$-form
on $M$, which we will also call $\theta\in\Omega^1(M)$.
To see that $\theta$ is a contact form, first relabel the forms on
${\mathcal F}$ (this will be useful later, as well)
$$
  \pi_i \stackrel{\mathit{def}}{=} \omega^0_i,
$$
and note the equation on ${\mathcal F}$
$$
  d\theta = -\sum \pi_i\wedge\omega^i.
$$
So on ${\mathcal F}$ we certainly have $\theta\wedge(d\theta)^n\neq 0$, and
because pullback of forms via the submersion ${\mathcal F}\to M$ is
injective, the same non-degeneracy holds on $M$.

To understand the Legendre submanifolds of $M$, consider an oriented immersion
$$
  N^n\stackrel{\iota}{\hookrightarrow} M^{2n+1},\quad y=(y^1,\ldots,y^n)\mapsto
     (x(y),e_0(y)).
$$
The Legendre condition is
$$
  (\iota^*\theta)_y(v) = 
  \langle dx_y(v),e_0(y)\rangle = 0,\qquad
    v\in T_yN.
$$
In the transverse\index{Legendre submanifold!non-transverse|(} case, when the composition $x\circ\iota:N^n\hookrightarrow
M^{2n+1}\to{\mathbf E}^{n+1}$ is a hypersurface immersion (equivalently,
$\iota^*(\bigwedge\omega^i)\neq 0$, suitably interpreted), this condition is
that $e_0(y)$ is a unit normal vector to the hypersurface $x\circ\iota(N)$.
These Legendre submanifolds may therefore be thought of as the graphs
of Gauss maps\index{Gauss map} of oriented hypersurfaces
$N^n\hookrightarrow{\mathbf E}^{n+1}$. 
Non-transverse Legendre submanifolds of $M$ are sometimes of interest.
To give some intuition for these, we exhibit two examples in the contact
manifold over ${\mathbf E}^3$.  First, over an immersed curve
$x\!:\!I\hookrightarrow\mathbf E^3$, one can define a cylinder
$N=S^1\times I\hookrightarrow M\cong{\mathbf E}^3\times S^2$ by
$$
  (v,w)\mapsto(x(w),R_v(\nu_x)),
$$
where $\nu$ is any normal vector field along the curve $x(w)$, and
$R_v$ is rotation through angle $v\in S^1$ about the tangent $x^\prime(w)$.
The image is just the unit normal bundle of the curve, and
it is easily verified that this is a Legendre submanifold.

Our second example corresponds to the {\em
  pseudosphere}\index{pseudosphere}, a singular
surface of revolution in ${\mathbf E}^3$ having constant Gauss curvature
$K=-1$ away from the singular locus.  The map $x:S^1\times\R\to{\mathbf E}^3$
given by
$$
  x:(v,w)\mapsto(\mbox{sech}\,w\cos v,-\mbox{sech}\,w\sin v,
    w-\mbox{tanh}\,w)
$$
fails to be an immersion where $w=0$.  However, the Gauss map of the
complement of this singular locus can be extended to a smooth map
$e_3:S^1\times\R\to S^2$ given by
$$
  e_3(v,w) =
  (-\mbox{tanh}\,w\cos v, \mbox{tanh}\,w\sin v, -\mbox{sech}\,w).
$$
The graph of the Gauss map is the product $(x,e_3):S^1\times
\R\hookrightarrow M$.  It is a
Legendre submanifold, giving a smooth surface in $M$ whose projection to
${\mathbf E}^3$ is one-to-one, is an immersion almost everywhere, and
has image equal to the singular pseudosphere.  We will discuss in
\S\ref{Subsection:Gauss}
the exterior differential system whose integral manifolds are graphs
of Gauss maps of $K=-1$ surfaces in ${\mathbf E}^3$.  In
\S\ref{Subsection:Backlund}, we will discuss the B\"acklund
transformation\index{B\"acklund transformation} for this system, which
relates this particular example
to a special case of the preceding example, the unit normal bundle of
a line.\index{Legendre submanifold!non-transverse|)}

\subsection{Euclidean-invariant Euler-Lagrange Systems}
\label{Subsection:Euclidean-invariant}
\index{area functional|(}
\index{minimal surface|(}

We can now introduce one of the most important of all variational
problems, that of finding minimal-area hypersurfaces in Euclidean
space.  Define the $n$-form
$$
  \Lambda = \omega^1\wedge\cdots\wedge\omega^n \in \Omega^n({\mathcal F}),
$$
and observe that it is basic over $M$; that is, it is the
pullback of
a well-defined $n$-form on $M$ (although its factors $\omega^i$
are not basic).  This defines a Lagrangian functional
$$
  {\mathcal F}_\Lambda(N) = \int_N\Lambda
$$
on compact Legendre submanifolds $N^n\hookrightarrow M^{2n+1}$, which in the
transverse case discussed earlier equals the area of $N$ induced by
the immersion $N\hookrightarrow {\mathbf E}^{n+1}$.  We calculate the
Poincar\'e-Cartan form up on ${\mathcal F}$ using the structure equations
(\ref{EucStructure}), as
$$
  d\Lambda = -\theta\wedge\sum\pi_i\wedge\omega_{(i)},
$$
so the Euler-Lagrange system ${\mathcal E}_\Lambda$ is generated
by $\I=\{\theta,d\theta\}$ and
$$
  \Psi = -\sum\pi_i\wedge\omega_{(i)},
$$
which is again well-defined on $M$.
A transverse Legendre submanifold $N\hookrightarrow M$
will locally have a basis of $1$-forms given by
pullbacks (by any section) of $\omega^1,\ldots,\omega^n$, so applying
the Cartan lemma\index{Cartan lemma} to
$$
  0 = d\theta|_N = -\pi_i\wedge\omega^i
$$
shows that restricted to $N$ there are expressions
$$
  \pi_i=\sum_j h_{ij}\omega^j
$$
for some functions $h_{ij}=h_{ji}$.  If $N\hookrightarrow M$ is also
an integral manifold of $\mathcal E_\Lambda\subset\Omega^*(M)$, then
additionally 
$$
  0 = \Psi|_N = -\left(\sum h_{ii}\right)\omega^1\wedge\cdots\wedge\omega^n.
$$
One can identify $h_{ij}$ with the second fundamental
form\index{second fundamental form} of
$N\hookrightarrow{\mathbf E}^{n+1}$ in this transverse case, and we then
have the usual criterion
that a hypersurface is stationary for the area functional if and only if 
its mean curvature \index{mean curvature} $\sum h_{ii}$ vanishes.  We
will return to the
study of this Euler-Lagrange system shortly.

Another natural $E(n+1)$-invariant PDE for hypersurfaces in Euclidean
space is that of prescribed constant mean 
curvature\index{mean curvature!constant} $H$, not
necessarily zero.  We first ask whether such an equation is even
Euler-Lagrange, and to answer this we apply our inverse
problem\index{inverse problem} test
to the Monge-Ampere system
$$
  {\mathcal E}_H = \{\theta,d\theta,\Psi_H\},\qquad
  \Psi_H = -\left(\sum\pi_i\wedge\omega_{(i)} - H\omega\right).
$$
Here, $H$ is the prescribed constant and $\omega =
\omega^1\wedge\cdots\wedge\omega^n$ is the induced volume form.  The
transverse integral manifolds of ${\mathcal E}_H$ correspond to the
desired Euclidean hypersurfaces.

To implement the test, we take the candidate Poincar\'e-Cartan form
$$
  \Pi_H = -\theta\wedge\left(\sum\pi_i\wedge\omega_{(i)} - H\omega\right)
$$
and differentiate; the derivative of the first term vanishes, as we
know from the preceding case of $H=0$, and we have
\begin{eqnarray*}
  d\Pi_H & = & H\,d(\theta\wedge\omega^1\wedge\cdots\wedge\omega^n) \\
         & = & H\,d(dx^0\wedge\cdots\wedge dx^n) \\
         & = & 0.
\end{eqnarray*}
So this ${\mathcal E}_H$ is at least locally the Euler-Lagrange system
for some functional $\Lambda_H$, which can be taken to be an
anti-derivative of $\Pi_H$.
One difficulty in finding $\Pi_H$ is that there is {\em no} such
$\Lambda_H$ that is invariant under the Euclidean group $E(n+1)$.
The next best thing would be to find a $\Lambda_H$ which is
invariant under the rotation subgroup $SO(n+1,\R)$, but not under
translations.  A little experimentation yields the Lagrangian
$$
  \Lambda_H = \omega + \sf{H}{n+1}x\innerprod\Omega,\qquad
     d\Lambda_H = \Pi_H,
$$
where $x=\sum x^a\frac{\partial}{\partial x^a}$ is the radial position
vector field, $\omega = \omega^1\wedge\cdots\wedge\omega^n$ is the
hypersurface area form, and $\Omega =
\omega^1\wedge\cdots\wedge\omega^{n+1}$ is the ambient volume form.
The choice of an origin from which to define the position vector $x$
reduces the symmetry 
group of $\Lambda_H$ from $E(n+1)$ to $SO(n+1,\R)$.  The functional
$\int_N\Lambda_H$ gives the area of the hypersurface $N$ plus a
scalar multiple of the signed volume of the cone on $N$ with vertex at
the origin.

It is actually possible to list all of the Euclidean-invariant
Poincar\'e-Cartan 
forms\index{Poincar\'e-Cartan form!Euclidean-invariant|(} 
on $M\to{\mathbf E}^{n+1}$.  Let
$$
  \Lambda_{-1} = -\sf{1}{n+1}x\innerprod\Omega,\quad
  \Lambda_k = \sum_{|I|=k}\pi_I\wedge\omega_{(I)}\ \ (0\leq k\leq n),
$$
and
$$
  \Pi_k = -\theta\wedge\Lambda_k,
$$
It is an exercise using the structure equations to show that
$$
  d\Lambda_k = \Pi_{k+1}.
$$
Although these forms are initially defined up on ${\mathcal F}$, it is easily
verified that they are pull-backs of forms on $M$, which we denote by
the same name.
It can be proved using the first fundamental theorem of orthogonal
invariants that any Euclidean-invariant Poincar\'e-Cartan form is a
linear combination of $\Pi_0,\ldots,\Pi_n$.  Note
that such a Poincar\'e-Cartan form is induced by a Euclidean-invariant
functional if and only if $\Pi_0$ is not involved.

We can geometrically interpret $\Lambda_k|_N$ for transverse Legendre
submanifolds $N$ as the sum 
of the $k\times k$ minor determinants of the second fundamental
form\index{second fundamental form}
$I\!I_N$, times the hypersurface area form of $N$.  In case $k=n$ we
have $d\Lambda_n = \Pi_{n+1} = 0$, reflecting the fact that the
functional
$$
  \int_N\Lambda_n = \int_NK\,dA
$$
is variationally trivial, where $K$ is the Gauss-Kronecker curvature.

\subsubsection{Contact Equivalence of Linear Weingarten Equations for
  Surfaces}
The Euclidean-invariant Poincar\'e-Cartan 
forms\index{Poincar\'e-Cartan form!Euclidean-invariant|)} for surfaces in
${\mathbf E}^3$ give rise to the {\em linear Weingarten
  equations},\index{Weingarten equation|(} of the
form
$$
  aK+bH+c=0
$$
for constants $a,b,c$.  Although these second-order PDEs are
inequivalent under 
point-transformations\index{transformation!point} 
for non-proportional choices of
$a,b,c$, we will show that under {\em contact} 
transformations\index{transformation!contact|(} there
are only five distinct equivalence classes of linear Weingarten equations.

To study surfaces, we work on the unit sphere bundle
$\pi:M^5\to{\mathbf E}^3$, and recall the formula for the contact form
$$
  \theta_{(x,e_0)}(v) = \langle \pi_*(v),e_0\rangle,
  \qquad v\in T_{(x,e_0)}M.
$$
We define two $1$-parameter groups of diffeomorphisms of $M$ as
follows:
\begin{eqnarray*}
   \varphi_t(x,e_0) & = & (x+te_0,e_0), \\
   \psi_s(x,e_0) & = & (\exp(s)x,e_0).
\end{eqnarray*}
It is not hard to see geometrically that these define contact
transformations on $M$, although this result will also come out of the
following calculations.  We will carry out
calculations on the full Euclidean frame 
bundle\index{Euclidean!frame bundle} ${\mathcal
  F}\to{\mathbf E}^3$, where
there is a basis of $1$-forms $\omega^1,\ \omega^2,\ \theta,\ \pi_1,\
\pi_2,\ \omega^1_2$ satisfying structure equations presented earlier.

To study $\varphi_t$ we use its generating vector field
$v = \frac{\partial}{\partial\theta}$, which is the 
dual of the $1$-form $\theta$ with respect to the preceding basis.  We 
can easily compute Lie derivatives
$$
  {\mathcal L}_v\omega^1 = -\pi_1,\quad
  {\mathcal L}_v\omega^2 = -\pi_2,\quad
  {\mathcal L}_v\theta = 0,\quad
  {\mathcal L}_v\pi_1 = 0,\quad
  {\mathcal L}_v\pi_2 = 0.
$$
Now, the fibers of ${\mathcal F}\to M$ have tangent spaces given by
$\{\omega^1,\omega^2,\theta,\pi_1,\pi_2\}^\perp$, and this
distribution is evidently preserved by the flow along $v$.  This
implies that $v$ induces a vector field downstairs on $M$, whose flow is
easily seen to be $\varphi_t$.  The fact that ${\mathcal L}_v\theta=0$ 
confirms that $\varphi_t$ is a contact transformation.

We can now examine the effect of $\varphi_t$ on the invariant
Euler-Lagrange systems corresponding to linear Weingarten equations
by introducing
$$
  \Psi_2 = \pi_1\wedge\pi_2,\quad
  \Psi_1 = \pi_1\wedge\omega^2-\pi_2\wedge\omega^1,\quad
  \Psi_0 = \omega^1\wedge\omega^2.
$$
Restricted to a transverse Legendre submanifold over a surface
$N\subset{\mathbf E}^3$, these give $K\,dA,\ H\,dA,$ and the area form $dA$ of
$N$, respectively.  Linear Weingarten surfaces are integral manifolds
of a Monge-Ampere system
$$
  \{\theta, d\theta, \Psi(a,b,c) \stackrel{\mathit{def}}{=}
  a\Psi_2+b\Psi_1+c\Psi_0\}.
$$
Our previous Lie derivative computations may be used to compute
$$
  {\mathcal L}_v\left(\begin{array}{c} \Psi_0 \\ \Psi_1 \\
       \Psi_2 \end{array}\right) =
  \left(\begin{array}{ccc} 0 & -1 & 0 \\ 0 & 0 & -2 \\ 0 & 0 & 0
    \end{array}\right)
  \left(\begin{array}{c} \Psi_0 \\ \Psi_1 \\ \Psi_2
       \end{array}\right).
$$
Exponentiate this to see
\begin{equation}
  \varphi_t^*\Psi(a,b,c) = \Psi(a-2bt+ct^2, b-ct, c).
\label{FirstContactTrans}
\end{equation}
This describes how the $1$-parameter group $\varphi_t$ acts on the
collection of linear Weingarten equations.
Similar calculations show that the $1$-parameter group $\psi_s$
introduced earlier consists of contact transformations, and acts on
linear Weingarten equations as
\begin{equation}
  \psi_s^*\Psi(a,b,c) = \Psi(a,\exp(s)b,\exp(2s)c).
\label{SecondContactTrans}
\end{equation}

It is reasonable to regard the coefficients $(a,b,c)$ which specify a
particular linear Weingarten equation as a point $[a:b:c]$ in the real 
projective plane $\mathbf{RP}^2$, and it is an easy exercise to
determine the orbits in $\mathbf{RP}^2$ of the group action generated
by (\ref{FirstContactTrans}) and (\ref{SecondContactTrans}). 
There are five orbits, represented by the points $[1:0:0]$, $[0:1:0]$, 
$[1:0:1]$, $[1:0:-1]$, $[0:0:1]$.  
The special case
$$
  \varphi_{\frac1A}^*\Psi(0,1,A) = \Psi(-\sf1A,0,A)
$$
gives the classically known fact that to every surface of non-zero
constant mean curvature\index{mean curvature!constant} $-A$, there is a (possibly singular) parallel
surface of constant positive Gauss curvature $A^2$.  Note finally that
the Monge-Ampere system corresponding to $[0:0:1]$ has for integral
manifolds those non-transverse Legendre 
submanifolds\index{Legendre submanifold!non-transverse} of $M$ which project
to curves in ${\mathbf E}^3$, instead of surfaces.
\index{transformation!contact|)}
\index{Weingarten equation|)}

\subsection{Conservation Laws for Minimal Hypersurfaces}
\label{Subsection:MinimalCLs}
\index{conservation law!for minimal hypersurfaces|(}

In Chapter 3, we will be concerned with conservation laws for various
Euler-Lagrange equations arising in conformal geometry.  We will
emphasize two questions:  how are conservation laws
found, and how can they be used?  In this section, we will explore
these two questions in the case of the minimal
hypersurface equation $H=0$, regarding conservation laws arising from
Euclidean symmetries.

We compute these conservation laws first for the
translations\index{translation}, and
then for the rotations.  The results of these computations will be
the two vector-valued conservation laws
$$
\boxed{
  d(*dx) = 0,\qquad d(*(x\wedge dx)) = 0.
}
$$
The notation will be explained in the course of the calculation.
These may be thought of as analogs of the conservation of linear
and angular momentum that are ubiquitous in physics.

To carry out the computation, note that the prescription for Noether's 
theorem\index{Noether's theorem} given in
(\ref{NPrescription1}, \ref{NPrescription2}, \ref{NPrescription3}) is
particularly simple for the case of the functional
$$
  \Lambda = \omega^1\wedge\cdots\wedge\omega^n
$$
on the contact manifold $M^{2n+1}$.  This is because first, $d\Lambda
= \Pi$ already, so no correction term is needed, and second, the
infinitesimal Euclidean symmetries (prolonged to act on $M$) actually
preserve $\Lambda$, and not merely the equivalence class $[\Lambda]$.
Consequently, the Noether prescription is (with a sign change)
$$
  \eta(v) = v\innerprod\Lambda.
$$
This $v\innerprod\Lambda$ is an $(n-1)$-form on $M$ which is closed
modulo the Monge-Ampere system ${\mathcal E}_\Lambda$.  

Proceeding, we can suppose that our translation vector field is
written up on the Euclidean frame bundle\index{Euclidean!frame bundle} as
$$
  v_{\mathcal F} = Ae_0 + A^ie_i,
$$
where the coefficients are such that the equation $dv=0$ holds; that
is, the functions $A$ and $A^i$ are the coefficients of a fixed vector
with respect to a varying oriented orthonormal frame.  We easily find
$$
  \varphi_v = v_{\mathcal F}\innerprod\Lambda = \sum_{i=1}^n A^i\omega_{(i)}.
$$
This, then, is the formula for an $(n-1)$-form on ${\mathcal F}$ which is
well-defined on the contact manifold $M$ and is
closed when restricted to integral manifolds of the Monge-Ampere
system ${\mathcal E}_\Lambda$.  To see it in another form, observe
that if we restrict our $(n-1)$-form to a transverse Legendre submanifold $N$,
$$
  \varphi_v|_N = \sum A^i(*\omega^i) = *\langle v,dx\rangle.
$$
Here and throughout, the star operator $*=*_N$ is defined with respect to
the induced metric and orientation on $N$, and
the last equality follows from the equation of ${\mathbf E}^{n+1}$-valued
$1$-forms $dx = e_0\theta + \sum e_i\omega^i$, where $\theta|_N=0$.
We now have a linear map from $\R^{n+1}$, regarded as the space of
translation vectors $v$, to the space of closed $(n-1)$-forms on any
minimal hypersurface $N$.  Tautologically, such a map may be regarded
as one closed $(\R^{n+1})^*$-valued $(n-1)$-form on $N$.  Using the
metric to identify $(\R^{n+1})^*\cong\R^{n+1}$, this may be
written as
$$
  \varphi_{trans} = *dx.
$$
This is the meaning of the conservation law stated at the beginning of 
this section.
Note that each component $d(*dx^a)=0$ of this conservation law
is equivalent to the claim that the coordinate function $x^a$ of the
immersion $x:N\hookrightarrow{\mathbf E}^{n+1}$ is a
harmonic\index{harmonic!function} function 
with respect to the induced metric on $N$.

Turning to the rotation\index{rotation} vector fields, we first write such a
vector field on ${\mathbf E}^{n+1}$ as
$$
  v = \sum_{a,b=1}^{n+1}x^aR^b_a\frac{\partial}{\partial x^b},\qquad
    R^a_b+R^b_a = 0.
$$
It is not hard to verify that this vector field lifts naturally to the 
frame bundle ${\mathcal F}$ as
$$
  v_{{\mathcal F}} = \sum x^aR^b_aA^b_c\frac{\partial}{\partial\omega^c}
    +\sum A^b_cR^b_aA^a_d\frac{\partial}{\partial\omega^c_d},
$$
where the coefficients $A^a_b$ are defined by the equation
$\frac{\partial}{\partial\omega^a} = \sum A^b_a\frac{\partial}{\partial
x^b}$, and the tangent
vectors $\frac{\partial}{\partial\omega^a}$,
$\frac{\partial}{\partial\omega^a_b}$ are dual to the canonical
coframing $\omega^a$, $\omega^a_b$ of ${\mathcal F}$.

We can now compute (restricted to $N$, for convenience)
\begin{eqnarray*}
  (v_{{\mathcal F}}\innerprod\Lambda)|_N & = &
    \sum x^aR^b_aA^b_i\omega_{(i)} \\
  & = & *(x^aR^b_aA^b_c\omega^c) \\
  & = & *\langle R\cdot x, dx\rangle.
\end{eqnarray*}
Reformulating the Noether map in a manner analogous to that used
previously, we can define a $\bigwedge^2\R^{n+1}\cong\lie{so}(n+1,\R)^*$-valued
$(n-1)$-form on $N$
$$
  \varphi_{rot} = *(x\wedge dx).
$$
Once again, $\varphi_{rot}$ is a conservation law by virtue of the
fact that it is closed if $N$ is a minimal hypersurface.

It is interesting to note that the conservation law for rotation
symmetry is a consequence of that for translation symmetry.  This is
because we have from $d(*dx)=0$ that
$$
  d(x\wedge*dx) = dx\wedge*dx = 0.
$$
The last equation holds because the exterior multiplication $\wedge$
refers to the ${\mathbf E}^{n+1}$ where the forms take values, {\em not} the
exterior algebra in which their components live.
It is an exercise to show that these translation conservation laws are
equivalent to minimality of $N$.

Another worthwhile exercise is to show that all of the classical
conservation laws for the $H=0$ system arise from infinitesimal
Euclidean symmetries.  In the next chapter, we will see directly that the
group of symmetries of the Poincar\'e-Cartan form for this system
equals the group of Euclidean motions, giving a more illuminating proof of
this fact.  At the end of this section, we will consider a dilation
vector field which preserves the minimal surface system $\mathcal E$,
but not the Poincar\'e-Cartan form, and use it to compute an
``almost-conservation law''.  

By contrast, in this case there is {\em
  no} discrepancy between $\lie{g}_{[\Lambda]}$ and $\lie{g}_{\Pi}$.
To see this, first note that by Noether's theorem\index{Noether's theorem}
\ref{NoetherTheorem}, $\lie{g}_{\Pi}$ 
is identified with $H^n({\mathcal E}_\Lambda)$, and
$\lie{g}_{[\Lambda]}\subseteq\lie{g}_{\Pi}$ is identified with the image
of the connecting map $\delta$ in the long exact sequence
$$
 \cdots\to H^{n-1}(\Omega^*/{\mathcal E}_\Lambda)\stackrel{\delta}{\to}
    H^n({\mathcal E}_\Lambda)\stackrel{\iota}{\to}H^n_{dR}(M)\to\cdots.
$$
With $M\cong {\mathbf E}^{n+1}\times S^n$, we have the isomorphism
$H^n_{dR}(M)\cong\R$ obtained by integrating an $n$-form along a fiber
of $M\to{\mathbf E}^{n+1}$, and it is not hard to see that any
$n$-form in ${\mathcal E}_\Lambda$ must vanish when restricted to such
a fiber.  Therefore the map $\iota$ is identically $0$, so $\delta$ is
onto, and that proves our claim.

\subsubsection{Interpreting the Conservation Laws for $H=0$}

To understand the meaning of the conservation law $\varphi_{trans}$, we
convert the equation $d\varphi_{trans}|_N = 0$ to integral form.  For a
smoothly bounded, oriented neighborhood $U\subset N\subset{\mathbf E}^{n+1}$ with $N$
minimal, we have by Stokes' theorem
$$
  \int_{\partial U}*dx = 0.
$$
To interpret this condition on $U$, we take an oriented orthonormal
frame field $(e_0,\ldots,e_n)$ along $U\cup\partial U$, such that
along the boundary $\partial U$ the following hold:
\begin{equation}
  \left\{\begin{array}{l}
    e_0 \mbox{ is the oriented normal to } N,\\
    e_n \mbox{ is the outward normal to }\partial U\mbox{ in } N, \\
    e_1,\ldots,e_{n-1}\mbox{ are tangent to }\partial U.
\end{array}\right.
\label{AdaptFrame}
\end{equation}
Calculations will be much easier in this adapted frame field.  
The dual coframe $\omega^a$ for ${\mathbf E}^{n+1}$ along $U\cup\partial U$ satisfies
$$
  dx = e_0\omega^0 + \sum_{i=1}^{n-1}e_i\omega^i + e_n\omega^n.
$$
Now, the first term vanishes when restricted to $N$.  The
last term vanishes when restricted to $\partial U$, but
cannot be discarded because it will affect $*_Ndx$, which we are trying
to compute.  Consequently,
$$
  *_Ndx = \sum_{i=1}^{n-1}e_i\omega_{(i)}+e_n\omega_{(n)}.
$$
{\em Now} we restrict to $\partial U$, and find
\begin{eqnarray*}
  *_Ndx|_{\partial U} & = &
      (-1)^{n-1}e_n\omega^1\wedge\cdots\wedge\omega^{n-1}\\
   & = & (-1)^{n-1}\mathbf{n}\ d\sigma.
\end{eqnarray*}
Here we use $\mathbf{n}$ to denote the normal to $\partial U$ in $N$
and $d\sigma$ to denote the area measure induced on $\partial U$.
Our conservation law therefore reads
$$
  \int_{\partial U}\mathbf{n}\ d\sigma = 0.
$$
In other words,
in a minimal hypersurface the average of the exterior
unit normal vectors over the smooth boundary of any oriented
neighborhood must vanish.  One consequence of this is that a minimal
surface can never be locally convex; that is, a neighborhood of a
point can never lie on one side of the tangent plane at that point.
This is intuitively reasonable from the notion of minimality.
Similar calculations give an analogous formulation for the rotation
conservation law:
$$
  \int_{\partial U}(\mathbf{x}\wedge\mathbf{n})d\sigma = 0.
$$

These interpretations have relevance to the classical Plateau
problem\index{Plateau problem}, 
which asks whether a given simple closed curve $\gamma$ in ${\mathbf
  E}^3$ bounds a
minimal surface.  The answer to this is affirmative, with the caveat
that such a surface is not necessarily unique and may not be smooth at
the boundary.  A more well-posed version 
gives not only a simple closed curve $\gamma\subset{\mathbf E}^3$, but a {\em
strip}, which is a curve $\gamma^{(1)}\subset M$ consisting of a base curve
$\gamma\subset{\mathbf E}^3$ along with a field of tangent planes along
$\gamma$ containing the tangent lines to $\gamma$.  Such a strip is
the same as a curve in $M$ along which the contact $1$-form vanishes.
Asking for a minimal surface whose boundary and boundary-tangent
planes are described by a given $\gamma^{(1)}$ is the same as asking
for a transverse integral manifold of ${\mathcal E}_\Lambda$ having boundary
$\gamma^{(1)}\subset M$.

The use of our two conservation laws in this context comes from the fact
that $\gamma^{(1)}$ determines the vector-valued form $*_Ndx$ along
$\partial N$ for {\em any} possible solution to this initial value
problem.  The conservation laws give integral constraints, often called
{\em moment conditions}\index{moment conditions|(}, on the values of $*_Ndx$, and hence constrain
the possible strips $\gamma^{(1)}$ for which our problem has an
affirmative answer.  However, the moment conditions on a strip
$\gamma^{(1)}$ are not sufficient for there to exist a minimal surface
with that boundary data.  We will discuss additional constraints
which have the feel of ``hidden conservation 
laws''\index{conservation law!hidden@``hidden''} after a
digression on similar moment conditions that arise for boundaries of
holomorphic curves.

\

It is natural to ask whether a given real, simple, closed curve
$\gamma_\C$ in complex space $\C^n$ (always $n\geq 2$) is the boundary
of some holomorphic disc\index{holomorphic curve|(}.  There is a
differential ideal ${\mathcal
J}\subset\Omega^*_\R(\C^n)$ whose integral manifolds are precisely
holomorphic curves, defined by
$$
  {\mathcal J} = \{(\Omega^{2,0}(\C^n)+\Omega^{0,2}(\C^n))\cap
       \Omega^2_\R(\C^n)\}.
$$
In other words, ${\mathcal J}$ is algebraically generated by real $2$-forms
which, when regarded as complex $2$-forms, have no part of type
$(1,1)$.  It is elementary to see that in degree $k\geq 3$, ${\mathcal
J}^k=\Omega^k_\R(\C^n)$, and that the integral $2$-planes in $T\C^n$
are exactly the complex $1$-dimensional subspaces.  This implies our
claim that integral manifolds of ${\mathcal J}$ are holomorphic
curves.

Now, ${\mathcal J}$ has many conservation laws.  Namely, for any
holomorphic $1$-form $\varphi\in\Omega^{1,0}_{hol}(\C^n)$, we find
that
$$
  d\varphi+d\bar\varphi \in{\mathcal J},
$$
so that $\varphi+\bar\varphi$ is a conservation law for ${\mathcal
J}$.  These give rise to infinitely many moment conditions
$$
  \int_{\gamma_\C}\varphi=0
$$
which must be satisfied by $\gamma_\C$, if it is to be the boundary of
a holomorphic disc.  

It is a fact which we shall not prove here that {\em every}
conservation law for ${\mathcal J}$ is of this form; trivial
conservation laws clearly arise when $\varphi = df$ for some
holomorphic function $f\in{\mathcal O}(\C^n)$.
Another fact, not to
be proved here, is that these moment conditions are
sufficient for $\gamma_\C$ to bound a (possibly branched) holomorphic
disc.

\

Returning to our discussion of minimal surfaces, suppose that
$x:U\to{\mathbf E}^3$ is a minimal immersion of a simply connected surface.
Then $*dx$ defines a closed, vector-valued $1$-form on $U$, so there
exists a vector-valued function $y\!:\!U\to{\mathbf E}^3$ satisfying
\begin{equation}
  dy = *dx.
\label{CREquation}
\end{equation}
Note that our ability to integrate the conservation law to obtain a
function relies essentially on the fact that we are in dimension
$n=2$. 

We can define
$$
  z=(x+iy)\!:\!U\to\C^3,
$$
and (\ref{CREquation}) is essentially the Cauchy-Riemann equations,
implying that $z$ is a holomorphic curve, with the conformal structure
induced from the immersion $z$.  Furthermore, the complex derivative
$z^\prime$ is at each point of $U$ a null vector for the complex
bilinear inner-product $\sum (dz^i)^2$.  This gives the classical
Weierstrass representation\index{Weierstrass representation} of a
minimal surface in ${\mathbf E}^3$ as locally
the real part of a holomorphic null curve in $\C^3$.

We can now incorporate the result of our digression on conservation
laws for holomorphic discs.  Namely, given a strip $\gamma^{(1)}$, the
Euclidean moment condition $\int_\gamma*dx=0$ implies that there
exists another real curve $y$ so that $dy=*dx$ (along $\gamma$).  Then we
can use $z=(x+iy)\!:\!\gamma\to\C^3$ as initial data for the holomorphic
disc problem, and all of the holomorphic moment conditions for that
problem come into play.  These are the additional hidden constraints
needed to fill the real curve $\gamma$ with a (possibly branched)
minimal surface.
\index{conservation law!for minimal hypersurfaces|)}
\index{holomorphic curve|)}

\subsubsection{Conservation Laws for Constant Mean Curvature}
\index{conservation law!for cmc hypersurfaces|(}

It is also a worthwhile exercise to determine the conservation laws
corresponding to Euclidean motions for the constant mean curvature
system when the constant $H$ is non-zero.  Recall that for that system the
Poincar\'e-Cartan form
$$
  \Pi_H = -\theta\wedge\left(\sum \pi_i\wedge\omega_{(i)} - H\omega\right)
$$
is invariant under the full Euclidean group, but that no particular
Lagrangian $\Lambda$ is so invariant; we will continue to work with
the $SO(n+1,\R)$-invariant Lagrangian
$$
  \Lambda_H = \omega + \sf{H}{n+1}x\innerprod\Omega,\qquad
     d\Lambda_H = \Pi_H.
$$
Fortunately, the equivalence class
$[\Lambda]\in H^n(\Omega^*(M)/\I)$ {\em is} invariant under the Euclidean
group, because as the reader can verify, the connecting map
$$
   \delta:H^n(\Omega^*/\I)\to H^{n+1}(\I)
$$
taking $[\Lambda]$ to $\Pi$ is an isomorphism for this contact
manifold.  This means that, as in the case $H=0$, we will find
conservation laws corresponding to the full Euclidean Lie algebra.

Computing the conservation laws corresponding to
translations\index{translation} requires
the more complicated form of the Noether prescription, because it is
the translation vector fields $v\in\R^{n+1}$ which fail to
preserve our $\Lambda_H$.  Instead, we have
\begin{eqnarray*}
  {\mathcal L}_v\Lambda_H & = &
    \sf{H}{n+1}{\mathcal L}_v(x\innerprod\Omega)
    \mbox{ (because ${\mathcal L}_v\omega=0$),} \\
  & = & \sf{H}{n+1}(({\mathcal L}_vx)\innerprod\Omega +
        x\innerprod({\mathcal L}_v\Omega)) \\
  & = & \sf{H}{n+1}(v\innerprod\Omega+0).
\end{eqnarray*}
In the last step, we have used ${\mathcal L}_vx=[v,x]=v$ (by a
simple calculation), and
${\mathcal L}_v\Omega=0$ (because the ambient volume $\Omega$ is
translation invariant).
To apply the Noether prescription, we need an anti-derivative of this
last term, which we find by experimenting:
\begin{eqnarray*}
  d(x\innerprod (v\innerprod\Omega)) & = &
    {\mathcal L}_x(v\innerprod\Omega)- 
        x\innerprod d(v\innerprod\Omega)\\
  & = & (({\mathcal L}_xv)\innerprod\Omega + 
     v\innerprod({\mathcal L}_x\Omega))-x\innerprod 0 \\
  & = & -v\innerprod\Omega+(n+1)v\innerprod\Omega,
\end{eqnarray*}
where we have again used ${\mathcal L}_xv=[x,v]=-v$,
and ${\mathcal L}_x\Omega=(n+1)\Omega$.
Combining these two calculations, we have
$$
  {\mathcal L}_v\Lambda_H = \sf{H}{n(n+1)}
    d(x\innerprod(v\innerprod\Omega)).
$$
The prescription
(\ref{NPrescription1}, \ref{NPrescription2}, \ref{NPrescription3}) now
gives
\begin{eqnarray*}
  \varphi_v & = & -v\innerprod\Lambda_H+
    \sf{H}{n(n+1)}x\innerprod(v\innerprod\Omega) \\
  & = & -v\innerprod\omega + \sf{H}{n}x\innerprod
       (v\innerprod\Omega).
\end{eqnarray*}
As in the case of minimal hypersurfaces, we consider the restriction
of $\varphi_v$ to an integral manifold $N$.
From the previous case, we know that $v\innerprod\omega$ restricts to
$*\langle v,dx\rangle$, where $*=*_N$ is the star operator of the
metric on $N$ and $\langle\cdot,\cdot\rangle$ denotes the ambient
inner-product.
To express the restriction of the other term of $\varphi_v$, decompose
$x=x_t+x_\nu\nu$ into tangential and normal parts along $N$ (so $x_t$
is a vector and $x_{\nu}$ is a scalar), and a
calculation gives
\begin{eqnarray*}
  \sf{H}{n}x\innerprod(v\innerprod\Omega)|_N & = &
    -\sf{H}{n}(x_\nu v\innerprod\omega -(v\innerprod\theta)
              (x_t\innerprod\omega)) \\
  & = & -\sf{H}{n}(x_\nu*\!\langle v,dx\rangle-
    (v\innerprod\theta)*\!\langle\cdot,x_t\rangle);
\end{eqnarray*}
the latter $*$ is being applied to the $1$-form on $N$ that is dual
via the metric to the
tangent vector $x_t$.  Again as in the $H=0$ case, we
can write these $(n-1)$-forms $\varphi_v$, which depend linearly on
$v\in\R^{n+1}$, as an $(\R^{n+1})^*$-valued $(n-1)$-form on $N$.  It is 
$$
  \varphi_{trans}=-(1+\sf{H}{n}x_\nu)*\!dx+
    \sf{H}{n}\nu*\!\langle\cdot,x_t\rangle.
$$
In the second term, the normal $\nu$ provides the ``vector-valued''
part (it replaced $\theta$, to which it is dual), and
$*\langle\cdot,x_t\rangle$ provides the ``$(n-1)$-form'' part.

Calculating the conservation laws for rotations\index{rotation} is a
similar process,
simplified somewhat by the fact that ${\mathcal L}_v\Lambda_H=0$; of
course, the lifted rotation vector fields $v$ are not so easy to work
with as the translations.  The resulting $\Lambda^2\R^{n+1}$-valued
$(n-1)$-form is
$$
  \varphi_{rot} = -(1+\sf{H}{n+1}x_\nu)*\!(x\wedge dx) +
    \sf{H}{n+1}(x\wedge \nu)*\!\langle\cdot,x_t\rangle.
$$
These can be used to produce moment conditions, just as in the $H=0$ case.
\index{moment conditions|)}
\index{Euclidean!motion|)}
\index{conservation law!for cmc hypersurfaces|)}

\

\index{conservation law!almost@``almost''|(}
We conclude with one more observation suggesting extensions of the
notion of a conservation law.  Recall that we showed in
(\ref{SymmBound}) that a Monge-Ampere system ${\mathcal E}_\Lambda$
might have an infinitesimal symmetry which scales the corresponding
Poincar\'e-Cartan form $\Pi$.  This is the case for the minimal
surface system, which is preserved by the dilation\index{dilation}
vector field on
${\mathbf E}^{n+1}$
$$
  x = \sum x^a\sf{\partial}{\partial x^a}.
$$
This induces a vector field $x$ on the contact manifold of
tangent hyperplanes to ${\mathbf E}^{n+1}$ where the functional $\Lambda$ and
Poincar\'e-Cartan form $\Pi$ are defined, and there are various ways to
calculate that
$$
  {\mathcal L}_x\Lambda = n\Lambda.
$$

If one tries to apply the Noether prescription to $x$ by writing
$$
  \varphi_{dil} = x\innerprod\Lambda,
$$
the resulting form satisfies
\begin{eqnarray*}
  d\varphi_{dil} & = &
        {\mathcal L}_x\Lambda - x\innerprod d\Lambda \\
    & = & n\Lambda - x\innerprod\Pi.
\end{eqnarray*}
Restricted to a minimal surface $N$, we will then have
$$
  d\varphi_{dil}|_N  =  n\Lambda.
$$
Because the right-hand side is not zero, we do not have a conservation
law, but it is still reasonable to look for consequences of
integrating on neighborhoods $U$ in $N$, where we find
\begin{equation}
  \int_{\partial U}\varphi = n\int_U\Lambda.
\label{AlmostConservation}
\end{equation}
The right-hand side equals $n$ times the area of $U$, and the
left-hand side can be investigated by choosing an oriented
orthonormal frame field $(e_0,\ldots,e_n)$ along $U\cup\partial U$
satisfying the conditions (\ref{AdaptFrame}) as before.  We write the
coefficients 
$$
   x = \sum x^a\sf{\partial}{\partial x^a} = \sum v^ae_a,
$$
and then restricted to $\partial U$, we have
$$
  \varphi|_{\partial U} = x\innerprod\Lambda = (-1)^{n-1}v^n\omega^1\wedge
     \cdots\wedge\omega^{n-1}.
$$
Up to sign, the form $\omega^1\wedge\cdots\wedge\omega^{n-1}$ along
$\partial U$ is exactly the $(n-1)$-dimensional area form for $\partial U$.

These interpretations of the two sides of (\ref{AlmostConservation})
can be exploited by taking for $U$ the family of neighborhoods $U_r$
for $r>0$, defined as the intersection of $N\subset{\mathbf E}^{n+1}$ with an
origin-centered ball of radius $r$.  In particular, along $\partial
U_r$ we will have $||x||=r$, so that $v^n\leq r$ and
\begin{equation}
  r\cdot\mbox{Area}(\partial U_r)\geq n\cdot\mbox{Vol}(U_r).
\label{EstimateArea}
\end{equation}
Observe that
$$
  \mbox{Area}(\partial U_r) = \frac{d}{dr}\mbox{Vol}(U_r),
$$
and (\ref{EstimateArea}) is now a differential inequality for
$\mbox{Vol}(U_r)$ which can be solved to give
$$
  \mbox{Vol}(U_r) \geq Cr^n
$$
for some constant $C$.  This is a remarkable result about minimal
hypersurfaces, and amply illustrates the power of
``almost-conservation laws'' like $\varphi_{dil}$.

\index{conservation law!almost@``almost''|)}
\index{Euclidean!space|)}
\index{area functional|)}
\index{minimal surface|)}

\chapter{The Geometry of Poincar\'e-Cartan Forms}

In this chapter, we will study some of the geometry associated to
Poincar\'e-Cartan forms using \'E.~Cartan's\index{Cartan, \'E.}
method of equivalence\index{equivalence method|(}.  The idea is
to identify such a
Poincar\'e-Cartan form with a 
$G$-structure\index{Gstructure@$G$-structure|(}---that is, a subbundle
of the 
principal coframe bundle of a manifold---and then attempt to find some 
canonically determined basis of
$1$-forms on the total space of that $G$-structure.  The differential 
structure equations of these $1$-forms will then exhibit 
associated geometric objects and invariants.

The pointwise linear algebra of a Poincar\'e-Cartan form in the case of $n=2$ 
``independent variables'' (that is, on a contact manifold of dimension $5$) 
is quite different from that of higher dimensional cases.  Therefore, in the 
first section we study only the former, which should serve as a good 
illustration of the method of equivalence for those not familiar with it.  
Actually, in case $n=2$ we will study the coarser equivalence of
Monge-Ampere systems rather than Poincar\'e-Cartan forms, and we will
do this without restricting to those systems which are locally
Euler-Lagrange.  An extensive study of the
geometry of Monge-Ampere systems in various low dimensions was carried out
in~\cite{Lychagin:Classification}, with a viewpoint somewhat similar
to ours.

In the succeeding sections, we will first identify in case $n\geq 3$ a 
narrower class of Poincar\'e-Cartan forms, called {\em
  neo-classical}\index{Poincar\'e-Cartan form!neo-classical}, which 
are of the same algebraic type as those arising from classical variational 
problems.  We will describe some of the geometry associated with 
neo-classical Poincar\'e-Cartan forms, consisting of a field of 
hypersurfaces in a vector bundle, well-defined up to fiberwise affine 
motions of the vector bundle.  A digression on the local geometry of
individual
hypersurfaces in affine\index{affine!hypersurface} space follows this.
We then turn to the very rich equivalence
problem for neo-classical Poincar\'e-Cartan forms; the differential 
invariants that this uncovers include those of the various associated affine 
hypersurfaces.  In the last section of this chapter, we use these
differential invariants to characterize systems locally
contact-equivalent to those for prescribed mean
curvature hypersurfaces\index{mean curvature!prescribed} in Riemannian
manifolds.

In the next chapter, we will specialize to the study of those neo-classical 
Poincar\'e-Cartan forms whose primary differential invariants all vanish.    
These correspond to interesting variational problems arising in conformal 
geometry.

\

We begin with a few elementary notions used in the method of equivalence.  On a
manifold $M$ of
dimension $n$, a {\em coframe}\index{coframe bundle|(} at a point $x\in M$ is
a linear isomorphism 
$$
  u_x : T_xM\stackrel{\sim}{\longrightarrow}\R^n.
$$
This is equivalent to a choice of basis for the cotangent space 
$T_x^*M$, and we will not maintain any distinction between these two
notions.  The set of all coframes for $M$ has the structure of a
principal $GL(n,\R)$-bundle $\pi:\F(M)\to M$, with right-action
$$
  u_x\cdot g \stackrel{\mathit{def}}{=} g^{-1}u_x,\qquad g\in GL(n,\R),
$$
where the right-hand side denotes composition of $u_x$ with multiplication by
$g^{-1}$.  A local section of $\pi:\F(M)\to M$ is called a {\em
  coframing},
or {\em coframe field}.  On the total space $\F(M)$, there is an 
$\R^n$-valued {\em tautological $1$-form}\index{tautological $1$-form|(}
$\omega$, given at $u\in\F(M)$ by
\begin{equation}
  \omega_{u}(v) = u(\pi_*v)\in\R^n,\qquad v\in T_u\F(M).
\label{DefTaut5}
\end{equation}
The $n$ components $\omega^i$ of this $\R^n$-valued $1$-form give a global 
basis for the semibasic $1$-forms of $\F(M)\to M$.

In terms of coordinates $x=(x^1,\ldots,x^n)$ on $M$, there is a trivialization
$M\times GL(n,\R)\cong\F(M)$ given by
$$
  (x,g)\leftrightarrow (x, g^{-1}dx),
$$
where on the right-hand side, $dx$ is a column of $1$-forms regarded
as a coframe at $x$, and $g^{-1}dx$ is the composition of that coframe
with multiplication by $g^{-1}\in GL(n,\R)$.  In this trivialization, we
can express the tautological $1$-form as
$$
  \omega = g^{-1}dx,
$$
where again the right-hand side represents the product of a
$GL(n,\R)$-valued fiber coordinate and an $\R^n$-valued
semibasic\index{semibasic form}
$1$-form.
\index{tautological $1$-form|)}

The geometric setting of the equivalence method is the following.
\begin{Definition}  Let $G\subset GL(n,\R)$ be a subgroup.  A
  {\em $G$-structure} on the $n$-manifold $M$ is a
  principal subbundle of the coframe bundle $\mathcal F(M)\to M$,
  having structure group $G$.
\end{Definition}
We will associate to a hyperbolic Monge-Ampere system (to be defined,
in case $n=2$), or to a
neo-classical Poincar\'e-Cartan form (in case $n\geq 3$), a succession
of $G$-structures on the contact manifold $M$, which carry
increasingly detailed information about the geometry of the system or
form, respectively.
\index{Gstructure@$G$-structure|)}
\index{equivalence method|)}
\index{coframe bundle|)}

\section{The Equivalence
  Problem for $n=2$}
\markright{2.1.  THE EQUIVALENCE PROBLEM FOR $n=2$}
\label{Section:SmallEquiv}
\index{equivalence!of Monge-Ampere systems|(}

In this section, we will study the equivalence problem for certain
Monge-Ampere systems on contact manifolds of dimension $5$.  We will
give criteria in terms of the differential invariants thus obtained
for a given system to be locally equivalent the system associated
to the linear homogeneous wave 
equation\index{wave equation!homogeneous}.  We will also give the
weaker criteria for a given system to be
locally equivalent to an Euler-Lagrange system, as in the previously
discussed inverse problem\index{inverse problem}.  Unless otherwise noted,
we use the index ranges $0\leq a,b,c\leq 4$, $1\leq i,j,k\leq 4$.

We assume given a $5$-dimensional contact manifold
$(M,I)$\index{contact!manifold} and a
Monge-Ampere system ${\mathcal E}$, locally algebraically generated as
$$
  {\mathcal E} = \{\theta,d\theta,\Psi\},
$$
where $0\neq\theta\in\Gamma(I)$ is a contact form, and $\Psi\in\Omega^2(M)$
is some $2$-form.  As noted previously, $\mathcal E$ determines $I$ and
$\mathcal I$.  We assume that $\Psi_x\notin\I_x$ for all $x\in M$.  
Recall from the discussion in \S\ref{Subsection:Inverse} that
given ${\mathcal E}$, the generator
$\Psi$ may be uniquely chosen modulo $\{I\}$ (and modulo
multiplication by functions) by the condition of
primitivity\index{primitive|(}; that is, we may assume
$$
  d\theta\wedge\Psi \equiv 0\pmod{\{I\}}.
$$
The assumption $\Psi_x\notin\I_x$ means that this primitive form is
non-zero everywhere.
We do {\em not} necessarily assume that ${\mathcal E}$ is
Euler-Lagrange.

On the contact manifold $M$, one can locally find a coframing 
$\eta=(\eta^a)$
such that
\begin{equation}
  \left\{\begin{array}{l} \eta^0\in\Gamma(I), \\
  d\eta^0\equiv \eta^1\wedge\eta^2 + \eta^3\wedge\eta^4 \pmod{\{I\}}.
  \end{array}\right.
\label{ContactDerivative}
\end{equation}
Then we can write $\Psi\equiv \frac12b_{ij}\eta^i\wedge\eta^j \mbox{ (mod 
$\{I\}$)}$, where the functions $b_{ij}$ depend on the choice of coframing
and on the choice of $\Psi$.  The assumption that $\Psi$ is primitive
means that in terms of a coframing satisfying (\ref{ContactDerivative}),
$$
  b_{12}+b_{34} = 0.
$$
We now ask what further conditions may be imposed on the coframing 
$\eta=(\eta^a)$ while preserving (\ref{ContactDerivative}).

To investigate this, we first consider changes of coframe that fix 
$\eta^0$; we will later take into account non-trivial rescalings of 
$\eta^0$.  In this case, an element of $GL(5,\R)$ preserves the condition 
(\ref{ContactDerivative}) if and only if it acts as a fiberwise {\em 
sympletic} transformation,\index{symplectic!transformation|(} modulo the
contact line bundle $I$.  Working  modulo $I$, we can split
$$
  \bw{2}(T^*M/I)\cong(\R\cdot d\eta^0)\oplus P^2(T^*M/I),
$$
where $P^2(T^*M/I)$ is the $5$-dimensional space of $2$-forms that are
primitive with 
respect to the symplectic structure on $I^\perp$ induced by
$d\eta^0$.  The
key observation is that the action of the symplectic group $Sp(2,\R)$ on 
$P^2(\R^4)$ is equivalent to the standard action of the group $SO(3,2)$ on 
$\R^5$.  This is because the symmetric bilinear form 
$\langle\cdot,\cdot\rangle$ on $P^2(T^*M/I)$ defined by
$$
  \psi_1\wedge\psi_2 = \langle\psi_1,\psi_2\rangle(d\eta^0)^2
$$
has signature $(3,2)$ and symmetry group $Sp(2,\R)$.  
Therefore, the orbit decomposition of the space of primitive forms
$\Psi$ modulo $\{I\}$ under
admissible changes of coframe will be a refinement of the standard
orbit decomposition under $SO(3,2)$.
\index{symplectic!transformation|)}
\index{primitive|)}

To incorporate rescaling of $\eta^0$ into our admissible changes of coframe, 
note that a rescaling of $\eta^0$ requires via (\ref{ContactDerivative}) the 
same rescaling of the symplectic form 
$\eta^1\wedge\eta^2+\eta^3\wedge\eta^4$, so we should 
actually allow changes by elements of $GL(5,\R)$ inducing the standard 
action of $CSp(2,\R)$; this is the group that preserves the standard 
symplectic form up to scale.  This in turn corresponds to the 
split-signature conformal group $CO(3,2)$,\index{conformal!group}
which acts on $\R^5$ with three
non-zero orbits:  a negative space, a null space, and a positive space.

The three orbits of this representation correspond to three types
of Monge-Ampere systems:
\begin{itemize}
  \item If $\Psi\wedge\Psi$ is a negative multiple of $d\eta^0\wedge 
d\eta^0$, then the local coframing $\eta$ may be chosen so that in
addition to (\ref{ContactDerivative}), 
$$
  \Psi\equiv \eta^1\wedge\eta^2-\eta^3\wedge\eta^4\pmod{\{I\}};
$$
for a classical variational problem, this occurs when the Euler-Lagrange PDE 
is hyperbolic.\index{Monge-Ampere system!hyperbolic|(}
  \item If $\Psi\wedge\Psi=0$, then $\eta$ may be chosen so that
$$
  \Psi\equiv \eta^1\wedge\eta^3\pmod{\{I\}};
$$
for a classical variational problem, this occurs when the Euler-Lagrange PDE 
is parabolic.
  \item If $\Psi\wedge\Psi$ is a positive multiple of $d\eta^0\wedge 
d\eta^0$, then $\eta$ may be chosen so that
$$
  \Psi\equiv \eta^1\wedge\eta^4-\eta^3\wedge\eta^2\pmod{\{I\}};
$$
for a classical variational problem, this occurs when the Euler-Lagrange PDE 
is elliptic.
\end{itemize}
The equivalence problem for elliptic Monge-Ampere systems in case 
$n=2$ develops in analogy with that for hyperbolic systems; we will
present the hyperbolic case.   The conclusion will be:
\begin{quote}
{\em
  Associated to a hyperbolic Monge-Ampere 
  system\index{Monge-Ampere system!hyperbolic} $(M^5,{\mathcal E})$
  is a canonical subbundle $B_1\to M$ of the coframe bundle of $M$
  carrying a pair of $2\times 2$-matrix-valued functions $S_1$ and
  $S_2$, involving up to second derivatives of the given system.  
  $(M,{\mathcal E})$ is locally of Euler-Lagrange
  type\index{Euler-Lagrange!system} if and only if $S_2$ vanishes
  identically, while it is
  equivalent to the system associated to the homogeneous wave equation
  $z_{xy}=0$ if and only if $S_1$ and $S_2$ both
  vanish identically.}
\end{quote}
An example of a hyperbolic Monge-Ampere
system, to be studied in more detail in Chapter 4, is the linear
Weingarten system for surfaces in $\mathbf E^3$ with Gauss
curvature\index{Gauss curvature} $K=-1$. 

To begin, assume that $(M^5,{\mathcal E})$ is a hyperbolic Monge-Ampere 
system.
A coframing $\eta=(\eta^a)$ of $M$ is said to be {\em $0$-adapted} to 
${\mathcal E}$ if
\begin{equation}
  {\mathcal E}=\{\eta^0,\ \eta^1\wedge\eta^2+\eta^3\wedge\eta^4,\
      \eta^1\wedge\eta^2-\eta^3\wedge\eta^4\}
\label{ZeroAdapt5.1}
\end{equation}
and also
\begin{equation}
  d\eta^0\equiv\eta^1\wedge\eta^2+\eta^3\wedge\eta^4\pmod{\{I\}}.
\label{ZeroAdapt5.2}
\end{equation}
\index{equivalence method|(}
\index{Gstructure@$G$-structure|(}
According to the following proposition, a hyperbolic Monge-Ampere system is 
equivalent to a certain type of $G$-structure, and it is the latter to which 
the equivalence method directly applies.
\begin{Proposition}
The $0$-adapted coframings for a hyperbolic Monge-Ampere system 
$(M^5,{\mathcal E})$ are the sections of a $G_0$-structure on $M$, where 
$G_0\subset GL(5,\R)$ is the (disconnected) subgroup generated by all 
matrices of the form (displayed in blocks of size $1,2,2$)
\begin{equation}
  g_0=\left(\begin{array}{ccc}
      a & 0 & 0 \\
      C & A & 0 \\
      D & 0 & B \end{array}\right),
\label{BlockForm}
\end{equation}
with $a=\mbox{det}(A)=\mbox{det}(B)\neq 0$, along with the matrix
\begin{equation}
  J=\left(\begin{array}{ccc} 1 & 0 & 0 \\ 0 & 0 & I_2 \\
    0 & I_2 & 0 \end{array}\right).
\label{DefJ}
\end{equation}
\label{GStrProp5}
\end{Proposition}
\begin{Proof}
The content of this proposition is that any two $0$-adapted coframes differ 
by multiplication by an element of $G_0$.  To see why this is so, note that 
the $2$-forms
$$
  \eta^1\wedge\eta^2+\eta^3\wedge\eta^4\quad\mbox{and}\quad
  \eta^1\wedge\eta^2-\eta^3\wedge\eta^4
$$
have, up to scaling, exactly $2$ decomposable\index{decomposable form}
linear combinations, $\eta^1\wedge\eta^2$ and 
$\eta^3\wedge\eta^4$.  These must be either preserved or exchanged by any 
change of coframe preserving their span modulo $\{I\}$, and this accounts for 
both the block form (\ref{BlockForm}) and the matrix $J$.  The condition on 
determinants then corresponds to (\ref{ZeroAdapt5.2}).
\end{Proof}

\

Although not every $G_0$-structure on a $5$-manifold $M$ is induced by a
hyperbolic Monge-Ampere system $\mathcal E$, it is easy to see that
those that do determine $\mathcal E$ uniquely.
We therefore make a digression to describe the first steps of the
equivalence
method, by which one investigates the local geometry of a general
$G$-structure.  This will be followed by application to the case at
hand of a $G_0$-structure, then a digression on
the next general steps, and application to the case at hand, and so
on.  One major step, that of 
{\em prolongation}\index{prolongation!of a $G$-structure},
will not appear in
this chapter but will be discussed in the study of conformal
geometry in Chapter 3.

\

Fix a subgroup $G\subset GL(n,\R)$.  Two
$G$-structures $B_i\to M_i$, $i=1,2$, are {\em
  equivalent}\index{equivalence!of $G$-structures} if there is
a diffeomorphism $M_1\to M_2$ such that under the induced isomorphism
of principal coframe bundles $\mathcal F(M_1)\to\mathcal F(M_2)$, the
subbundle $B_1\subset\mathcal F(M_1)$ is mapped to $B_2\subset\mathcal
F(M_2)$.  One is typically interested only in those properties of a
$G$-structure which are preserved under this notion of equivalence.
For instance, if one has a pair of $5$-manifolds with hyperbolic
Monge-Ampere systems, then a diffeomorphism of the $5$-manifolds
carries one of these systems to the other if and only if it induces an
equivalence of the associated $G_0$-structures.

It is easy to see that a diffeomorphism $F:B_1\to B_2$ between the
total spaces of two $G$-structures $B_i\to M_i$ is an equivalence in
the above sense if and only if $F^*(\omega_2)=\omega_1$, where $\omega_i$
is the restriction of the tautological $\R^n$-valued
form\index{tautological $1$-form|(} (\ref{DefTaut5}) on $\mathcal
F(M_i)\supseteq B_i$.  The first step in investigating the
geometry of a $G$-structure $B\to M$ is therefore to understand the
local behavior of this tautological form.  To do this, we seek an
expression for its exterior derivative, and to understand what such an
expression should look like, we proceed as follows.  

Consider a local trivialization $B\cong M\times G$, induced by a
choice of section $\eta$ of $B\to M$ whose image is identified with
$M\times\{e\}\subset M\times G$.  The section $\eta$ is in particular
an $\R^n$-valued $1$-form on $M$, and the tautological $1$-form is
$$
  \omega = g^{-1}\eta\in\Omega^1(B)\otimes\R^n.
$$
The exterior derivative of this equation is
\begin{equation}
  d\omega = -g^{-1}dg\wedge\omega+g^{-1}d\eta.
\label{CoordStreqn5}
\end{equation}
Note that the last term in this equation is semibasic
for $B\to M$,
and that the matrix $1$-form $g^{-1}dg$ takes values in the Lie
algebra $\lie{g}$ of $G$.  Of course, these pieces 
$g^{-1}d\eta$ and $g^{-1}dg$
each depend on the choice of trivialization.  To better understand the
pointwise linear algebra of (\ref{CoordStreqn5}), we introduce the following
notion. 
\begin{Definition}
A {\em pseudo-connection}\index{pseudo-connection|(} in the
$G$-structure $B\to M$ is a
$\lie{g}$-valued $1$-form on $B$ whose restriction to the fiber
tangent spaces $\mathcal V_b\subset T_bB$ equals the identification
$\mathcal V_b\cong\lie{g}$ induced by the right $G$-action on $B$.
\end{Definition}
This differs from the definition of a {\em
  connection}\index{connection} in the 
principal bundle $B\to M$ by omission of an equivariance requirement.
In terms of our trivialization above, a pseudo-connection on $M\times
G$ is any $\lie{g}$-valued $1$-form of the form
$$
  g^{-1}dg+\mbox{(semibasic $\lie{g}$-valued $1$-form)};
$$
in particular, every $G$-structure carries a pseudo-connection.
A consequence of (\ref{CoordStreqn5}) is that any pseudo-connection
$\varphi\in\Omega^1(B)\otimes\lie{g}$ satisfies a structure equation
that is fundamental for the equivalence method:
\begin{equation}
\boxed{
  d\omega = -\varphi\wedge\omega+\tau,
}
\label{StructureEqn1.1}
\end{equation}
where $\tau=(\sf12T^i_{jk}\omega^j\wedge\omega^k)$ is
a semibasic\index{semibasic form} $\R^n$-valued $2$-form on
$B$, called the {\em torsion}\index{torsion!of a pseudo-connection|(} of
the pseudo-connection $\varphi$.  It
is natural to consider exactly how a different choice of
pseudo-connection---remember that any two differ by an arbitrary
semibasic $\lie{g}$-valued $1$-form---yields a different torsion
form.  We will pursue this after considering the situation for our
hyperbolic Monge-Ampere systems.

\

Let $B_0\subset\F(M)$ be the $G_0$-bundle of $0$-adapted coframes for
a hyperbolic Monge-Ampere system ${\mathcal E}$.  A local
section $\eta$ corresponds to an $\R^5$-valued $1$-form $(\eta^a)$
satisfying (\ref{ZeroAdapt5.1}, \ref{ZeroAdapt5.2}).  In terms of the
trivialization $B_0\cong M\times G_0$ induced by $\eta$, the
tautological $\R^5$-valued $1$-form is $\omega = g_0^{-1}\eta$.
Locally (over neighborhoods in $M$), there is a structure equation
(\ref{StructureEqn1.1}), in which
$$
  \omega = \left(\begin{array}{c} \omega^0 \\ \omega^1 \\ \omega^2 \\
    \omega^3 \\ \omega^4 \end{array}\right) \quad \mbox{and} \quad
  \varphi = \left(\begin{array}{ccccc}
    \varphi^0_0 & 0 & 0 & 0 & 0 \\
    \varphi^1_0 & \varphi^1_1 & \varphi^1_2 & 0 & 0 \\
    \varphi^2_0 & \varphi^2_1 & \varphi^2_2 & 0 & 0 \\
    \varphi^3_0 & 0 & 0 & \varphi^3_3 & \varphi^3_4 \\
    \varphi^4_0 & 0 & 0 & \varphi^4_3 & \varphi^4_4
  \end{array}\right)
$$
are the tautological $\R^5$-valued $1$-form and the pseudo-connection form, 
respectively; note that the condition for $\varphi$ to be $\lie{g}_0$-valued 
includes the condition
$$
  \varphi^0_0 = \varphi^1_1+\varphi^2_2 =\varphi^3_3+\varphi^4_4.
$$
The torsion $\tau$ of $\varphi$ is an $\R^5$-valued $2$-form, 
semibasic for $B_0\to M$ and depending on a choice of pseudo-connection.

\

Returning to the general situation of a $G$-structure $B\to M$, our
goal is to understand how different choices of pseudo-connection in
(\ref{StructureEqn1.1}) yield different torsion forms.  We will use
this to restrict attention to those pseudo-connections whose torsion
is in some normal form.

The linear-algebraic machinery for this is as follows.  Associated to
the linear Lie algebra $\lie{g}\subset\lie{gl}(n,\R)$ is a map of
$G$-modules
$$
  \delta:\lie{g}\otimes(\R^n)^*\to
    \R^n\otimes\bw{2}(\R^n)^*,
$$
defined as the restriction to
\begin{equation}
  \lie{g}\otimes(\R^n)^* \subset (\R^n\otimes (\R^n)^*)\otimes(\R^n)^*
\label{LieInclusion5}
\end{equation}
of the surjective skew-symmetrization map
$$
  \R^n\otimes (\R^n)^*\otimes(\R^n)^* \to \R^n\otimes\bw{2}(\R^n)^*.
$$
The cokernel of $\delta$
\begin{equation}
  H^{0,1}(\lie{g}) \stackrel{\mathit{def}}{=} 
    (\R^n\otimes\bw{2}(\R^n)^*)/\delta(\lie{g}\otimes(\R^n)^*)
\label{DefCokernel5}
\end{equation}
is one of the {\em Spencer cohomology 
groups}\index{Spencer cohomology} of
$\lie{g}\subset\lie{gl}(n,\R)$.  Note that to each $b\in B$ is
associated an isomorphism $T_{\pi(b)}M\stackrel{\sim}{\to}\R^n$, and
consequently an identification of semibasic $1$-forms at $b\in B$ with
$(\R^n)^*$.  Now, given a pseudo-connection in the
$G$-structure, the semibasic $\R^n$-valued torsion $2$-form
$(\frac12T^i_{jk}\omega^j\wedge\omega^k)$ at $b\in B$ can be identified
with an element
$\tau_b\in\R^n\otimes\Lambda^2(\R^n)^*$.  Similarly, a permissible change at
$b\in B$ of the pseudo-connection---that is, a semibasic $\lie{g}$-valued
$1$-form---can be identified with an element
$\varphi_b^\prime\in\lie{g}\otimes(\R^n)^*$.  Under these two
identifications, the map $\delta$ associates to a
change $\varphi_b^\prime$ the corresponding change in the torsion
$\varphi_b^\prime\wedge\omega_b$, where in this expression we have
contracted the middle factor of $(\R^n)^*$ in $\varphi_b^\prime$
(see (\ref{LieInclusion5})) with the
values of the $\R^n$-valued $1$-form $\omega_b$.
Therefore, different
choices of pseudo-connection yield torsion maps differing by elements of
$Im(\delta)$, so what is determined by the $G$-structure alone,
independent of a choice of pseudo-connection, is a
map $\bar\tau\!:\!B\to H^{0,1}(\lie{g})$, called the {\it intrinsic
  torsion}\index{torsion!intrinsic, of a $G$-structure} of $B\to M$.

This suggests a major step in the equivalence method, called {\em
  absorption of torsion}\index{torsion!absorption of|(}, which one
implements by choosing a (vector space) splitting of the projection
\begin{equation}
  \R^n\otimes\bw{2}(\R^n)^*\to H^{0,1}(\lie{g})\to 0.
\label{Surjection5}
\end{equation}
As there may be no $G$-equivariant splitting, one is merely choosing
some vector subspace $T\subset \R^n\otimes\bigwedge^2(\R^n)^*$ which
complements the kernel $\delta(\lie{g}\otimes(\R^n)^*)$.  Fixing a
choice of $T$, it holds by construction that any
$G$-structure $B\to M$ locally has pseudo-connections whose torsion at
each $b\in B$ corresponds to a tensor lying in $T$.  

\

We will see from our example of hyperbolic Monge-Ampere systems that
this is not as complicated as it may seem.
Denote the semibasic $2$-form components of the $\R^5$-valued torsion by
$$
  \tau = \left(\begin{array}{c}
    \tau^0 \\
    \tau^1 \\
    \tau^2 \\
    \tau^3 \\
    \tau^4 \end{array}\right).
$$
We know from the condition (\ref{ZeroAdapt5.2}) in the
definition of $0$-adapted that
$$
  \tau^0 \stackrel{\mathit{def}}{=}
   d\omega^0 + \varphi^0_0\wedge\omega^0 =
   \omega^1\wedge\omega^2+\omega^3\wedge\omega^4 
     +\sigma\wedge\omega^0
$$
for some semibasic $1$-form $\sigma$.  We may now replace $\varphi^0_0$ by 
$\varphi^0_0-\sigma$ in our pseudo-connection, eliminating the term 
$\sigma\wedge\omega^0$ from the torsion.  We then rename this altered 
pseudo-connection entry again as $\varphi^0_0$; to keep the
pseudo-connection $\lie{g}_0$-valued, we have to make a similar change
in $\varphi^1_1+\varphi^2_2$ and $\varphi^3_3+\varphi^4_4$.  
What we have just shown is that given an arbitrary pseudo-connection
in a $G_0$-structure $B_0\to M$, there is another pseudo-connection
whose torsion satisfies (using obvious coordinates on
$\R^n\otimes\bigwedge^2(\R^n)^*$) $T^0_{0a}=T^0_{a0}=0$.  By choosing this
latter pseudo-connection, we are absorbing the
corresponding torsion components into $\varphi$.  
Furthermore, the fact that our
$G_0$-structure is not arbitrary, but comes from a hyperbolic
Monge-Ampere system, gave us the additional information that
$T^0_{12}=T^0_{34}=1$, and all other independent $T^0_{ij}=0$.  
Note incidentally that our 
decision to use a pseudo-connection giving $\sigma=0$ determines
$\varphi^0_0$ uniquely, up to addition of multiples of $\omega^0$;
this uniqueness applies also to $\varphi^1_1+\varphi^2_2$ and  
$\varphi^3_3+\varphi^4_4$.  The effort to uniquely determine
pseudo-connection forms should guide the choices one makes in the
equivalence method.

Other torsion terms may be absorbed using similar methods. 
Using the index range $1\leq i,j,k\leq 4$, we write
$$
  \tau^i = T^i_{j0}\omega^j\wedge\omega^0 +
    \sf12T^i_{jk}\omega^j\wedge\omega^k
$$
for functions $T^i_{j0}$ and $T^i_{jk}=-T^i_{kj}$.  First, by altering the 
nilpotent part $\varphi^i_0$, we can arrange that all $T^i_{j0}=0$.  
Second, by altering the off-diagonal terms $\varphi^1_2$, $\varphi^2_1$, 
$\varphi^4_3$, $\varphi^3_4$, we can arrange that
$$
   T^1_{2j}=T^2_{1j}=T^3_{4j}=T^4_{3j}=0.
$$
Third, by altering the traceless diagonal parts $\varphi^1_1-\varphi^2_2$ 
and $\varphi^3_3-\varphi^4_4$, we can arrange that
$$
  T^1_{13}=T^2_{23},\ T^1_{14}=T^2_{24},\ T^3_{13}=T^4_{14},\
    T^3_{23}=T^4_{24}.
$$
We summarize this by renaming
\begin{eqnarray*}
  \tau^1 &=& (V_3\omega^3+V_4\omega^4)\wedge\omega^1+
     U^1\omega^3\wedge\omega^4, \\
  \tau^2 &=& (V_3\omega^3+V_4\omega^4)\wedge\omega^2+
     U^2\omega^3\wedge\omega^4, \\
  \tau^3 &=& (V_1\omega^1+V_2\omega^2)\wedge\omega^3+
     U^3\omega^1\wedge\omega^2, \\
  \tau^4 &=& (V_1\omega^1+V_2\omega^2)\wedge\omega^4+
     U^4\omega^1\wedge\omega^2,
\end{eqnarray*}
for $8$ torsion functions $U_i$, $V_i$ on $B_0$.  The collection of
torsion tensors $(T^a_{bc})$ taking this form, and satisfying
$T^0_{0a}=T^0_{a0}=0$, constitutes the
splitting of (\ref{Surjection5}) given in the general discussion, to
which we will return shortly.

At this point, we can uncover more consequences of the fact
that we are dealing not with an arbitrary $G_0$-structure on a
$5$-manifold, but a special one induced by a hyperbolic Monge-Ampere
system.  We already found as one consequence the fact that 
$$
  \tau^0 \equiv \omega^1\wedge\omega^2+\omega^3\wedge\omega^4
    \pmod{\{\omega^0\}},
$$
which has nothing to do with our choices in absorbing torsion;
absorbing torsion allowed us to render this congruence into an
equality.\index{torsion!absorption of|)}  Similarly, we now
obtain pointwise relations among other torsion coefficients 
by computing, modulo $\{I\}$ (which 
in this case means ignoring all $\omega^0$ terms after differentiating),
\begin{eqnarray*}
  0 & \equiv & d(d\omega^0) \\ & \equiv &
    \varphi^0_0\wedge(\omega^1\wedge\omega^2+\omega^3\wedge\omega^4) \\
    & & \qquad + ((-\varphi^1_1+V_3\omega^3+V_4\omega^4)\wedge\omega^1
                    +U^1\omega^3\wedge\omega^4)\wedge\omega^2 \\
    & & \qquad -\omega^1\wedge((-\varphi^2_2+V_3\omega^3+V_4\omega^4)
              \wedge\omega^2+U^2\omega^3\wedge\omega^4) \\
    & & \qquad + ((-\varphi^3_3+V_1\omega^1+V_2\omega^2)\wedge\omega^3
                    +U^3\omega^1\wedge\omega^2)\wedge\omega^4 \\
    & & \qquad -\omega^3\wedge((-\varphi^4_4+V_1\omega^1+V_2\omega^2)
              \wedge\omega^4+U^4\omega^1\wedge\omega^2) \\
    & \equiv & (U^1+2V_2)\omega^2\wedge\omega^3\wedge\omega^4
              -(U^2-2V_1)\omega^1\wedge\omega^3\wedge\omega^4 \\
    & & \qquad
              +(U^3+2V_4)\omega^1\wedge\omega^2\wedge\omega^4
              -(U^4-2V_3)\omega^1\wedge\omega^2\wedge\omega^3,
\end{eqnarray*}
so that
$$
  U^1=-2V_2,\ U^2=2V_1,\ U^3=-2V_4,\ U^4=2V_3.
$$
These are pointwise linear-algebraic relation among our $8$ torsion
functions.  

\

In the general study of $G$-structures $B\to M$, we now have to
consider the group action in more detail.  Specifically,
$H^{0,1}(\lie{g})$ is the cokernel of a map of $G$-modules, so it
inherits a $G$-action as well, and it is easy to see that the
intrinsic torsion\index{torsion!intrinsic, of a $G$-structure|(} $\bar\tau:B\to
H^{0,1}(\lie{g})$ is equivariant for
this action.  Therefore, there is an induced map
$$
  [\bar\tau]:M\to H^{0,1}(\lie{g})/G,
$$
which is an invariant of the equivalence class of the $G$-structure
$B\to M$; that is, under a diffeomorphism $M_1\to M_2$ inducing an
equivalence of $G$-structures, $[\tau_2]$ must pull back to
$[\tau_1]$.  Now, $H^{0,1}(\lie{g})/G$ typically has a
complicated topology, and is rarely a manifold.  However, in many
cases of
interest one can find a {\em slice} $W\subset H^{0,1}(\lie{g})$,
a submanifold whose points all have the same stabilizer $G_1\subset G$,
and which is a cross-section of the orbits which $W$ itself intersects.
If the intrinsic torsion $\bar\tau:B\to H^{0,1}(\lie{g})$ of a
$G$-structure takes values in a union of orbits represented by such a
slice, then the set
$$
  B_1 \stackrel{\mathit{def}}{=} \bar\tau^{-1}(W)
$$
is a smooth principal subbundle of $B\to M$ having structure group
$G_1\subset G$.

The process of reducing to a
subbundle defined as the locus where intrinsic torsion lies in a slice
is called {\em normalizing}\index{torsion!normalization of} the
torsion.  If $G_1$ is a proper
subgroup of $G$, then we can essentially start the process over,
starting with an arbitrary pseudo-connection, absorbing torsion, and
so on.  Typically, one inherits from $B\to M$ some information about
the torsion of the subbundle $B_1\to M$, because the original
structure equations restrict to the submanifold $B_1\subset B$.  We
will see an example of this below.

In practice, one typically studies the $G$-action on
$H^{0,1}(\lie{g})$ by transporting it to the representing vector space
$T\subset\R^n\otimes\bigwedge^2(\R^n)^*$.  If $T$ is not an
invariant subspace of $\R^n\otimes\bigwedge^2(\R^n)^*$, then
typically $G$ will act by {\em affine-linear} motions on $T$.  This
is the case in the next step of our equivalence problem for
hyperbolic Monge-Ampere systems.

\

We have represented the intrinsic 
torsion\index{torsion!intrinsic, of a $G$-structure|)}
of a $G_0$-structure $B_0\to
M$ corresponding to a hyperbolic Monge-Ampere system by $4$
independent functions on $B$; that is, our torsion takes values in a
$4$-dimensional subspace of the lift $T$ of $H^{0,1}(\lie{g}_0)$.
The next step is to determine how the independent torsion functions vary 
along the fibers of $B_0\to M$.  This will be expressed 
infinitesimally, in an equation for the exterior derivative of the torsion 
functions, modulo the space of forms that are semibasic for $B_0\to M$; the
expressions will be in
terms of the pseudo-connection forms which parallelize the fibers.  They are 
obtained as follows.

We first consider the equations for $d\omega^1$, $d\omega^2$.  Taking the 
exterior derivative of each, modulo the algebraic ideal 
$\{\omega^0,\omega^1,\omega^2\}$, yields equivalences of $3$-forms that 
{\em do not involve derivatives of any psuedo-connection forms}, but do 
involve $dU^1$, $dU^2$.  From each of these can be factored the $2$-form 
$\omega^3\wedge\omega^4$, yielding a pair of equivalences modulo
$\{\omega^0,\ldots,\omega^4\}$, expressible in 
matrix form as
$$
  0\equiv d\left(\begin{array}{c} U^1 \\ U^2\end{array}\right)
    + \left(\begin{array}{c}\varphi^1_0 \\ \varphi^2_0\end{array}\right)
    + \left(\begin{array}{cc}\varphi^1_1 & \varphi^1_2 \\
        \varphi^2_1 & \varphi^2_2 \end{array}\right)\cdot
      \left(\begin{array}{c} U^1 \\ U^2\end{array}\right) -
      \varphi^0_0\cdot\left(\begin{array}{c} U^1\\
         U^2\end{array}\right). 
$$
A similar procedure applied to the equations for $d\omega^3$, $d\omega^4$ 
yields the pair
$$
  0\equiv d\left(\begin{array}{c} U^3 \\ U^4\end{array}\right)
    + \left(\begin{array}{c}\varphi^3_0 \\ \varphi^4_0\end{array}\right)
    + \left(\begin{array}{cc}\varphi^3_3 & \varphi^3_4 \\
        \varphi^4_3 & \varphi^4_4 \end{array}\right)\cdot
      \left(\begin{array}{c} U^3 \\ U^4\end{array}\right) -
      \varphi^0_0\cdot\left(\begin{array}{c} U^3\\
         U^4\end{array}\right). 
$$
These describe the derivatives of the functions $U^i$ along the fibers of 
$B_0\to M$.  They are to be interpreted as giving
$$
  \left.\sf{d}{dt}\right|_{t=0}U_i(u\cdot g_t),
$$
where $g_t$ is a path in $G_0$ passing through the identity matrix at $t=0$. 
Exponentiated, we see that the vector-valued functions $(U_1, U_2)$ and 
$(U_3, U_4)$ on $B_0$
each transform by an {\em affine-linear action} of $G_0$ along 
the fibers; that is, they vary by a linear representation composed with a 
translation.\footnote{Strictly speaking, we have only shown that the
  torsion function $(U^i)$ varies by an affine-linear action under the
  {\em identity component} of $G_0$.  What will be important, however,
  is that if $u\in B_0$ satisfies $U^i(u)=0$, then $U^i(u\cdot J)=0$
  as well, and likewise for some matrix in each component where
  $a<0$.  These claims can be verified directly.} 
It is the ``nilpotent'' part of the group, with components 
$g^i_0$, which gives rise to the translation.  Specifically, we have
for $g_0$ as in (\ref{BlockForm})
\begin{eqnarray}
  \left(\begin{array}{c} U^1(u\cdot g_0) \\ U^2(u\cdot g_0)\end{array}\right)
     & = & a A^{-1}
  \left(\begin{array}{c} U^1(u) \\ U^2(u)\end{array}\right)
  - A^{-1}C,
\label{TransRule1} \\
  \left(\begin{array}{c} U^3(u\cdot g_0) \\ U^4(u\cdot g_0)\end{array}\right)
     & = & a B^{-1}
  \left(\begin{array}{c} U^3(u) \\ U^4(u)\end{array}\right)
  - B^{-1}D.
\label{TransRule2}
\end{eqnarray}

Now define a {\em $1$-adapted coframe} to be a $0$-adapted coframe
$u\in B_0$ satisfying $U^i(u)=0$ for $1\leq i\leq 4$.  It then
follows from the above reasoning that the subset $B_1\subset B_0$ of
$1$-adapted coframes 
is a $G_1$-subbudle of $B_0$, where the subgroup $G_1\subset G_0$ is 
generated by the matrix $J$ of (\ref{DefJ}), and by matrices of the
form (again, in
blocks of size $1,2,2$) 
\begin{equation}
  g_1 = \left(\begin{array}{ccc} a & 0 & 0 \\ 0 & A & 0 \\
    0 & 0 & B\end{array}\right),
\label{BlockG1}
\end{equation}
with $a = \mbox{det}(A) = \mbox{det}(B) \neq 0$.  The structure equation 
(\ref{StructureEqn1.1}) on $B_0$ still holds when restricted to $B_1$, with 
$\tau^i|_{B_1}=0$; but the pseudo-connection forms 
$\varphi^i_0|_{B_1}$ are semibasic
for $B_1\to M$, and their contribution 
should be regarded as torsion.\index{torsion!of a pseudo-connection}
With everything now restricted to $B_1$, we write
$$
  \varphi^i_0 = P^i_0\omega^0 + P^i_j\omega^j
$$
and then have
$$
  d\omega = -\varphi\wedge\omega+\tau
$$
with
$$
  \varphi =
  \left(\begin{array}{ccccc} \varphi^0_0 & 0 & 0 & 0 & 0 \\
        0 & \varphi^1_1 & \varphi^1_2 & 0 & 0 \\
        0 & \varphi^2_1 & \varphi^2_2 & 0 & 0 \\
        0 & 0 & 0 & \varphi^3_3 & \varphi^3_4 \\
        0 & 0 & 0 & \varphi^4_3 & \varphi^4_4 \end{array}\right)
  \mbox{ and }\tau =
  \left(\begin{array}{c} \omega^1\wedge\omega^2+\omega^3\wedge\omega^4 \\
       -P^1_j\omega^j\wedge\omega^0 \\
       -P^2_j\omega^j\wedge\omega^0 \\
       -P^3_j\omega^j\wedge\omega^0 \\
       -P^4_j\omega^j\wedge\omega^0 \end{array}\right).
$$
As before, we can absorb some of this torsion into the pseudo-connection 
form, respecting the constraint $\varphi^0_0=\varphi^1_1+\varphi^2_2 = 
\varphi^3_3+\varphi^4_4$, until the torsion is of the form
\begin{equation}
  d\omega+\varphi\wedge\omega = \left(\begin{array}{c}
    \omega^1\wedge\omega^2+\omega^3\wedge\omega^4 \\
    -(P\omega^1+P^1_3\omega^3+P^1_4\omega^4)\wedge\omega^0 \\
    -(P\omega^2+P^2_3\omega^3+P^2_4\omega^4)\wedge\omega^0 \\
    -(Q\omega^3+P^3_1\omega^1+P^3_2\omega^2)\wedge\omega^0 \\
    -(Q\omega^4+P^4_1\omega^1+P^4_2\omega^2)\wedge\omega^0
  \end{array}\right).
\label{Torsion1}
\end{equation}
We can go further:  recall that $\varphi^0_0$ was uniquely determined up to 
addition of a multiple of $\omega^0$.  We now exploit this, and take the 
unique choice of 
$\varphi^0_0=\varphi^1_1+\varphi^2_2=\varphi^3_3+\varphi^4_4$ that yields a 
torsion vector of the form (\ref{Torsion1}), with
$$
  P+Q = 0.
$$
Now that $\varphi^0_0$ is uniquely determined, it is reasonable to try to 
get information about its exterior derivative.  To do this, we differentiate 
the equation
$$
  d\omega^0 = -\varphi^0_0\wedge\omega^0 + \omega^1\wedge\omega^2 +
    \omega^3\wedge\omega^4,
$$
which simplifies to
\begin{eqnarray*}
  0 &=& (-d\varphi^0_0 + 2P\omega^1\wedge\omega^2 +
       2Q\omega^3\wedge\omega^4 \\
    &\quad& \quad + (P^2_3-P^4_1)\omega^1\wedge\omega^3
          -(P^1_3+P^4_2)\omega^2\wedge\omega^3 \\
    &\quad& \quad  +(P^2_4+P^3_1)\omega^1\wedge\omega^4
          -(P^1_4-P^3_2)\omega^2\wedge\omega^4)\wedge\omega^0.
\end{eqnarray*}
This tells us the derivative of $\varphi^0_0$ modulo the algebraic ideal 
$\{\omega^0\}=\{I\}$, which we now use in a somewhat unintuitive way.

We easily compute that
$$
  d(\omega^0\wedge\omega^1\wedge\omega^2) =
    -2\varphi^0_0\wedge\omega^0\wedge\omega^1\wedge\omega^2
        +\omega^1\wedge\omega^2\wedge\omega^3\wedge\omega^4.
$$
With knowledge of $d\varphi^0_0\wedge\omega^0$ from above, we can 
differentiate this equation to find
$$
  0 = 2(P-Q)\omega^0\wedge\omega^1\wedge
    \omega^2\wedge\omega^3\wedge\omega^4.
$$
This implies that $P-Q=0$, and combined with our normalization $P+Q=0$, we 
have
$$
  P=Q=0,
$$
which somewhat simplifies our structure equations (\ref{Torsion1}).
In particular, we have modulo $\{I\}$
\begin{equation}
  d\varphi^0_0 \equiv 
     (P^2_3-P^4_1)\omega^1\wedge\omega^3
          -(P^1_3+P^4_2)\omega^2\wedge\omega^3
   +(P^2_4+P^3_1)\omega^1\wedge\omega^4
          -(P^1_4-P^3_2)\omega^2\wedge\omega^4.
\label{short}
\end{equation}

As before, the next step is to study the $8$ torsion coefficients
$P^1_3$, $P^1_4$, $P^2_3$, $P^2_4$, $P^3_1$, $P^3_2$, $P^4_1$,
$P^4_2$.  We can again obtain a 
description of how they vary along the connected components of
the fibers using infinitesimal methods, and then get a
full description of their variation along fibers by explicitly
calculating how they transform under one representative of each
component of the structure group $G_1$.

We state only the result of this calculation.  The torsion in each
fiber transforms by an $8$-dimensional linear representation of the
group $G_1$, which decomposes as the direct sum of two $4$-dimensional
representations.  
Motivated by (\ref{short}), we define a pair of $2\times 2$ matrix-valued
functions on $B_1$
$$
  S_1(u) = \left(\begin{array}{cc}
    P^1_3-P^4_2 & P^1_4+P^3_2 \\ P^2_3+P^4_1 & P^2_4-P^3_1
  \end{array}\right)(u),\quad
  S_2(u) = \left(\begin{array}{cc}
    P^1_3+P^4_2 & P^1_4-P^3_2 \\ P^2_3-P^4_1 & P^2_4+P^3_1
  \end{array}\right)(u).
$$
Now, for $g_1\in G_1$ as in (\ref{BlockG1}), one finds that
$$
  S_1(u\cdot g_1) = aA^{-1}S_1(u)B,\quad
  S_2(u\cdot g_1) = aA^{-1}S_2(u)B.
$$
In particular, the two summand representations for our torsion are the
same, when restricted to the components of $G_1$ of
(\ref{BlockG1}).  However, one may also verify that
\begin{eqnarray*}
  S_1(u\cdot J) & = &  -\left(\begin{array}{cc}
    0 & 1 \\ -1 & 0 \end{array}\right)S_1^t(u)\left(\begin{array}{cc}
    0 & -1 \\ 1 & 0 \end{array}\right), \\
  S_2(u\cdot J) & = & \left(\begin{array}{cc}
    0 & 1 \\ -1 & 0 \end{array}\right)S_2^t(u)\left(\begin{array}{cc}
    0 & -1 \\ 1 & 0 \end{array}\right).
\end{eqnarray*}
An immediate conclusion to be drawn from this is that if $S_1(u)=0$ at
some point, then $S_1(u)=0$ everywhere on the same fiber of $B_1\to
M$, and likewise for $S_2$.

If the torsion vector takes its values in a union of non-trivial
orbits having conjugate stabilizers, then we can try to make a further
reduction to the subbundle consisting of those coframes on which the
torsion lies in a family of normal forms.  However, it is usually
interesting in equivalence problems to consider the case when no
further reduction is possible; in the present situation, this occurs
when all of the invariants vanish identically.
\index{Gstructure@$G$-structure|)}
\index{tautological $1$-form|)}

We first claim that $S_2=0$ identically if and only if the uniquely
determined form $\varphi^0_0$ is closed.  To see this, note first that from
(\ref{short}) we have $S_2=0$ if and only if
$$
  d\varphi^0_0 = \mu\wedge\omega^0
$$
for some $1$-form $\mu$.  We differentiate modulo $\{\omega^0\}$ to obtain
$$
  0 \equiv -\mu\wedge d\omega^0 \pmod{\{\omega^0\}}
$$
which by symplectic linear algebra\index{symplectic!linear algebra}
implies that
$$
  \mu\equiv 0\pmod{\{\omega^0\}}.
$$
But then $d\varphi^0_0 = 0$, as claimed.  Conversely, if
$d\varphi^0_0=0$, then obviously $S_2 = 0$.

Now suppose that $S_1=S_2=0$ identically.  Then because 
$d\varphi^0_0=0$, we can locally find a function $\lambda > 0$ satisfying
$$
  \varphi^0_0 = \lambda^{-1}d\lambda.
$$
We can also compute in case $S_1=S_2=0$ that
$$
   d(\omega^1\wedge\omega^2) = -\varphi^0_0\wedge\omega^1\wedge\omega^2,
$$
so that
$$
  d(\lambda\,\omega^1\wedge\omega^2) = 0.
$$
Now, by a variant of the Darboux theorem,\index{Darboux theorem}
this implies that there are locally defined functions $p$, $x$ such that
$$
  -dp\wedge dx = \lambda\,\omega^1\wedge\omega^2.
$$
Similar reasoning gives locally defined functions $q$, $y$ such that
$$
  -dq\wedge dy = \lambda\,\omega^3\wedge\omega^4.
$$
In terms of these functions, note that
$$
  d(\lambda\,\omega^0) = \lambda(\omega^1\wedge\omega^2
     +\omega^3\wedge\omega^4) = -dp\wedge dx - dq\wedge dy,
$$
which by the Poincar\'e lemma\index{Poincar\'e lemma} implies that
there is another locally defined function $z$ such that
$$
  \lambda\omega^0 = dz-p\,dx - q\,dy.
$$
The linear independence of $\omega^0,\ldots,\omega^4$ implies that pulled 
back by any $1$-adapted coframe (that is, any section of $B_1$), the 
functions $x,y,z,p,q$ form local coordinates on $M$.  In terms of these 
local coordinates, our hyperbolic Monge-Ampere system is
\begin{eqnarray}
  {\mathcal E} &=& \{\omega^0,\ \omega^1\wedge\omega^2 + 
\omega^3\wedge\omega^4,
     \omega^1\wedge\omega^2 - \omega^3\wedge\omega^4\} \\
  &=& \{dz-p\,dx-q\,dy,\ dp\wedge dx + dq\wedge dy,\ dp\wedge dx - dq\wedge 
dy\}.
\label{WaveSystem}
\end{eqnarray}
In an obvious way, transverse local integral surfaces of ${\mathcal E}$ are 
in one-to-one correspondence with solutions to the wave
equation\index{wave equation!homogeneous} for $z(x,y)$
$$
  \frac{\partial^2 z}{\partial x\,\partial y}=0.
$$
This establishes the following.

\begin{Theorem} A hyperbolic Monge-Ampere system $(M^5,{\mathcal E})$ 
satisfies $S_1 = S_2 = 0$ if and only if it is locally equivalent to the 
Monge-Ampere system (\ref{WaveSystem}) for the linear homogeneous wave
equation.
\end{Theorem}
This gives us an easily computable method for determining when a given 
second-order scalar Monge-Ampere equation in two variables is {\em 
contact-equivalent} to this wave equation.

Looking at the equation (\ref{short}) for $d\varphi^0_0$ (mod $\{I\}$), 
it is natural to ask about the situation in which $S_2=0$, but 
possibly $S_1\neq 0$.  This gives an alternative version of the
solution to the inverse
problem\index{inverse problem} discussed in the previous chapter.

\begin{Theorem}
A hyperbolic Monge-Ampere system $(M^5,{\mathcal E})$ is locally
equivalent to an
Euler-Lagrange system\index{Euler-Lagrange!system} if and only if its invariant $S_2$ vanishes identically.
\end{Theorem}
\begin{Proof}
The condition for our ${\mathcal E}$ to contain a Poincar\'e-Cartan form 
$$
  \Pi=\lambda\,\omega^0\wedge
    (\omega^1\wedge\omega^2-\omega^3\wedge\omega^4)
$$
is that this $\Pi$ be {\em closed} for some function $\lambda$ on $B_1$, 
which we can assume satisfies $\lambda>0$.  Differentiating then gives
$$
  0 = (d\lambda-2\lambda\varphi^0_0)\wedge\omega^0\wedge
       (\omega^1\wedge\omega^2-\omega^3\wedge\omega^4).
$$
Exterior algebra shows that this is equivalent to 
$d\lambda-2\lambda\varphi^0_0$ being a multiple of $\omega^0$, say
$$
  d\lambda-2\lambda\varphi^0_0 = \sigma\,\lambda\,\omega^0
$$
for some function $\sigma$,
or in other words,
$$
  d(\log\lambda) - 2\varphi^0_0 = \sigma\,\omega^0.
$$
Such an equation can be satisfied if and only if $d\varphi^0_0$ is
equivalent modulo $\{I\}$ to a multiple of $d\omega^0$.  But we know
that 
$$
  d\omega^0 \equiv \omega^1\wedge\omega^2 + \omega^3\wedge\omega^4,
$$
and from (\ref{short}) we see that $d\varphi^0_0$ is a multiple of
this just in case $S_2=0$.  
\end{Proof}

\

This result may be thought of as follows.  For any hyperbolic
Monge-Ampere system, $d\varphi^0_0\in\Omega^2(B_1)$ is both closed and
semibasic for $B_1\to M$.  This means that it is the pullback of a
$2$-form on $M$ canonically associated with ${\mathcal
  E}$.\footnote{This statement also requires one to verify that
$d\varphi^0_0$ is invariant under the action of some element
of each connected component of $G_1$; this is easily done.}  
We showed
that this $2$-form vanishes if and only if $S_2=0$, which is
equivalent to ${\mathcal E}$ being locally Euler-Lagrange.  This
condition is reminiscent of the vanishing of a curvature, when
$\varphi^0_0$ is viewed as a connection\index{connection} in the
contact line bundle\index{contact!line bundle} $I$.
\index{pseudo-connection|)}
\index{Monge-Ampere system!hyperbolic|)}
\index{torsion!of a pseudo-connection|)}
\index{equivalence!of Monge-Ampere systems|)}
\index{equivalence method|)}

\section{Neo-Classical Poincar\'e-Cartan Forms}
\index{Poincar\'e-Cartan form!neo-classical|(}

We now turn to the geometry of Poincar\'e-Cartan forms in case $n\geq 3$.  
In the preceding section, we emphasized the corresponding Monge-Ampere 
system; from now on, we will instead emphasize the more specialized 
Poincar\'e-Cartan form.

Let $M^{2n+1}$ be a manifold with contact line bundle $I$, locally
generated by a $1$-form $\theta$.  Let $\Pi\in\Omega^{n+1}(M)$ be a
closed $(n+1)$-form locally expressible as
$$
  \Pi = \theta\wedge\Psi,
$$
where $\bar\Psi\in P^n(T^*M/I)$ is primitive modulo $\{I\}$.  As in
the preceding section, the pointwise linear algebra of this data
involves the action of the conformal symplectic group $CSp(n,\R)$ on
the space $P^n(\R^{2n})\subset\bw{n}\R^{2n}$.  When $n=2$, there are
four orbits (including $\{0\}$) for this action, but for $n>2$, the
situation is more complicated.  For example, when $n=3$, the space of
primitive $3$-forms on $\R^6$ has two open orbits and many degenerate
orbits, while for $n=4$ there are no open orbits.

Which orbits contain the Poincar\'e-Cartan forms of most interest to us?
Consider the classical case, in which $M=J^1(\R^n,\R)$, $\theta=dz-p_idx^i$, 
and $\Lambda = L(x,z,p)dx$.  We have already seen that
\begin{eqnarray}
  \Pi & = & d(L\,dx + \theta\wedge L_{p_i}dx_{(i)}) \\
      & = & -\theta\wedge(d(L_{p_i})\wedge dx_{(i)}-L_zdx) \\
      & = & -\theta\wedge(L_{p_ip_j}dp_j\wedge dx_{(i)} +
          (L_{p_iz}p_i+L_{p_ix^i}-L_z)dx).
\label{PCCoords}
\end{eqnarray}
This suggests the following definition, which singles out Poincar\'e-Cartan 
forms of a particular algebraic type; it is these---with a slight refinement 
in the case $n=3$, to be introduced below---whose geometry we will study.
Note that non-degeneracy of the functional is built in to the definition.

\begin{Definition}  A closed $(n+1)$-form $\Pi$ on a contact manifold 
$(M^{2n+1},I)$ is 
{\em almost-classical}\index{Poincar\'e-Cartan form!almost-classical|(}  
if it can locally be expressed as
\begin{equation}
  \Pi = -\theta\wedge(H^{ij}\pi_i\wedge\omega_{(j)}-K\omega)
\label{AlmostClassical}
\end{equation}
for some coframing $(\theta,\omega^i,\pi_i)$ of $M$ with $\theta\in\Gamma(I)$, 
some invertible matrix of functions $(H^{ij})$, and some function $K$.
\label{DefAlmostClassical}
\end{Definition}
Later, we will see the extent to which this definition generalizes the 
classical case.  We remark that the almost-classical forms
$\Pi = \theta\wedge\Psi$ are those for which the primitive $\bar\Psi$ lies
in the tangent variety of the cone of totally
decomposable\index{decomposable form}\footnote{A $n$-form
  is {\em totally decomposable} if it is equal to the exterior product of $n$
  $1$-forms.} $n$-forms in
$P^n(T^*M/I)$, but not in the cone itself.

Applying the equivalence 
method\index{equivalence method|(} will yield differential invariants and 
geometric structures intrinsically associated to our Poincar\'e-Cartan forms. 
  This will be carried out in \S\ref{Section:Bigequiv}, but prior to
  this, it is
best to directly look for some naturally associated geometry.  The preview 
that this provides will make easier the task of interpreting the results of 
the equivalence method.

First note that
the local coframings and functions appearing in the definition of an 
almost-classical form are not uniquely determined by $\Pi$.  The extent of 
the non-uniqueness of the coframings is described in the following
lemma, which prepares us for the equivalence 
method\index{equivalence method|)}.

\begin{Lemma}
If $(\theta,\omega^i,\pi_i)$ is a coframing adapted to an almost-classical 
form $\Pi$ as in Definition~\ref{DefAlmostClassical}, then $(\bar\theta,\bar\omega^i,\bar\pi_i)$ 
is another if and only if the transition matrix is of the form (in blocks of 
size $1,n,n$)
$$
  \left(\begin{array}{c} \bar\theta\\ \bar\omega^i\\ \bar\pi_i
     \end{array}\right) =
  \left(\begin{array}{ccc} a & 0 & 0 \\ C^i & A^i_j & 0 \\
     D_i & E_{ij} & B^j_i\end{array}\right)
  \left(\begin{array}{c} \theta\\ \omega^j\\ \pi_j\end{array}\right).
$$
\label{AlmostClassicalLemma}
\end{Lemma}
\begin{Proof}
That the first row of the matrix must be as shown is clear from the 
requirement that $\theta,\bar\theta\in\Gamma(I)$.  The real content of the 
lemma is that Pfaffian system
$$
  J_\Pi\stackrel{\mathit{def}}{=}
    \mbox{Span}\{\theta,\omega^1,\ldots,\omega^n\} 
$$
is uniquely determined by $\Pi$.  This follows from the claim that $J_\Pi$ 
is characterized as the set of $1$-forms $\xi$ such that $\xi\wedge\Pi$ is 
totally decomposable; this claim we leave as an exercise for the reader.
\end{Proof}

\

The Pfaffian system 
$J_\Pi=\{\theta,\omega^1,\ldots,\omega^n\}$ associated to $\Pi$ 
is crucial for all that follows.  It is canonical in the sense that
any local diffeomorphism of $M$ preserving $\Pi$ also preserves
$J_\Pi$.  In the classical case described
previously we have $J_\Pi=\{dz,dx^1,\ldots,dx^n\}$, which is
integrable and has leaf space $J^0(\R^n,\R)$.

\begin{Proposition}  If $n\geq 4$, then for any almost-classical form $\Pi$ 
on a contact manifold $(M^{2n+1},I)$, the Pfaffian system $J_\Pi$ is 
integrable.
\label{Prop:Neoclassical6}
\end{Proposition}
\begin{Proof}
We need to show that $d\theta$, $d\omega^i\equiv 0\mbox{ (mod
  $\{J_\Pi\}$)}$, for some (equivalently, any) coframing
$(\theta,\omega^i,\pi_i)$ adapted to $\Pi$ as in the definition.  
We write
$$
    \Pi = -\theta\wedge( H^{ij}\pi_i\wedge\omega_{(j)}-K\omega),
$$
and
$$
  d\theta\equiv  a^{ij}\pi_i\wedge\pi_j \pmod{\{J_\Pi\}}.
$$
Then taking those terms of the equation $d\Pi\equiv0\ (\mbox{mod
  }{\{I\}})$ that are 
cubic in $\pi_i$, we find
$$
  a^{ij}\pi_i\wedge\pi_j\wedge H^{kl}\pi_l = 0.
$$
Then the $2$-form $a^{ij}\pi_i\wedge\pi_j$ has at least $n\geq 3$ linearly 
independent $1$-forms as divisors, which is impossible unless 
$a^{ij}\pi_i\wedge\pi_j=0$.  Therefore, $d\theta\equiv 0\mbox{ (mod 
$\{J_\Pi\}$)}$ (with only the hypothesis $n\geq 3$).

For the next step, it is useful to work with the $1$-forms
$$
  \pi^i \stackrel{\mathit{def}}{=}  H^{ij}\pi_j,
$$
and write
$$
  d\omega^i \equiv P^i_{jk}\pi^j\wedge\pi^k\pmod{\{J_\Pi\}},\qquad
    P^i_{jk}+P^i_{kj}=0.
$$
From the form of $\Pi$ (\ref{AlmostClassical}), we have
$$
  0 = \omega^i\wedge\omega^j\wedge\Pi
$$
for any pair of indices $1\leq i,j\leq n$.  Differentiating, we obtain
\begin{eqnarray*}
  0 & = & (d\omega^i\wedge\omega^j - \omega^i\wedge d\omega^j)\wedge\Pi \\
    & = & (P^i_{kl}\pi^k\wedge\pi^l\wedge\pi^j -
             P^j_{kl}\pi^k\wedge\pi^l\wedge\pi^i)\wedge\theta\wedge\omega\\
    & = & (\delta^j_mP^i_{kl}-\delta^i_mP^j_{kl})
      \pi^k\wedge\pi^l\wedge\pi^m\wedge\theta\wedge\omega.
\end{eqnarray*}
It is now an exercise in linear algebra to show that if $n\geq 4$, then this 
implies $P^i_{jk}=0$.  The hypotheses are that $P^i_{jk}=-P^i_{kj}$ and
\begin{equation}
  (\delta^j_mP^i_{kl}-\delta^i_mP^j_{kl}) +
  (\delta^j_kP^i_{lm}-\delta^i_kP^j_{lm}) +
  (\delta^j_lP^i_{mk}-\delta^i_lP^j_{mk}) = 0.
\label{LinAlg6}
\end{equation}
By contracting first on $jk$ and then on $il$, one finds that for $n\neq
2$ the contraction $P^i_{ik}$ vanishes.  Contracting (\ref{LinAlg6})
only on $jk$ and using $P^i_{ik}=0$, one finds that for $n\neq 3$, all
$P^i_{jk}$ vanish.
\end{Proof}

\

There do exist counterexamples in case $n=3$, for which
(\ref{LinAlg6}) implies only that
$$
  d\omega^i = P^{ij}\pi_{(j)},\quad P^{ij}=P^{ji}.
$$
For example, if we fix constants $P^{ij}=P^{ji}$ also satisfying
$P^{ii}=0$, then there is a
unique simply connected, $7$-dimensional Lie group $G$ having a basis
of left-invariant $1$-forms $(\omega^i,\theta,\pi^i)$ satisfying
structure equations
$$
   d\omega^i = P^{ij}\pi_{(j)},\quad
   d\theta = -\pi^i\wedge\omega^i,\quad
   d\pi^i = 0.
$$
In this case, $\theta$ generates a homogeneous contact structure on
$G$, and the form
$$
  \Pi \stackrel{\mathit{def}}{=} -\theta\wedge\pi^i\wedge\omega_{(i)}
$$
is closed, giving an almost-classical form for which
$J_\Pi$ is not integrable.  

These counterexamples cannot arise from classical cases, however, and
this suggests that
we consider the following narrower class of Poincar\'e-Cartan forms.
\begin{Definition}
An almost-classical Poincar\'e-Cartan form $\Pi$ is {\em neo-classical} if 
its associated Pfaffian system $J_\Pi$ is integrable.
\end{Definition}
So the preceding Proposition states that in case $n\geq 4$, every
almost-classical Poincar\'e-Cartan form is neo-classical, and we have 
narrowed the definition only in case $n=3$. 

The foliation corresponding to the integrable Pfaffian system $J_\Pi$ is the 
beginning of the very rich geometry associated to a neo-classical 
Poincar\'e-Cartan form.  Before investigating it further,
we justify the study of this class of objects with the 
following.
\begin{Proposition}
Every neo-classical Poincar\'e-Cartan form $\Pi$ on a contact manifold
$(M,I)$ is locally
equivalent to that arising from some classical variational problem.
More precisely, given such $(M,I,\Pi)$, there are local coordinates
$(x^i,z,p_i)$ on $M$ with respect to which the contact system $I$ is
generated by $dz-p_idx^i$, and there is a Lagrangian of the form
$L(x^i,z,p_i)dx$ whose Poincar\'e-Cartan form is $\Pi$.
\label{ClassicalProp}
\end{Proposition}

Note that we have already observed the converse, that those
non-degenerate Poincar\'e-Cartan forms arising form classical
variational problems (in case $n\geq 3$) are neo-classical.

\

\begin{Proof}
We fix a coframing $(\theta,\omega^i,\pi_i)$ as in the definition of
an almost-classical form.
Using the Frobenius theorem\index{Frobenius theorem}, we take
independent functions $(x^i,z)$ on $M$ so that
$$
  J_\Pi = \{\omega^i,\theta\} = \{dx^i, dz\}.
$$
By relabelling if necessary, we may assume
$\theta\notin\{dx^i\}$, and we find that there are functions $p_i$ 
so that
$$
  \theta \in \R\cdot(dz- p_idx^i).
$$
The fact that $\theta\wedge(d\theta)^n\neq 0$ implies that
$(x^i,z,p_i)$ are local coordinates on $M$.

We now introduce a technical device that is often useful in the study of 
exterior differential systems.  Let
$$
  {\mathcal F}^p\Omega^q\subset\Omega^{p+q}(M)
$$
be the collection of $(p+q)$-forms with at least $p$ factors in $J_\Pi$; 
this is well-defined.  With this notation, the fact that $J_\Pi$ is 
integrable may be expressed as
$$
  d({\mathcal F}^p\Omega^q)\subset {\mathcal F}^p\Omega^{q+1}.
$$
There is a version of the Poincar\'e lemma\index{Poincar\'e lemma}
that can be applied to each leaf of the foliation determined by
$J_\Pi$, with
smooth dependence on the leaves' parameters; it says precisely that the
complex  
$$
  {\mathcal F}^p\Omega^0\stackrel{d}{\longrightarrow}
  {\mathcal F}^p\Omega^1\stackrel{d}{\longrightarrow}\cdots
$$
is locally exact for each $p$.
Now, any almost-classical form $\Pi$ lies in ${\mathcal F}^n\Omega^1$; so 
not only is the closed form $\Pi$ locally equal to $d\Lambda$ for some 
$\Lambda\in\Omega^n(M)$, we can actually choose $\Lambda$ to lie in 
${\mathcal F}^n\Omega^0$.  In other words, we can locally find a Lagrangian 
$\Lambda$ of the form
$$
  \Lambda = L^0(x,z,p)dx + L^i(x,z,p)dz\wedge dx_{(i)}
$$
for some functions $L^0$, $L^i$.  This may be rewritten as
$$
  \Lambda = (L^0+p_iL^i)dx + \theta\wedge(L^idx_{(i)}),
$$
and then the condition $\theta\wedge d\Lambda = 0$ (recall that this
was part of the construction of the Poincar\'e-Cartan form associated
to any class in $H^n(\Omega^*/\I)$) gives the relation
$$
  L^i(x,z,p) = \frac{\partial L}{\partial p_i}(x,z,p).
$$
This is exactly the condition for $\Lambda$ to locally be a classical 
Lagrangian.
\end{Proof}
\index{Poincar\'e-Cartan form!almost-classical|)}

\

Returning to the geometry associated to a neo-classical Poincar\'e-Cartan 
form $\Pi$, we have found (or in case $n=3$, postulated) an integrable 
Pfaffian system $J_\Pi$ which is invariant under contact
transformations\index{transformation!contact|(}
preserving $\Pi$.  Locally in $M$, the induced foliation has 
a smooth ``leaf-space'' $Q$ 
of dimension $n+1$, and  there is a smooth submersion $q:M\to Q$ 
whose fibers are 
$n$-dimensional integral manifolds of $J_\Pi$.  On such a neighborhood,
the foliation will be called {\em simple}\index{simple foliation},
and as we are only going to
consider the local geometry of $\Pi$ in this section, we assume that
the foliation is simple on all of $M$.
We may restrict to smaller neighborhoods
as needed in the following.

To explore the geometry of the situation, we ask what the data 
$(M^{2n+1},I,\Pi)$ look like from the point of view of $Q^{n+1}$.  The 
first observation is that we can locally identify $M$, as a contact 
manifold, with the standard contact manifold $G_n(TQ)$, the
Grassmannian\index{Grassmannian} bundle parameterizing 
$n$-dimensional subspaces of fibers of $TQ$.  This is easily 
seen in coordinates as follows.  If, as in the preceding proof, we
integrate $J_\Pi$ as
$$
  J_\Pi = \{dz,dx^i\}
$$
for some local functions $z,x^i$ on $M$, then the same functions $z,x^i$ may 
be regarded as coordinates on $Q$.  With the assumption that 
$\theta\notin\{dx^i\}$ (on $M$, again), we must have 
$dz-p_idx^i\in\Gamma(I)$ for some local functions $p_i$ on $M$, which by the 
non-degeneracy condition for $I$ make $(x^i,z,p_i)$ local coordinates on 
$M$.  These $p_i$ can also thought of as local fiber coordinates for $M\to 
Q$, and we can map $M\to G_n(TQ)$ by
$$
  (x^i,z,p_i) \mapsto ((x^i,z);\{dz- p_idx^i\}^\perp).
$$
The latter notation refers to a hyperplane in the tangent space of $Q$ at 
$(x^i,z)$.  Under this map, the standard contact system on $G_n(TQ)$ 
evidently pulls back to $I$, so we have a 
local contact diffeomorphism commuting with projections to $Q$.
Every point transformation\index{transformation!point|(} of $Q$
prolongs to give a contact transformation  
of $G_n(TQ)$, hence of $M$ as well.  Conversely, every contact 
transformation of $M$ that preserves $\Pi$ 
is the prolongation of a 
point transformation of $Q$, because the foliation by integral
manifolds of $J_\Pi$ defining $Q$ is associated to $\Pi$ in a
contact-invariant manner.\footnote{This statement is
  only valid in case the foliation  by integral manifolds of 
  $J_\Pi$ is simple; in other
  cases, only a cumbersome local version of the statement holds.}  In
this sense, studying the geometry of a 
neo-classical Poincar\'e-Cartan form (in case $n\geq 3$) under contact 
transformations is locally no different than studying the geometry of
an equivalence class of classical non-degenerate first-order scalar
Lagrangians under point transformations.
\index{transformation!contact|)}
\index{transformation!point|)}

We have now interpreted $(M,I)$ as a natural object in terms of $Q$,
but our real interest lies in $\Pi$.  What kind of geometry does $\Pi$
define in terms of $Q$?  We will answer this question in terms of the
following notion.
\begin{Definition}  A 
{\em Lagrangian potential}\index{Lagrangian potential|(} for a
  neo-classical Poincar\'e-Cartan form $\Pi$ on $M$ is an $n$-form
  $\Lambda\in\mathcal F^n\Omega^0$ (that is, $\Lambda$ is
  semibasic\index{semibasic form} for $M\to Q$) such that $d\Lambda = \Pi$.
\end{Definition}
We saw in the proof of Proposition
\ref{ClassicalProp} that locally a Lagrangian potential $\Lambda$ exists.
Such $\Lambda$ are not
unique, but are determined only up to addition of closed forms in
$\F^n\Omega^0$.  It will be important below to note that a closed form
in $\F^n\Omega^0$ must actually be basic for $M\to Q$;
that is, it
must be locally the pull-back of a (closed) $n$-form on $Q$.  In
particular, the difference between any two Lagrangian potentials for a
give neo-classical form $\Pi$ must be basic.\index{basic form}

Consider one such Lagrangian potential $\Lambda$, semibasic over $Q$.
Then at each point $m\in M$, one may
regard $\Lambda_m$ as an element of $\bw{n}(T^*_{q(m)}Q)$, an $n$-form at 
the corresponding point of $Q$.  This defines a map
$$
  \nu:M\to\bw{n}(T^*Q),
$$
commuting with the natural projections to $Q$.  Counting dimensions shows 
that if $\nu$ is an immersion, then we actually obtain a hypersurface in 
$\bw{n}(T^*Q)$; to be more precise, we have a smoothly varying field of 
hypersurfaces in the vector bundle $\bw{n}(T^*Q)\to Q$.  It is not
hard to see
that $\nu$ is an immersion if the Poincar\'e-Cartan form $\Pi$ is 
non-degenerate, which is a standing hypothesis.  We can work 
backwards, as well:  given a hypersurface
$M\hookrightarrow\bw{n}(T^*Q)$ over an $(n+1)$-dimensional manifold $Q$,
we may restrict to $M$ the tautological $n$-form on $\bw{n}(T^*Q)$ to 
obtain a form $\Lambda\in\Omega^n(M)$.  Under mild technical hypotheses on 
the hypersurface $M$, the form $d\Lambda\in\Omega^{n+1}(M)$ will be a 
neo-classical Poincar\'e-Cartan form.

\index{affine!hypersurface|(}
So we have associated to a Poincar\'e-Cartan form $\Pi$, and a choice
of Lagrangian potential $\Lambda\in\F^n\Omega^0$, a field of 
hypersurfaces in $\bw{n}(T^*Q)\to Q$.  However, we noted that $\Lambda$ 
was not canonically defined in terms of $\Pi$, so neither are these 
hypersurfaces.  As we have seen, the ambiguity in $\Lambda$ is that
another
admissible $\tilde\Lambda$ may differ from $\Lambda$ by a form that is
basic over $Q$.  This means that $\Lambda-\tilde\Lambda$ does not depend on 
the fiber-coordinate for $M\to Q$, and therefore the two corresponding 
immersions $\nu,\tilde\nu$ differ in each fiber $M_q$ ($q\in Q$)
only by a translation in $\bw{n}(T^*_qQ)$.  Consequently, we have
in each $\bw{n}(T^*_qQ)$ a hypersurface well-defined up to
translation.  A contact transformation of $M$ which preserves $\Pi$
will therefore carry the field of hypersurfaces for a particular
choice of $\Lambda$ to a field of hypersurfaces differing by (a field
of) {\em affine transformations}.\index{affine!transformation}

To summarize, 
\begin{quote}{\em
one can canonically associate to any neo-classical
Poincar\'e-Cartan form $(M,\Pi)$ a field of hypersurfaces in the bundle 
$\bw{n}(T^*Q)\to Q$, regarded as a bundle of affine spaces.  We expect 
the differential invariants of $\Pi$ to include information about the geometry 
of each of these affine hypersurfaces, and this will turn out to be the 
case.
}
\end{quote}
\index{Poincar\'e-Cartan form!neo-classical|)}
\index{affine!hypersurface|)}
\index{Lagrangian potential|)}

\section{Digression on Affine Geometry of Hypersurfaces}
\index{affine!hypersurface|(}
\newcommand{\A}{{\mathbf A}}

Let $\A^{n+1}$ denote $(n+1)$-dimensional affine space, which is
simply $\R^{n+1}$ regarded as a homogeneous space of the group
$A(n+1)$ of affine transformations\index{affine!transformation|(} 
$$
  x\mapsto g\cdot x+v,\qquad g\in GL(n+1,\R),\ 
     v\in\R^{n+1}.
$$
Let $x:{\mathbf F}\to\A^{n+1}$ denote the principal
$GL(n+1,\R)$-bundle\index{affine!frame bundle|(} of affine
frames; that is,
$$
  {\mathbf F} = \{f=(x,(e_0,\ldots,e_n))\},
$$
where $x\in\A^{n+1}$ is a point, and $(e_0,\ldots,e_n)$ is a basis for
the tangent space $T_x\A^{n+1}$.
The action is given by
\begin{equation}
  (x,(e_0,\ldots,e_n))\cdot(g^a_b) \stackrel{\mathit{def}}{=}
  (x,(e_bg^b_0,\ldots,e_bg^b_n)).
\label{AffineFrameAction}
\end{equation}
For this section, we adopt the index ranges $0\leq a,b,c\leq n$,
$1\leq i,j,k\leq n$, and always assume $n\geq 2$.

There is a basis of
$1$-forms $\omega^a$, $\varphi^a_b$ on ${\mathbf F}$ defined by decomposing
the $\A^{n+1}$-valued $1$-forms
$$
  dx= e_a\cdot\omega^a,\quad de_a= e_b\cdot\varphi^b_a.
$$
These equations implicitly use a trivialization of $T\A^{n+1}$
that commutes with affine transformations.
Differentiating, we obtain the structure equations for ${\mathbf F}$:
\begin{equation}
  d\omega^a = -\varphi^a_b\wedge\omega^b,\quad
  d\varphi^a_b = -\varphi^a_c\wedge\varphi^c_b.
\label{AffineGpStructure}
\end{equation}
Choosing a reference frame $f_0\in{\mathbf F}$ determines an 
identification ${\mathbf F}\cong A(n+1)$, and under this identification the
$1$-forms $\omega^a$, $\varphi^a_b$ on ${\mathbf F}$ correspond to a basis of
left-invariant $1$-forms on the Lie group $A(n+1)$.  The structure equations 
(\ref{AffineGpStructure}) on ${\mathbf F}$ then correspond to the
usual Maurer-Cartan\index{Maurer-Cartan!equation} 
structure equations for left-invariant $1$-forms on a Lie group.

In this section, we will study the geometry of smooth hypersurfaces 
$M^n\subset\A^{n+1}$, to be called {\em affine hypersurfaces}, using the 
method of moving frames\index{moving frames}; 
no previous knowledge of this method is assumed.  
In particular, we give constructions that associate to $M$ geometric 
objects in a manner invariant under affine transformations of the ambient 
$\A^{n+1}$.  Among these objects are tensor fields $H_{ij}$, $U^{ij}$, and 
$T_{ijk}$ on $M$, called the affine first and second fundamental
forms and the affine cubic form\index{affine!fundamental forms} 
of the hypersurface.  We will classify 
those non-degenerate\index{non-degenerate!affine hypersurface} 
(to be defined) hypersurfaces for which $T_{ijk}=0$ everywhere.
This is of interest because the particular neo-classical 
Poincar\'e-Cartan forms\index{Poincar\'e-Cartan form!neo-classical}
that we study later induce fields of affine hypersurfaces of this
type.

\

Suppose given a smooth affine hypersurface $M\subset\A^{n+1}$.
We define the collection of {\em $0$-adapted frames} along $M$ by
$$
  {\mathbf F}_0(M) = \{(x,(e_0,\ldots,e_n))\in{\mathbf F}: x\in M,\
      e_1,\ldots,e_n\mbox{ span }T_xM\}\subset{\mathbf F}.
$$
This is a principal subbundle of ${\mathbf F}|_M$ whose structure
group is\footnote{Here and throughout, $\R^*$ denotes the connected
  group of positive real numbers under multiplication.}
\begin{equation}
  G_0 \stackrel{\mathit{def}}{=} \left\{g_0=\left(\begin{array}{cc} a
    & 0 \\ v & A\end{array}\right):
    a\in\R^*,\ A\in GL(n,\R),\ v\in\R^n\right\}.
\label{DefG0Affine}
\end{equation}
Restricting forms on ${\mathbf F}$ to 
${\mathbf F}_0(M)$ (but supressing notation), we have
$$
  \omega^0 = 0,\qquad \omega^1\wedge\cdots\wedge\omega^n\neq 0.
$$
Differentiating the first of these gives
$$
   0 = d\omega^0 = -\varphi^0_i\wedge\omega^i,
$$
and we apply the Cartan lemma\index{Cartan lemma} to obtain
$$
  \varphi^0_i =  H_{ij}\omega^j
    \mbox{ for some functions }H_{ij}=H_{ji}.
$$

One way to understand the meaning of these functions $H_{ij}$, which 
constitute the 
{\em first fundamental form}\index{affine!fundamental forms|(} 
of $M\subset\A^{n+1}$, is
as follows.  At any given point of $M\subset\A^{n+1}$, one can find
an affine frame and associated coordinates with respect to which 
$M$ is locally a graph
$$
  x^0 = \sf12 \bar H_{ij}(x^1,\ldots,x^n)x^ix^j
$$
for some functions $\bar H_{ij}$.  Restricted to the $0$-adapted frame
field defined by
$$
  \bar e_0(x) = \frac{\partial}{\partial x^0},\
  \bar e_i(x) = (\bar H_{ij}(x)x^j + \sf12\partial_i\bar
      H_{jk}(x)x^jx^k)\frac{\partial}{\partial x^0} +
      \frac{\partial}{\partial x^i},
$$
one finds that the values over $0\in M$ of the functions $H_{ij}$ equal 
$\bar H_{ij}(0)$.  Loosely speaking, the functions $H_{ij}$ express the second 
derivatives of a defining function for $M$.

Returning to the general situation, we calculate as follows.
We substitute the expression $\varphi^0_i=H_{ij}\omega^j$ into the
structure equation $d\varphi^0_i=-\varphi^0_b\wedge\varphi^b_i$,
collect terms, and conclude
$$
  0 = (dH_{ij}+H_{ij}\varphi^0_0 -H_{kj}\varphi^k_i
         -H_{ik}\varphi^k_j)\wedge\omega^j.
$$
Using the Cartan lemma\index{Cartan lemma}, we have
$$
  dH_{ij} = -H_{ij}\varphi^0_0 + H_{kj}\varphi^k_i
    +H_{ik}\varphi^k_j + T_{ijk}\omega^k
$$
for some functions $T_{ijk}=T_{ikj}=T_{kji}$.  This infinitesimally
describes how the functions $H_{ij}$ vary along the fibers of
${\mathbf F}_0(M)$, on which $\omega^j=0$.  In particular, as a
matrix-valued function $H=(H_{ij})$ on ${\mathbf F}_0(M)$, it
transforms by a linear representation of the structure group:
$$
   H(f\cdot g_0) = (a^{-1})\,^t\!AH(f)A,
$$
where $g_0\in G_0$ is as in (\ref{DefG0Affine}).\footnote{As usual,
  our argument only proves this claim for $g_0$ in the identity
  component of $G_0$, but it may be checked directly for
  representative elements of each of the other components.}
Now we consider the quantity
$$
  \Delta(f)\stackrel{\mathit{def}}{=}\mbox{det}(H_{ij}(f)),
$$ 
which vanishes at
some point of ${\mathbf F}_0(M)$ if and only if it vanishes on the
entire fiber containing that point.  We will say that
$M\subset\A^{n+1}$ is {\em
  non-degenerate}\index{non-degenerate!affine hypersurface|(}
if $\Delta\neq 0$
everywhere on ${\mathbf F}_0(M)$.  Also note that the absolute
signature of $H_{ij}$ is well-defined at each point of $M$.  It is
easy to see that $H_{ij}$ is definite if and only if
$M\subset\A^{n+1}$ is convex.  In what follows, we will assume that
$M$ is a non-degenerate 
hypersurface\index{non-degenerate!affine hypersurface|)}, but not
necessarily that it is
convex.

It turns out that $T=(T_{ijk})$, which one would like to regard as a
sort of covariant derivative of $H=(H_{ij})$, is not a tensor; that
is, it does not transform by a linear representation along the fibers
of ${\mathbf F}_0(M)\to M$.  We will exploit this below to reduce the
principal bundle ${\mathbf F}_0(M)\to M$ to a subbundle of frames
satisfying a higher-order adaptivity condition.  Namely, ${\mathbf
  F}_1(M)\subset{\mathbf F}_0(M)$ will consist of those frames where
$T_{ijk}$ is traceless with respect to the non-degenerate symmetric
bilinear form $H_{ij}$, meaning $H^{jk}T_{ijk}=0$, where $(H^{ij})$ is
the matrix inverse of $(H_{ij})$.  Geometrically, the reduction will
amount to a canonical choice of line field $\R e_0$ transverse to $M$,
which we will think of as giving at each point of $M$ a canonical 
{\em affine normal line}.\index{affine!normal line|(}

To justify this, we let $(H^{ij})$ denote the matrix inverse of
$(H_{ij})$, and let 
$$
  C_i \stackrel{\mathit{def}}{=} H^{jk}T_{ijk}
$$
be the vector of traces of $T$ with respect to $H$.  We compute
\begin{eqnarray*}
  d(\mbox{log }\Delta) & = & \Delta^{-1}d\Delta \\
       & = & \mbox{Tr}(H^{-1}dH) \\
       & = & H^{ij}dH_{ij} \\
       & = & -n\varphi^0_0+2\varphi^i_i+C_i\omega^i.
\end{eqnarray*}
Now differentiate again and collect terms to find
\begin{equation}
   0 = (dC_i- C_j\varphi^j_i-(n+2) H_{ij}\varphi^j_0)
           \wedge\omega^i.
\label{USymmetry7}
\end{equation}
Therefore, we have
\begin{equation}
  dC_i \equiv C_j\varphi^j_i+(n+2)H_{ij}\varphi^j_0
     \pmod{\{\omega^1,\ldots,\omega^n\}},
\label{VarnOfTrace}
\end{equation}
which expresses how the traces $C_i$ vary along the fibers of
${\mathbf F}_0(M)\to M$.
In particular, if the matrix $(H_{ij})$ is non-singular,
as we are assuming, then the action of the structure group on the
values of the vector $(C_i)\in\R^n$ is transitive; that is, every
value in $\R^n$ is taken by $(C_i)$ in each fiber.  Therefore,
the set of $0$-adapted frames $f\in{\mathbf F}_0(M)$ where each
$C_i(f)=0$ is a
principal subbundle ${\mathbf F}_1(M)\subset{\mathbf F}_0(M)$, whose
structure group is the
stabilizer of $0\in\R^n$ under the action.  This stabilizer is
$$
  G_1 \stackrel{\mathit{def}}{=} \left\{g_1 = 
    \left(\begin{array}{cc} a & 0 \\ 0 & A\end{array}\right):
    a\in\R^*,\ A\in GL(n,\R)\right\}.
$$
Comparing to the full action (\ref{AffineFrameAction}) of the affine group 
$A(n+1)$ on ${\mathbf F}$, we see that along each fiber of ${\mathbf
  F}_1(M)$, the direction $\R e_0$ is fixed.
Thus, we have uniquely chosen the direction of $e_0$ at each point of $M$ by 
the condition $C_i=0$ for  $i=1,\ldots,n$.

A more concrete explanation of what we have done is seen by locally
presenting our hypersurface in the form
$$
  x^0 = \frac12\bar H_{ij}(0)x^ix^j+\frac16\bar
    T_{ijk}(x^1,\ldots,x^n)x^ix^jx^k.
$$
An affine change of coordinates that will preserve this form is the 
addition of a multiple of $x^0$ to each $x^i$; the $n$ choices that this 
entails can be uniquely made so that $\bar T_{ijk}(0)$ is traceless with 
respect to $\bar H_{ij}(0)$.  Once such choices are fixed, then so is the 
direction of $\frac{\partial}{\partial x^0}$, and this gives the canonical 
affine normal line at $x=0$.

There is a remarkable interpretation of the affine normal direction
at a point where $H_{ij}$ is positive-definite (see
\cite{Blaschke:Vorlesungen}\index{Blaschke, W.}).
Consider the $1$-parameter family of hyperplanes parallel to the
tangent plane at the given point.  For those planes sufficiently near
the tangent plane, the
intersection with a fixed neighborhood in the surface is a closed submanifold
of dimension $n-2$ in $M$,
having an affine-invariant center-of-mass.  These centers-of-mass form
a curve in affine space, passing through the point of interest; this
curve's tangent line at that point is the affine normal direction.

We can see from (\ref{VarnOfTrace}) that on ${\mathbf F}_1(M)$, where
$T_{ijk}$ is 
traceless, the forms $\varphi^j_0$ are semibasic over $M$.  It is less
convenient to express these in terms of the basis $\omega^i$ than
to instead use $\varphi^0_i = H_{ij}\omega^j$, assuming that
$M\subset\A^{n+1}$ is non-degenerate.  On ${\mathbf F}_1(M)$ we write
$$
  \varphi^j_0 = U^{jk}\varphi^0_k.
$$
Now (\ref{USymmetry7}), restricted to ${\mathbf F}_1(M)$ where
$C_i=0$, implies that $U^{ij}=U^{ji}$.

The reader may carry out computations similar to those above to show that the
$(U^{ij})$ and $(T_{ijk})$ are tensors; that is, they transform along
the fibers of ${\mathbf F}_1(M)$ by a linear representation of $G_1$.
For example, if $T_{ijk}=0$ at some point of ${\mathbf
  F}_1(M)$, then $T_{ijk}=0$ everywhere along the same fiber  
of ${\mathbf F}_1(M)\to M$.  Furthermore, the transformation law for
$U^{ij}$ is such 
that if $U^{ij}=\lambda H^{ij}$ at some point, for some $\lambda$, 
then the same is true---with possibly varying $\lambda$---everywhere on the 
same fiber.  We can now give some additional interpretations of the
simplest cases of the affine second 
fundamental form $U^{ij}$ and the affine cubic form
$T_{ijk}$.  
The following theorem is the main purpose of this digression.

\begin{Theorem} (1) If $U^{ij}=\lambda H^{ij}$ everywhere on 
${\mathbf F}_1(M)$---that is, if the second fundamental form is a
scalar multiple of the first fundamental
form---then either $\lambda= 0$ everywhere or $\lambda\neq0$ everywhere.  In 
the first case, the affine normal lines of $M$ are all parallel, and in the 
second case, the affine normal lines of $M$ are all concurrent.

(2) If $T_{ijk}=0$ everywhere on ${\mathbf F}_1(M)$, then $U^{ij}=\lambda H^{ij}$ 
everywhere.  In this case, if $\lambda=0$, then $M$ is a paraboloid,
while if $\lambda\neq 0$, then $M$ is a non-degenerate quadric.
\end{Theorem}
\begin{Proof}
Suppose first that $U^{ij}=\lambda H^{ij}$ on ${\mathbf F}_1(M)$ for
some function $\lambda$.  This is same as writing
$$
  \varphi^j_0 = \lambda H^{jl}\varphi^0_l = \lambda\omega^j.
$$
We differentiate this equation (substituting itself), and obtain
$$
   (d\lambda - \lambda\varphi^0_0)\wedge\omega^j = 0
   \mbox{ for each }j.
$$
Under the standing assumption $n>1$, this means that
$$
  d\lambda = \lambda\varphi^0_0.
$$
So assuming that $M$ is connected, we have the first statement of (1).
We will describe the geometric consequences of each of the two 
possibilities.

First, suppose that $\lambda = 0$, so that $U^{ij}=0$, and then 
$\varphi^j_0=0$ throughout ${\mathbf F}_1(M)$.  Then the definition of our original 
basis of $1$-forms gives
$$
  de_0 =  e_a\varphi^a_0 = e_0\varphi^0_0,
$$
meaning that the direction in $\A^{n+1}$ of $e_0$ is fixed throughout 
${\mathbf F}_1(M)$, or equivalently, all of the affine normals of $M$ are 
parallel.

Next, suppose $\lambda\neq 0$, and assume for simplicity that
$\lambda<0$.  The differential equation $d\lambda =
\lambda\varphi^0_0$ implies that we can restrict to the principal
subbundle ${\mathbf F}_2(M)$ where $\lambda = -1$.  This amounts to a
choice of a particular vector field $e_0$ along the affine normal line
field already defined.  Note that on ${\mathbf F}_2(M)$, we have
\begin{equation}
  \varphi^j_0 = -\omega^j = -H^{jk}\varphi^0_k,\qquad \varphi^0_0=0.
\label{VeryRestricted}
\end{equation}
As a result, the structure equations $dx = e_i\omega^i$ and 
$de_0=e_a\varphi^a_0 = -e_i\omega^i$ imply that
$$
  d(x+e_0) = 0,
$$
so that $x+e_0$ is a constant element of $\A^{n+1}$. In particular, all 
of the affine normal lines of $M$ pass through this point.  This
completes the proof of (1).
\index{affine!normal line|)}

Now assume that $T_{ijk}=0$ identically; this will be satisfied by
each member of the fields of affine hypersurfaces associated to
certain neo-classical Poincar\'e-Cartan forms of interest.  Our first
claim is that $U^{ij}=\lambda H^{ij}$ for some function $\lambda$ on
${\mathbf F}_1(M)$.  To see this, note that our hypothesis means
$$
  dH_{ij} = -H_{ij}\varphi^0_0 + H_{kj}\varphi^k_i + H_{ik}\varphi^k_j.
$$
We differentiate this, using the structure equations in the simplified form 
that defined the reduction to ${\mathbf F}_1(M)$, and obtain
$$
  0 = -H_{kj}\varphi^k_0\wedge\varphi^0_i - 
H_{ik}\varphi^k_0\wedge\varphi^0_j.
$$
If we use $H_{ij}$ to raise and lower indices and define
$$
  U_{ij} = H_{ik}H_{lj}U^{kl},
$$
then the preceding equation may be written as
$$
  0 = -(U_{jl}H_{ik}+U_{il}H_{jk})\omega^l\wedge\omega^k.
$$
The coefficients of this vanishing $2$-form then satisfy
$$
  0 = U_{jl}H_{ik}+U_{il}H_{jk}-U_{jk}H_{il}-U_{ik}H_{jl};
$$
we multiply by $H^{ik}$ (and sum over $i,k$) to conclude
$$
  U_{jl}=\frac{1}{n}(H^{ik}U_{ik})H_{jl}.
$$
This proves that
$$
  U^{ij} = \lambda H^{ij},
$$
with $\lambda = \frac{1}{n}H^{kl}U_{kl}$.

We now return to the possibilities $\lambda=0$, $\lambda\neq 0$ under the 
stronger hypothesis $T_{ijk}=0$.

In the first case, note that with the condition $\varphi^j_0=0$ on 
${\mathbf F}_1(M)$, we have that the Pfaffian system generated by $\varphi^0_0$ and 
$\varphi^i_j$ (for $1\leq i,j\leq n$) is integrable.  Let $\tilde M$ be any 
leaf of this system.  Restricted to $\tilde M$, we have
$$
  dH_{ij}=0,
$$
so that the functions $H_{ij}$ are constants.  Furthermore, the linearly 
independent $1$-forms $\omega^i$ on $\tilde M$ are each closed, so that (at 
least locally, or else on a simply connected cover) there are coordinates 
$u^i$ on $\tilde M$ with
$$
  \omega^i = du^i.
$$
Substituting all of this into the structure equations, we have:
\begin{itemize}
\item $de_0 = 0$, so that $e_0$ is a constant element of $\A^{n+1}$ on 
$\tilde M$;
\item
$de_i = e_0\varphi^0_i = e_0H_{ij}\omega^j = d(e_0H_{ij}u^j),$
so that
$$
  e_i = \bar e_i + e_0H_{ij}u^j
$$
for some constant $\bar e_i\in\A^{n+1}$;
\item
$dx = e_i\omega^i = (\bar e_i+e_0H_{ij}u^j)du^i =
  d(u^i\bar e_i+\frac12 e_0H_{ij}u^iu^j)$, so that
$$
  x = \bar x + u^i\bar e_i + \frac12 H_{ij}u^iu^je_0
$$
for some constant $\bar x\in\A^{n+1}$.
\end{itemize}
The conclusion is that as the coordinates $u^i$ vary on $\tilde M$, the $\A^{n+1}$-valued 
function $x$ on $\tilde M$ traces out a paraboloid, with vertex at $\bar x$ 
and axis along the direction of $e_0$.

Turning to the case $\lambda\neq 0$, recall that under the assumption 
$\lambda < 0$, we can reduce to a subbundle ${\mathbf F}_2(M)\subset{\mathbf F}_1(M)$ on 
which $\lambda = -1$.  We use the differential equation
$$
  dH_{ij} = H_{ik}\varphi^k_j + H_{kj}\varphi^k_i
$$
to reduce again to a subbundle ${\mathbf F}_3(M)\subset{\mathbf F}_2(M)$ on which $H_{ij}=\bar 
H_{ij}$ is some constant matrix.  On ${\mathbf F}_3(M)$, the forms
$\varphi^i_j$ satisfy linear algebraic relations
$$
  0 = \bar H_{ik}\varphi^k_j+\bar H_{kj}\varphi^k_i.
$$
Our assumption $\lambda=-1$ allows us to
combine these with the relations (\ref{VeryRestricted}) by defining
$$
  \Phi = \left(\begin{array}{cc} 0 & \varphi^0_j \\ \varphi^i_0
      & \varphi^i_j \end{array}\right),\qquad
  {\mathbf H} = \left(\begin{array}{cc} 1 & 0 \\ 0 & \bar 
H\end{array}\right),
$$
and then
$$
  {\mathbf H}\Phi +\, ^t\Phi{\mathbf H} = 0.
$$
In other words, the matrix-valued $1$-form $\Phi$ on ${\mathbf F}_3(M)$ takes values 
in the Lie algebra of the stabilizer of the bilinear form ${\mathbf H}$.  
For instance, if our hypersurface $M$ is convex, so that $(H_{ij})$ is 
definite everywhere, then we could have chosen $\bar 
H_{ij}=\delta_{ij}$, and then $\Phi$ would take values in the Lie algebra 
$\lie{so}(n+1,\R)$.  Whatever the signature of $H_{ij}$, let the stabilizer of 
${\mathbf H}$ be denoted by $O({\mathbf H})\subset GL(n+1,\R)$, with Lie algebra 
$\lie{so}({\mathbf H})$.
Then the structure equation
$$
  d\Phi+\Phi\wedge\Phi = 0
$$
implies that there is locally (alternatively, on a simply connected
cover) a map
$$
  g:{\mathbf F}_3(M)\to O({\mathbf H})
$$
such that
$$
  \Phi = g^{-1}dg.
$$
Using the structure equations $de_a = e_b\varphi^b_a$, this implies
$$
  d(e_a\cdot(g^{-1})^a_b) = 0,
$$
so that
$$
  e_a = \bar e_bg^b_a
$$
for some fixed affine frame $(\bar e_b)$.
In particular, the $\A^{n+1}$-valued function $e_0$ on ${\mathbf
  F}_3(M)$ takes as its
values precisely the points of a level surface of a non-degenerate
quadratic form,  
defined by ${\mathbf H}$.  Recalling from the first part of the proof
that $x+e_0$ is constant on $\A^{n+1}$,  
this means that the hypersurface $M$, thought of as the image of the map 
$x:{\mathbf F}_3(M)\to\A^{n+1}$, is a constant translate of a
non-degenerate quadric  
hypersurface.  The signature of the quadric is $(p,q)$, where $(p-1,q)$ is 
the signature of the first fundamental form $(H_{ij})$.

The case $\lambda>0$ instead of $\lambda <0$ is quite similar, but $M$ is a 
quadric of signature $(p,q)$ when $(H_{ij})$ has signature $(p,q-1)$.
\end{Proof}

\index{affine!hypersurface|)}
\index{affine!fundamental forms|)}
\index{affine!transformation|)}
\index{affine!frame bundle|)}

\section{The Equivalence Problem for $n\geq 3$}
\markright{2.4.  THE EQUIVALENCE PROBLEM FOR $n\geq 3$}
\label{Section:Bigequiv}

\index{equivalence method|(}
\index{equivalence!of Poincar\'e-Cartan forms|(}
We now consider a contact manifold $(M,I)$ with a closed,
almost-classical\index{Poincar\'e-Cartan form!almost-classical|(}
form
\begin{equation}
  \Pi = -\theta\wedge(H^{ij}\pi_i\wedge\omega_{(j)}-K\omega).
\label{AlmostClassicalAgain}
\end{equation}
We will shortly specialize to the case in which $\Pi$ is
neo-classical\index{Poincar\'e-Cartan form!neo-classical|(}.
The coframes in which $\Pi$ takes the form
(\ref{AlmostClassicalAgain}), for some
functions $H^{ij}$ and $K$, constitute a $G$-structure as described in
Lemma~\ref{AlmostClassicalLemma}.
The purpose of this section is to describe a canonical reduction of this
$G$-structure to one carrying a pseudo-connection satisfying structure
equations of a prescribed form, as summarized in
(\ref{SummaryGp8}--\ref{MasterDRho}), at
least in case the matrix $(H^{ij})$ is either positive- or
negative-definite everywhere.  This application of the equivalence
method involves no techniques beyond those introduced in
\S\ref{Section:SmallEquiv}, but some of the linear-algebraic
computations are more involved.

We begin by refining our initial $G$-structure as follows.
\begin{Lemma}
Let $(M,I)$ be a contact manifold with almost-classical form $\Pi$.

\noindent
(1)
There exist local coframings $(\theta,\omega^i,\pi_i)$ on $M$ such
that $\Pi$ has the form (\ref{AlmostClassicalAgain}) and such that
$$
  d\theta\equiv -\pi_i\wedge\omega^i\pmod{\{I\}}.
$$

\noindent
(2)
Local coframings as in (1) are the sections of a $G_0$-structure
$B_0\to M$, where $G_0$ is the group of matrices of the form (in blocks of
size $1,n,n$)
\begin{equation}
  g_0=\left(\begin{array}{ccc}
    a & 0 & 0 \\ C^i & A^i_j & 0 \\ D_i & S_{ik}A^k_j &
      a(A^{-1})^j_i
  \end{array}\right),\quad A\in GL(n,\R),\ S_{ij}=S_{ji}.
\label{DefG0Big}
\end{equation}

\noindent
(3)
If two local coframings as in (1) are related as
$$
  \left(\begin{array}{c} \theta \\ \omega^i \\ \pi_i
    \end{array}\right) = g_0^{-1}\cdot
  \left(\begin{array}{c}
    \bar\theta \\ \bar\omega^j \\ \bar\pi_j
  \end{array}\right),
$$
and if $\Pi = -\theta\wedge(H^{ij}\pi_i\wedge\omega_{(j)}-K\omega) =
  -\bar\theta\wedge(\bar H^{ij}\bar\pi_i\wedge\bar\omega_{(j)}-\bar
  K\bar\omega)$ are the expressions for $\Pi$ with respect to
  these coframings, 
then
\begin{eqnarray}
    H & = & a^2(\mbox{det }A)A^{-1}\bar H\,^t\!A^{-1},
\label{TransH8} \\
    K & = & a(\mbox{det }A)(\bar K - Tr(\bar HS)).
\end{eqnarray}
\end{Lemma}
\begin{Proof}
(1) 
First observe that in any coframing, we may write
$$
  d\theta \equiv a^{ij}\pi_i\wedge\pi_j + b^j_i\pi_j\wedge\omega^i +
     c_{ij}\omega^i\wedge\omega^j \pmod{\{I\}}.
$$
We will deal with each of the three coefficient matrices $(a^{ij})$, 
$(b^j_i)$,
$(c_{ij})$ to obtain the desired condition $d\theta\equiv -
\sum\pi_i\wedge\omega^i$.
\begin{itemize}
  \item  The proof of Proposition~\ref{Prop:Neoclassical6} showed for
    $n\geq 3$ that
    $$
       0 \equiv d\theta \equiv a^{ij}\pi_i\wedge\pi_j
            \pmod{\{J_\Pi\}},
    $$
    which implies $a^{ij}\pi_i\wedge\pi_j=0$.  This followed from
    calculating $0=d\Pi$ modulo $\{I\}$.
  \item  From the fact that $\theta$ is a contact form, we have
    $$
      0\neq\theta\wedge(d\theta)^n =
         \pm\mbox{det}(b^j_i)\theta\wedge\omega\wedge\pi,
    $$
    so that $(b^j_i)$ is an invertible matrix.  Therefore, we may apply
    the matrix $-(b^j_i)^{-1}$ to the $1$-forms $\pi_j$ to obtain a new
    basis in which we have $b^j_i=-\delta^j_i$, so that
    $$
      d\theta\equiv -\pi_i\wedge\omega^i + c_{ij}\omega^i\wedge\omega^j
          \pmod{\{I\}}.
    $$
    Note that this coframe change is of the type admitted by
    Lemma~\ref{AlmostClassicalLemma}, preserving the form
    (\ref{AlmostClassicalAgain}).
  \item  Finally, we can replace $\pi_i$ by $\pi_i+c_{ij}\omega^j$ to have
    the desired $d\theta\equiv-\pi_i\wedge\omega^i$.  This coframe
    change also preserves the form (\ref{AlmostClassicalAgain}).
\end{itemize}

\noindent
(2)
We already know that any matrix as in Lemma~\ref{AlmostClassicalLemma}
will preserve the form (\ref{AlmostClassicalAgain}).  We write the
action of such a matrix as
$$
  \left\{\begin{array}{ccc}
     \bar\theta & = & a\theta \\
     \bar\omega^i & = & C^i\theta + A^i_j\omega^j \\
     \bar\pi_i & = & D_i\theta + S_{ik}A^k_j\omega^j+B^j_i\pi_j.
  \end{array}\right.
$$
It is easily verified that the condition $d\theta\equiv-\pi_i\wedge\omega^i$ 
implies the analogous condition 
$d\bar\theta\equiv-\bar\pi_i\wedge\bar\omega^i$ if and 
only if
$$
  \left\{\begin{array}{l}
    B^j_iA^i_k=a\delta^j_k, \\
    S_{jk} = S_{kj}.
  \end{array}\right.
$$
This is what we wanted to prove.

\noindent
(3)
These formulae are seen by
substituting the formulae for $(\bar\theta,\bar\omega^i,\bar\pi_i)$
into the equation for the two expressions for $\Pi$, and comparing
terms.  One uses the following fact from linear algebra: if
$$
  \bar\omega^i \equiv A^i_j\omega^j \pmod{\{I\}},
$$
then
$$
  \bar\omega_{(j)} \equiv (\mbox{det }A)(A^{-1})^i_j\omega_{(i)}
  \pmod{\{I\}}; 
$$
that is, the coefficients of $\bar\omega_{(j)}$ in terms of
$\omega_{(i)}$ are the cofactors of the coefficient matrix of
$\bar\omega^i$ in terms of $\omega^j$.
\end{Proof}

\

We can see from (\ref{TransH8}) that the matrix $H=(H^{ij})$ transforms under
coframe changes like a bilinear form, up to scaling, and in
particular that its absolute signature is fixed at each point of
$M$.  To proceed, we have to assume that this signature is constant
throughout $M$.  In particular, we shall from now on assume that $H$
is positive or negative definite everywhere, and refer to
almost-classical forms $\Pi$ with this property as {\em
  definite}.\index{Poincar\'e-Cartan form!neo-classical!definite|(}
Cases of different constant
signature are of interest, but can be easily reconstructed by the
reader in analogy with the definite case examined below.
\index{Poincar\'e-Cartan form!neo-classical|)}

Once we assume that the matrix-valued function $H$ on $B_0$ is
definite, the following is an easy consequence of the
preceding lemma.
\begin{Lemma}
Given a definite, almost-classical Poincar\'e-Cartan form $\Pi$ on a
contact manifold $(M,I)$, there are $0$-adapted local coframings
$(\theta,\omega^i,\pi_i)$ for which
$$
  \Pi = -\theta\wedge(\delta^{ij}\pi_i\wedge\omega_{(j)}),
$$
and these form a $G_1$-structure $B_1\subset B_0\to M$, where $G_1$ is
the group of matrices $g_1$ of the form (\ref{DefG0Big}) with
$$
  \mbox{det }A>0,\quad  a(\mbox{det }A)^{\frac12}A\in O(n,\R),\quad
    S_{ii}=0.
$$
\label{OneLemma8}
\end{Lemma}
This follows from imposing the conditions $\bar H=H=I_n$, $\bar K=K=0$
in the previous lemma.  Unfortunately, it is difficult to give a
general expression in coordinates for such a $1$-adapted coframing in
the classical case, because such an expression requires that we
normalize the Hessian matrix $(L_{p_ip_j})$.  In practice, however,
such a coframing is usually easy to compute.
\index{Poincar\'e-Cartan form!almost-classical|)}

It is convenient for later purposes to use a different
parameterization of our group $G_1$.  Namely, an arbitrary element
will be written as
\begin{equation}
  g_1 = \left(\begin{array}{ccc} \pm r^{n-2} & 0 & 0 \\
     C^i & r^{-2}A^i_j & 0 \\ D_i & S_{ik}A^k_j & \pm r^n(A^{-1})^j_i
  \end{array}\right),
\label{DefG1Big}
\end{equation}
where $A=(A^i_j)\in SO(n,\R)$, $r>0$, $S_{ij}=S_{ji}$, $S_{ii}=0$.
Also, now that the orthogonal group has appeared,
some of the representations occuring in the sequel are isomorphic to
their duals, for which it may be unuseful and sometimes
confusing to maintain the usual summation 
convention\index{summation convention}, in which one
only contracts a pair of indices in which one index is raised and the
other lowered.  Therefore, we will now sum any index
occuring twice in a single term, regardless of its positions.

We now assume that we have a definite, neo-classical Poincar\'e-Cartan
form $\Pi$ with associated $G_1$-structure $B_1\to M$, and we begin
searching for differential invariants.
There are local pseudo-connection $1$-forms 
$\rho,\gamma^i,\delta_i,\alpha^i_j,\sigma_{ij}$ defined so that
equations of the following form hold:
$$
  d\left(\begin{array}{c} \theta \\ \omega^i \\ \pi_i \end{array}\right)
  = -\left(\begin{array}{ccc}
      (n-2)\rho & 0 & 0 \\
      \gamma^i & -2\rho\delta^i_j+\alpha^i_j & 0 \\
      \delta_i & \sigma_{ij} & n\rho\delta^j_i - \alpha^j_i
    \end{array}\right) \wedge
  \left(\begin{array}{c} \theta \\ \omega^j \\ \pi_j \end{array}\right)
  + \left(\begin{array}{c} \Theta \\ \Omega^i \\ \Pi_i\end{array}\right),
$$
where $\theta,\omega^i,\pi_i$ are the tautological $1$-forms on $B_1$,
the torsion $2$-forms $\Theta,\Omega^i,\Pi_i$ are semibasic for 
$B_1\to M$, and the psuedo-connection $1$-forms satisfy
$$
  \alpha^i_j+\alpha^j_i = 0,\quad \sigma_{ij}=\sigma_{ji},\quad
    \sigma_{ii}=0.
$$
These last conditions mean that the psuedo-connection matrix takes
values in the Lie algebra $\lie{g}_1\subset\lie{gl}(2n+1,\R)$ of $G_1$.

The psuedo-connection $1$-forms are not uniquely determined, and our 
next step is to exploit this indeterminacy to try to absorb components
of the torsion.

First, we know that $d\theta\equiv-\pi_i\wedge\omega^i\mbox{ (mod 
$\{I\}$)}$.
The difference between $\Theta = d\theta+(n-2)\rho\wedge\theta$ and
$-\pi_i\wedge\omega^i$ is therefore a semibasic multiple of $\theta$,
which can be absorbed by a semibasic change in $\rho$.  We can
therefore simply assume that
$$
  d\theta = -(n-2)\rho\wedge\theta - \pi_i\wedge\omega^i,
$$
or equivalently, $\Theta = -\pi_i\wedge\omega^i$.

Second, our assumption that $\Pi$ is neo-classical means that the Pfaffian 
system $J_\Pi = \{\theta,\omega^i\}$ is integrable (even up on $B_1$).  In 
the structure equation
\begin{equation}
  d\omega^i = -\gamma^i\wedge\theta - 
    (-2\rho\delta^i_j+\alpha^i_j)\wedge\omega^j + \Omega^i,
\label{SecondStructure}
\end{equation}
this means that $\Omega^i\equiv 0\mbox{ (mod $\{J_\Pi\}$)}$.
Also, $\Omega^i$ is semibasic over $M$, so we can write
\begin{equation}
  \Omega^i \equiv T^{ijk}\pi_j\wedge\omega^k + \frac12P^i_{jk}
    \omega^j\wedge\omega^k\pmod{\{I\}}.
\label{SecondTorsion}
\end{equation}
Now, adding semibasic $1$-forms to $\gamma^i$ allows us to preserve the 
equation (\ref{SecondStructure}) while also making (\ref{SecondTorsion})
an equality, and not merely a congruence.  A little linear algebra shows 
that there is a {\em unique} linear combination of the $\omega^i$ that can 
be added to $\alpha^i_j$, preserving $\alpha^i_j+\alpha^j_i=0$, to 
absorb the term $\frac12P^i_{jk}\omega^j\wedge\omega^k$.  This leaves us 
only with
$$
  \Omega^i = T^{ijk}\pi_j\wedge\omega^k.
$$
As in the elimination of the $P^i_{jk}$, we can add a combination of the 
$\pi_i$ to $\alpha^i_j$ to arrange
$$
  T^{ijk}=T^{kji}.
$$

To investigate the third torsion term $\Pi_i$, we use an alternate
derivation of the equation for $d\pi_i$. 
Namely, we differentiate the equation
$$
  d\theta = -(n-2)\rho\wedge\theta-\pi_i\wedge\omega^i,
$$
and take the result only modulo $\{I\}$ to avoid the unknown quantity
$d\rho$.  This eventually yields
$$
  0\equiv -(\Pi_k - T^{ijk}\pi_i\wedge\pi_j)\wedge\omega^k
    \pmod{\{I\}}.
$$
As before, multiples of $\theta$ may be absorbed by redefining $\delta_i$, 
so that we can assume this congruence is an equality.  Reasoning
similar to that which proves the Cartan lemma gives
$$
  \Pi_k - T^{ijk}\pi_i\wedge\pi_j = \nu_{kl}\wedge\omega^l
$$
for some semibasic $1$-forms $\nu_{kl}=\nu_{lk}$.  Now, most of these forms 
$\nu_{kl}$ can be subtracted from the psuedo-connection forms $\sigma_{kl}$, 
simplifying the torsion; but the condition $\sigma_{ii}=0$ prevents us from 
completely absorbing them.  Instead, the trace remains, and we have
$$
  \Pi_k = \delta_{kl}\nu\wedge\omega^l + T^{ijk}\pi_i\wedge\pi_j.
$$
We can learn more about $\nu$ using the integrability condition $d\Pi=0$, 
taken modulo terms quadratic in the $\pi_i$:
$$
  0 = d\Pi \equiv n\theta\wedge\nu\wedge\omega.
$$
A consequence is that $\nu\equiv 0 \mbox{ (mod $\{\theta,\omega^i\}$)}$; in 
other words, $\nu$ has no $\pi_i$-terms, and may be written (using
again a change in $\delta_i$) as
$$
  \nu = \sum N_i\omega^i.
$$
Then replacing $\sigma_{ij}$ by
$$
  \sigma_{ij} + \sf{n}{n+2}
    (\delta_{ik}N_j+\delta_{jk}N_i-\sf2n\delta_{ij}N_k)\omega^k
$$
yields new psuedo-connection forms, for which the third torsion term is simply
$$
  \Pi_k = T^{ijk}\pi_i\wedge\pi_j.
$$
This completes the major step of absorbing torsion by altering the
pseudo-connection.

Before proceeding to the next major step, we look for linear-algebraic
conditions on the torsion which may simplify later calculations.  In
particular, 
we made only very coarse use of $d\Pi=0$ above.  Now we compute more 
carefully
$$
  0 = d\Pi = -\theta\wedge(2T^{ijk}+\delta^{ij}T^{lkl})\pi_k\wedge\pi_j
    \wedge\omega_{(i)},
$$
so we must have
\begin{equation}
  2T^{ijk}+\delta^{ij}T^{lkl} = 2T^{ikj}+\delta^{ik}T^{ljl}.
\label{ThirdTorsion}
\end{equation}

The next major step is a reduction of our $G_1$-structure.  We will 
examine the variation of the functions $T^j \stackrel{\mathit{def}}{=}
T^{iji}$ along fibers  
of $B_1\to M$, and observe that the zero-locus $\{T^j=0\}$ defines a 
$G_2$-structure for a certain codimension-$n$ subgroup $G_2\subset G_1$.  

As usual, the variation of $T^j$ will be described
infinitesimally.  To study $dT^j$ without knowledge of the traceless part 
of $dT^{ijk}$, we exploit the exterior algebra, writing
\begin{equation}
  d(\theta\wedge\omega^1\wedge\cdots\wedge\omega^n) =
    ((n+2)\rho + T^k\pi_k)\wedge\theta\wedge\omega.
\label{CleverTrace}
\end{equation}
We will differentiate this for information about $dT^k$, but in doing
so we will need information about $d\rho$ as well.  Fortunately, this
is available by differentiating the first structure equation
$$
  d\theta = -(n-2)\rho\wedge\theta - \pi_k\wedge\omega^k,
$$
yielding
$$
  (n-2)d\rho \equiv \gamma^k\wedge\pi_k \pmod{\{\theta,\omega^i\}}.
$$
Now we return to differentiating (\ref{CleverTrace}) and eventually find
$$
  dT^k \equiv -\sf{n+2}{n-2}\gamma^k +
    (n\delta^k_j\rho - \alpha^k_j)T^j
    \pmod{\{\theta,\omega^i,\pi_i\}}.
$$
This means that along fibers of $B_1\to M$, the vector-valued function 
$T(u)=(T^j(u))$, $u\in B_1$, is orthogonally rotated (infinitesimally, by 
$\alpha^k_j$), scaled (by $\rho$), and translated (by $\gamma^i$).  In fact,
for $g_1\in G_1$ as in (\ref{DefG1Big}),
$$
  T(u\cdot g_1) = \pm r^2A^{-1}(r^{n-2}T(u) - \sf{n+2}{n-2}C).
$$
Now the set
$$
  B_2 \stackrel{\mathit{def}}{=} \{u\in B_1: T(u) = 0\}\subset B_1
$$
is a $G_2$-subbundle of $B_1\to M$, where $G_2$ consists of matrices as in 
(\ref{DefG1Big}) with $T^i = 0$.

On the submanifold $B_2\subset B_1$, we have from (\ref{ThirdTorsion})
the symmetry
$$
  T^{ijk}=T^{kji} = T^{ikj}.
$$
As a consequence, the torsion $\Pi_k$ restricts to
$$
  \Pi_k = T^{ijk}\pi_j\wedge\pi_k = 0.
$$

The previous structure equations continue to hold, but the forms 
$\gamma^i|_{B_2}$ should not be regarded as part of the psuedo-connection, as 
they are now semibasic over $M$.  We therefore write
\begin{equation}
  d\left(\begin{array}{c} \theta \\ \omega^i \\ \pi_i \end{array}\right)
  = -\left(\begin{array}{ccc}
      (n-2)\rho & 0 & 0 \\
      0 & -2\rho\delta^i_j+\alpha^i_j & 0 \\
      \delta_i & \sigma_{ij} & n\rho\delta^j_i - \alpha^j_i
    \end{array}\right) \wedge
  \left(\begin{array}{c} \theta \\ \omega^j \\ \pi_j \end{array}\right)
  + \left(\begin{array}{c} \Theta \\ \Omega^i \\ \Pi_i\end{array}\right),
\label{Str8.1}
\end{equation}
where still 
\begin{equation}
  \alpha^i_j+\alpha^j_i=0,\ \sigma_{ij}=\sigma_{ji},\ \sigma_{ii}=0,
\label{Str8.5}
\end{equation}
and now
\begin{equation}
  \left\{\begin{array}{l}
    \Theta = -\pi_i\wedge\omega^i, \\
    \Omega^i = -(S^i_j\omega^j+U^{ij}\pi_j)\wedge\theta +
      T^{ijk}\pi_j\wedge\omega^k, \\
    \Pi_i = 0.  \end{array}\right.
\label{Str8.2}
\end{equation}
Here we have denoted $\gamma^i\equiv S^i_j\omega^j+U^{ij}\pi_j\mbox{ (mod 
$\{I\}$)}$.  Also, we still have
\begin{equation}
  T^{ijk}=T^{kji} = T^{ikj},\quad T^{iik}=0.
\label{Str8.3}
\end{equation}
Notice that we can alter $\alpha^i_j$ and $\rho$ to assume that
\begin{equation}
  S^i_j = S^j_i,\qquad S^i_i = 0,
\label{Str8.4}
\end{equation}
where we also have to add combinations of $\omega^i$ to $\delta_i$ 
to preserve $\Pi_i = 0$.  In fact, these assumptions uniquely
determine $\alpha^i_j$ and 
$\rho$, although $\delta_i$ and $\sigma_{ij}$ still admit some 
ambiguity.

Equations (\ref{Str8.1}--\ref{Str8.4}) summarize the results of the
equivalence method carried out to this point.
We have uncovered the primary differential invariants of a definite 
neo-classical Poincar\'e-Cartan form: they are the functions $T^{ijk}$, 
$S^i_j$ and $U^{ij}$.  Their properties are central in what follows.

For example, note that the rank-$n$ Pfaffian system $\{\omega^i\}$ on $B_2$ 
is invariant under the action of the structure group $G_2$, and therefore it 
is the pullback of a Pfaffian system (also to be denoted $\{\omega^i\}$) 
down on $M$.  Testing its integrability, we find
\begin{equation}
  d\omega^i \equiv -U^{ij}\pi_j\wedge\theta \pmod{\{\omega^i\}}.
\label{SomeEqn}
\end{equation}
We will see shortly that the matrix-valued function $(U^{ij})$ varies along 
the fibers of $B_2\to M$ by a linear representation of $G_2$, so that it is 
plausible to ask about those Poincar\'e-Cartan forms for which $U^{ij}=0$; 
(\ref{SomeEqn}) shows that this is equivalent to the integrability of 
$\{\omega^i\}$.  In this case, in addition to the local fibration $M\to Q$ 
whose fibers are leaves of $J_\Pi$, we have
$$
  Q^{n+1}\to N^n,
$$
where $N$ is the locally-defined $n$-dimensional ``leaf space'' for
$\{\omega^i\}$.  Coordinates on $N$---equivalently, functions on $M$
whose differentials lie in $\{\omega^i\}$---may be thought of as
``preferred independent variables'' for the contact-equivalence class
of our Euler-Lagrange
equation\index{Euler-Lagrange!equation}, canonical in the sense that 
every symmetry of $M$ preserving the Poincar\'e-Cartan form preserves the 
fibration $M\to N$ and therefore acts on $N$.  Note that even if an 
$(M,\Pi)$ satisfying $U^{ij}=0$ came to us from a classical Lagrangian
with independent variables $(x^i)$, we need not have 
$\{\omega^i\} = \{dx^i\}$.  

This is not to say that the case $U^{ij}\neq 0$ is
uninteresting.  In the next section, we will see an important family
of examples from Riemannian geometry with $U^{ij}=\lambda\delta^i_j.$  To 
obtain preliminary information about $U^{ij}$ in a manner that will not
require much knowledge of $S^i_j$ or $T^{ijk}$, we start with the equation
\begin{equation}
  d(\omega^1\wedge\cdots\wedge\omega^n) = 2n\rho\wedge(\omega^1\wedge
    \cdots\wedge\omega^n) + U^{ij}\theta\wedge\pi_j\wedge\omega_{(i)}.
\label{dOmega}
\end{equation}
We will differentiate again, but we need more refined information about 
$d\rho$; this is obtained from
$$
  0 = d^2\theta = -((n-2)d\rho+\delta_i\wedge\omega^i+\pi_i\wedge\gamma^i)
    \wedge\theta.
$$
Keep in mind that $\gamma^i=S^i_j\omega^j+U^{ij}\pi_j\mbox{ (mod $\{I\}$)}$
on this reduced bundle.
We can now write
\begin{equation}
  (n-2)d\rho + \delta_i\wedge\omega^i + \pi_i\wedge\gamma^i = \tau\wedge
     \theta
\label{CrudeDRho}
\end{equation}
for some unknown $1$-form $\tau$.  Returning to the derivative of
(\ref{dOmega}), we find
$$
  0\equiv \sf{2n}{n-2}(-\pi_i\wedge U^{ij}\pi_j)\wedge\omega
    +U^{ij}\pi_i\wedge\pi_j\wedge\omega \pmod{\{I\}}.
$$
This implies that $U^{ij}\pi_i\wedge\pi_j\wedge\omega = 0$, so that we
have
$$
  U^{ij}=U^{ji}.
$$

We will need an even more refined version of the equation
(\ref{CrudeDRho}) for $d\rho$.  In the preceding paragraph, we
substituted that equation into the equation for $0\equiv
d^2(\omega^1\wedge\cdots\wedge\omega^n)\mbox{ (mod $\{I\}$)}$.  Now, we
substitute it instead into
\begin{eqnarray*}
  0 & \equiv & d^2(\omega^1\wedge\cdots\wedge\omega^n) 
    \pmod{\{\pi_1,\ldots,\pi_n\}} \\
    & \equiv & \left(\sf{2n}{n-2}\tau-U^{ij}\sigma_{ij}\right)\wedge\theta
        \wedge\omega^1\wedge\cdots\wedge\omega^n.
\end{eqnarray*}
This means that $\frac{2n}{n-2}\tau-U^{ij}\sigma_{ij}$ lies in
$\{\theta,\omega^i,\pi_i\}$.  
Recall that also $\gamma^i = S^i_j\omega^j+U^{ij}\pi_j+V^i\theta$ for
some functions $V^i$, and we can put this back into (\ref{CrudeDRho})
to finally obtain
$$
  (n-2)d\rho = -\delta_i\wedge\omega^i-S^i_j\pi_i\wedge\omega^j
    +\left(\sf{n-2}{2n}\right)U^{ij}\sigma_{ij}\wedge\theta+
    (s_i\omega^i-t^i\pi_i)\wedge\theta,
$$
for some functions $s_i$, $t^i$.
Furthermore, we can replace each $\delta_i$ by $\delta_i-s_i\theta$,
preserving previous equations, to assume that $s_i=0$.  This gives
$$
  (n-2)d\rho = -\delta_i\wedge\omega^i-S^i_j\pi_i\wedge\omega^j
    +\left(\sf{n-2}{2n}\right)U^{ij}\sigma_{ij}\wedge\theta
    -t^i\pi_i\wedge\theta,
$$
which will be used in later sections.

The last formulae that we will need are those for the
transformation rules for $T^{ijk}$, $U^{ij}$, $S^i_j$ along fibers of
$B\to M$.  These are obtained by computations quite similar to those
carried out above, and we only state the results here, which are:
\begin{itemize}
\item
  $dT^{ijk}\equiv n\rho T^{ijk}-\alpha^i_lT^{ljk}-\alpha^j_lT^{ilk}
    -\alpha^k_lT^{ijl},$
\item
  $dU^{ij}\equiv 2n\rho U^{ij}-\alpha^i_lU^{lj}-\alpha^j_lU^{il},$
\item
  $dS^i_j\equiv (n-2)\rho S^i_j-\alpha^i_lS^l_j+S^i_l\alpha^l_j
    +\frac12(U^{il}\sigma_{lj}+U^{jl}\sigma_{li})
    -\frac1n\delta^i_jU^{kl}\sigma_{kl}+T^{ijk}\delta_k,$
\end{itemize}
all modulo $\{\theta,\omega^i,\pi_i\}$.
Notice in particular that $T^{ijk}$ and $U^{ij}$ transform by a
combination of rescaling and a standard representation of $SO(n)$.
However, $(S^i_j)$ is only a tensor when the tensors $(T^{ijk})$ and
$(U^{ij})$ both vanish.  We will consider this situation in the next
chapter.

An interpretation of the first two transformation rules is that the
objects
\begin{eqnarray*}
  {\mathbf T} & = & T^{ijk}(\pi_i\circ\pi_j\circ\pi_k)\otimes
    |\pi_1\wedge\cdots\wedge\pi_n|^{-\frac{2}{n}}, \\
  {\mathbf U} & = & U^{ij}\pi_i\circ\pi_j
\end{eqnarray*}
are invariant modulo $J_\Pi=\{\theta,\omega^i\}$ under flows along
fibers over $M$; that is, when
restricted to a fiber of $B_2\to Q$, they actually descend to
well-defined objects on the smaller fiber of $M\to Q$.  
The restriction to fibers suggests our next result, which nicely
relates the differential invariants of
the Poincar\'e-Cartan form with the affine\index{affine!hypersurface}
geometry of hypersurfaces discussed in the preceding section.  
\begin{Theorem}  The functions $T^{ijk}$ and $U^{ij}$ are coefficients of 
the affine cubic form and affine second fundamental
form\index{affine!fundamental forms} for the
fiberwise affine
hypersurfaces in $\bw{n}(T^*Q)$ induced by a semibasic Lagrangian
potential\index{Lagrangian potential} $\Lambda$ of $\Pi$.
\end{Theorem}
Proving this is a matter of identifying the bundles where the two
sets of invariants are defined, and unwinding the definitions.

In the next section, we will briefly build on the preceding results in
the case where $T^{ijk}=0$ and $U^{ij}\neq 0$, showing that these
conditions roughly characterize those definite neo-classical
Poincar\'e-Cartan forms appearing in the problem of finding prescribed
mean curvature hypersurfaces, in Riemannian or Lorentzian manifolds.
In the next chapter, we will extensively consider the case
$T^{ijk}=0$, $U^{ij}=0$, which includes remarkable Poincar\'e-Cartan
forms arising in conformal geometry.
About the case for which $T^{ijk}\neq 0$, nothing is known. 

For reference, we summarize the results of the equivalence
method that will be used below.  Associated to a definite,
neo-classical Poincar\'e-Cartan form $\Pi$ on a contact manifold
$(M,I)$ is a $G$-structure $B\to M$, where
\begin{equation}
  G = \left\{\left(\begin{array}{ccc}
    \pm r^{n-2} & 0 & 0 \\ 0 & r^{-2}A^i_j & 0 \\
      D_i & S_{ik}A^k_j & \pm r^n(A^{-1})^j_i\end{array}
    \right):
  \begin{array}{l} (A^i_j)\in SO(n,\R),\ r>0, \\ S_{ij}=S_{ji},\
    S_{ii}=0
    \end{array}
  \right\}.
\label{SummaryGp8}
\end{equation}
$B\to M$ supports a pseudo-connection (not uniquely determined)
$$
  \varphi = 
     -\left(\begin{array}{ccc}
      (n-2)\rho & 0 & 0 \\
      0 & -2\rho\delta^i_j+\alpha^i_j & 0 \\
      \delta_i & \sigma_{ij} & n\rho\delta^j_i - \alpha^j_i
    \end{array}\right),
$$
with $\alpha^i_j+\alpha^j_i = 0$, $\sigma_{ij}=\sigma_{ji}$,
$\sigma_{ii}=0$, such that in the structure equation
$$
  d\left(\begin{array}{c} \theta \\ \omega^i \\ \pi_i
    \end{array}\right) = -\varphi\wedge
  \left(\begin{array}{c} \theta \\ \omega^j \\ \pi_j
    \end{array}\right) + \tau,
$$
the torsion is of the form
$$
  \tau = \left(\begin{array}{c}
    -\pi_i\wedge\omega^i \\ 
 -(S^i_j\omega^j+U^{ij}\pi_j)\wedge\theta +
      T^{ijk}\pi_j\wedge\omega^k \\
    0 \end{array}\right),
$$
with
$$
  T^{ijk}=T^{jik}=T^{kji},\ T^{iik}=0;\ U^{ij}=U^{ji};\
    S^i_j=S^j_i,\ S^i_i=0.
$$
In terms of any section of $B\to M$, the Poincar\'e-Cartan form is
$$
  \Pi = -\theta\wedge\pi_i\wedge\omega_{(i)}.
$$
One further structure equation is 
\begin{equation}
  (n-2)d\rho = -\delta_i\wedge\omega^i-S^i_j\pi_i\wedge\omega^j
    +\left(\sf{n-2}{2n}\right)U^{ij}\sigma_{ij}\wedge\theta
    -t^i\pi_i\wedge\theta.
\label{MasterDRho}
\end{equation}
\index{Poincar\'e-Cartan form!neo-classical!definite|)}
\index{equivalence method|)}

\section{The Prescribed Mean Curvature System}
\label{Section:PrescribedH}

\index{Poincar\'e-Cartan form!neo-classical|(}
In this section, we will give an application of the part of the
equivalence method completed so far.  We will show that a
definite, neo-classical Poincar\'e-Cartan form with $T^{ijk}=0$, and
satisfying an additional open condition specified below, is
locally equivalent to that which arises in the problem of finding in
a given Riemannian manifold\index{Riemannian!manifold} a hypersurface
whose mean curvature\index{mean curvature!prescribed}
coincides with a prescribed background function.  This conclusion is
presented as Theorem~\ref{Theorem:PrescribedH}.

To obtain this result, we continue applying the equivalence
method\index{equivalence method} where we left off in
the preceding section, and take up the case $T^{ijk}=0$.  From our
calculations in affine
hypersurface\index{affine!hypersurface} geometry, we know that this
implies that
$$
  U^{ij}=\lambda\delta^{ij},
$$
for some function $\lambda$ on the principal bundle $B\to M$;
alternatively, this can be shown
by computations continuing those of the preceding section.  We will
show that under the
hypothesis $\lambda < 0$, the Poincar\'e-Cartan form $\Pi$ is locally
equivalent to that occuring in a prescribed mean curvature system.

We have in general on $B$ that
$$
  dU^{ij}\equiv 2n\rho U^{ij}-\alpha^i_kU^{kj}-
   \alpha^j_kU^{ik} \pmod{\{\theta,\omega^i,\pi_i\}}.
$$
Then for our $U^{ij}=\lambda\delta^{ij}$, the function $\lambda$ scales
positively along fibers of $B\to M$, so under our assumption $\lambda
<0$ we may make a reduction to
$$
  B_1 = \{u\in B: \lambda(u) = -1\}\subset B;
$$
this defines a subbundle of $B$ of codimension $1$, on which $\rho$
is semibasic over $M$.\footnote{In this section, we will denote by $B_1$,
  $B_2$, etc., successive reductions of the $G$-structure
  $B\to M$ which was constructed in the preceding section.  These are
  {\em not} the same as the bundles of the same name used in
  constructing $B$, which are no longer needed.}
  In particular, on $B_1$ we may write
$$
  \rho = -\frac{H}{2n}\theta+E_i\omega^i + F^i\pi_i
$$
for some functions $H$, $E_i$, $F^i$.
The reason for the normalization of the $\theta$-coefficient will
appear shortly.  

We claim that $F^i = 0$.  To see this,
start from the equation (\ref{MasterDRho}) for $d\rho$,
which on $B_1$ reads
$$
  (n-2)d\rho = -\delta_i\wedge\omega^i - t^i\pi_i\wedge\theta
    -S^i_j\pi_i\wedge\omega^j.
$$
Then, as we have done so often, we compute $d^2\omega$, where $\omega =
\omega^1\wedge\cdots\wedge\omega^n$ and
$$
  d\omega^i = 2\rho\wedge\omega^i - \alpha^i_j\wedge\omega^j
    +\pi_i\wedge\theta-S^i_j\omega^j\wedge\theta.
$$
We find
$$
  d\omega = 2n\rho\wedge\omega - \theta\wedge\pi_i\wedge\omega_{(i)}
    = 2n\rho\wedge\omega + \Pi,
$$
and the next step is simplified by knowing $d\Pi = 0$:
\begin{eqnarray*}
  0 & = & d^2\omega \\
    & = & 2n\,d\rho\wedge\omega - 2n\rho\wedge d\omega \\
    & = & -\left(\sf{2n}{n-2}\right)(t^i\pi_i\wedge\theta)
               \wedge\omega\\
    & \quad & \qquad + 2n(E_j\omega^j+F^j\pi_j)\wedge
               \theta\wedge\pi_i\wedge\omega_{(i)} \\
    & = & 2n\theta\wedge\left(F^j\pi_i\wedge\pi_j\wedge\omega_{(i)}
               +\left(\frac{t^i}{n-2}+E^i\right)\pi_i\wedge\omega
               \right).
\end{eqnarray*}
This gives our claim $F^i=0$, as well as
\begin{equation*}
  (n-2)E^i = -t^i.
\end{equation*}

For our next reduction, we will show that we can define a
principal subbundle
$$
  B_2 = \{u\in B_1: E^i(u) = 0 \}\subset B_1,
$$
having structure group defined by the condition $D_i=0$, $r=1$ in
(\ref{SummaryGp8}).  This follows by computing modulo
$\bigwedge^2\{\theta,\omega^i,\pi_i\}$:
$$
  (n-2)d\rho\equiv - \delta_i\wedge\omega^i,
$$
and also
$$
  d\rho \equiv -\sf{1}{2n}dH\wedge\theta + dE_i\wedge\omega^i -
     E_j\alpha^j_i\wedge\omega^i.
$$
Comparing these, we obtain
$$
  \left(dE_i - E_j\alpha^j_i+\sf{1}{n-2}\delta_i\right)
    \wedge\omega^i - \sf{1}{2n}dH\wedge\theta \equiv 0.
$$
This implies that
$$
  dE_i - E_j\alpha^j_i + \sf{1}{n-2}\delta_i\equiv 0
   \pmod{\{\theta,\omega^i,\pi_i\}},
$$
justifying the described reduction to $B_2\to M$, on which $\rho$ and
$\delta_i$ are semibasic.

Finally, a third reduction is made possible by the general equation
$$
  dS^i_j \equiv (n-2)\rho S^i_j - \alpha^i_kS^k_j + S^i_k\alpha^k_j
    +\sf12(U^{il}\sigma_{lj} + U^{jl}\sigma_{li}) - \sf1n\delta^i_j
    U^{kl}\sigma_{kl},
$$
modulo $\{\theta,\omega^i,\pi_i\}$.  On $B_2$, where in particular
$\lambda=-1$ and $\rho$ is semibasic, we have
$$
  dS^i_j \equiv -\alpha^i_kS^k_j + S^i_k\alpha^k_j - \sigma_{ij}
    \pmod{\{\theta,\omega^i,\pi_i\}}.
$$
This means that the torsion matrix $(S^i_j)$ can undergo translation by
an arbitrary traceless symmetric matrix along the fibers of $B_2\to M$,
so the locus
$$
  B_3 = \{u\in B_2: S^i_j(u) = 0\}\subset B_2
$$
is a subbundle, whose structure group is $SO(n,\R)$ with Lie
algebra represented by matrices of the form
$$
   a_2 = \left(\begin{array}{ccc} 0 & 0 & 0 \\
     0 & \alpha^i_j & 0 \\ 0 & 0 & -\alpha^j_i \end{array}\right),\qquad
   \alpha^i_j+\alpha^j_i = 0.
$$

This is all the reduction that we shall need.  On $B_3$, we have
equations
$$
  \left\{\begin{array}{l}
  \rho  =  -\frac{H}{2n}\theta, \\
  d\theta  =  -\pi_i\wedge\omega^i \quad
    \mbox{(because $\rho\wedge\theta = 0$ on $B_3$)}, \\
  (n-2)d\rho  =  -\delta_i\wedge\omega^i \quad
    \mbox{(because $t^i = -(n-2)E^i=0$ on $B_3$)}.\end{array}\right.
$$
The $\delta_i$ appearing the third equation are semibasic over $M$, and
the three equations together imply that
$$
  dH \equiv 0 \pmod{\{\theta,\omega^i\}}.
$$
This last observation is quite important.  Recall the integrable
Pfaffian system $J_\Pi =\{\theta,\omega^i\}$, assumed to have a
well-defined leaf-space $Q^{n+1}$ with submersion $M\to Q$.  The last
equation shows that $H$ is locally constant along the
fibers of $M\to Q$, and may therefore be thought of as a
function on $Q$.

Now, considering the two structure equations
$$
  \left\{\begin{array}{l}
    d\theta = -\pi_i\wedge\omega^i, \\
    d\omega^i = 2\rho\wedge\omega^i-\alpha^i_j\wedge\omega^j +
      \pi_i\wedge\theta,
  \end{array}\right.
$$
it is tempting to define
$$
  \tilde\pi_i = \pi_i+\sf{H}{n}\omega^i,
$$
and rewrite them as
$$
  d\left(\begin{array}{c}\theta\\ \omega^i\end{array}\right) = -
    \left(\begin{array}{cc} 0 & \tilde\pi_j\\
        -\tilde\pi_i & \alpha^i_j \end{array}\right)\wedge
    \left(\begin{array}{c}\theta\\ \omega^j\end{array}\right).
$$
Observe that this looks exactly like the structure equation
characterizing the Levi-Civita 
connection\index{connection!Levi-Civita|(} of a Riemannian
metric\index{Riemannian!metric|(}.  We
justify and use this as follows.  

Consider the quadratic form on $B_3$
$$
  \theta^2+\sum(\omega^i)^2.
$$
An easy computation shows that
for any vertical vector field $v\in Ker(\pi_*)$ for $\pi:B_3\to Q$,
$$
  {\mathcal L}_v\left(\theta^2+\sum(\omega^i)^2\right) = 0.
$$
This means that our quadratic form is the pullback of a quadratic form on 
$Q$, which defines there a Riemannian metric $ds^2$.  There is locally 
a bundle isomorphism over $Q$
$$
  B_3\to{\mathcal F}(Q,ds^2)
$$
\index{Riemannian!frame bundle|(}
from $B_3$, which was constructed from the neo-classical Poincar\'e-Cartan 
form $\Pi$, to the orthonormal frame bundle of this Riemannian metric. 
Under this isomorphism, the $Q$-semibasic forms $\theta,\omega^i$
correspond to the tautological semibasic forms on ${\mathcal
  F}(Q,ds^2)$\index{Riemannian!frame bundle|)}, while the matrix
$$
  \left(\begin{array}{cc} 0 & \tilde\pi_j\\
    -\tilde\pi_i & \alpha^i_j \end{array}\right)
$$
corresponds to the Levi-Civita 
connection\index{connection!Levi-Civita|)} matrix.  The contact
manifold $M$\index{contact!manifold}, 
as a quotient of $B_3$, may be then identified with the manifold of tangent 
hyperplanes to $Q$; and the Poincar\'e-Cartan form is
\begin{eqnarray*}
  \Pi & = & -\theta\wedge(\pi_i\wedge\omega_{(i)}) \\
      & = & -\theta\wedge(\tilde\pi_i\wedge\omega_{(i)} - H\omega).
\end{eqnarray*}
We recognize this as exactly the Poincar\'e-Cartan form for the prescribed 
mean curvature $H=H(q)$ system, in an arbitrary $(n+1)$-dimensional 
Riemannian manifold.  The following is what we have shown.

\begin{Theorem}
A definite neo-classical Poincar\'e-Cartan form $(M,\Pi)$ whose differential 
invariants satisfy $T^{ijk}=0$ and $U^{ij}=\lambda\delta^i_j$ with 
$\lambda<0$ is locally equivalent to the Poincar\'e-Cartan of the
prescribed mean curvature\index{mean curvature!prescribed} system on
some Riemannian manifold $(Q^{n+1}, ds^2)$.
\label{Theorem:PrescribedH}
\end{Theorem}
We will consider these Poincar\'e-Cartan forms further in
\S\ref{Section:SecondVarn}, when we
discuss the formula for the second variation\index{second variation}
of a Lagrangian
functional $\mathcal F_\Lambda$.  At that time, we will also see an
interpretation of the partial reduction $B_2\supset B_3$ in terms of
the Riemannian geometry.
Note that it is easy, given $(M,\Pi)$ as in the proposition, to determine 
the prescribed function $H(q)$ by carrying out the reductions described above, 
and to determine the Riemann curvature\index{Riemann curvature tensor}
of the ambient $(n+1)$-manifold in terms of the connection $1$-forms
$\tilde\pi_i$, $\alpha^i_j$.  The Euclidean minimal
surface\index{minimal surface} system discussed 
in \S\ref{Section:Hypersurface} is the case $H=0$, $R_{ijkl}=0$.
\index{Riemannian!metric|)}

The fact that such an $(M,\Pi)$ canonically determines $(Q,ds^2)$ implies 
the following.\footnote{As usual, this assumes that the foliation
  associated to $J_\Pi$ is
  simple\index{simple foliation|nn}; 
otherwise, only a local reformulation holds.}
\begin{Corollary}
The symmetry\index{symmetry} group of $(M,\Pi)$ is 
equal to the group of isometries of $(Q,ds^2)$ that preserve the 
function $H$.
\end{Corollary}
A consequence of this is the fact, claimed in
\S\ref{Section:Hypersurface}, that all symmetries 
of the minimal surface Poincar\'e-Cartan form---and hence, all classical 
conservation laws for the Euler-Lagrange equation---are induced by Euclidean 
motions\index{Euclidean!motion}.

Finally, in case $T^{ijk}=0$ and $U^{ij}=\lambda\delta^i_j$ with $\lambda>0$ 
instead of $\lambda <0$, one can carry out similar reductions, eventually 
producing on the quotient space $Q^{n+1}$ a Lorentz
metric\index{Lorentz!metric} $ds^2 = 
-\theta^2+\sum(\omega^i)^2$; the Poincar\'e-Cartan form is then equivalent 
to that for prescribed mean curvature of space-like
hypersurfaces\index{space-like hypersurface}.
\index{Poincar\'e-Cartan form!neo-classical|)}
\index{equivalence!of Poincar\'e-Cartan forms|)}

\chapter{Conformally Invariant Systems}

Among non-linear Euler-Lagrange equations on $\R^n$, 
the largest symmetry\index{symmetry|(} group that seems to occur is the
$\frac{(n+1)(n+2)}{2}$-dimensional conformal
group\index{conformal!group}.  This consists
of diffeomorphisms of the $n$-sphere that preserve its standard
conformal structure\index{conformal!structure|(}, represented by the
Euclidean metric under stereographic 
projection\index{stereographic projection} to $\R^n$.  These
maximally symmetric equations have a number of special properties,
including of course an abundance of classical conservation
laws\index{conservation law} as
predicted by Noether's theorem\index{Noether's theorem}.  
This chapter concerns the geometry of the
Poincar\'e-Cartan forms associated to these equations, and that of the
corresponding conservation laws.
\label{Chapter:Conformal}

We will begin by presenting background material on conformal
geometry.  This includes a discussion of the flat conformal
structure\index{conformal!structure!flat model|(}
on the $n$-sphere and its symmetry group, a construction of a canonical
parallelized principal bundle over a manifold with conformal
structure, and the
definition of the {\em conformal
  Laplacian}\index{conformal!Laplacian}, a second-order
differential operator associated to a conformal structure.  This
material will provide the framework for understanding the geometry of
non-linear Poisson equations\index{Poisson equation}, 
in particular the maximally symmetric non-linear example
$$
  \Delta u = Cu^{\frac{n+2}{n-2}},\qquad C\neq 0.
$$
After developing the geometric context for this equation, we will
continue the equivalence problem\index{equivalence method} 
for Poincar\'e-Cartan forms, pursuing
the branch in which these Euler-Lagrange equations occur.

We then turn to conservation laws for these conformally invariant
equations.  The elaborate geometric structure allows several
approaches to computing these conservation laws, and we will carry out
one of them in detail.  The analogous development for non-linear wave
equations\index{wave equation!non-linear} 
involves conformal structures with Lorentz signature, and
the conserved quantities for maximally symmetric Euler-Lagrange
equations in this case give rise to integral identities that have been very
useful in analysis.

\section{Background Material on Conformal Geometry}
\label{Section:Conformal}

In this section, we discuss some of the less widely known aspects
of conformal geometry.  In the first subsection, we define a flat
model for conformal geometry which is characterized by its large
symmetry\index{symmetry} group, and we give structure equations in
terms of the Maurer-Cartan\index{Maurer-Cartan!form} 
form of this group.  In the second
subsection, we give Cartan's solution to the local equivalence
problem for general conformal structures on manifolds.  This consists
of an algorithm by which one associates to any conformal structure 
$(N, [ds^2])$ a parallelized principal bundle $P\to N$
having structure equations of a specific algebraic form.  
In the third subsection, we introduce a second-order differential
operator $\Delta$, called the {\em conformal
  Laplacian}\index{conformal!Laplacian}, which is
associated to any conformal structure and which appears in the
Euler-Lagrange equations of conformal geometry that we study in the
remainder of the chapter.  The fundamental definition is the following.
\begin{Definition}  A 
{\em conformal inner-product}\index{conformal!inner-product} 
at a point $p\in
  N$ is an equivalence class of positive inner-products on $T_pN$,
  where two such inner-products are equivalent if one is a positive
  scalar multiple of the other.  A {\em conformal structure}, or {\em
    conformal metric}, on $N$ consists of a conformal inner-product at
  each point $p\in N$, varying smoothly in an obvious sense.
\end{Definition}
Note that this emphasizes the pointwise data of the conformal
structure, unlike the usual definition of a conformal structure as an
equivalence class of global Riemannian
metrics\index{Riemannian!metric}.  
An easy topological argument shows that these notions are equivalent.

\subsection{Flat Conformal Space}
\label{Subsection:FlatConformal}

We start with oriented Lorentz space\index{Lorentz!space} ${\mathbf
  L}^{n+2}$, with coordinates $x=(x^0,\ldots,x^{n+1})$, orientation 
$$
  dx^0\wedge\cdots\wedge dx^{n+1} > 0,
$$
and inner-product
$$
  \langle x,y\rangle = -(x^0y^{n+1}+x^{n+1}y^0)+\sum_ix^iy^i.
$$
Throughout this section, we use the index ranges
$ 0\leq a,b\leq n+1$ and $1\leq i,j\leq n$.

A non-zero vector $x\in{\mathbf L}^{n+2}$ 
is {\em null}\index{null vector} if $\langle x,x\rangle=0$.  A
null vector $x$ is {\em positive}\index{null vector!positive|(} if $x^0
> 0$ or $x^{n+1} >
0$; this designation is often called a
``time-orientation''\index{time-orientation} for
$\mathbf L^{n+2}$.  The symmetries of Lorentz space are the linear
transformations of $\mathbf L^{n+2}$ preserving the inner-product, the
orientation, and the time-orientation, and they constitute a connected
Lie group $SO^o(n+1,1)$.
We denote the space of positive null 
vectors\index{null vector!positive|)} by
$$
  Q = \{x\in{\mathbf L}^{n+2}:\langle x,x\rangle=0,\mbox{
    and } x^0 > 0 \mbox{ or } x^{n+1} > 0 \},
$$
which is one half of the familiar light-cone, with axis
$\{x^i=x^0-x^{n+1}=0\}$. 

We now define {\em flat conformal space} $R$ to be the space of null
lines in ${\mathbf L}^{n+2}$.  As a manifold, $R$ is a non-singular
quadric in the projective space $\mathbf{P}(\mathbf{L}^{n+2})$, which is
preserved by the natural action of the symmetry group $SO^o(n+1,1)$ of
${\mathbf L}^{n+2}$.  We will describe the flat conformal structure on
$R$ below, in terms of the Maurer-Cartan
form\index{Maurer-Cartan!form|(} 
of the group.  Note that the obvious map $Q\to R$, which we will write as
$x\mapsto [x]$, gives a principal bundle
with structure group $\R^*$.

In the literature, $R$ is usually defined as $\R^n$ with a point added at
infinity to form a topological sphere.  To make this identification, note that for
$x,y\in Q$, we have $\langle x,y\rangle\leq 0$, with
equality if and only if $[x]=[y]$.  We then claim that
$$
  H_y \stackrel{\mathit{def}}{=} \{x\in Q:\langle x,y\rangle = -1\}
$$
is diffeomorphic to both $\R^n$ and $R\backslash[y]$; this is
easily proved for $y=(0,\ldots,0,1)$, for instance, where the
map $\R^n\to H_y$ is given by
\begin{equation}
  (x^1,\ldots,x^n)\mapsto (1,x^1,\ldots,x^n,\sf12||x||^2).
\label{StandardCoords10}
\end{equation}  
The classical
description of the conformal structure on $R$ is obtained by
transporting the Euclidean metric on $\R^n$ to $H_y$, and noting
that for $y\neq y^\prime$ with $[y]=[y^\prime]$, this gives
unequal but conformally equivalent metrics on $R\backslash[y]$.
The fact that $SO^o(n+1,1)$ acts transitively on $R$ then implies that
for $[x]\neq[y]$ the conformal structures obtained on 
$R\backslash[x]$ and $R\backslash[y]$ are the same.

A {\em Lorentz frame}\index{Lorentz!frame} is a positively oriented basis
$f=(e_0,\ldots,e_{n+1})$ of ${\mathbf L}^{n+2}$, in which $e_0$ and
$e_{n+1}$ positive null vectors, and for which the
inner-product is (in blocks of size $1,n,1$, like most matrices in
this section) 
$$
  \langle e_a,e_b\rangle = \left(\begin{array}{ccc}
    0 & 0 & -1 \\ 0 & I_n & 0 \\ -1 & 0 & 0
    \end{array}\right).
$$
We let $P$ denote the set of all Lorentz frames.
There is a standard simply transitive right-action of $SO^o(n+1,1)$ on $P$,
by which we can identify the two spaces in a way that depends on a
choice of basepoint in $P$; this gives $P$ the structure of a smooth
manifold.  Because we have used the {\em right}-action,
the pullback to $P$ of any {\em left}-invariant $1$-form on $SO^o(n+1,1)$ is
independent of this choice of basepoint.  These pullbacks can be
intrinsically described on $P$ as follows.
We view each $e_a$ as a map $P\to{\mathbf L}^{n+2}$, and
we define $1$-forms $\rho$, $\omega^i$, $\beta_j$, $\alpha^i_j$ on $P$ by
decomposing the ${\mathbf L}^{n+2}$-valued $1$-forms $de_a$ in terms
of the bases $\{e_b\}$:
$$
  \left\{\begin{array}{l}
    de_0 = 2e_0\rho + e_i\omega^i, \\
    de_j = e_0\beta_j + e_i\alpha^i_j + e_{n+1}\omega^j, \\
    de_{n+1} = e_i\beta_i - 2e_{n+1}\rho.
  \end{array}\right.
$$
Equivalently,
$$
  d\left(\begin{array}{lcr} e_0 & e_j & e_{n+1}\end{array}\right)
    = \left(\begin{array}{lcr} e_0 & e_i &
        e_{n+1}\end{array}\right)
    \left(\begin{array}{ccc}  2\rho & \beta_j & 0 \\
      \omega^i & \alpha^i_j & \beta_i \\
      0 & \omega^j & -2\rho \end{array}\right).
$$
These forms satisfy $\alpha^i_j+\alpha^j_i=0$ but are otherwise
linearly independent, and they span the left-invariant $1$-forms on
$SO^o(n+1,1)$ under the preceding identification with $P$.  Decomposing the
exterior derivatives of these equations gives the Maurer-Cartan
equations\index{Maurer-Cartan!equation}, expressed in matrix form as
\begin{equation}
     d\left(\begin{array}{ccc}  2\rho & \beta_j & 0 \\
      \omega^i & \alpha^i_j & \beta_i \\
      0 & \omega^j & -2\rho \end{array}\right) +
    \left(\begin{array}{ccc}  2\rho & \beta_k & 0 \\
      \omega^i & \alpha^i_k & \beta_i \\
      0 & \omega^k & -2\rho \end{array}\right)\wedge
    \left(\begin{array}{ccc}  2\rho & \beta_j & 0 \\
      \omega^k & \alpha^k_j & \beta_k \\
      0 & \omega^j & -2\rho \end{array}\right) = 0.
\label{FlatConfStreqns}
\end{equation}

All of the local geometry of $R$ that is invariant under $SO^o(n+1,1)$
can be expressed in terms of these Maurer-Cartan
forms\index{Maurer-Cartan!form|)}.  In
particular, the fibers of the map $\pi_R:P\to R$ given by
$$
  \pi_R:(e_0,\ldots,e_{n+1})\mapsto [e_0]
$$
are the integral manifolds of the integrable Pfaffian
system
$$
  I_R = \{\omega^1,\ldots,\omega^n\}.
$$
This fibration has the structure of a principal bundle, whose
structure group consists of matrices in $SO^o(n+1,1)$ of the form
\begin{equation}
 g= \left(\begin{array}{ccc}
    r^2 & b_j & \sf12r^{-2}\ss b_j^2 \\
    0 & a^i_j & r^{-2}a^i_kb_k \\
    0 & 0 & r^{-2} \end{array}\right),
\label{ConfPrincipalGroup}
\end{equation}
where $r>0$, $a^i_ka^j_k=\delta^{ij}$.
Now, the symmetric differential form on $P$ given by
$$
  q = \sum(\omega^i)^2
$$
is semibasic for $\pi_R:P\to R$,
and a Lie derivative computation using the structure equations
(\ref{FlatConfStreqns}) gives, for any vertical vector field
$v\in\mbox{Ker }(\pi_R)_*$,
$$
  {\mathcal L}_vq = 4(v\innerprod\rho)q.
$$
This implies that there is a unique conformal structure $[ds^2]$ on $R$
whose representative metrics pull back under $\pi_R^*$ to multiples
of $q$.  By construction, this conformal structure is invariant under
the action of $SO^o(n+1,1)$, and one can verify that it gives the same
structure as the classical construction described above.

In \S\ref{Subsection:ConfEquiv}, we will follow Cartan in showing that
associated to any conformal
structure $(N,[ds^2])$ is a principal bundle $P\to N$ with $1$-forms
$\alpha^i_j=-\alpha^j_i$, $\rho$, $\omega^i$, and $\beta_j$, satisfying
structure equations like (\ref{FlatConfStreqns}) but with generally
non-zero curvature\index{conformal!curvature} terms on the right-hand side.

Before doing this, however, we point out a few more structures in the
flat model which will have useful generalizations.  These correspond to
Pfaffian systems
$$
  I_R = \{\omega^i\},\
  I_Q = \{\omega^i,\rho\},\ 
  I_M = \{\omega^i,\rho,\beta_j\},\ I_{P_0} = \{\omega^i,\rho,\alpha^i_j\},
$$
each of which is integrable, and in fact has a global quotient; that
is, there are manifolds $R$, $Q$, $M$, and $P_0$, and surjective
submersions from $P$ to each of these, whose leaves are the integral
manifolds of $I_R$, $I_Q$, $I_M$, and $I_{P_0}$, respectively:
$$
  \begin{array}{ccccc}
    & & P & & \\
    & \swarrow & & \searrow & \\
    P_0 & & & & M \\
    & \searrow & & \swarrow & \\
    &  & Q  &  & \\
    & & \downarrow & & \\
    & & R. & &
  \end{array}
$$
We have already seen that the leaves of the
system $I_R$ are fibers of the map $\pi_R:P\to R$.  Similarly, the
leaves of $I_Q$ are fibers of the map $\pi_Q:P\to
Q$ given by
$$
  \pi_Q:(e_0,\ldots,e_{n+1})\mapsto e_0.
$$

To understand the leaves of $I_M$, we let $M$ be the set of ordered pairs
$(e,e^\prime)$ of positive null vectors satisfying
$\langle e,e^\prime\rangle= -1$.
We then have a surjective submersion $\pi_M:P\to M$
defined by
$$
  \pi_M:(e_0,\ldots,e_{n+1})\mapsto(e_0, e_{n+1}),
$$
and the fibers of this map are the leaves of the Pfaffian system $I_M$.
Note that the $1$-form $\rho$ and its exterior derivative are
semibasic
for $\pi_M:P\to M$, and this means that there is a $1$-form
(also called $\rho$) on $M$ which pulls back to $\rho\in\Omega^1(P)$.
In fact, the equation for $d\rho$ in (\ref{FlatConfStreqns}) shows that on $P$,
$$
  \rho\wedge(d\rho)^n\neq 0,
$$
so the same is true on $M$.  Therefore, $\rho$ defines an
$SO^o(n+1,1)$-invariant contact structure\index{contact!manifold} 
on $M$.  The reader can
verify that $M$ has the structure of an $\R^*$-bundle over the
space $G^{(1,1)}(\mathbf L^{n+2})$ parameterizing those oriented
$2$-planes in $\mathbf
L^{n+2}$ on which the Lorentz metric has signature $(1,1)$.  In this
context, $2\rho\in\Omega^1(M)$ can be interpreted as a connection
$1$-form.  

Finally, to understand the leaves of $I_{P_0}$, we proceed as follows.
Define a {\em conformal frame}\index{conformal!frame bundle|(} 
for $(R,[ds^2])$
at a point $[x]\in R$ to be a positive basis $(v_1,\ldots,v_n)$
for $T_{[x]}R$ normalizing the conformal
inner-product\index{conformal!inner-product} as
$$
  ds^2(v_i,v_j) = \lambda\delta_{ij},
$$
for some $\lambda\in\R^*$ not depending on $i,j$.
The set of conformal frames for $(R,[ds^2])$ is the total space of a
principal bundle $P_0\to R$\index{conformal!frame bundle|)}, and there
is a surjective submersion
$P\to P_0$.  This last is induced by the maps $\overline{e_i}:P\to TR$
associating to a Lorentz frame\index{Lorentz!frame}
$f=(e_0,\ldots,e_{n+1})$ an obvious tangent vector
$\overline{e_i}$ to $R$ at $[e_0]$.  The reader can verify that the fibers of
the map $P\to P_0$ are the leaves of the Pfaffian system $I_{P_0}$.

Each of the surjective submersions $P\to
R$, $P\to Q$, $P\to M$, $P\to P_0$ has the structure of a principal
bundle, defined as a quotient of $P$ by a subgroup of $SO^o(n+1,1)$.
Additionally, the spaces $P$, $R$, $Q$, $M$, and $P_0$ are homogeneous
spaces of $SO^o(n+1,1)$, induced by the standard left-action on
${\mathbf L}^{n+2}$. 

We conclude with a brief description of the geometry of
$SO^0(n+1,1)$ acting on flat conformal space
$R$.\index{conformal!structure!flat model|)}  This will be useful later
in understanding the space of conservation laws of conformally
invariant Euler-Lagrange equations.  There are four main types of
motions.
\begin{itemize}
\item
The {\em
  translations}\index{translation} are defined as motions of $R$ induced by
left-multiplication by matrices of the
form
\begin{equation}
  \left(\begin{array}{ccc} 1 & 0 & 0 \\ w^i & I_n & 0 \\ 
    \sf{||w||^2}{2} & w^j & 1
    \end{array}\right).
\label{Translation10}
\end{equation}
In the standard coordinates on $R\backslash\{\infty\}$ described in
(\ref{StandardCoords10}), this is simply translation by the vector
$(w^i)$.
\item
The {\em rotations}\index{rotation} are defined as motions of $R$ induced by matrices
of the form
$$
  \left(\begin{array}{ccc} 1 & 0 & 0 \\ 0 & a^i_j & 0 \\ 0 & 0 & 1
    \end{array}\right),
$$
where $(a^i_j)\in SO(n,\R)$.  In the standard coordinates, this is the
usual rotation action of the matrix $(a^i_j)$.
\item
The {\em dilations}\index{dilation} are defined as motions of $R$ induced by matrices
of the form
$$
  \left(\begin{array}{ccc} r^2 & 0 & 0 \\ 0 & I & 0 \\ 0 & 0 & r^{-2}
  \end{array}\right).
$$
In the standard coordinates, this is dilation about the origin by a
factor of $r^{-2}$.
\item
The {\em inversions}\index{inversion} are defined as motions of $R$ induced by matrices
of the form
$$
  \left(\begin{array}{ccc} 1 & b_j & \sf{||b||^2}{2} \\ 0 & I & b_i \\
    0 & 0 & 1\end{array}\right).
$$
Note that these are exactly conjugates of the translation matrices
(\ref{Translation10}) by the matrix
$$
  J = \left(\begin{array}{ccc} 0 & 0 & 1 \\ 0 & I & 0 \\ 1 & 0 & 0
    \end{array}\right).
$$
Now, $J$ itself is not in $SO^o(n+1,1)$, but it still acts in an obvious way
on $R$; in standard coordinates, it gives the familiar inversion
in the sphere of radius $\sqrt{2}$.  So the inversions can be thought
of as conjugates of translation by the standard sphere-inversion, or
alternatively, as ``translations with the origin fixed''. 
\end{itemize}
These four subgroups generate  $SO^o(n+1,1)$.  Although the conformal
isometry group of $R$ has more that this one component, the others do not
appear in the Lie algebra, so they do not play a role in
calculating conservation laws for conformally invariant Euler-Lagrange
equations.
\index{symmetry|)}

\subsection{The Conformal Equivalence Problem}
\label{Subsection:ConfEquiv}

\index{equivalence!of conformal structures|(}
\index{equivalence method|(}
We will now apply the method of equivalence to conformal structures of
dimension $n\geq 3$.
This will involve some
of the ideas used in the equivalence problem for definite
Poincar\'e-Cartan forms discussed in the preceding chapter, but we
we will also encounter the new concept of {\em
  prolongation}\index{prolongation!of a $G$-structure|(}.  This is
the step that one takes when the usual process of
absorbing\index{torsion!absorption of} and
normalizing\index{torsion!normalization of} the
torsion in a $G$-structure\index{Gstructure@$G$-structure|(}
\index{torsion!of a pseudo-connection|(} does not
uniquely determine a pseudo-connection\index{pseudo-connection|(}. 

Let $(N,[ds^2])$ be an oriented conformal manifold of
dimension $n\geq 3$, and let $P_0\to N$ be the bundle
of {\em $0^{th}$-order oriented conformal coframes}
\index{conformal!frame bundle}
$\omega=(\omega^1,\ldots,\omega^n)$, which by definition satisfy
$$
  [ds^2] = \left[\ss(\omega^i)^2\right],\quad
  \omega^1\wedge\cdots\wedge\omega^n > 0.
$$
This is a principal bundle with structure group 
$$
  CO(n,\R)=\{A\in GL^+(n,\R):A\,^t\!A=\lambda I, 
    \mbox{ for some }\lambda\in\R^*\},
$$
having Lie algebra
\begin{eqnarray*}
  \lie{co}(n,\R) & = & \{a\in\lie{gl}(n,\R): a+\,^t\!a=\lambda I,
    \mbox{ for some }\lambda\in\R\} \\
  & = & \{(-2r\delta^i_j+a^i_j):a^i_j+a^j_i=0,\ a^i_j,r\in\R\}.
\end{eqnarray*}
We will describe a principal bundle $P\to P_0$, called the 
{\em prolongation} of $P_0\to N$,
whose sections correspond to torsion-free pseudo-connections in $P_0\to
N$, and construct a canonical parallelism of $P$ which defines a {\em
  Cartan connection}\index{connection!Cartan} in $P\to N$.  In
case $(N,[ds^2])$ is isomorphic to an open subset of flat conformal
space\index{conformal!structure!flat model}, this
will correspond to the restriction of the Lorentz frame bundle $P\to
R$ to that open subset, with parallelism given by the Maurer-Cartan
forms\index{Maurer-Cartan!form} of $SO^o(n+1,1)\cong P$.  

Recall that a pseudo-connection in $P_0\to N$ is a
$\lie{co}(n,\R)$-valued $1$-form
$$
  \varphi = (\varphi^i_j) = (-2\rho\delta^i_j+\alpha^i_j),\qquad
   \alpha^i_j+\alpha^j_i=0,
$$
whose restriction to each tangent space of a fiber of $P_0\to N$ gives
the canonical identification with $\lie{co}(n,\R)$ induced by the
group action.  As discussed previously (see \S\ref{Section:SmallEquiv}),
this last requirement means that $\varphi$
satisfies a structure equation
\begin{equation}
  d\omega^i = -\varphi^i_j\wedge\omega^j + 
    \sf12T^i_{jk}\omega^j\wedge\omega^k, \quad T^i_{jk}+T^i_{kj}=0,
\label{GenericStreqn}
\end{equation}
where $\omega^i$ are the components of the tautological $\R^n$-valued
$1$-form\index{tautological $1$-form} 
on $P_0$, and $\sf12T^i_{jk}\omega^j\wedge\omega^k$ is the semibasic
$\R^n$-valued torsion $2$-form.  We also noted previously that a
psuedo-connection $\varphi$ is a
genuine connection if and only if it is $Ad$-equivariant for the
action of $CO(n,\R)$ on $P_0$, meaning that
$$
  R_g^*\varphi = Ad_{g^{-1}}(\varphi),
$$
where $R_g:P_0\to P_0$ is the right-action of $g\in CO(n,\R)$ and
$Ad_{g^{-1}}$ is the adjoint action on $\lie{co}(n,\R)$, where
$\varphi$ takes its values.  However, completing this equivalence
problem requires us to consider the more general notion of a
pseudo-connection.  Although the parallelism that we eventually
construct is sometimes called the ``conformal connection'', there
is no canonical way (that is, no
way that is invariant under all conformal automorphisms) to associate
to a conformal structure a linear connection in the usual sense.

What we seek instead is a psuedo-connection $\varphi^i_j$ for which
the torsion vanishes, $T^i_{jk}=0$.  We know from the
fundamental lemma of Riemannian geometry, which guarantees a unique
torsion-free connection in the orthonormal frame bundle of any
Riemannian manifold, that whatever structure equation
(\ref{GenericStreqn}) we have with some initial pseudo-connection, we
can alter the pseudo-connection-forms $\alpha^i_j=-\alpha^j_i$ to
arrange that $T^i_{jk}=0$.  Specifically, we replace
$$
  \alpha^i_j\leadsto\alpha^i_j + 
    \sf12(T^i_{jk}-T^j_{ik}-T^k_{ij})\omega^k.
$$
So we can assume that $T^i_{jk}=0$, and we have simply
$$
  d\omega^i = -\varphi^i_j\wedge\omega^j =
    -(-2\rho\delta^i_j+\alpha^i_j)\wedge\omega^j,
$$
with $\alpha^i_j+\alpha^j_i=0$.  However, in contrast to
Riemannian geometry, this condition on the torsion does not uniquely
determine the psuedo-connection forms $\rho,\alpha^i_j$.  If we write
down an undetermined semibasic change in psuedo-connection
$$
  \left\{\begin{array}{l}
    \rho\leadsto \rho+t_k\omega^k, \\
    \alpha^i_j\leadsto \alpha^i_j + t^i_{jk}\omega^k, \quad 
      t^i_{jk}+t^j_{ik}=0,
  \end{array}\right.
$$
then the condition that the new pseudo-connection
be torsion-free is that
$$
  (2\delta^i_jt_k-t^i_{jk})\omega^j\wedge\omega^k = 0.
$$
This boils down eventually to the condition
$$
  t^i_{jk}=2(\delta^j_kt_i-\delta^i_kt_j).
$$
Therefore, given one torsion-free pseudo-connection $\varphi^i_j$ in
$P_0\to N$, the most general is obtained by adding
\begin{equation}
  2(\delta^i_jt_k -\delta^j_kt_i +
       \delta^i_kt_j)\omega^k,
\label{ChangeConn10}
\end{equation}
where $t=(t_k)\in\R^n$ is arbitrary.  This fact is needed for the next
step of the equivalence method, which consists of prolonging our
$CO(n,\R)$-structure.  We now digress to explain this general concept,
starting with the abstract machinery underlying the preceding
calculation.

\

We begin by amplifying the discussion of
normalizing\index{torsion!normalization of|(} torsion in
\S\ref{Section:SmallEquiv}. Associated to any linear Lie algebra
$\lie{g}\subset\lie{gl}(n,\R)$ is
an exact sequence of $\lie{g}$-modules
\begin{equation}
  0\to\lie{g}^{(1)}\to \lie{g}\otimes(\R^n)^* \stackrel{\delta}{\to}
    \R^n\otimes\bw{2}(\R^n)^* \to H^{0,1}(\lie{g})\to 0.
\label{SpencerSequence}
\end{equation}
Here, the map $\delta$ is the restriction to the subspace 
$$
  \lie{g}\otimes(\R^n)^* \subset (\R^n\otimes (\R^n)^*)\otimes(\R^n)^*
$$
of the surjective skew-symmetrization map
$$
  \R^n\otimes (\R^n)^*\otimes(\R^n)^* \to \R^n\otimes\bw{2}(\R^n)^*.
$$
The space $\lie{g}^{(1)}$ is the kernel of this restriction, and is
called the 
{\em prolongation}\index{prolongation!of a linear Lie algebra} 
of $\lie{g}$;
the cokernel $H^{0,1}(\lie{g})$, a {\em Spencer cohomology
  group}\index{Spencer cohomology} of
$\lie{g}$, was encountered in \S\ref{Section:SmallEquiv}.  Note that
$\lie{g}^{(1)}$ and $H^{0,1}(\lie{g})$ depend on the representation
$\lie{g}\hookrightarrow\lie{gl}(n,\R)$, and not just on the abstract Lie
algebra $\lie{g}$.

Recall from \S\ref{Section:SmallEquiv} that the intrinsic
torsion\index{torsion!intrinsic, of a $G$-structure} of a
$G$-structure vanishes if and only if there exist (locally)
torsion-free pseudo-connections in that $G$-structure.  This is a
situation in which further canonical reduction of the structure
group is not generally possible.  In particular, this will always
occur for $G$-structures with $H^{0,1}(\lie{g})=0$.

In this situation, the torsion-free pseudo-connection is unique if and
only if $\lie{g}^{(1)}=0$.  For example, when
$\lie{g}=\lie{so}(n,\R)$, both $\lie{g}^{(1)}=0$ and
$H^{0,1}(\lie{g})=0$, which accounts for the existence and uniqueness
of a torsion-free, metric-preserving 
connection\index{connection!Levi-Civita} on any Riemannian
manifold\index{Riemannian!manifold}.  In this favorable situation, we
have essentially completed
the method of equivalence, because the tautological form and the
unique torsion-free pseudo-connection constitute a canonical, global
coframing for the total space of our $G$-structure.  Equivalences of
$G$-structures correspond to isomorphisms of the associated
coframings, and there is a systematic procedure for determining when
two parallelized manifolds are locally isomorphic.

However, one frequently works with a structure group for which
$\lie{g}^{(1)}\neq 0$.
The observation that allows us to proceed in this case is that any
pseudo-connection
$\varphi$ in a $G$-structure $P\to N$ defines a particular type of
$\lie{g}\oplus\R^n$-valued 
coframing 
\begin{equation}
  \varphi\oplus\omega:TP\to\lie{g}\oplus\R^n
\label{ProlConn10}
\end{equation}
of the total space $P$.  Our previous discussion implies
that given some torsion-free pseudo-connection
$\varphi$, any change $\varphi^\prime$ lying in
$\lie{g}^{(1)}\subset\lie{g}\otimes(\R^n)^*$ yields a pseudo-connection
$\varphi+\varphi^\prime$ which is also torsion-free.  This means that
the coframings of $P$ as in (\ref{ProlConn10}), with $\varphi$
torsion-free, are
exactly the sections of a $\lie{g}^{(1)}$-structure $P^{(1)}\to P$,
where we regard $\lie{g}^{(1)}$ as an abelian Lie group.  This
$P^{(1)}\to P$ is by definition the {\em prolongation} of the
$G$-structure $P\to N$, and differential invariants of the former are
also differential invariants of the latter.\footnote{Situations with
  non-unique torsion-free pseudo-connections are not the only ones
  that call for prolongation; sometimes one finds intrinsic torsion
  lying in the fixed set of $H^{0,1}(\lie{g})$, and essentially
  the same process being described here must be used.  However, we
  will not face such a situation.}  The next natural step in
studying $P\to N$ is therefore to start over with $P^{(1)}\to P$, by
choosing a pseudo-connection, absorbing and normalizing its torsion,
and so forth.
\index{torsion!normalization of|)}

In practice, completely starting over would be wasteful.  The total
space $P^{(1)}$ supports tautological forms $\varphi$ and $\omega$,
valued in $\lie{g}$ and $\R^n$, respectively; and the equation
$d\omega+\varphi\wedge\omega=0$ satisfied by any particular torsion-free
psuedo-connection $\varphi$ on $P$ still holds on $P^{(1)}$ with
$\varphi$ replaced by a tautological form.  We can therefore
differentiate this equation and try to extract results about
the algebraic form of $d\varphi$.  These results can be interpreted as
statements about the intrinsic torsion of $P^{(1)}\to P$.
Only then do we return to the
usual normalization process.  We will now illustrate this,
returning to our situation in the conformal structure equivalence
problem.
\index{Gstructure@$G$-structure|)}
\index{prolongation!of a $G$-structure|)}
\index{equivalence method|)}

\

We have shown the existence of torsion-free pseudo-connections $\varphi^i_j$
in the $CO(n,\R)$-structure $P_0\to N$, so the intrinsic torsion
of $P_0\to N$ vanishes.\footnote{In fact, what we proved is that $\delta$ is
  surjective for $\lie{g}=\lie{co}(n,\R)$, so $H^{0,1}(\lie{co}(n,\R))=0$.}
We also have that such
$\varphi^i_j$ are unique modulo addition of a semibasic
$\lie{co}(n,\R)$-valued $1$-form linearly depending on an arbitrary
choice of $(t_k)\in\R^n$.  Therefore $\lie{co}(n,\R)^{(1)}\cong\R^n$,
and the inclusion
$\lie{co}(n,\R)^{(1)}\hookrightarrow\lie{co}(n,\R)\otimes(\R^n)^*$ is
described by (\ref{ChangeConn10}).
As explained above, we have an $\R^n$-structure
$P\stackrel{\mathit{def}}{=}(P_0)^{(1)}\to P_0$,
whose sections correspond to torsion-free pseudo-connections in
$P_0\to N$.  Any choice of the latter trivializes $P\to P_0$, and then
$(t_k)\in\R^n$ is a fiber coordinate.  We now search for structure
equations on $P$, with the goal of identifying a canonical
$\R^n$-valued pseudo-connection form for $P\to P_0$.

The first structure equation is still
$$
  d\omega^i = -\varphi^i_j\wedge\omega^j,
$$
where
$$
    \varphi^i_j = -2\delta^i_j\rho + \alpha^i_j,\quad 
      \alpha^i_j + \alpha^j_i = 0,
$$
and $\alpha^i_j$, $\rho$ are tautological forms on $P$.
Differentiating this gives
\begin{equation}
  (d\varphi^i_j + \varphi^i_k\wedge\varphi^k_j)\wedge\omega^j = 0,
\label{ProlongMaster}
\end{equation}
so
$$
  d\varphi^i_j + \varphi^i_k\wedge\varphi^k_j \equiv 0
    \pmod{\{\omega^1,\ldots,\omega^n\}}.
$$
Taking the trace of this equation of matrix $2$-forms shows that
$d\rho\equiv 0$, so guided
by the flat model (\ref{FlatConfStreqns}), we write
\begin{equation}
  d\rho = -\sf12\beta_i\wedge\omega^i
\label{DRho10}
\end{equation}
for some $1$-forms $\beta_i$ which are not uniquely determined.  We
will recognize these below as pseudo-connection forms in $P\to P_0$,
to be uniquely
determined by conditions on the torsion which we will uncover
shortly.  Substituting (\ref{DRho10}) back into
(\ref{ProlongMaster}), we have 
$$
  (d\alpha^i_j+\alpha^i_k\wedge\alpha^k_j-\beta_j\wedge\omega^i+
    \beta_i\wedge\omega^j)\wedge\omega^j = 0,
$$
and we set
$$
  A^i_j = d\alpha^i_j+\alpha^i_k\wedge\alpha^k_j -
    \beta_j\wedge\omega^i + \beta_i\wedge\omega^j.
$$
Note that $A^i_j + A^j_i = 0$.   We can write
$$
  A^i_j = \psi^i_{jk}\wedge\omega^k
$$
for some $1$-forms $\psi^i_{jk}=\psi^i_{kj}$,
in terms of which the condition $A^i_j+A^j_i=0$ is
$$
  (\psi^i_{jk}+\psi^j_{ik})\wedge\omega^k = 0,
$$
which implies
$$
  \psi^i_{jk} + \psi^j_{ik} \equiv 0 \pmod{\{\omega^1,\ldots,\omega^n\}}.
$$
Now computing modulo $\{\omega^1,\ldots,\omega^n\}$ as in Riemannian
geometry, we have
\begin{equation}
  \psi^i_{jk}\equiv - \psi^j_{ik} \equiv -\psi^j_{ki}\equiv
    \psi^k_{ji}\equiv \psi^k_{ij} \equiv -\psi^i_{kj}\equiv
    -\psi^i_{jk},
\label{StandardArgument}
\end{equation}
so $\psi^i_{jk}\equiv 0$.
We can now write
$$
  A^i_j = \psi^i_{jk}\wedge\omega^k =
    \sf12A^i_{jkl}\omega^k\wedge\omega^l,
$$
and forget about the $\psi^i_{jk}$, as our real interest is in
$d\alpha^i_j$.  We can assume that $A^i_{jkl}+A^i_{jlk}=0$, and we
necessarily have $A^i_{jkl}+A^j_{ikl}=0$.  Substituting once more into
$A^i_j\wedge\omega^j = 0$, we find that
$$
  A^i_{jkl}+A^i_{klj} + A^i_{ljk}=0.
$$
In summary, we have
$$
  d\alpha^i_j + \alpha^i_k\wedge\alpha^k_j -\beta_j\wedge\omega^i
    +\beta_i\wedge\omega^j = \sf12A^i_{jkl}\omega^k\wedge\omega^l,
$$
where $A^i_{jkl}$ has the symmetries of the Riemann curvature
  tensor\index{Riemann curvature tensor}.

In particular, we need only $n$ new $1$-forms $\beta_i$ to express the
derivatives of $d\rho$,
$d\alpha^i_j$.  The $\beta_i$ are pseudo-connection forms for
the prolonged $\lie{co}(n,\R)^{(1)}$-bundle $P\to P_0$, chosen to eliminate
torsion in the equation for $d\rho$, while the functions $A^i_{jkl}$
constitute the torsion in the equations for $d\alpha^i_j$.  Some of this
torsion will now be absorbed in the usual manner, by making a
uniquely determined choice of $\beta_i$.

Notice that the equation (\ref{DRho10}) for $d\rho$ is preserved
exactly under substitutions of the form
$$
  \beta_i \leadsto \beta_i + s_{ij}\omega^j,\qquad s_{ij}=s_{ji}.
$$
This substitution will induce a change
$$
  A^i_{jkl}\leadsto A^i_{jkl} + (-\delta^i_ls_{jk}+\delta^j_ls_{ik}
    +\delta^i_ks_{jl}-\delta^j_ks_{il}).
$$
Now, we know from the symmetries of the Riemann curvature tensor
that
$$
  A^l_{jkl} = A^l_{kjl},
$$
and on this contraction (the ``Ricci'' component) our substitution
will induce the change
$$
  A^l_{jkl}\leadsto A^l_{jkl}-(n-2)s_{jk}-\delta_{jk}s_{ll}.
$$
As we are assuming $n\geq 3$, there is a unique choice of
$s_{ij}$ which yields
$$
  A^l_{jkl}=0.
$$
It is not difficult to compute that the appropriate $s_{ij}$ is given
by
$$
  s_{ij} = \sf{1}{n-2}\left(A^l_{ijl}-\sf{1}{2n-2}
     \delta_{ij}A^l_{kkl}\right).
$$
In summary,
\begin{quote}
{\em
  On $P$, there is a unique coframing
  $\omega^i,\rho,\beta_j,\alpha^i_j=-\alpha^j_i$, where $\omega^i$ are the
  tautological forms over $N$, and such that the following structure equations
  are satisfied:
}
\end{quote}
\begin{eqnarray*}
   d\omega^i  & = &  (2\delta^i_j\rho-\alpha^i_j)\wedge\omega^j, \\
   d\rho  & = &  -\sf12\beta_i\wedge\omega^i, \\
   d\alpha^i_j  & = &  -\alpha^i_k\wedge\alpha^k_j+\beta_j\wedge\omega^i
     -\beta_i\wedge\omega^j + \sf12A^i_{jkl}\omega^k\wedge\omega^l, \\
   & & \qquad \mbox{\em with }A^l_{jkl}=0.
\end{eqnarray*}
\index{pseudo-connection|)}
\index{torsion!of a pseudo-connection|)}

We now seek structure equations for $d\beta_j$.  We start by
differentiating the simplest equation in which $\beta_j$ appears,
which is $d\rho = -\sf12\beta_j\wedge\omega^j$, and this gives
$$
  (d\beta_j+2\rho\wedge\beta_j+\beta_k\wedge\alpha^k_j)\wedge\omega^j
    = 0.
$$
We write
\begin{equation}
  d\beta_j+2\rho\wedge\beta_j+\beta_k\wedge\alpha^k_j = B_{jk}\wedge\omega^k
\label{Dbeta}
\end{equation}
for some $1$-forms $B_{jk}=B_{kj}$.  Because the equation for $d\rho$
did not determine $\beta_j$ uniquely, we cannot expect to use it to
completely 
determine expressions for $d\beta_j$; we need to
differentiate the equations for $d\alpha^i_j$, substituting
(\ref{Dbeta}).  This gives
$$
  (DA^i_{jkl}-B_{ik}\delta^j_l+B_{jk}\delta^i_l+B_{il}\delta^j_k
    -B_{jl}\delta^i_k)\wedge\omega^k\wedge\omega^l=0.
$$
Here we have defined for convenience the ``covariant
derivative''\index{covariant derivative}
\begin{equation}
  DA^i_{jkl}=dA^i_{jkl}+4\rho A^i_{jkl}+\alpha^i_mA^m_{jkl}
    -A^i_{mkl}\alpha^m_j - A^i_{jml}\alpha^m_k-A^i_{jkm}\alpha^m_l.
\label{DefDA10}
\end{equation}
Now we can write
$$
  DA^i_{jkl}-B_{ik}\delta^j_l+B_{jk}\delta^i_l+B_{il}\delta^j_k
    -B_{jl}\delta^i_k \equiv 0 \pmod{\{\omega^1,\ldots,\omega^n\}},
$$
and contracting on $il$ gives
$$
  B_{jk}\equiv 0 \pmod{\{\omega^1,\ldots,\omega^n\}}.
$$
This allows us to write simply
$$
  d\beta_i+2\rho\wedge\beta_i + \beta_j\wedge\alpha^j_i = 
    \sf12B_{ijk}\omega^j\wedge\omega^k,
$$
for some functions $B_{ijk}=-B_{ikj}$.  Returning to the equation
$$
  0 = d^2\rho = -\sf12d(\beta_j\wedge\omega^j)
$$
now yields the cyclic symmetry
$$
  B_{ijk}+B_{jki}+B_{kij}=0.
$$

We now have complete structure equations, which can be summarized in
the matrix form suggested by the flat
model\index{conformal!structure!flat model} (\ref{FlatConfStreqns}):
\begin{equation}
\boxed{
  \phi \stackrel{\mathit{def}}{=} \left(\begin{array}{ccc}
    2\rho & \beta_j & 0 \\ \omega^i & \alpha^i_j & \beta_i \\
      0 & \omega^j & -2\rho\end{array}\right),\quad
  \Phi\stackrel{\mathit{def}}{=}d\phi+\phi\wedge\phi =
    \left(\begin{array}{ccc} 0 & B_j & 0 \\ 0 & A^i_j & B_i \\
      0 & 0 & 0 \end{array}\right),
}
\label{ConfConn10}
\end{equation}
where
\begin{eqnarray*}
  A^i_j  & = &  \sf12A^i_{jkl}\omega^k\wedge\omega^l, \\
   & & \quad A^i_{jkl}  +A^j_{ikl}=A^i_{jkl}+A^i_{jlk} = 0, \\
   & & \quad A^i_{jkl}+A^i_{klj}+A^i_{ljk} = A^l_{jkl}=0, \\
  B_j & = &  \sf12B_{jkl}\omega^k\wedge\omega^l, \\
    & & \quad B_{jkl}  +B_{jlk}=B_{jkl}+B_{klj}+B_{ljk}=0.
\end{eqnarray*}
Furthermore, the action of $\R^n$ on $P\to P_0$ and that of
$CO(n,\R)$ on $P_0\to N$ may be combined, to realize $P\to
N$ as a principal bundle having structure group $G\subset SO^o(n+1,1)$
consisting of matrices
of the form (\ref{ConfPrincipalGroup}).  The matrix
$1$-form $\phi$ in (\ref{ConfConn10}) defines an $\lie{so}(n+1,1)$-valued
parallelism on $P$, under which the tangent spaces of fibers of $P\to
N$ are carried to the Lie algebra $\lie{g}\subset\lie{so}(n+1,1)$ of
$G$, and $\phi$ is equivariant with respect to the adjoint action of
$G$ on $\lie{so}(n+1,1)$.  The data of $(P\to N,\phi)$ is often called
a {\em Cartan connection}\index{connection!Cartan} modelled on
$\lie{g}\hookrightarrow\lie{so}(n+1,1)$.

We conclude this discussion by describing some properties of the
functions $A^i_{jkl}$, $B_{jkl}$ on $P$.  Differentiating the
definition of $\Phi$ (\ref{ConfConn10}) yields the {\em Bianchi
  identity}\index{Bianchi identity|(}
$$
  d\Phi = \Phi\wedge\varphi-\varphi\wedge\Phi.
$$
The components of this matrix equation yield linear-algebraic
consequences about the derivatives of $A^i_{jkl}$, $B_{jkl}$.  First,
one finds that
\begin{equation}
  \sf12DA^i_{jkl}\wedge\omega^k\wedge\omega^l =
    \sf12B_{ikl}\omega^j\wedge\omega^k\wedge\omega^l -
    \sf12B_{jkl}\omega^i\wedge\omega^k\wedge\omega^l.
\label{Dependence10}
\end{equation}
Detailed
information can be obtained from this equation, but note immediately
the fact that
$$
  DA^i_{jkl}\equiv 0 \pmod{\{\omega^i\}}.
$$
In particular, referring to the definition (\ref{DefDA10}), this shows
that the collection of functions $(A^i_{jkl})$ vary along the fibers
of $P\to N$ by a linear representation of the structure group $G$.  In
other words, they correspond to a section of an associated vector
bundle over $N$.  Specifically, we can see that the expression
$$
  A \stackrel{\mathit{def}}{=} \sf14A^i_{jkl}(\omega^i\wedge\omega^j\otimes
    \omega^k\wedge\omega^l)\otimes
      (\omega^1\wedge\cdots\wedge\omega^n)^{-2/n}
$$
on $P$ is invariant under the group action, so $A$ defines a
section of
$$
  \mbox{Sym}^2(\bw{2}T^*N)\otimes D^{-2/n},
$$
where $D$ is the {\em density line bundle}\index{density line bundle} 
for the conformal
structure, to be defined shortly.  This section is called the {\em Weyl
  tensor}\index{Weyl tensor|(} of the conformal structure.

Something different happens with $B_{jkl}$.  Namely, the Bianchi
identity\index{Bianchi identity|)} for $dB_{jkl}$ yields
\begin{eqnarray*}
  DB_{jkl} & \stackrel{\mathit{def}}{=} & dB_{jkl} +
    6\rho B_{jkl} - B_{mkl}\alpha^m_j-B_{jml}\alpha^m_k-
      B_{jkm}\alpha^m_l \\
  & \equiv & - \beta_i A^i_{jkl} \pmod{\{\omega^i\}}.
\end{eqnarray*}
In particular, the collection $(B_{jkl})$ transforms by a
representation of $G$ if and only if the Weyl tensor $A=0$.  In
case $n=3$, the symmetry identities of $A^i_{jkl}$ imply that $A=0$
automatically; there is no Weyl tensor in $3$-dimensional conformal
geometry.  In this case, $(B_{jkl})$ defines a section of the
vector bundle $T^*N\otimes\bigwedge^2T^*N$, which actually lies in a
subbundle, consisting of traceless elements of the kernel of
$$
  T^*N\otimes\bw{2}T^*N\to\bw{3}T^*N\to 0.
$$
This section is called the {\em Cotten tensor}\index{Cotten tensor} 
of the $3$-dimensional
conformal structure.  If the Cotten tensor vanishes, then the
conformal structure is locally equivalent to the flat conformal
structure on the $3$-sphere.

In case $n>3$, from (\ref{Dependence10}) one can show that the
functions $B_{jkl}$ can be expressed as linear combinations of the
covariant derivatives\index{covariant derivative} 
of $A^i_{jkl}$.  In particular, if the Weyl
tensor\index{Weyl tensor|)} $A$ vanishes, then so do all of the
$B_{jkl}$, and the
conformal structure of $N$ is locally equivalent to the flat conformal
structure on the $n$-sphere.
\index{equivalence!of conformal structures|)}

\subsection{The Conformal Laplacian}
\label{Subsection:ConfLaplacian}
\index{conformal!Laplacian|(}
\index{density line bundle|(}

To every conformal manifold $(N^n,[ds^2])$ is canonically associated a
linear differential operator $\Delta$, called the {\em conformal Laplacian}.
In this section, we define this operator and discuss its elementary
properties.  One subtlety is that $\Delta$ does not act on
functions, but on sections of a certain {\em density line bundle}, and
our first task is to define this.  We will use the parallelized
principal bundle $\pi:P\to N$ canonically associated to $[ds^2]$
as in the preceding discussion.

To begin, note that any $n$-form $\sigma$ on $N$ pulls back to $P$ to give a
closed $n$-form
$$
  \pi^*\sigma = u\,\omega^1\wedge\cdots\wedge\omega^n\in\Omega^n(P),
$$
where $u$ is a function on $P$ whose values on a fiber $\pi^{-1}(x)$
give the coefficient of
$\sigma_x\in\bw{n}(T_x^*N)$ with respect to various conformal coframes
at $x\in N$.  Among all $n$-forms on $P$ of the form
$u\,\omega^1\wedge\cdots\wedge\omega^n$, those that are locally pullbacks from
$N$ are characterized by the property of being closed.  Using the
structure equations, we find that this is equivalent to
$$
  (du + 2nu\rho)\wedge\omega^1\wedge\cdots\wedge\omega^n = 0,
$$
or
$$
  du \equiv -2nu\rho \pmod{\{\omega^1,\ldots,\omega^n\}}.
$$
This is the infinitesimal form of the relation
\begin{equation}
  u(p\cdot g) = r^{-2n}u(p),
\label{DensityEquivar}
\end{equation}
for $p\in P$ and $g\in G$ as in (\ref{ConfPrincipalGroup}).
This is in turn the same as saying that the function $u$ on $P$
defines a section of the oriented line bundle $D\to N$ associated to the
$1$-dimensional representation $g\mapsto r^{2n}$ of the structure
group.\footnote{That is, $D$ is the quotient of $P\times\R$ by the
  equivalence relation $(p,u)\sim(p\cdot g, r^{-2n}u)$ for $p\in P$,
  $u\in\R$, $g\in G$; a series of elementary exercises shows that this
  is naturally a line bundle over $N$, whose sections correspond to
  functions $u(p)$ satisfying (\ref{DensityEquivar}).}
Positive sections of $D$ correspond to oriented volume forms on $N$,
which in an obvious way correspond to Riemannian metrics representing the
conformal class $[ds^2]$.  Because so many of the PDEs studied
in the conformal geometry literature describe conditions on such
a metric, we should expect our study of Euler-Lagrange equations in
conformal geometry to involve this density bundle.

In analogy with this, we define for any positive real number $s$ the
degree-$\frac{s}{n}$ density bundle $D^{s/n}$ associated to the $1$-dimensional
representation $g\mapsto r^{2s}$; the degree-$1$ density bundle is
the preceding $D$.  Sections are represented by functions $u$ on $P$
satisfying 
\begin{equation}
  u(p\cdot g) = r^{-2s}u(p),
\label{DensityEquivar2}
\end{equation}
or infinitesimally,
\begin{equation}
  du \equiv -2su\rho\pmod{\{\omega^1,\ldots,\omega^n\}}.
\label{GenDensityInfEquivar}
\end{equation}
Summarizing, we will say that any function $u$ on $P$ satisfying
(\ref{GenDensityInfEquivar}) defines a section of the
degree-$\frac{s}{n}$ density bundle, and write
$$
  u\in\Gamma(D^{s/n}).
$$

We further investigate the local behavior of $u\in\Gamma(D^{s/n})$,
writing
$$
  du + 2su\rho = u_i\omega^i
$$
for some ``first covariant derivative''\index{covariant derivative} 
functions $u_i$.  Differentiating again, and applying the Cartan
lemma\index{Cartan lemma}, we obtain
\begin{equation}
  du_i + su\beta_i + 2(s+1)u_i\rho - u_j\alpha^j_i =
    u_{ij}\omega^j,
\label{SecondDeriv10}
\end{equation}
for some ``second covariant derivatives'' $u_{ij}=u_{ji}$; this is the
infinitesimal form of the transformation rule
\begin{equation}
  u_i(p\cdot g) = r^{-2(s+1)}(u_j(p)a^j_i-sb_iu(p)).
\label{JetDensityEquivar}
\end{equation}
Note that unless $s=0$ (so that $u$ is actually a function on $N$),
the vector-valued function $(u_i)$ on $P$ does not represent a section
of any associated vector bundle.

Differentiating again, and factoring out $\omega^k$, we obtain
modulo $\{\omega^1,\ldots,\omega^n\}$
$$
  du_{ij} \equiv \delta_{ij}u_k\beta_k - (s+1)(u_j\beta_i+u_i\beta_j) 
    -2(s+2)u_{ij}\rho+u_{kj}\alpha^k_i+u_{ik}\alpha^k_j,
$$
so once again, $u_{ij}$ is not a section of any associated vector
bundle.  However, we can take the trace
$$
  du_{ii} \equiv (n-2s-2)\beta_ku_k -2(s+2)\rho u_{ii}
    \pmod{\{\omega^1,\ldots,\omega^n\}},
$$
and we see that in case $s=\frac{n-2}{2}$, the function
$u_{ii}$ on $P$ is a section of $D^{\frac{s+2}{n}}$.  To summarize,
\begin{quote}
{\em
  the map $u\mapsto u_{ii}$ defines a second-order linear
    differential operator, called the conformal Laplacian,
}
$$
  \Delta:\Gamma(D^{\frac{n-2}{2n}})\to \Gamma(D^{\frac{n+2}{2n}}).
$$
\end{quote}
Note that for sections $u,v\in\Gamma(D^{\sf{n-2}{2n}})$, the quantity
$u\Delta v\in\Gamma(D^1)$ can be thought of as an $n$-form on $N$, and
integrated.  Furthermore, the reader can compute that
$$
  (u\Delta v - v\Delta u)\omega =
    d((uv_i-vu_i)\omega_{(i)}).
$$
We interpret this as saying that $u\Delta v-v\Delta u$ is
canonically a divergence, so that $(\cdot,\Delta\cdot)$ is a symmetric
bilinear form on $\Gamma_o(D^{\frac{n-2}{2n}})$, where
$$
  (\cdot,\cdot):\Gamma_o(D^{\sf{n-2}{2n}})\times
    \Gamma_o(D^{\sf{n+2}{2n}})\to\R 
$$
is given by integration on $N$ of the product.

To clarify the meaning of $\Delta$, we can choose a particular Riemannian
metric\index{Riemannian!metric} 
$g$ representing the conformal structure, and compare the
second covariant derivatives\index{covariant derivative} 
of an $\sf{s}{n}$-density $u$ taken in the
conformal sense with those derivatives taken in the usual sense of
Riemannian geometry.  By
construction of $P$, the pulled-back quadratic form
$\pi^*g\in\mbox{Sym}^2(T^*P)$ may be expressed as
$$
  \pi^*g = \lambda((\omega^1)^2+\cdots+(\omega^n)^2)
$$
for some function $\lambda>0$ on $P$.  Proceeding in a manner similar
to the preceding, we note that
$$
  {\mathcal L}_v(\pi^*g) = 0
$$
for any vector field $v$ that is vertical for $P\to N$.  Knowing the
derivatives of $\omega^i$ quite explicitly, we can then calculate that
\begin{equation}
  d\lambda = -4\lambda\rho + \lambda_i\omega^i
\label{ConfNorm1}
\end{equation}
for some functions $\lambda_i$.  Differentiating again, we find
\begin{equation}
  d\lambda_i = -2\lambda\beta_i - 6\lambda_i\rho + \lambda_j\alpha^j_i
    +\lambda_{ij}\omega^j,
\label{ConfNorm2}
\end{equation}
for some functions $\lambda_{ij}=\lambda_{ji}$.  Now we can reduce
our bundle $P\to N$ to a subbundle $P_g\subset P$, defined by
$$
  P_g = \{p\in P : \lambda(p) = 1,\ \lambda_1(p)=\cdots=
    \lambda_n(p)=0\}.
$$
Equations (\ref{ConfNorm1}, \ref{ConfNorm2}) imply that $P_g$ has
structure group $SO(n,\R)\subset G$, and using bars to denote
restrictions to $P_g$, we have for the pseudo-connection forms
$$
  \bar\rho = 0,\quad 
  \bar\beta_i = \sf12\bar\lambda_{ij}\bar\omega^j,\quad
  d\bar\omega^i = -\bar\alpha^i_j\wedge\bar\omega^j.
$$
The last of these means that if we identify $P_g$ with the usual
orthonormal frame\index{Riemannian!frame bundle}
bundle of $(N,g)$, then $\bar\alpha^i_j$ gives the Levi-Civita
connection\index{connection!Levi-Civita}.  The curvature is by definition
$$
  d\bar\alpha^i_j+\bar\alpha^i_k\wedge\bar\alpha^k_j = 
    \sf12R^i_{jkl}\bar\omega^k\wedge\bar\omega^l,
$$
but we have an expression for the left-hand side coming from the
conformal geometry; namely,
$$
  d\bar\alpha^i_j+\bar\alpha^i_k\wedge\bar\alpha^k_j
    = -\bar\beta_i\wedge\bar\omega^j + \bar\beta_j\wedge\bar\omega^i
      +\sf12A^i_{jkl}\bar\omega^k\wedge\bar\omega^l.
$$
Substituting $\bar\beta_i=\frac12\bar\lambda_{ij}\bar\omega^j$ and comparing these
two expressions gives
$$
  R^i_{jkl} = \sf12(\delta^i_l\bar\lambda_{jk} - \delta^j_l\bar\lambda_{ik}
    -\delta^i_k\bar\lambda_{jl}+\delta^j_k\bar\lambda_{il}) + A^i_{jkl}.
$$
From this we find the other components of curvature
\index{Ricci curvature}
\index{scalar curvature}
\begin{eqnarray*}
  & &  \mbox{Ric}_{jl} = R^i_{jil}
      = \sf12((2-n)\bar\lambda_{jl}-\delta^j_l\bar\lambda_{ii}), \\
  & & R = \mbox{Ric}_{ll}= (1-n)\bar\lambda_{ii}.
 \end{eqnarray*}

Now we will compute the conformal Laplacian of
$u\in\Gamma(D^{\frac{n-2}{2n}})$, but restrict the computation to
$P_g$.  Note that the choice of $g$ amounts to a trivialization of
$D$ and of all of its powers, so in this setting it is correct to
think of $u$ as a function.  We have
\begin{eqnarray*}
  d\bar u & = & -(n-2)\bar\rho \bar u +\bar u_i\bar\omega^i \\
    & = & \bar u_i\bar\omega^i, \\
  d\bar u_i & = & -\sf{n-2}2\bar u\bar\beta_i-n\bar u_i\bar\rho +
    \bar u_j\bar\alpha^j_i+\bar u_{ij}\bar\omega^j \\
   & = & \bar u_j\bar\alpha^j_i + 
     (\bar u_{ij}-\sf{n-2}4\bar u\bar\lambda_{ij})\bar\omega^j.
\end{eqnarray*}
Denoting by $\Delta_g$ the Riemannian
Laplacian\index{Riemannian!Laplacian}, we now have
\begin{eqnarray*}
  \Delta_g \bar u & = & \bar u_{ii}-\sf{n-2}4 \bar u\bar\lambda_{ii} \\
    & = & \Delta \bar u+\sf{n-2}{4(n-1)}R \bar u.
\end{eqnarray*}
This is the more familiar expression for the conformal Laplacian,
defined in terms of the Riemannian Laplacian of some representative
metric.  In the case of the flat 
model\index{conformal!structure!flat model} of conformal geometry, if one
uses standard coordinates on $\R^n=R\backslash\{\infty\}$, then the
Euclidean metric represents the conformal class, and we can use the
ordinary Laplacian $\Delta=\sum(\sf{\p}{\p x^i})^2$.
Its transformation properties, often stated and proved with
tedious calculations, can be easily derived from the present
viewpoint.

Of particular interest to us will be non-linear Poisson
equations\index{Poisson equation!non-linear|(}, of
the form
\begin{equation}
  \Delta u = f(x^i,u),
\label{Poisson10}
\end{equation}
where we will have an interpretation of $\Delta$ as the conformal
Laplacian on a conformal manifold with coordinates $x^i$.  We will
therefore want to interpret the unknown $u$ as a section of the
density bundle $D^{\frac{n-2}{2n}}$, and we will want to
interpret $f(x,u)$ as a ($0^{th}$-order) bundle map
$$
  f:D^{\frac{n-2}{2n}}\to D^{\frac{n+2}{2n}}.
$$
Certain obvious bundle maps $f$ come to mind.  One
kind is given by multiplication by any section
$\lambda\in\Gamma(D^{2/n})$; this would make (\ref{Poisson10}) a
linear equation.  Another is the appropriate power map
$$
  u\mapsto u^{\frac{n+2}{n-2}}.
$$
This yields a non-linear Poisson 
equation\index{Poisson equation!non-linear|)}, and we will examine it
quite closely in what follows.

We conclude this discussion with an alternate perspective on
the density bundles $D^{s/n}$.  First, note that for any
conformal manifold $(N,[ds^2])$ with its associated
parallelized  bundle $P\to N$, the Pfaffian system
$$
  I_Q = \{\rho,\omega^1,\ldots,\omega^n\}
$$
is integrable, and its associated foliation is 
simple\index{simple foliation}.
The leaf space of this foliation is just the quotient $Q$ of $P$ by the
action of a subgroup of its structure group, and this $Q$ is also a
fiber-bundle over $N$, with fiber $\R^*$.  This generalizes the
space $Q$ of positive null vectors\index{null vector!positive} 
in ${\mathbf L}^{n+2}$ which
appeared in the discussion of the flat model.  Now, the density
bundles $D^{s/n}$ are all canonically oriented, and we claim that
$Q$ is canonically identified with the {\em positive} elements of
$D^{s/n}$, for any $s$.

To see this, note that any positive $u\in D^{s/n}$, over $x\in N$, is
defined as a positive function on the fiber $P_x\subset P$ satisfying
(\ref{DensityEquivar2}).  It is not hard to see that the locus
$\{p\in P_x:u(p)=1\}\subset P_x$ is a leaf of the foliation defined by
$I_Q$.  Conversely, let $L_Q\subset P$ be a leaf of the foliation
defined by $I_Q$.  Then $L_Q$ lies completely in some fiber $P_x$ of $P\to N$,
and we can define a function $u$ on $P_x$ by setting $u=1$ on $L_Q$, and
extending to $P_x$ by the rule (\ref{DensityEquivar2}).  These
are clearly inverse processes.

We can extend the identification as follows.  Let
$J^1(N,D_+^{s/n})$ be the space of $1$-jets of positive sections
of $D$; it is a contact manifold\index{contact!manifold|(}, 
in the usual manner.  Let $M$ be the
leaf space of the simple foliation associated to the 
integrable Pfaffian system on $P$
$$
  I_M = \{\rho,\omega^1,\ldots,\omega^n,\beta_1,\ldots,\beta_n\}.
$$
This $M$ is also a contact manifold, a with global contact form pulling back
to $\rho\in\Omega^1(P)$, and it generalizes the contact manifold $M$
mentioned in our discussion of the flat model.  
We claim that there is a canonical contact isomorphism between
$J^1(N,D_+^{s/n})$ and $M$.

To see this, note that a $1$-jet at $x\in N$ of a positive section of
$D^{s/n}$ is specified by $n+1$ functions $(u,u_1,\ldots,u_n)$
on the fiber $P_x$ satisfying
(\ref{DensityEquivar2}, \ref{JetDensityEquivar}).  It is then not hard to
see that the locus $\{p\in P_x:u(p)=1,\ u_i(p)=0\}\subset P_x$ is a leaf of the foliation defined by $I_M$.
Conversely, let $L_M\subset P$ be a leaf of the foliation
defined by $I_M$.  Then $L_M$ lies completely in some fiber $P_x$ of
$P\to N$, and we can define $n+1$ functions $(u,u_1,\ldots,u_n)$ on
$P_x$ by setting $u=1$ and $u_i=0$ on $L_M$, and extending to $P_x$ by
the rules (\ref{DensityEquivar2}, \ref{JetDensityEquivar}).  
These are again inverse processes, and
we leave it to the reader to investigate the correspondence between
contact structures.
\index{contact!manifold|)}
\index{conformal!Laplacian|)}
\index{density line bundle|)}

\section{Conformally Invariant Poincar\'e-Cartan\\ Forms}
\markright{3.2. CONFORMALLY INVARIANT POINCAR\'E-CARTAN FORMS}
\index{Poincar\'e-Cartan form!conformally invariant|(}

In this section, we identify the Poincar\'e-Cartan forms on the
contact manifold $M$ over flat conformal
space\index{conformal!structure!flat model|(} $R$ that are
invariant under the action of the conformal
group\index{conformal!group} $SO^o(n+1,1)$.  We
then specialize to one that is neo-classical, and
determine expressions for the corresponding Euler-Lagrange equation
in coordinates; it turns out to be the non-linear Poisson equation
\index{Poisson equation!non-linear} with critical exponent
$$
  \Delta u = Cu^{\frac{n+2}{n-2}}.
$$ 
The calculation should clarify some of the more abstract
constructions of the preceding section.  It will also be helpful in
understanding the branch of the equivalence problem in which this
Poincar\'e-Cartan form appears, which is the topic of the next section.

We denote by $P$ the set of Lorentz frames\index{Lorentz!frame} 
for ${\mathbf L}^{n+2}$, by $M$ the set 
of pairs $(e,e^\prime)$ of positive null vectors with $\langle
e,e^\prime\rangle=-1$,
by $Q$ the space of positive
null vectors, and by $R$ the flat conformal
space of null lines.  There are $SO^o(n+1,1)$-equivariant maps
$$
  \left\{\begin{array}{ll}
    \pi_M:P\to M, & (e_0,\ldots,e_{n+1})\mapsto (e_0,e_{n+1}),
 \\ \pi_Q:P\to Q, & (e_0,\ldots,e_{n+1})\mapsto e_0,
 \\ \pi_R:P\to R, & (e_0,\ldots,e_{n+1})\mapsto [e_0].
  \end{array}\right.
$$
For easy reference we recall the structure equations for
Lorentz frames
\begin{eqnarray}
 & &  \left\{\begin{array}{l}
     de_0=2e_0\rho+e_i\omega^i, \\
     de_j=e_0\beta_j+e_i\alpha_j^i+e_{n+1}\omega^j,
        \qquad \alpha_j^i+\alpha_i^j=0, \\
     de_{n+1}=e_i\beta_i-2e_{n+1}\rho;
  \end{array}\right.
\label{DefConformalMC}
\\
 & & \left\{\begin{array}{l}
    d\rho+\sf12\beta_i\wedge\omega^i=0, \\
    d\omega^i-2\rho\wedge\omega^i+\alpha_j^i\wedge\omega^j=0, \\
    d\beta_i+2\rho\wedge\beta_i+\beta_j\wedge\alpha_i^j=0, \\
    d\alpha_j^i+\alpha_k^i\wedge\alpha_j^k+\beta_i\wedge\omega^j-
      \beta_j\wedge\omega^i=0.
  \end{array}\right.
\end{eqnarray}
We noted in the previous section that
$M$ has a contact $1$-form\index{contact!form} which pulls back to
$\rho$, and this is the setting for our Poincar\'e-Cartan forms.

\begin{Proposition}
The $SO^o(n+1,1)$-invariant Poincar\'e-Cartan forms on $M$, pulled back
to $P$, are constant linear combinations of
$$
  \Pi_k\stackrel{\mathit{def}}{=}\rho\wedge\sum\limits_{|I|=k}
    \beta_I\wedge\omega_{(I)},
$$
where $0\leq k\leq n$.  Those that are
neo-classical\index{Poincar\'e-Cartan form!neo-classical} 
with respect to $Q$ are of the form
\begin{equation}
  \Pi=c_1\Pi_1+c_0\Pi_0,\quad c_0,c_1\in\R.
\label{NeoClassical11}
\end{equation}
\end{Proposition}

\begin{Proof}
In this setting, an invariant Poincar\'e-Cartan form on $M$,
pulled back to $P$, is an $(n+1)$-form that is a 
multiple of $\rho$, semibasic over $M$,
invariant under the left-action of $SO^o(n+1,1)$, invariant under
the right-action of the isotropy subgroup $SO(n,\R)$ of $M$, and
closed.  That $\Pi$ must be semibasic and $SO^o(n+1,1)$-invariant
forces it to be a constant linear combination of exterior products of
$\rho,\beta_i,\omega^i$.
It is then a consequence of the Weyl's theory of vector invariants
that the further conditions of being a multiple of $\rho$ and
$SO(n,\R)$-invariant force $\Pi$ to be a linear combination of the
given $\Pi_k$.  
It follows from the structure equations of $P$ that $d\Pi_k=0$,
so each $\Pi_k$ is in fact the pullback of a Poincar\'e-Cartan form.
\end{Proof}

\

We note that for the $n$-form $\Lambda_k$ defined by
$$
  \Lambda_k \stackrel{\mathit{def}}{=} 
    \sum_{|I|=k}\beta_I\wedge\omega_{(I)}
$$
we have
$$
  d\Lambda_k = 2(n-2k)\Pi_k.
$$
This means that for $k\neq\frac{n}{2}$, the Poincar\'e-Cartan form
$\Pi_k$ is associated to an $SO^o(n+1,1)$-invariant functional, which in
the standard coordinates discussed below is second-order.
For the exceptional case $n=2k$, there is no
invariant functional corresponding to $\Pi_k$, but in the
neo-classical case $k\leq 1$ with $n\geq 3$, this is not an issue.

We now focus on the neo-classical case (\ref{NeoClassical11}), for
which it will be convenient to rescale and study
\begin{equation}
\boxed{
   \Pi = \rho\wedge\left(\beta_i\wedge\omega_{(i)} - 
    \sf{2C}{n-2}\omega\right),
}
\label{FlatNeoConfPCForm}
\end{equation}
where $C$ is a constant.  This is the exterior derivative of the
Lagrangian
$$
  \Lambda = \sf{1}{2(n-2)}\beta_i\wedge\omega_{(i)}
    -\sf{C}{n(n-2)}\omega,
$$
and our Monge-Ampere differential system\index{Monge-Ampere system} 
is generated by $\rho$ and the $n$-form
$$
  \Psi = \beta_i\wedge\omega_{(i)}-\sf{2C}{n-2}\omega.
$$

\begin{Proposition}
The Euler-Lagrange equation\index{Euler-Lagrange!equation} 
corresponding to the Poincar\'e-Cartan
form (\ref{FlatNeoConfPCForm}) is locally equivalent to
\begin{equation}
  \Delta u = Cu^{\frac{n+2}{n-2}}.
\label{FlatPoisson}
\end{equation}
\label{Proposition:FlatPoisson}
\end{Proposition}
The meaning of ``locally equivalent'' will come out in the proof.  It
includes an explicit and computable
correspondence between integral manifolds of the
Monge-Ampere system and solutions to the PDE.

We remark that the PDEs corresponding to
higher Poincar\'e-Cartan forms $\Pi_k$, with $k>1$,
have been computed and analyzed by J. Viaclovsky
in~\cite{Viaclovsky:Conformal}.

\

\begin{Proof} We begin by defining a map
  $\sigma: J^1(\R^n,\R)\hookrightarrow P$,\index{jets|(} 
  which can be projected to
  $M$ to give an open inclusion of contact manifolds with dense image.
  This map will
  be expressed in terms of the usual contact coordinates $(x^i,z,p_i)$
  on $J^1(\R^n,\R)$, except that $z$ is replaced by $u=e^{\lambda z}$
  for some undetermined constant $\lambda\neq 0$, so that in
  particular,
$$
  dz-p_idx^i = (\lambda u)^{-1}du-p_idx^i.
$$
  When using coordinates $(x^i,u,p_i)$, we denote our jet
  space by $J^1(\R^n,\R^+)$.
  We then
  pull back $\Psi$ via $\sigma$, and consider its restriction to a
  transverse Legendre submanifold.  With a convenient choice of $\lambda$,
  we will obtain a non-zero multiple of
  $\Delta u-Cu^{\frac{n+2}{n-2}}$, implying the Proposition.

We define $\sigma$ as a lift of the following map $\R^n\hookrightarrow
P$, to be
extended to $J^1(\R^n,\R^+)$ shortly:
\begin{equation}
  \bar e_0(x) = \left(\begin{array}{c}
    1 \\ x^1 \\ \vdots \\ x^n \\ \frac{||x||^2}{2}
  \end{array}\right),\
  \bar e_i(x) = \left(\begin{array}{c}
    0 \\ \vdots \\ 1_i \\ \vdots \\ x^i
  \end{array}\right),\
  \bar e_{n+1}(x) = \left(\begin{array}{c}
    0 \\ 0 \\ \vdots \\ 0 \\ 1
  \end{array}\right).
\label{PrelimBasis11}
\end{equation}
It is easy to verify that this does take values in $P$.  Also, note
that the composition $\R^n\hookrightarrow P\to R$ gives standard
(stereographic)
coordinates on 
$R\backslash\{\infty\}$\index{conformal!structure!flat model|)}.  
This partly indicates the
notion of ``locally equivalent'' used in this Proposition.  
We now let
\begin{equation}
  \begin{array}{rcl}
  e_0(x,u,p) & = & u^{2k}\bar e_0(x), \\
  e_i(x,u,p) & = & \bar e_i(x) + p_i\bar e_0(x), \\
  e_{n+1}(x,u,p) & = & u^{-2k}(\bar e_{n+1}(x) +
    p_j\bar e_j(x) + \sf{||p||^2}{2}\bar e_0(x)),
  \end{array}
\label{FinalBasis11}
\end{equation}
for some constant $k\neq 0$ to be determined shortly.
Our use of the dependent variable $u$ as a scaling factor for $e_0$
reflects the fact that we expect $u$ to represent a section of some
density line bundle.  The formula for $e_{n+1}$ is
chosen just so that our map takes values in $P$.
\index{jets|)}

Now we can compute directly 
\begin{eqnarray*}
  de_0 & = & 2ku^{-1}e_0du+u^{2k}\bar e_i dx^i  \\
    & = & 2(ku^{-1}du-\sf{1}{2}p_i dx^i)e_0+(u^{2k}dx^i)e_i,
\end{eqnarray*}
so by comparison with the expression in (\ref{DefConformalMC}) we
obtain some of the pulled-back
Maurer-Cartan forms\index{Maurer-Cartan!form}:
\begin{eqnarray*}
  \sigma^*\omega^i & = & u^{2k}dx^i,\\
  \sigma^*\rho & = & ku^{-1}du-\sf{1}{2}p_idx^i.
\end{eqnarray*}
Similarly, we have
\begin{eqnarray*}
  \beta_i & = & -\langle e_{n+1}, de_i\rangle \\
    & = & -\langle u^{-2k}(\bar e_{n+1}+p_j\bar e_j + 
        \sf{||p||^2}{2}\bar e_0), \\ & & \qquad\quad
      \bar e_{n+1}dx^i+p_i\bar e_kdx^k
        + \bar e_0dp_i\rangle \\
    & = & u^{-2k}\left(dp_i-p_ip_jdx^j+\sf{||p||^2}{2}
      dx^i\right).
\end{eqnarray*}

Because we want the projection to $M$ of $\sigma:J^1(\R^n,\R^+)\to P$ to
be a contact mapping, we need $\sigma^*\rho$ to be a multiple of
$dz-p_idx^i = (\lambda u)^{-1}du-p_idx^i$,
which holds if we choose
$$
  k=\sf{1}{2\lambda}.
$$
Now, $\lambda$ is still undetermined, but it will shortly be chosen to
simplify the expression for the restriction of $\Pi$ to a transverse Legendre
submanifold.  Namely, we find that
$$
  \beta_i\wedge\omega_{(i)} = u^{\frac{n-2}{\lambda}}
    \left(dp_i\wedge dx_{(i)}+\sf{n-2}{2}||p||^2dx\right),
$$
and also
$$
  \omega = u^{\frac{n}{\lambda}}dx.
$$

On transverse Legendre submanifolds, we have
$$
  du = \lambda e^{\lambda z}dz = \lambda up_idx^i,
$$
so that
$$
  p_i = \frac{1}{\lambda u}\frac{\partial u}{\partial x^i}.
$$
Differentiating, we obtain
$$
  dp_i = \frac{1}{\lambda}
    \left(\frac{1}{u}\frac{\partial^2u}{\partial x^i\partial x^j}
      -\frac{1}{u^2}\frac{\partial u}{\partial x^i}
        \frac{\partial u}{\partial x^j}\right)dx^j,
$$
so that on transverse Legendre submanifolds,
$$
  \Psi = u^{\sf{n-2}{\lambda}}\left(\frac{1}{\lambda}\frac{\Delta u}{u} + 
     \frac{n-2-2\lambda}{2\lambda^2}
     \frac{||\nabla u||^2}{u^2} - 
       \left(\frac{2C}{n-2}\right)u^{2/\lambda}\right)dx=0.
$$
We can eliminate the first-order term by choosing
$$
  \lambda = \sf{n-2}{2},
$$
and then
$$
  \Psi = \frac{2u}{n-2}\left(\Delta u-Cu^{\sf{n+2}{n-2}}\right)dx,
$$
which is the desired result.
\end{Proof}

\

Note that $z=\lambda^{-1}\log u$ satisfies a PDE that is slightly more
complicated, but equivalent under a classical
transformation\index{transformation!classical}.
Also, note that (\ref{FlatPoisson}) is usually given
as the Euler-Lagrange equation of the functional
$$
  \int\left(\sf12||\nabla u||^2+\sf{n-2}{2n}Cu^{\frac{2n}{n-2}}\right)dx,
$$
which has the advantage of being first-order, but the disadvantage of
not being preserved by the full conformal group\index{conformal!group}
$SO^o(n+1,1)$.  
In contrast, our Lagrangian $\Lambda$ restricts to transverse Legendre
submanifolds (in the coordinates of the preceding proof) as the
variationally equivalent integrand
$$
  \Lambda = \left(\sf{1}{(n-2)^2}u\Delta u - \sf{C}{n(n-2)}
                 u^{\frac{2n}{n-2}}\right)dx.
$$

\section{The Conformal Branch of the Equivalence Problem}
\markright{3.3. CONFORMAL BRANCH OF THE EQUIVALENCE PROBLEM}

\label{Section:ConfEquiv}

\index{Poisson equation!non-linear|(}
\index{equivalence method|(}
\index{equivalence!of Poincar\'e-Cartan forms|(}
\index{Poincar\'e-Cartan form!neo-classical!definite|(}

Let $(M^{2n+1},\Pi)$ be a manifold with a non-degenerate
Poincar\'e-Cartan form; that
is, $\Pi\in\Omega^{n+1}(M)$ is closed, and has a linear divisor that is
unique modulo scaling and defines a contact structure\index{contact!form}.
We also assume that $n\geq 3$ and that $\Pi$ is neo-classical and definite.
Then as discussed in \S\ref{Section:Bigequiv} we may associate to
$(M,\Pi)$ a $G$-structure\index{Gstructure@$G$-structure}
$B\to M$, where $G$ is a subgroup of $GL(2n+1,\R)$ whose Lie
algebra consists of matrices of the form
\begin{equation}
  \left(\begin{array}{ccc}
    (n-2)r & 0 & 0 \\
    0 & -2r\delta_j^i+a_j^i & 0 \\
    d_i & s_{ij} & nr\delta^j_i-a^j_i
  \end{array}\right),
\label{LieAlg12}
\end{equation}
where $a_j^i+a_i^j=0$ and $s_{ij}=s_{ji}$, $s_{ii}=0$.
In this section, we show how to uniquely characterize
in terms of the invariants of the $G$-structure
those $(M,\Pi)$ which are locally equivalent to the Poincar\'e-Cartan form for
the equation
\begin{equation}
 \Delta u= Cu^{\frac{n+2}{n-2}},\qquad C\neq 0.
\label{FlatPoissonEqn}
\end{equation}
on flat conformal space\index{conformal!structure!flat model}.
The result may be loosely summarized as follows.
\begin{quote}{\em
  The vanishing of the primary invariants $T^{ijk}$, $U^{ij}$, $S^i_j$
  is equivalent to the existence of a foliation $B\to N$ over a
  conformal manifold $(N,[ds^2])$, for which $[ds^2]$ pulls back to
  the invariant
  $[\sum(\omega^i)^2]$.  In this case, under open conditions on
  further invariants, three successive reductions of $B\to M$ yield a
  subbundle which is naturally identified with the conformal
  bundle over $N$.  The Poincar\'e-Cartan form can then be identified
  with that associated to a non-linear Poisson equation.  In case a
  further invariant is constant, this equation is equivalent to
  (\ref{FlatPoissonEqn}).}
\end{quote}

We find these conditions by continuing to apply the equivalence method
begun in \S\ref{Section:Bigequiv}, pursuing the case in which all of the
non-constant torsion vanishes.
One corollary of the discussion is a characterization of
Poincar\'e-Cartan forms locally equivalent to those for general
non-linear Poisson equations of the form
\begin{equation}
  \Delta u=f(x,u), \qquad x\in N,
\label{PoissonEqn}
\end{equation}
on an $n$-dimensional conformal manifold $(N,[ds^2])$; here and in the
following, $\Delta$ is the conformal
Laplacian\index{conformal!Laplacian}.
The condition that (\ref{PoissonEqn}) be non-linear can be
characterized in terms of the geometric invariants associated to
$(M,\Pi)$, as can the condition that $(N,[ds^2])$ be
conformally flat.  The characterization of (\ref{FlatPoissonEqn}) will
imply that this equation has
maximal symmetry\index{symmetry} group among non-linear Euler-Lagrange
equations satisfying certain geometric conditions on the
torsion.  We will not actually prove the characterization result for
general
Poisson equations (\ref{PoissonEqn}), but we will use these equations
(in the conformally flat case, with $\Delta =
\sum\left(\frac{\partial}{\partial x^i}\right)^2$) as an example at
each stage of the following calculations.  

We first recall the structure equations
of the $G$-structure $B\to M$, associated to a neo-classical, definite
Poincar\'e-Cartan form
$$
  \Pi = -\theta\wedge(\pi_i\wedge\omega_{(i)}).
$$
There is a pseudo-connection\index{pseudo-connection}
\begin{equation}
  \varphi =  \left(\begin{array}{ccc}
    (n-2)\rho &0 &0\\
    0 &-2\rho\delta_j^i+\alpha_j^i &0\\
    \delta_i &\sigma_{ij} &n\rho\delta_i^j-\alpha_i^j
  \end{array}\right),\mbox{ with }
  \left\{\begin{array}{l}
    \alpha^i_j+\alpha^j_i=0, \\
    \sigma_{ij}=\sigma_{ji},\ \sigma_{ii}=0,
  \end{array}\right.
\label{PCPseudoConn}
\end{equation}
having torsion\index{torsion!of a pseudo-connection}
\begin{equation}
  d\left(\begin{array}{c} \theta\\ \omega^i\\ \pi_i
     \end{array}\right) +
  \varphi\wedge
  \left(\begin{array}{c} \theta\\ \omega^j\\ \pi_j
    \end{array}\right) =
  \left(\begin{array}{c}
    -\pi_i\wedge\omega^i\\ 
 -(S^i_j\omega^j+U^{ij}\pi_j)\wedge\theta +
      T^{ijk}\pi_j\wedge\omega^k \\
     0 \end{array}\right),
\label{GeneralStructure}
\end{equation}
where enough torsion has been absorbed so that
\begin{equation}
  T^{ijk}=T^{jik}=T^{kji},\ T^{iik}=0;\ U^{ij}=U^{ji};\
    S^i_j=S^j_i,\ S^i_i=0.
\label{TorsionSymmetry}
\end{equation}
We also recall the structure equation (\ref{MasterDRho})
\begin{equation}
  (n-2)d\rho=-\delta_i\wedge\omega^i-S_j^i\pi_i\wedge
    \omega^j+\left(\sf{n-2}{2n}\right)U^{ij}\sigma_{ij}\wedge
    \theta-t^i\pi_i\wedge\theta.
\label{DRho}
\end{equation}
The equations
(\ref{PCPseudoConn}, \ref{GeneralStructure}, \ref{TorsionSymmetry},
\ref{DRho}) uniquely determine the forms $\rho$,
$\alpha_j^i$, and we are still free to alter our pseudo-connection by
\begin{equation}
  \left\{\begin{array}{l}
    \delta_i\leadsto\delta_i+ b_i\theta+t_{ij}\omega^j,
       \mbox{ with }t_{ij}=t_{ji}\mbox{ and }t_{ii}=0, \\
    \sigma_{ij}\leadsto\sigma_{ij}+t_{ij}\theta+t_{ijk}\omega^k,
      \mbox{ with } t_{ijk}=t_{jik}=t_{kji}\mbox{ and }t_{iik}=0,
  \end{array}\right.
\label{Ambiguity12}
\end{equation}
requiring also 
$$
  2nb_i+(n-2)U^{jk}t_{ijk}=0.
$$

We set up our example (\ref{PoissonEqn}) by taking coordinates
$(x^i,u,q_i)$ on $M=J^1(\R^n,\R)$,\index{jets} with contact
form\index{contact!form}
$$ \tilde\theta \stackrel{\mathit{def}}{=} du-q_idx^i. $$  
Then transverse Legendre submanifolds
which are also integral manifolds of
$$
  \tilde\Psi \stackrel{\mathit{def}}{=} -dq_i\wedge dx_{(i)}+f(x,u)dx
$$
correspond locally to solutions of (\ref{PoissonEqn}).  One can verify
that the form 
$$ \tilde\Pi \stackrel{\mathit{def}}{=} \tilde\theta\wedge\tilde\Psi $$ 
is closed, so in
particular our Poisson equation is an Euler-Lagrange equation.  We
find a particular $1$-adapted coframing of $J^1(\R^n,\R)$ as in Lemma
\ref{OneLemma8} by writing
$$
  \tilde\Pi = -\tilde\theta\wedge\left((dq_i-\sf{f}{n}dx^i)\wedge
  dx_{(i)}\right), 
$$
and then setting
$$
  \left(\begin{array}{c} \tilde\theta\\ \tilde\omega^i \\ \tilde\pi_i
    \end{array}\right) = \left(\begin{array}{c}
    du-q_idx^i \\ dx^i \\ dq_i-\sf{f}{n}dx^i\end{array}\right).
$$
It turns out that this coframing is actually a section of $B\to
J^1(\R^n,\R)$, as
one discovers by setting
$$
  \tilde\rho = 0,\ \tilde\alpha^i_j=0,\
  \tilde\delta_i=-\sf{1}{n}f_u\tilde\omega^i,
$$
and noting that the structure equations (\ref{PCPseudoConn},
\ref{GeneralStructure}) hold (with some complicated choice of
$\tilde\sigma_{ij}$ which we will not need).  In fact,
(\ref{GeneralStructure}) holds with torsion coefficients $S^i_j$,
$U^{ij}$, $T^{ijk}$ all vanishing, and we will see the significance of
this presently. 

In the general setting,
we seek conditions under which the quadratic form on $B$
$$
  q\stackrel{\mathit{def}}{=} \sum(\omega^i)^2
$$
can be regarded as defining a conformal
structure\index{conformal!structure} on some quotient of
$B$.  For the appropriate quotient to exist, at least locally, the
necessary and sufficient condition is that the Pfaffian system
$I=\{\omega^1,\ldots,\omega^n\}$ be integrable; it is easily seen from
the structure equations (and we noted in
\S\ref{Section:Bigequiv}) that this is equivalent to the condition
$$
  U^{ij}=0.
$$
We assume this in what follows, and for convenience assume further
that the foliation of $B$ by leaves of $I$ is
simple\index{simple foliation}; that is, there
is a smooth manifold $N$ and a surjective submersion $B\to N$ whose
fibers are the leaves of $I$.  Coordinates on $N$ may be thought of as
``preferred independent variables'' for the contact-equivalence class
of our Euler-Lagrange PDE, as indicated in \S\ref{Section:Bigequiv}.

We can now compute the Lie derivative of $q$ under a vector field $v$
which is vertical for $B\to N$, satisfying $v\innerprod\omega^i=0$;
using the hypothesis $U^{ij}=0$ and 
the structure equations, we find
$$
  {\mathcal L}_vq
     = 2\left(T^{ijk}(v\innerprod\pi_j)\omega^i\omega^k
       + S_j^i(v\innerprod\theta)\omega^i\omega^j\right)
       +4(v\innerprod\rho)q.
$$
It follows that if $T^{ijk}=0$ and $S^i_j = 0$,
then there is a quadratic form on $N$ which pulls back to a
non-zero multiple of $q$ on $B$.  A short calculation shows that the
converse as true as well, so we have the following.
\begin{Proposition}
The conditions $U^{ij}=T^{ijk}=S^i_j=0$ are
necessary and sufficient for there to exist (locally) a
conformal manifold $(N,[ds^2])$ and a map $B\to N$ such that the pullback to
$B$ of $[ds^2]$ is equal to $[q]=[\sum(\omega^i)^2]$.
\end{Proposition}
From now on, we assume $U^{ij}=S^i_j=T^{ijk}=0$.

From the discussion of the conformal equivalence problem in
\S\ref{Subsection:ConfEquiv}, we know that associated to
$(N,[ds^2])$ is the second-order conformal frame
bundle\index{conformal!bundle} $P\to N$ with
global coframing $\bar\omega^i$,
$\bar\rho$, $\bar\alpha_j^i$, $\bar\beta_i$ satisfying structure
equations
\begin{equation}
  \left\{  \begin{array}{l}
    d\bar\omega^i - 2\bar\rho\wedge\bar\omega^i + 
      \bar\alpha_j^i\wedge\bar\omega^j  =  0, \\
    d\bar\rho + \frac{1}{2}\bar\beta_i\wedge\bar\omega^i  =  0,\\ 
    d\bar\alpha_j^i +\bar\alpha_k^i\wedge\bar\alpha_j^k +
     \bar\beta_i\wedge\bar\omega^j - \bar\beta_j\wedge\bar\omega^i
        =  \sf12\bar A_{jkl}^i\bar\omega^k\wedge\bar\omega^l, \\
    d\bar\beta_i + 2\bar\rho\wedge\bar\beta_i +
      \bar\beta_j\wedge\bar\alpha_i^j
       =  \sf12 \bar B_{ijk}\bar\omega^j\wedge\bar\omega^k.
  \end{array}\right.
\label{InducedConfEqns}
\end{equation}
Our goal is to directly relate the principal bundle
$B\to M$ associated to the Poincar\'e-Cartan form $\Pi$ on
$M$ to the principal bundle $P\to N$ associated to the induced conformal geometry on $N$.
We shall eventually find that under some further conditions stated below,
the main one of which reflects the non-linearity of
the Euler-Lagrange system associated to $\Pi$, there is
a canonical reduction $B_3\to M$ of the $G$-structure $B\to M$
such that locally $B_3\cong P$ as parallelized manifolds.\footnote{As
  in the characterization in \S\ref{Section:PrescribedH} of
  prescribed mean
  curvature\index{mean curvature!prescribed|nn} 
  systems, we will denote by $B_1$, $B_2$, etc., reductions
  of the bundle $B\to M$ associated to $\Pi$, and these are unrelated
  to the bundles of the same names used in the construction of $B$.}
Because the canonical coframings on $B_3$ and $P$ determine the bundle
structure of each, we will then have shown that the subbundle $B_3\to N$ of
$B\to N$ can be locally 
identified with the bundle $P\to N$ associated to the conformal
structure $(N,[ds^2])$.

In the special case of our Poisson equation, we have
$\tilde\omega^i=dx^i$ as part of a section of $B\to M$, so we can
already see that our quotient space $N\cong\R^n$ is conformally
flat. 
This reflects the fact
that the differential operator $\Delta$ in (\ref{PoissonEqn}) is the
conformal Laplacian for flat conformal space.

We return to the general case, and make the
simplifying observation that under our hypotheses,
$$
  0 = d^2(\omega^1\wedge\cdots\wedge\omega^n) =
     -\sf{2n}{n-2}t^i\pi_i\wedge\theta\wedge\omega^1\wedge
       \cdots\wedge\omega^n,
$$
so that $t^i=0$ in the equation (\ref{DRho}) for $d\rho$.
We now have on $B$ the equations
\begin{eqnarray}
    d\omega^i &=& 2\rho\wedge\omega^i - \alpha^i_j\wedge\omega^j, 
\label{domega12}
\\
    d\rho &=& -\sf{1}{n-2}\delta_i\wedge\omega^i.
\label{AlmostConf12}
\end{eqnarray}
With the goal of making our structure equations on $B$ resemble the
conformal structure equations (\ref{InducedConfEqns}), we define
$$
  \beta_i \stackrel{\mathit{def}}{=}
    \left(\sf{2}{n-2}\right)\delta_i.
$$
The equations for $d\omega^i$ and $d\rho$ are now formally identical
to those for $d\bar\omega^i$ and $d\bar\rho$, and
computing exactly as in the conformal
equivalence problem, we find that
$$
  d\alpha_j^i+\alpha_k^i\wedge\alpha_j^k
    +\beta_i\wedge\omega^j-\beta_j\wedge\omega^i=\sf12
    A_{jkl}^i\omega^k\wedge\omega^l,
$$
for some functions $A^i_{jkl}$ on $B$ having the symmetries of the
Riemann curvature tensor\index{Riemann curvature tensor}.

Of course, we want $A^i_{jkl}$ to correspond to the Weyl
tensor\index{Weyl tensor} $\bar A^i_{jkl}$ of
$(N,[ds^2])$, so we would like to alter our pseudo-connection forms
(\ref{PCPseudoConn}) in a way that will give
$$
  A^l_{jkl}=0.
$$
Again, reasoning exactly as we did in the conformal equivalence
problem, we know that there are uniquely determined functions
$t_{ij}=t_{ji}$ such that replacing
$$
  \beta_i\leadsto\beta_i+t_{ij}\omega^j
$$
accomplishes this goal.  However, these may have $t_{ii}\neq 0$,
meaning that we cannot make the compensating change in $\sigma_{ij}$ (see
(\ref{Ambiguity12})) without introducing torsion in the equation for $d\pi_i$.
We proceed anyway, and now have structure equations
\begin{equation}
  d\left(\begin{array}{c}\theta \\ \omega^i \\ \pi_i\end{array}\right)
    + \left(\begin{array}{ccc}
      (n-2)\rho & 0 & 0 \\ 0 & \alpha^i_j\!-\!2\rho\delta^i_j & 0 \\
      \left(\sf{n-2}{2}\right)\beta_i & \sigma_{ij} &
        n\rho\delta^j_i\!-\!\alpha^j_i\end{array}\right)\wedge
   \left(\begin{array}{c}\theta \\ \omega^j \\ \pi_j\end{array}\right)
    = \left(\begin{array}{c} -\pi_i\wedge\omega^i \\ 0 \\
      A\omega^i\wedge\theta \end{array}\right)
\label{MoreStructure12}
\end{equation}
where $A=\frac{n-2}{2n}t_{ii}$ is a component of the original
$A^i_{jkl}$, analogous to scalar curvature\index{scalar curvature} 
in the Riemannian setting.  Also, we have
\begin{equation}
  d\rho = -\sf12\beta_i\wedge\omega^i,
\label{drho12}
\end{equation}
\begin{equation}
  d\alpha_j^i +\alpha_k^i\wedge\alpha_j^k
    +\beta_i\wedge\omega^j-\beta_j\wedge\omega^i  = \sf12
    A_{jkl}^i\omega^k\wedge\omega^l,
\label{dalpha12}
\end{equation}
with $A^l_{jkl}=0$.
These uniquely determine the pseudo-connection forms $\rho$,
$\alpha^i_j$, $\beta_i$, and leave $\sigma_{ij}$ determined only up to
addition of terms of the form $t_{ijk}\omega^k$, with $t_{ijk}$
totally symmetric and trace-free.

Now that $\beta_i$ is uniquely determined, we can once again mimic
calculations from the conformal equivalence problem, deducing from
(\ref{domega12}, \ref{drho12}, \ref{dalpha12}) that
\begin{equation}
  d\beta_i+2\rho\wedge\beta_i + \beta_j\wedge\alpha^j_i =
    \sf12B_{ijk}\omega^j\wedge\omega^k,
\label{dbeta12}
\end{equation}
with $B_{ijk}+B_{ikj}=0$, $B_{ijk}+B_{jki}+B_{kij}=0$.

In the case of our non-linear Poisson equation (\ref{PoissonEqn}), a
calculation shows that the modification of
$\tilde\beta_i=\sf{2}{n-2}\tilde\delta_i=0$ is not necessary, and that
with everything defined as before, we have not only
(\ref{MoreStructure12}), but also
(\ref{drho12}, \ref{dalpha12}, \ref{dbeta12}) with
$A^i_{jkl}=B_{ijk}=0$.  This gives us another way of seeing that the
conformal structure associated to (\ref{PoissonEqn}) is flat.  What
will be
important for us, however, is the fact that along this section of
$B\to J^1(\R^n,\R)$, the torsion function $A$ is
$$
  \tilde A=\sf{1}{n}f_u(x,u).
$$
This comes out of the calculations alluded to above.

We now begin to reduce $B\to M$, as promised.  
To get information about the derivative of the torsion coefficient
$A$ without knowing anything about $d\sigma_{ij}$,
we consider
\begin{eqnarray*}
  0 & = & d^2(\pi_i\wedge\omega_{(i)}) \\
    & = & d((n-2)\rho\wedge\pi_i\wedge\omega_{(i)} +
       \sf{n-2}{2}\theta\wedge\beta_i\wedge\omega_{(i)}
       +A\theta\wedge\omega) \\
    & = & (dA+4\rho A)\wedge\theta\wedge\omega.
\end{eqnarray*}
This describes the variation of the function $A$ along the fibers of
$B\to Q$, where we recall that $Q$ is the leaf space of the integrable
Pfaffian system $J_\Pi = \{\theta,\omega^i\}$.
In particular, we can write
\begin{equation}
  dA + 4\rho A = A_0\theta+A_i\omega^i,
\label{dA12}
\end{equation}
for some functions $A_0$, $A_i$ on $B$.  We see that on each fiber of
$B\to Q$, either $A$ vanishes
identically or $A$ never vanishes, and we assume that the latter holds
throughout $B$.  This is motivated by the case of the Poisson equation
(\ref{PoissonEqn}), for which $\tilde A=\sf{1}{n}f_u$ (so we are assuming
in particular that the zero-order term $f(x,u)$ depends on $u$).
Because the sign of $A$ is fixed, we assume $A>0$ in
what follows.  The case $A<0$ is similar, but the case $A=0$ is
quite different.

For the first reduction of $B\to M$, we define
$$
   B_1 = \{b\in B: A(b)=\sf14\}\subset B.
$$
From equation (\ref{dA12}) with the assumption $A>0$ everywhere, it
is clear that $B_1\to M$ is a principal subbundle of $B$, whose
structure group's Lie algebra consists of matrices (\ref{LieAlg12}) with
$r=0$.  Furthermore, restricted to $B_1$ there is a relation
\begin{equation}
  \rho = A_0\theta+A_i\omega^i.
\label{RhoReln12}
\end{equation}

In the case of a Poisson equation (\ref{PoissonEqn}), our section
$(\tilde\theta,\tilde\omega^i,\tilde\pi_i)$ of $B\to J^1(\R^n,\R)$ is
generally not a section of $B_1\subset B$, because we have along this
section that $\tilde A=\sf{f_u}n$.  However, (\ref{dA12}) guides us in
finding a section of $B_1$.  Namely, we define a function $r(x,u,q)>0$
on $M$ by
\begin{equation}
  r^4 = 4\tilde A = \sf{4}{n}f_u,
\label{rDef12}
\end{equation}
and then one can verify that for the coframing
$$
  \left(\begin{array}{c}\hat\theta\\ \hat\omega^i\\ \hat\pi_i
    \end{array}\right) \stackrel{\mathit{def}}{=} \left(\begin{array}{ccc}
      r^{2-n} & 0 & 0 \\ 0 & r^2\delta^i_j & 0 \\ 0 & 0 & r^{-n}\delta^j_i
     \end{array}\right)\left(\begin{array}{c}
    \tilde\theta\\ \tilde\omega^j \\ \tilde\pi_j
    \end{array}\right),
$$
one has the structure equation (\ref{MoreStructure12}), with
$$
  \hat\rho=r^{-1}dr=\sf14f_u^{-1}df_u,\quad
  \hat\beta_i=\hat\alpha^i_j=0,\quad \hat A=\sf14.
$$
Again, we won't have any need for $\hat\sigma_{ij}$.  Observe that
along this section of $B_1$, $\hat\rho=\sf14f_u^{-1}df_u=\hat
A_0\hat\theta+\hat A_i\hat\omega^i$, so that
\begin{equation}
  \hat A_0 = \sf14r^{n-2}f_u^{-1}f_{uu},\quad
  \hat A_i = \sf14r^{-2}f_u^{-1}(f_{ux^i}+f_{uu}q_i),
\label{Nonlinearity12}
\end{equation}
with $r$ given by (\ref{rDef12}).

Returning to the general situation on $B_1$, we differentiate
(\ref{RhoReln12}) and find
\begin{equation*}
  (dA_0 - (n-2)\rho A_0)\wedge\theta + 
    (dA_i+2\rho A_i-A_j\alpha^j_i+\sf12\beta_i-A_0\pi_i)
    \wedge\omega^i = 0,
\end{equation*}
and the Cartan lemma\index{Cartan lemma} then gives
\begin{eqnarray}
    & & dA_0-(n-2)\rho A_0  =  A_{00}\theta+A_{0i}\omega^i,
\label{Reducer1} \\
    & & dA_i + 2\rho A_i -A_j\alpha^j_i  +  
      \sf12\beta_i-A_0\pi_i = A_{i0}\theta+A_{ij}\omega^j,
\label{Reducer2}
\end{eqnarray}
with $A_{0i}=A_{i0}$ and $A_{ij}=A_{ji}$.

We interpret (\ref{Reducer1}) as saying that if $A_0$ vanishes at one
point of a fiber of $B_1\to Q$, then it vanishes everywhere on that
fiber.  We make the assumption that $A_0\neq 0$; the other extreme
case, where $A_0=0$ everywhere, gives a different branch of the
equivalence problem.  Note that in the case of a Poisson equation
(\ref{PoissonEqn}), the condition $A_0\neq 0$ implies by
(\ref{Nonlinearity12}) that the equation is everywhere non-linear.  This
justifies our decision to pursue, among the many branches of the
equivalence problem within the larger conformal branch, the case
$A>0$, $A_0\neq 0$.  This justification was our main reason to carry
along the example of the Poisson equation, and we will not mention it
again.  General calculations involving it become rather messy at this
stage, but how to continue should be clear from the preceding.

Returning to the general setting,
our second reduction uses (\ref{Reducer2}), which tells us that the
locus
$$
  B_2 = \{b\in B_1:A_i(b)=0\}\subset B_1
$$
is a principal subbundle of $B_1\to M$, whose structure group's Lie
algebra consists of matrices (\ref{LieAlg12}) with $r=d_i=0$.
Furthermore, restricted to $B_2$ there are relations
$$
  \beta_i = 2(A_{i0}\theta+A_{ij}\omega^j+A_0\pi_i),
$$
and also
$$
  \rho = A_0\theta.
$$

With the $A_i$ out of the way we differentiate once more, and
applying the Cartan lemma\index{Cartan lemma} find that on $B_2$, modulo
$\{\theta,\omega^i,\pi_i\}$,
\begin{eqnarray}
  dA_{00} & \equiv & 0,
\label{Trans1} \\
  dA_{0i} & \equiv & A_{0j} \alpha^j_i, 
\label{Trans2} \\
  dA_{ij} & \equiv & A_{kj}\alpha^k_i
    +A_{ik}\alpha^k_j + A_0\sigma_{ij}.
\label{Trans3}
\end{eqnarray}
We interpret (\ref{Trans1}) as saying that $A_{00}$ descends to a
well-defined function on $M$.
We interpret
(\ref{Trans2}) as saying that the vector valued function $(A_{0i})$
represents a section a vector bundle associated to $B_2\to M$.  
We interpret (\ref{Trans3}) as saying that if $A_0 =0$, then the
matrix $(A_{ij})$ represents a section of a vector bundle associated
to $B_2\to M$.  However, we have already made the assumption that
$A_0\neq 0$ everywhere.  In some examples of interest, most notably
for the equation $\Delta u=Cu^{\frac{n+2}{n-2}}$, the section
$(A_{0i})$ vanishes; for a general non-linear Poisson equation, this
vanishing loosely corresponds to the non-linearity being
translation-invariant on flat conformal space.  We will not need to
make any assumptions about this quantity.

This allows us to make a third reduction.  With $A_0\neq 0$,
(\ref{Trans3}) tells us that
the locus where the trace-free part of $A_{ij}$ vanishes,
$$
  B_3 = \{b\in B_2:
  A^0_{ij}(b) \stackrel{\mathit{def}}{=} A_{ij}(b)-
    \sf1n\delta_{ij}A_{kk}(b)=0\},
$$
is a subbundle $B_3\to
M$ of $B_2\to M$.  In terms of (\ref{LieAlg12}), the Lie algebra of the
structure group of $B_3$ is defined by $r=d_i=s_{ij}=0$.

Let us summarize what we have done.  Starting from the structure
equations 
(\ref{PCPseudoConn}, \ref{GeneralStructure}, \ref{TorsionSymmetry},
\ref{DRho})
on $B\to M$ for a definite, neo-classical
Poincar\'e-Cartan form with $n\geq 3$, we specialized to the case
where the torsion satisfies
$$
  U^{ij}=S^i_j=T^{ijk}=0.
$$
In this case, we found that the leaf space $N$ of the Pfaffian system
$\{\omega^1,\ldots,\omega^n\}$ has a conformal structure
pulling back to $\left[\sum(\omega^i)^2\right]$.  We
replaced each pseudo-connection form $\delta_i$ by its multiple
$\beta_i$, and guided by computations from conformal geometry, we
determined the torsion in the equation for $d\alpha^i_j$, which
resembled a Riemann curvature tensor.  This torsion's analog of scalar
curvature provided our fundamental invariant $A$, which had first
``covariant derivatives''\index{covariant derivative} 
$A_0$, $A_i$, and second ``covariant
derivatives'' $A_{00}$, $A_{i0}=A_{0i}$, $A_{ij}=A_{ji}$.  With the
assumptions
$$
  A\neq 0,\quad A_0\neq 0,
$$
we were able to make successive reductions by passing to the loci where
$$
  A=\sf14,\quad A_i=0,\quad A_{ij}=\sf1n\delta_{ij}A_{kk}.
$$

This leaves us on a bundle $B_3\to M$ with a coframing
$\omega^i,\rho,\beta_i,\alpha^i_j$, satisfying structure equations
exactly like those on the conformal bundle $P\to N$ associated with
$(N,[ds^2])$.  From here, a standard result shows that there is a
local diffeomorphism $B_3\to P$ under which the two coframings
correspond.  In particular, the invariants $A^i_{jkl}$ and $B_{jkl}$
remaining in the bundle $B_3$ equal the invariants named similarly in
the conformal structure, so we can tell for example if the conformal
structure associated to our Poincar\'e-Cartan form is flat.

We now write the restricted Poincar\'e-Cartan form, 
\begin{eqnarray*}
  \Pi & = & -\theta\wedge(\pi_i\wedge\omega_{(i)}) \\
      & = & -\sf{1}{A_0}\rho\wedge\sf{1}{A_0}
              \left(\sf12\beta_i-\sf1nA_{kk}\omega^i\right)
              \wedge\omega_{(i)} \\
      & = & -\sf{1}{2A_0^2}\rho\wedge(\beta_i\wedge\omega_{(i)}-2A_{kk}\omega).
\end{eqnarray*}

We can see from previous equations that $A_{kk}$ is constant on fibers
of $B_3\to M$.  Therefore, it makes sense to say that $A_{kk}$ is or
is not constant on $B_3$.  If it is constant, and if the conformal
structure on $N$ is flat (that is, $A^i_{jkl}=0$ if $n\geq 4$, or
$B_{jkl}=0$ if $n=3$), then our Poincar\'e-Cartan form is equivalent
to that associated to the non-linear Poisson equation
$$
  \Delta u = Cu^\frac{n+2}{n-2},
$$
where $C=(n-2)A_{kk}$.
This completes the characterization of Poincar\'e-Cartan
forms equivalent to that of this equation.  Our next goal is to
determine the conservation laws associated to this Poincar\'e-Cartan
form.
\index{conformal!structure|)}
\index{Poisson equation!non-linear|)}
\index{equivalence method|)}
\index{equivalence!of Poincar\'e-Cartan forms|)}
\index{Poincar\'e-Cartan form!neo-classical!definite|)}

\section{Conservation Laws for $\Delta u=Cu^{\frac{n+2}{n-2}}$}
\markright{3.4. CONSERVATION LAWS FOR $\Delta u=Cu^{\frac{n+2}{n-2}}$}
\index{conservation law!for x@for $\Delta u = Cu^{\frac{n+2}{n-2}}$|(}
In this section, we will determine the classical conservation laws for
the conformally invariant non-linear Poisson equation
\index{Poisson equation!conformally invariant|(}
\begin{equation}
  \Delta u = Cu^\frac{n+2}{n-2}.
\label{NLPoisson13}
\end{equation}
Recall that from $\Lambda_0=\omega$ and $\Lambda_1 = \beta_i\wedge
\omega_{(i)}$ we constructed the functional
$$
  \Lambda = \sf{1}{2(n-2)}\Lambda_1-\sf{C}{n(n-2)}\Lambda_0
$$
having the Poincare-Cartan form
$$
  \Pi = d\Lambda = \rho\wedge(\beta_i\wedge\omega_{(i)}
    -\sf{2C}{n-2}\omega),
$$
and that under a certain embedding
$\sigma:J^1(\R^n,\R^+)\hookrightarrow P$, the Euler-Lagrange system of
$\Pi$ restricted to a transverse Legendre submanifold is generated by
$$
  \Psi = \left(\Delta u - Cu^\frac{n+2}{n-2}\right)dx,
$$
for coordinates on $J^1(\R^n,\R^+)$ described in the proof of
Proposition~\ref{Proposition:FlatPoisson}.  We also proved that the
composition of $\sigma:J^1(\R^n,\R^+)\hookrightarrow P$ with the
projection $P\to M$ gives an open contact 
embedding of $J^1(\R^n,\R^+)$ as a dense subset of $M$.
Our invariant forms on $P$ pull back via $\sigma$ to give the
following forms on $J^1(\R^n,\R^+)$, expressed in terms of the
canonical coordinates $(x^i,u,p_i)$:
\begin{equation}
  \left\{\begin{array}{l}
    \rho = \sf{1}{n-2}u^{-1}du -\sf12p_idx^i, \\
    \omega^i = u^{\frac{2}{n-2}}dx^i,  \\
    \beta_i = u^{-\frac{2}{n-2}}\left(
      dp_i-p_ip_jdx^j+\sf{||p||^2}{2}dx^i\right), \\
    \alpha^i_j = p_jdx^i-p_idx^j.
  \end{array}\right.
\label{RestrictedMCForms}
\end{equation}
To describe the conservation laws, we first calculate for symmetry
vector fields $V\in\lie{g}_\Pi$ the expression
$$
  \varphi_V = V\innerprod\Lambda\in\Omega^{n-1}(P)
$$
at points of $J^1(\R^n,\R^+)\subset P$, and then restrict this $(n-1)$-form to
that submanifold, where it will be a conserved integrand for the
equation. 

\subsection{The Lie Algebra of Infinitesimal Symmetries}
\index{symmetry|(}

We know that the Poincar\'e-Cartan forms
$$
  \Pi_k = \rho\wedge\left(\sum_{|I|=k}\beta_I\wedge\omega_{(I)}\right)
$$
on $P$ are invariant under the simple, transitive left-action of the
conformal group\index{conformal!group} 
$SO^o(n+1,1)$.  The infinitesimal generators of this
action are the vector fields on $P$ corresponding under the
identification $P\cong SO^o(n+1,1)$ to {\it right-invariant} vector
fields.  Our first task is to determine the
right-invariant vector fields in terms of the basis
$$
  \left\{\frac{\p}{\p\rho},\ \frac{\p}{\p\omega^i},\ \frac{\p}{\p\beta_i},\
    \frac{\p}{\p\alpha^i_j}\right\}
$$
of left-invariant vector fields dual to the basis of left-invariant
$1$-forms used previously; this is because the Maurer-Cartan equation
in our setup
only allows us to compute in terms of left-invariant objects.

For an unknown vector field
\begin{equation}
  V = g\frac{\p}{\p\rho} + V^i\frac{\p}{\p\omega^i} +
    V_i\frac{\p}{\p\beta_i}+V^i_j\frac{\p}{\p\alpha^i_j}
  \quad (V^i_j+V^j_i=0)
\label{RightFieldCoefs}
\end{equation}
to be right-invariant is equivalent to the conditions
\begin{equation}
  {\mathcal L}_V\rho = {\mathcal L}_V\omega^i = 
    {\mathcal L}_V\beta_i = {\mathcal L}_V\alpha^i_j = 0;
\label{RightFieldEqns}
\end{equation}
that is, the flow of $V$ should preserve all left-invariant
$1$-forms.  We will solve the system
(\ref{RightFieldEqns}) of first-order differential equations 
for $V$ along the submanifold
$J^1(\R^n,\R^+)\subset P$.  Such $V$ are not
generally tangent to $J^1(\R^n,\R^+)$, but the calculation of
conservation laws as $V\innerprod\Lambda$ is still valid, as
$J^1(\R^n,\R^+)$ is being used only as a
slice of the foliation $P\to M$.  The solution will give the
coefficient functions $g$, $V^i$, $V_i$ of
(\ref{RightFieldCoefs}) in terms of the coordinates $(x^i,u,p_i)$ of
$J^1(\R^n,\R^+)$.  We will not need the coefficients
$V^i_j$, because they do not appear in
$\varphi_V=V\innerprod\Lambda$; in fact, we compute
$g=V\innerprod\rho$ only because it simplifies the rest of the solution.

First, we use the equation ${\mathcal L}_V\rho=0$, which gives
\begin{eqnarray*}
  0 & = & d(V\innerprod\rho) + V\innerprod d\rho \\
    & = & dg -\sf12(V_i\omega^i-V^i\beta_i).
\end{eqnarray*}
We have the formulae (\ref{RestrictedMCForms})
for the restrictions of $\omega^i$ and $\beta_i$ to $J^1(\R^n,\R^+)$,
by which the last condition becomes
$$
  dg = \frac12\left(V_iu^{\frac{2}{n-2}}dx^i -
    V^iu^{-\frac{2}{n-2}}(dp_i-p_ip_jdx^j+\sf{||p||^2}{2}dx^i)\right).
$$
This suggests that we replace the unknowns $V_i$, $V^i$ in our PDE system
(\ref{RightFieldEqns}) with
$$
  v^i \stackrel{\mathit{def}}{=} \sf12V^iu^{-\frac{2}{n-2}},\quad
  v_i \stackrel{\mathit{def}}{=} \sf12V_iu^{\frac{2}{n-2}}-\sf12V^ju^{-\frac{2}{n-2}}
          (-p_jp_i+\delta_{ij}\sf{||p||^2}{2}).
$$
Then we have the result
\begin{equation}
  \frac{\p g}{\p x^i} = v_i,\quad \frac{\p g}{\p p_i} = -v^i,
  \quad \frac{\p g}{\p u} = 0.
\label{Diffeqs13.1}
\end{equation}
In particular, we now need to determine only the function $g$.

For this, we use the equation ${\mathcal L}_V\omega^i=0$, which gives
\begin{eqnarray*}
  0 & = & d(V\innerprod\omega^i) + V\innerprod d\omega^i \\
    & = & dV^i-2\rho V^i + \alpha^i_jV^j + 2g\omega^i - 
      V^i_j\omega^j.
\end{eqnarray*}
When we restrict to $J^1(\R^n,\R^+)$ using (\ref{RestrictedMCForms})
and use our new dependent variables $v_i$, $v^i$, this gives
\begin{equation}
  dv^i = (p_idx^j-p_jdx^i)v^j-(p_jdx^j)v^i-g\,dx^i + \sf12V^i_jdx^j.
\label{Diffeqs13.2}
\end{equation}
This says in particular that $v^i(x,u,p)$ is a function of the variables $x^i$
alone, so along with (\ref{Diffeqs13.1}) we find that
$$
  g(x,u,p) = f(x) + f^i(x)p_i,
$$
for some functions $f(x)$, $f^i(x)$.  Substituting this back into
(\ref{Diffeqs13.2}), we have
$$
  df^i = (p_if^j-p_jf^i-\sf12V^i_j+\delta^i_jf)dx^j.
$$

This is a PDE system
$$
  \frac{\p f^i}{\p x^j} = p_if^j-p_jf^i-\sf12V^i_j+\delta^i_jf
$$
for the unknowns $f^i(x)$, and it can be solved in the
following elementary way.  We first let
$$
  h^i_j = p_if^j-p_jf^i-\sf12V^i_j = -h^j_i
$$
so that our equation is
\begin{equation}
  \frac{\p f^i}{\p x^j} = h^i_j +\delta^i_jf.
\label{Diffeq13.11}
\end{equation}
Differentiating this with respect to $x^k$ and equating mixed partials
implies that the expression
\begin{equation}
  \frac{\p h^i_k}{\p x^j}-\delta^i_j\frac{\p f}{\p x^k}
    +\delta^k_j\frac{\p f}{\p x^i}
\label{Diffeqs13.3}
\end{equation}
is symmetric in $j,k$.  It is also clearly skew-symmetric in $i,k$,
and therefore equals zero (as in (\ref{StandardArgument})).  
Now we can equate mixed partials of $h^i_k$ to obtain
$$
  \delta^i_k\frac{\p^2f}{\p x^jx^l}  -
  \delta^j_k\frac{\p^2f}{\p x^ix^l} =
  \delta^i_l\frac{\p^2f}{\p x^jx^k}  -
  \delta^j_l\frac{\p^2f}{\p x^ix^k}.
$$
With the standing assumption $n\geq 3$, this implies that all of
these second partial derivatives of $f$ are zero, and we can finally write
$$
  f(x) = r+\sf12b_kx^k,
$$
for some constants $r$, $b_k$.  The reasons for our labelling of these
and the following constants of integration will be indicated below.
Because the expressions
(\ref{Diffeqs13.3}) vanish, we can integrate to obtain
$$
  h^i_j = -\sf12a^i_j+\sf12b_jx^i-\sf12b_ix^j
$$
for some constants $a^i_j=-a^j_i$, and then integrate (\ref{Diffeq13.11}) to find
$$
  f^i(x) = -\sf12w^i + (\delta^i_jr-\sf12a^i_j)x^j - \sf14b_i||x||^2
    +\sf12\langle b,x\rangle x^i,
$$
where we have written $\langle b,x\rangle=\sum b_kx^k$ and $||x||^2 =
\sum (x^k)^2$.
We summarize the discussion in the following.

\begin{Proposition}
The coefficients of the vector fields (\ref{RightFieldCoefs}) on $P$
preserving the
left-invariant $1$-forms $\rho$, $\omega^i$ along $J^1(\R^n,\R^+)$ are
of the form
\begin{eqnarray*}
  g & = & r+\sf12\langle b,x\rangle(1+\langle p,x\rangle)
    +\left(-\sf12w^i+\left(\delta^i_jr-\sf12a^i_j\right)x^j
        -\sf14b_i||x||^2\right)p_i,
  \\ v^i &\stackrel{\mathit{def}}{=}& \sf12V^iu^{-\frac{2}{n-2}}\ = \
       -\frac{\p g}{\p p_i}, \\
    v_i &\stackrel{\mathit{def}}{=}& \sf12V_iu^{\frac{2}{n-2}}
          -\sf12V^ju^{-\frac{2}{n-2}}
          (-p_jp_i+\delta_{ij}\sf{||p||^2}{2})\ =\ \frac{\p g}{\p x^i},
\end{eqnarray*}
where $r$, $b_i$, $w^i$, $a^i_j=-a^j_i$ are constants.
\label{SymmetryFieldPropn}
\end{Proposition}
It is easy to verify that such $g$, $V^i$, $V_i$ uniquely determine
$V^i_j=-V^j_i$ such that the vector field (\ref{RightFieldCoefs})
preserves $\beta_i$ and $\alpha^i_j$ as well, but we will not need
this fact.  
Note that the number of constants in the Proposition
equals the dimension of the Lie algebra $\lie{so}(n+1,1)$, as
expected.

The reader may be aware that one should not have to solve differential
equations to determine right-invariant vector fields in terms of
left-invariant vector fields.  In fact, an algebraic calculation will
suffice, which in this case would consist of writing an arbitrary Lie
algebra element
$$
  g_L = \left(\begin{array}{ccc} 2r & b_j & 0 \\ w^i & a^i_j & b_i \\
    0 & w^j & -2r\end{array}\right)
$$
interpreted as a left-invariant vector field, and conjugating by
$\sigma(x,u,p)\in P$  regarded as a matrix with columns $e_0(x,u,p)$,
$e_j(x,u,p)$, $e_{n+1}(x,u,p)$ given by
(\ref{PrelimBasis11}, \ref{FinalBasis11}).  The resulting
$\lie{so}(n+1,1)$-valued function on $J^1(\R^n,\R^+)$ then has entries
which are the coefficients of a right-invariant vector field $V$.  The
calculation is tedious, but of course the vector fields so obtained
are as in Proposition~\ref{SymmetryFieldPropn}.
\index{symmetry|)}

\subsection{Calculation of Conservation Laws}

We can now use the formulae for the infinitesimal symmetries derived above
to calculate the conservation laws for $\Pi$, which are
$(n-1)$-forms on $J^1(\R^n,\R^+)$ that are closed when restricted to
integral submanifolds of the Euler-Lagrange system.

The Noether prescription\index{Noether's theorem} 
is particularly simple in this case, because the equations
$$
  {\mathcal L}_V\Lambda = 0 \qquad \mbox{and} \qquad
    d\Lambda = \Pi
$$
mean that there are no compensating terms, and we can take for the
conserved integrand just
$$
  \varphi_V = V\innerprod\Lambda.
$$
This is straightforward in principle, but there are some
delicate issues of signs and constants.  We find that for $V$ as in
Proposition \ref{SymmetryFieldPropn},
$$
  V\innerprod\Lambda_0 = V^i\omega_{(i)},
$$
and restricting to $J^1(\R^n,\R^+)$, using $v^i$ instead of $V^i$, we obtain
$$
  (V\innerprod\Lambda_0)|_{J^1(\R^n,\R^+)} =
  2u^{\frac{2n}{n-2}}v^idx_{(i)}.
$$
The analogous computation for $V\innerprod\Lambda_1$ is a little more
complicated and gives
$$
  (V\innerprod\Lambda_1)_{J^1(\R^n,\R^+)} = 2u^2(-v^jdp_i\wedge dx_{(ij)} + 
    (v_i +\sf{n-2}{2}v^i||p||^2)dx_{(i)}).
$$
On a transverse Legendre submanifold $S$ of $J^1(\R^n,\R^+)$, we can use
the condition $\rho=0$ from (\ref{RestrictedMCForms}) to write
\begin{equation}
  p_i = \sf{2}{n-2}u^{-1}\sf{\p u}{\p x^i},
\label{LegendreDerivs}
\end{equation}
and if we compute $dp_i$ and $||p||^2$ for
such a submanifold, then we can substitute and obtain
$$
  (V\innerprod \Lambda_1)|_S = 
    \sf{4}{n-2}\left(-uu_{x^ix^j}v^j+uu_{x^jx^j}v^i + 
      u_{x^i}u_{x^j}v^j+\sf{n-2}{2}u^2v_i\right)dx_{(i)}.
$$
We summarize with the following.
\begin{Proposition}
The restriction of $V\innerprod\Lambda$ to the $1$-jet graph of
$u(x^1,\ldots,x^n)$ equals
$$
\boxed{
  \varphi_V = \left(\sf{2}{n-2}\left(
    uu_{x^jx^j}v^i - uu_{x^ix^j}v^j
      +u_{x^i}u_{x^j}v^j\right) -\sf{2C}{n}u^{\frac{2n}{n-2}}v^i
      +u^2v_i\right)dx_{(i)}.
}
$$
\label{CLPropn}
\end{Proposition}

We now have a representative for each of the classical conservation
laws corresponding to a conformal symmetry of our equation
\begin{equation}
  \Delta u = Cu^{\frac{n+2}{n-2}}.
\label{NLP13}
\end{equation}
We say ``representative'' because a conservation law is actually an
equivalence class of $(n-1)$-forms.  In fact, our $\varphi_V$ is not
the $(n-1)$-form classically taken to represent the conservation law
corresponding to $V$; our $\varphi_V$ involves second
derivatives of the unknown $u(x)$, while the classical expressions are
all first-order.  We can find the first-order expressions by adding to
$\varphi_V$ a suitable exact $(n-1)$-form, obtaining
\begin{eqnarray*}
  \varphi_g & \stackrel{\mathit{def}}{=} & \varphi_V + 
     \sf{2}{n-2}d(uu_{x^i}v^jdx_{(ij)}) \\  & = &
  \left(\sf4{n-2}u_{x^i}u_{x^j}v^j - 
    \left(\sf2{n-2}||\nabla u||^2+\sf{2C}nu^\frac{2n}{n-2}\right)v^i
    \right. \\
   &\ &\qquad \left. +
   u^2v_i+\sf{2}{n-2}u(u_{x^i}v^j_{x^j}-u_{x^j}v^i_{x^j})\right)dx_{(i)}.
\end{eqnarray*}

This turns out to give the classical expressions for the conservation laws
associated to our equation (\ref{NLP13}), up to multiplicative
constants.  It could have been obtained more directly using the methods of
Section~\ref{Section:Noether}.  For this, one would work on the usual
$J^1(\R^n,\R)$, with standard coordinates $(x^i,u,q_i)$ in which the
contact structure is generated by
$$
  \theta = du-q_idx^i,
$$
and then consider the Monge-Ampere system generated by $\theta$ and
$$
  \Psi = -dq_i\wedge dx_{(i)}+Cu^\frac{n+2}{n-2}dx.
$$
A little experimenting yields a Lagrangian density
$$
  L\,dx = \left(\sf{||q||^2}{2}+\sf{n-2}{2n}Cu^{\sf{2n}{n-2}}\right)dx,
$$
so the functional
$$
  \Lambda = L\,dy + \theta\wedge L_{q_i}dy_{(i)}
$$
induces the Poincar\'e-Cartan form
$$
  \Pi = \theta\wedge\Psi = d\Lambda.
$$
One can then determine the Lie algebra of the symmetry group of $\Pi$
by solving an elementary PDE system, with a result closely resembling
that of Proposition~\ref{SymmetryFieldPropn}.  Applying the Noether
prescription\index{Noether's theorem} 
to these vector fields and this $\Lambda$ yields
$(n-1)$-forms which restrict to transverse Legendre submanifolds to give
$\varphi_g$ above.

Returning to our original situation,
we now compute $\varphi_g$ explicitly for various choices of $g$ as in
Proposition \ref{SymmetryFieldPropn}.
These choices of $g$ correspond
to subgroups of the conformal group.

\subsubsection{Translation: $g = w^ip_i$.}
\index{translation}
\index{conformal!structure!flat model|(}
In this case, we have $v^i = -w^i$, $v_i=0$, so we find on a transverse
Legendre submanifold of $J^1(\R^n,\R^+)$ that
$$
  \varphi_g = \left(\sf2{n-2}||\nabla u||^2w^i - 
    \sf4{n-2}u_{x^i}u_{x^j}w^j +
    \sf{2C}{n}u^{\frac{2n}{n-2}}w^i\right)dx_{(i)}.
$$
The typical use of a conservation law involves its integration along
the smooth $(n-1)$-dimensional boundary of a region
$\Omega\subset\R^n$.  To make more sense of the preceding expression,
we take such a region to have unit normal $\nu$ and area element
$d\sigma$ (with respect to the Euclidean metric), and using the fact that $q^idx_{(i)}|_{\p\Omega}=\langle
q,\nu\rangle d\sigma$ for a vector $q=q^i\sf{\p}{\p x^i}$, we have
$$
  \varphi_g|_{\p\Omega} = \left\langle\sf2{n-2}||\nabla u||^2w -
     \sf4{n-2} \langle\nabla u,w\rangle\nabla u +
      \sf{2C}nu^{\frac{2n}{n-2}}w, \nu\right\rangle d\sigma.
$$
Here, we have let $w=w^i\sf{\p}{\p x^i}$ be the translation vector
field induced
on flat conformal space $R=\R^n\cup\{\infty\}$ 
by the right-invariant vector field on $P$ which gives this conservation law.

\subsubsection{Rotation: $g=a^i_jp_ix^j,\ a^i_j+a^j_i=0$.}
\index{rotation}
In this case, we have $v^i=-a^i_jx^j$, $v_i=a^j_ip_j$.  On a transverse
Legendre submanifold of $J^1(\R^n,\R^+)$, we have from
(\ref{LegendreDerivs}) that $p_i = \sf{2}{n-2}u^{-1}u_{x^i}$, and we find that
$$
  \varphi_g = \left(\left(\sf2{n-2}||\nabla u||^2+\sf{2C}n
    u^\frac{2n}{n-2}\right)a^i_jx^j-\sf4{n-2}u_{x^i}u_{x^k}a^k_jx^j
    +\sf2{n-2}uu_{x^j}a^j_i\right)dx_{(i)}.
$$
In this formula, the last term represents a trivial conservation
law\index{conservation law!trivial}; that is,
$d(uu_{x^j}a^j_idx_{(i)})=0$ on any transverse Legendre submanifold, so it
will be ignored below.
Restricting as in the preceding case to the smooth boundary of
$\Omega\subset\R^n$
with unit normal $\nu$ and area element $d\sigma$, this is
$$
  \varphi_g|_{\p\Omega} =
    \left\langle\left(\sf2{n-2}||\nabla u||^2+\sf{2C}nu^\frac{2n}{n-2}\right)a
      -\sf4{n-2}\langle\nabla u,a\rangle\nabla u,\nu\right\rangle
      d\sigma.
$$
Here, we have let $a =a^i_jx^j\sf{\p}{\p x^i}$ be the rotation vector
field induced on flat conformal space $R$ by the right-invariant
vector field on $P$ which gives this conservation law.

\subsubsection{Dilation: $g=1+x^ip_i$.}
\index{dilation}
This generating function gives the right-invariant vector field whose
value at the identity is the Lie algebra element (in blocks of size $1,n,1$)
$$
  \left(\begin{array}{ccc} 2 & 0 & 0 \\ 0 & 0 & 0 \\ 0 & 0 & -2
    \end{array}\right),
$$
which generates a $1$-parameter group of dilations about the origin in
flat conformal space $R$.
In this case, we have $v^i = -x^i$, $v_i=p_i$, and on a transverse
Legendre submanifold with $p_i=\sf{2}{n-2}u^{-1}u_{x^i}$, we find that
\begin{equation}
  \varphi_g = \left(\left(\sf2{n-2}||\nabla u||^2
    + \sf{2C}nu^\frac{2n}{n-2}\right)x^i-\sf4{n-2}u_{x^i}u_{x^j}x^j
    -2uu_{x^i}\right)dx_{(i)}.
\label{GenDil13}
\end{equation}
For this conservation law, it is instructive to take for
$\Omega\subset\R^n$ the open ball of radius $r>0$ centered at the
origin, and then
\begin{equation}
  \varphi_g|_{\p\Omega} =\left(
    r\left(\sf2{n-2}||\nabla u||^2+\sf{2C}nu^\frac{2n}{n-2}
      -\sf4{n-2}\langle\nabla u,\nu\rangle^2\right)-2u
      \langle\nabla u,\nu\rangle\right)d\sigma.
\label{DilCons13}
\end{equation}
A simple consequence of this conservation law is the following uniqueness
theorem.\footnote{See \cite{Pohozaev:Eigenfunctions},
where a non-existence theorem is
  proved for a more general class of equations, for which dilation gives
  an integral identity instead of a conservation law.}
\index{Poho{\v z}aev's theorem}
\begin{Theorem}[Poho{\v z}aev]  
  If $u(x)\in C^2(\bar\Omega)$ is a solution to $\Delta u = 
  Cu^{\sf{n+2}{n-2}}$ in the ball $\Omega$ of radius $r$, 
  with $u\geq 0$ in $\Omega$ and $u=0$ on $\p\Omega$, then $u=0$.
\end{Theorem}
\begin{Proof}
We will first use the conservation law to show that $\nabla u=0$ everywhere
on $\p\Omega$.
If we decompose $\nabla u = u_\tau+u_\nu\nu$ into tangential and
normal components along $\p\Omega$, so that in particular $u_\tau=0$
by hypothesis, then the conserved integrand (\ref{DilCons13}) is
$$
  \varphi_g|_{\p\Omega} = -\sf{2r}{n-2}u_\nu^2d\sigma,
$$
so the conservation law $\int_{\p\Omega}\varphi_g=0$ implies that
$u_\nu=0$ on $\p\Omega$.  

Now with $\nabla u=0$ on $\p\Omega$, we
can compute
\begin{eqnarray*}
  0 & = & \int_{\p\Omega}*du \\ & = & \int_\Omega d*du \\
    & = & \int_\Omega\Delta u\, dx.
\end{eqnarray*}
But it is clear from the PDE that $\Delta u$ cannot change sign,
so it must vanish identically, and this implies that $u=0$ throughout
$\Omega$. 
\end{Proof}

\

\noindent
In fact, looking at the expression (\ref{GenDil13}) for $\varphi_g$ for
a more general region, it is not hard to see that the same proof applies
whenever $\Omega\subset\R^n$ is bounded and star-shaped.

\subsubsection{Inversion: $g=-\sf12p_jb_j||x||^2+b_jx^j
    (1+p_ix^i)$.}
\index{inversion}
This is the generating function for the vector field in
$R=\R^n\cup\{\infty\}$ which is the conjugate of a translation
vector field by inversion in an origin-centered sphere.

In this case, we have $v^i=\sf12b_i||x||^2-b_jx^jx^i$,
$v_i = b_ix^jp_j-b_jx^ip_j+b_jx^jp_i+b_i$, and on a
transverse Legendre submanifold, we find after some tedious calculation
that
\begin{eqnarray*}
  \varphi_g & = & \left[\left((\sf{2}{n-2}||\nabla u||^2+\sf{2C}{n}
       u^{\sf{2n}{n-2}})\delta_{ij}
    -\sf4{n-2}u_{x^i}u_{x^j}\right)
      (b_kx^kx^j-\sf12b_j||x||^2)\right. \\
  & \ &\left.\qquad -2ub_jx^ju_{x^i}+u^2b_i
    \right]dx_{(i)}.
\end{eqnarray*}
Again taking $\Omega\subset\R^n$ to be the open ball of radius $r>0$
centered at the origin, we have
\begin{eqnarray*}
  (n-2)\varphi_g|_{\p\Omega} & = & \left\langle
    (r^2(-4\langle\nabla u,\nu\rangle^2
    +||\nabla u||^2
    +\sf{C}{n}u^{\sf{2n}{n-2}})
    +u^2)b\right. \\
 & \  & \qquad + 2(r^2\langle b,\nabla u\rangle - 
         (n-2)ru\langle b,\nu\rangle)\nabla u,
   \nu\bigg\rangle\,d\sigma,
\end{eqnarray*}
where $b=b_i\sf{\p}{\p x^i}$ is the vector field whose conjugate by a
sphere-inversion is the vector field generating the conservation law.
\index{Poincar\'e-Cartan form!conformally invariant|)}
\index{conservation law!for x@for $\Delta u = Cu^{\frac{n+2}{n-2}}$|)}
\index{Poisson equation!conformally invariant|)}
\index{conformal!structure!flat model|)}

\section{Conservation Laws for Wave Equations}
\index{conservation law!for wave equations|(}

\def\L{\mathbf L}

In this section, we will consider {\em non-linear wave equations}
\index{wave equation!non-linear|(}
\begin{equation}
  \square z = f(z),
\label{NLWave14}
\end{equation}
which are hyperbolic analogs of the non-linear Poisson
equations{\index{Poisson equation!non-linear}
considered previously.  Here, we are working in Minkowski
space\index{Minkowski space|(}
$\L^{n+1}$ with coordinates $(t,y^1,\ldots,y^n)$, and the wave
operator\index{wave operator} is
$$
  \square = -\left(\frac\p{\p t}\right)^2+\sum\left(\frac\p{\p
    y^i}\right)^2.
$$
It is in this hyperbolic case that conservation laws have been
most effectively used.  Everything developed previously in this
chapter for the
Laplace operator and Poisson equations on Riemannian manifolds has an
analog for the wave operator and wave equations on {\em Lorentzian}
manifolds, which by definition carry a metric of signature $(n,1)$.
Indeed, even the coordinate formulae for conservation laws that we
derived in the preceding section are easily altered by a sign change
to give corresponding conservation laws for the wave equation.  Our
goal in this section is to see how certain analytic conclusions can be
drawn from these conservation laws.

Before doing this, we will illustrate the usefulness of understanding
the wave operator geometrically, by presenting a
result of Christodoulou\index{Christodoulou, D.} 
asserting that the Cauchy problem\index{Cauchy problem|(} 
for the non-linear hyperbolic 
equation\index{wave equation!conformally invariant|(}
\begin{equation}
  \square z = z^{\frac{n+3}{n-1}}
\label{InvarWave14}
\end{equation}
has solutions for all time, given sufficiently small initial
data.\footnote{See \cite{Christodoulou:Global};
what is proved there is somewhat more general.}
The proof exploits conformal invariance of the equation in an
interesting way, and this is what we want to explain.
Note that (\ref{InvarWave14}) is the hyperbolic analog of the
maximally symmetric non-linear 
Poisson\index{Poisson equation!conformally invariant}
equation $\Delta z=z^{\frac{n+2}{n-2}}$ considered previously; the
change in exponent reflects the fact that the number of
independent variables is now $n+1$, instead of $n$.  This equation
will be of special interest in our discussion of conservation laws, as
well.

The idea for proving the long-time existence result is to map
Minkowski space $\L^{n+1}$, which is the domain
for the unknown $z$ in (\ref{InvarWave14}), to a {\em bounded} domain,
in such a way that the equation (\ref{InvarWave14}) corresponds to an
equation for which short-time existence of the Cauchy problem is
already known.  With sufficiently small initial data, the
``short-time'' will cover the bounded domain, and back on $\L^{n+1}$
we will have a global solution.

The domain to which we will map $\L^{n+1}$ is actually part of a {\em
  conformal compactification} of Minkowski 
space\index{Minkowski space!compactified}, analogous to the
conformal compactification of Euclidean space constructed in
\S\ref{Subsection:FlatConformal}.  This compactification is diffeomorphic to a
product $S^1\times S^n$, and topologically may be thought of as the
result of adding a point at spatial-infinity for each time, and a
time-at-infinity for each spatial point.  Formally, one can begin with
a vector space with inner-product of signature $(n+1,2)$, and consider
the projectivized null-cone; it is a smooth, real quadric
hypersurface in $\mathbf P^{n+2}$, which in certain homogeneous
coordinates is given by
$$
  \xi_1^2+\xi_2^2 = \eta_1^2 + \cdots + \eta_n^2,
$$
evidently diffeomorphic to $S^1\times S^n$.
The $(n+1,2)$-inner-product induces a
Lorentz metric\index{Lorentz!metric} 
on this hypersurface, well-defined up to
scaling, and its conformal isometry 
group\index{Lorentz!conformal group} has identity component
$SO^o(n+1,2)$, which we will revisit in considering
conservation laws.  What is important for us now is that among the
representative Lorentz metrics for this conformal structure
one finds
$$
  g = -dT^2 + dS^2,
$$
where $T$ is a coordinate on $S^1$, and $dS^2$ is the standard
metric on $S^n$.
In certain spherical coordinates (``usual'' spherical
coordinates applied to $\R^n$, after stereographic 
projection\index{stereographic projection}, this
may be written
$$
  g = -dT^2 + dR^2 + (\sin^2R)\,dZ^2,
$$
where $R\in[0,\pi]$, and $dZ^2$ is the standard metric on the unit
$(n-1)$-sphere.

Now we will conformally embed Minkowski space $\L^{n+1}$ as a bounded
domain in the finite part $\R\times\R^n$ of $S^1\times S^n$, the
latter having coordinates $(T,R,Z)$.  The map $\varphi(t,r,z) =  (T,
R, Z)$ is given by
$$
  \left(\begin{array}{c} T \\ R \\ Z \end{array}\right) =
  \left(\begin{array}{c}
    \arctan(t+r)+\arctan(t-r) \\
      \arctan(t+r)-\arctan(t-r) \\
      z \end{array}\right),
$$
and one can easily check that
$$
  \varphi^*(-dT^2 + dR^2 + (\sin^2R)\,dZ^2) =
    \Omega^2(-dt^2+dr^2+r^2dz^2),
$$
where $dZ^2$ and $dz^2$ are both the standard metric on the unit
$(n-1)$-sphere; the conformal factor is
$$
  \Omega = 2(1+(t+r)^2)^{-\frac12}(1+(t-r)^2)^{-\frac12},
$$
and the right-hand side is a multiple of the flat Minkowski metric.
The image of $\varphi$ is the ``diamond''
$$
  {\mathcal D} = \{(T,R,Z):R-\pi < T < \pi-R,\ R\geq 0\}.
$$
Note that the initial hyperplane $\{t=0\}$ corresponds to $\{T=0\}$,
and that with fixed $(R,Z)$, as $T\to \pi-R$, $t\to\infty$.
Consequently, the long-time Cauchy problem for the invariant wave
equation (\ref{InvarWave14}) corresponds to a short-time Cauchy
problem on the bounded domain ${\mathcal D}$ for some other equation.

We can see what this other equation is without carrying out tedious
calculations by considering the {\em conformally invariant wave
  operator}\index{wave operator!conformal}, 
an analog of the conformal Laplacian\index{conformal!Laplacian} 
discussed in \S\ref{Subsection:ConfLaplacian}.
This is a differential operator
$$
  \square_c : \Gamma(D^\frac{n-1}{2(n+1)})\to
  \Gamma(D^\frac{n+3}{2(n+1)}) 
$$
between certain density line bundles\index{density line bundle} 
over a manifold with Lorentz metric. 
With a choice of Lorentz metric $g$ representing the conformal class,
$u\in\Gamma(D^\frac{n-1}{2(n+1)})$ is represented
a function $u_g$, and the density $\square_c u$ is represented by the
function
$$
  (\square_c u)_g = \square_g (u_g) + \sf{n-1}{4n}R_gu_g,
$$
where $R_g$ is the scalar curvature and $\square_g$ is the wave
operator associated to $g$.  We interpret our wave
equation (\ref{InvarWave14})
as a condition on a density represented
in the flat (Minkowski) metric $g_0$ by the function $u_0$, and the
equation 
transformed by the map $\varphi$ introduced above should express the
same condition represented in the new metric $g$.  The representative
functions are related by
$$
  u_g = \Omega^{-\frac{n-1}{2}}u_0, \qquad
  (\square_cu)_g = \Omega^{-\frac{n+3}{2}}(\square_cu)_0,
$$
so the condition (\ref{InvarWave14}) becomes
\begin{eqnarray*}
  \square_g(u_g) + \sf{n-1}{4n}R_gu_g & = & (\square_cu)_g \\
    & = & \Omega^{-\frac{n+3}{2}}(\square_cu)_0 \\
    & = & \Omega^{-\frac{n+3}{2}}u_0^\frac{n+3}{n-1} \\
    & = & u_g^\frac{n+3}{n-1}.
\end{eqnarray*}
The scalar curvature is just that of the round metric on the
$n$-sphere, $R_g = n(n-1)$, so letting $u=u_0$, $U=u_g$, the equation
(\ref{InvarWave14}) is transformed into
\begin{equation}
  \square_g U + \sf{(n-1)^2}{4}U = U^{\frac{n+3}{n-1}}.
\label{ChristoWave14}
\end{equation}

Finally, suppose given compactly supported initial data $u(0,x) =
u_0(x)$ and $u_t(0,x) = u_1(x)$ for (\ref{InvarWave14}).  These
correspond to initial data 
$U_0(X)$ and $U_1(X)$ for (\ref{ChristoWave14}), supported in the ball
of radius $\pi$.  The standard result on local existence implies that the
latter Cauchy problem\index{Cauchy problem|)} 
can be solved for all $X$, in some time interval
$T\in[0,T_0]$, with a lower-bound for $T_0$ determined by the size of
the initial data.  Therefore, with sufficiently small initial data, we
can arrange $T_0\geq \pi$, and translated back to the original
coordinates, this corresponds to a global solution of
(\ref{InvarWave14}).

\

We now turn to more general wave equations
(\ref{NLWave14}), where conservation laws have been most effectively
used.\footnote{This 
  material and much more may be found in
  \cite{Strauss:Nonlinear}\index{Strauss, W.|nn}.}
Equation (\ref{NLWave14}) is the Euler-Lagrange equation for the
action functional
$$
  \int_\R\int_{\R^n}\left(\sf12(-z_t^2+||\nabla
     z||^2)+F(z)\right)dy\,dt, 
$$
where $F^\prime(z)=f(z)$, and the gradient $\nabla z$ is with respect
to the ``space'' variables $y^1,\ldots,y^n$.

Rather than redevelop the machinery of conformal geometry in the
Lorentz case,
we work in the classical setting, on $J^1(\L^{n+1},\R)$ with
coordinates $t,y^i,z,p_a$ (as usual, $1\leq i\leq n$, $0\leq a\leq n$),
contact form $\theta = dz-p_0dt-p_idy^i$, Lorentz inner-product
$ds^2=-dt^2+\sum(dy^i)^2$, and Lagrangian
$$
  L(t,y,z,p) = \sf12(-p_0^2+\ss|p_i|^2)+F(z).
$$
A normalized representative functional is then
\begin{equation}
  \Lambda = L\,dt\wedge dy+\theta\wedge(-p_0dy-p_idt\wedge dy_{(i)}),
\label{WaveLagrangian14}
\end{equation}
satisfying 
$$
  d\Lambda = \Pi =\theta\wedge
    (dp_0\wedge dy+dp_i\wedge dt\wedge dy_{(i)} + f(z)dt\wedge dy).
$$
This is the example discussed at the end of
\S\ref{Section:Noether}.  As mentioned there, the invariance of
the equation under time-translation gives an important conservation
law, and its uses will be our first topic below.  In fact, there are
conservation laws associated to 
space-translations and Lorentz rotations, the latter generated by
ordinary spatial rotations $b^i_jy^j\frac{\partial}{\partial y^j}$
($b^i_j+b^j_i=0$) and {\em Lorentz boosts}\index{Lorentz boost}
$y^i\frac{\partial}{\partial t}+t\frac{\partial}{\partial y^i}$;
however, these seem to have been used less widely in the analysis of
(\ref{NLWave14}).

Especially interesting is the case of (\ref{InvarWave14}), which is
preserved under a certain action of the {\em conformal Lorentz
  group}\index{Lorentz!conformal group}
$SO^0(n+1,2)$ on $J^1(\L^{n+1},\R)$.  In particular, there are
extensions to $J^1(\L^{n+1},\R)$ of the dilation\index{dilation} and
inversion\index{inversion} vector
fields on $\L^{n+1}$, and these give rise to more conservation laws.
We will consider these after discussing uses of the time-translation
conservation law for the more general wave equations (\ref{NLWave14}).

\subsection{Energy Density}

The time-translation vector field $\frac{\partial}{\partial t}$ on
$\mathbf L^{n+1}$ lifts to $J^1(\L^{n+1},\R)$ to a symmetry of
$\Lambda$ having the same expression, $V=\frac{\partial}{\partial t}$.
The Noether prescription gives
$$
  \varphi_t = V\innerprod\Lambda =
     \left(\sf12(p_0^2+\textstyle\sum p_i^2)+F(z)\right)dy+
        p_0p_idt\wedge dy_{(i)},
$$
as calculated in Section~\ref{Section:Noether}.
The coefficient of $dy$ here is the {\it energy density}
\index{energy|(}
$$
  e \stackrel{\mathit{def}}{=} \sf12(p_0^2+|p_i|^2) + F(z),
$$
and it appears whenever we integrate $\varphi_t$ along
a constant-time level surface $\R^n_t = \{t\}\times\R^n$.  The
energy function
$$
  E(t) \stackrel{\mathit{def}}{=} \int_{\R^n_t}e\,dy\geq 0,
$$
is constant by virtue of the equation (\ref{NLWave14}), assuming
sufficient decay of $z$ and its derivatives in the space variables for
the integral to make sense.

A more substantial application involves a region $\Omega\subset\L^{n+1}$ of the form
$$
  \Omega = \bigcup_{t\in(t_0,t_1)}\{||y||<r_0-(t-t_0)\},
$$
a union of open balls in space, with initial radius $r_0$ decreasing with
speed $1$.
The boundary $\p\Omega$ is $T-B+K$, where
\begin{itemize}
\item
$B=\{t_0\}\times\{||y||\leq r_0\}$ is the initial disc, 
\item
$T=\{t_1\}\times\{||y||\leq r_0-(t_1-t_0)\}$ is the final disc, and
\item
$K=\cup_{t\in[t_0,t_1]}\{||y||=r_0-(t-t_0)\}$ is part of a null
cone. 
\end{itemize}
The conservation of $\varphi_t$ on $\p\Omega$ reads
\begin{equation}
  0 = \int_{\p\Omega}\varphi_t = \int_{T-B}e\,dy+\int_K\varphi_t.
\label{CLTrans14}
\end{equation}
The term $\int_K\varphi_t$ describes the flow of energy across part of
the null cone; we will compute the integrand more explicitly in terms
of the area form $dK$ induced from an ambient {\em Euclidean}
metric $dt^2+\sum (dy^i)^2$, with the goal of showing that
$\int_K\varphi_t\geq 0$.  This area form is the contraction of
the outward unit normal $\frac{1}{\sqrt{2}}(\frac{\p}{\p t}+
\frac{y^i}{||y||}\frac{\p}{\p y^i})$ with the ambient Euclidean volume
form $dt\wedge dy$, giving
$$
  dK = \sf{1}{\sqrt{2}}(dy-\sf{y^i}{||y||}dt\wedge dy_{(i)})|_K
    = \sqrt{2}\,dy|_K.
$$
It is easy to calculate that the restriction to $K$ of $\varphi_t$ is
$$
  \varphi_t|_K  =  \sf{1}{\sqrt{2}}(e-\sf{y^ip_i}{||y||}p_0)dK.
$$
Separating the radial and tangential space derivatives
$$
  p_r \stackrel{\mathit{def}}{=} \frac{y^ip_i}{||y||},\quad 
  p_\tau \stackrel{\mathit{def}}{=} \sqrt{\sum p_i^2-p_r^2},
$$
we can rewrite this as
\begin{eqnarray}
  \varphi_t|_K  & = & \sf{1}{\sqrt{2}}(\sf12(p_0^2+p_r^2+p_\tau^2)
                        +F(z)-p_rp_0)dK \\
    & = & \sf{1}{\sqrt{2}}(\sf12(p_0-p_r)^2+\sf12p_\tau^2)+F(z))
                    dK.
\label{SideEnergy}
\end{eqnarray} 
In the region of $J^1(\L^{n+1},\R)$ where $F(z)\geq 0$ this integrand
is positive, and from (\ref{CLTrans14}) we obtain the bound
$$
  \int_Te\,dy\leq\int_Be\,dy.
$$
This says that ``energy travels with at most
unit speed''---no more energy can end up in $T$ than was already
present in $B$.  If $F(z)\geq 0$ everywhere, then one can obtain
another consequence of the expression (\ref{SideEnergy}) by writing
$$
  \int_K\varphi_t = \int_B\varphi_t-\int_T\varphi_t \leq
    \int_B\varphi_t\leq E(t_0) = E.
$$
This gives an upper bound for
$$
  ||(dz|_K)||^2_{L^2} = \int_K((p_0-p_r)^2+p_\tau^2)
     dK \leq 2\sqrt{2}E
$$
which holds for the entire backward null cone; that is, our bound is
independent of $t_0$.  Here the $L^2$-norm is with respect
to Euclidean measure.

We should also mention that the spatial-translation and Lorentz
rotations give rise to conserved quantities that may be thought of as
linear and angular momenta, respectively.  The uses of these are
similar to, though not as extensive as, the uses of the conserved
energy.
\index{energy|)}

\subsection{The Conformally Invariant Wave Equation}

\index{conservation law!for wave equations!for $\square z = Cz^{\frac{n+3}{n-1}}$|(}
We now determine some additional conservation laws for the conformally
invariant wave equation (\ref{InvarWave14}).  Again, we could
duplicate the process used for Poisson equations by calculating
restrictions of the right-invariant vector fields of $SO^o(n+1,2)$ to
the image of an embedding $J^1(\L^{n+1},\R^+)\hookrightarrow
SO^o(n+1,2)$, and contracting with the left-invariant Lagrangian.
Instead, we will illustrate the more concrete, coordinate-based
approach, though we will still make some use of the geometry.

\subsubsection{The Dilation Conservation Law}
\index{dilation}

To find the conservation law corresponding to dilation symmetry
of (\ref{InvarWave14}), we have to first determine a formula for this
symmetry on $J^1(\L^{n+1},\R^+)$, and then apply the Noether
prescription.  For this, we will first determine the vector field's
action on $J^0(\L^{n+1},\R^+)$; the lift of this action to
$J^1(\L^{n+1},\R^+)$ is determined by the requirement that it preserve
the contact line bundle\index{contact!line bundle}.

By analogy with (\ref{FinalBasis11}), we have an embedding of
$J^0(\L^{n+1},\R^+)$ into the null-cone of $\L^{n+1,2}$ given by
$$
  (t,y^i,z)\mapsto z^{\frac{2}{n-1}}\left(\begin{array}{c}
    1 \\ t \\ y^i \\ \sf12(-t^2+|y|^2)\end{array}\right),
$$
and the dilation matrix (in blocks of size $1,1,n,1$) acts
projectively on this slice of the null-cone by
$$
   \left[z_r^{\frac{2}{n-1}}\left(\begin{array}{c}
      1 \\ t_r \\ y^i_r \\ \sf12||(t_r,y_r)||^2)\end{array}\right)
     \right] \stackrel{\mathit{def}}{=}
  \left(\begin{array}{cccc} r^{-1} & 0 & 0 & 0 \\ 0 & 1 & 0 & 0 \\
     0 & 0 & I_n & 0 \\ 0 & 0 & 0 & r\end{array}\right)
   \cdot\left[z^{\frac{2}{n-1}}\left(\begin{array}{c}
      1 \\ t \\ y^i \\ \sf12||(t,y)||^2\end{array}\right)\right].
$$
Taking the derivative with respect to $r$ and setting $r=1$ gives the vector
field
$$
  \bar{V}_{\mathit{dil}} = 
    -\sf{n-1}2z\sf{\partial}{\partial z}+t\sf{\partial}
      {\partial t}+y^i\sf{\partial}{\partial y^i}.
$$
The scaling in the $z$-coordinate reflects an interpretation of the
unknown $z(t,y)$ as a section of a certain density line
bundle\index{density line bundle}.  We
then find the lift from $\bar{V}_{\mathit{dil}}\in\mathcal
V(J^0(\L^{n+1},\R^+))$ to $V_{\mathit{dil}}\in\mathcal
V(J^1(\L^{n+1},\R^+))$ by the requirement that the contact
form\index{contact!form}
$\theta = dz - p_0dt - p_idy^i$ be preserved up to scaling; that is,
$$
  V_{\mathit{dil}} = \bar{V}_{\mathit{dil}} + v^0p_0 + v^ip_i
$$
must satisfy
$$
  \mathcal L_{V_{\mathit{dil}}}\theta \equiv
    0 \pmod{\{\theta\}},
$$
where $v^0$ and $v^i$ are the unknown coefficients of the lift.
This simple calculation yields
$$
  V_{dil} = -\sf{n-1}2z\sf{\partial}{\partial z}+t\sf{\partial}
    {\partial t}+y^i\sf{\partial}{\partial y^i}-\sf{n+1}2
      p_a\sf{\partial}{\partial p_a}.
$$
Then one can compute even for the general wave equation (\ref{NLWave14}) that
$$
  \mathcal L_{V_{dil}}(L(t,y,z,p)dt\wedge dy) = 
    ((n+1)F(z)-\sf{n-1}{2}zf(z))dt\wedge dy.
$$
Tentatively following the Noether prescription for the general wave
equation, we set
$$
  \varphi_{dil} =V_{dil}\innerprod\Lambda,
$$
where $\Lambda$ is given by (\ref{WaveLagrangian14}),
and because $\Lambda\equiv L\,dt\wedge dy$ modulo $\{I\}$, we can
calculate
\begin{eqnarray*}
  d\varphi_{\mathit{dil}} & = & \mathcal L_{V_\mathit{dil}}\Lambda - 
    V_\mathit{dil}\innerprod\Pi \\
   & \equiv & \mathcal L_{V_\mathit{dil}}(L\,dt\wedge dy) - 0
      \pmod{\mathcal E_\Lambda} \\
   & = &  ((n+1)F(z)-\sf{n-1}{2}zf(z))dt\wedge dy
      \pmod{\mathcal E_\Lambda}.
\end{eqnarray*}
The condition on the equation (\ref{NLWave14}) that $\varphi_{dil}$ be a
conservation law is therefore
$$
   F(z) = Cz^{\frac{2(n+1)}{n-1}},
$$
so the PDE is
$$
  \square z = C^\prime z^{\frac{n+3}{n-1}},
$$
as expected (cf.~(\ref{InvarWave14})); we work with $C^\prime = 1$, $C
= \frac{n-1}{2(n+1)}$.
Now, one can calculate that restricted to any Legendre submanifold
the conserved density is
$$
  \varphi_\mathit{dil} = L(t,y,z,p)(t\,dy-y^idt\wedge dy_{(i)})
     +(\sf{n-1}{2}z+tp_0+y^ip_i)(p_0dy + p_jdt\wedge dy_{(j)}).
$$
Typically, one considers the restriction of this form to the
constant-time hyperplanes $\R^n_t = \{t\}\times\R^n$, which is
$$
  \varphi_{dil} \equiv (te+rp_0p_r+\sf{n-1}2zp_0)dy\pmod
    {\{dt\}}.
$$
For example, we find that for solutions to (\ref{InvarWave14}) with compact
support in $y^i$,
\begin{equation}
    \frac{d}{dt}\int_{\R_t^n}(te+rp_0p_r+\sf{n-1}2zp_0)dy = 0.
\label{DilCL14}
\end{equation}

For more general wave equations (\ref{NLWave14}), an identity like
(\ref{DilCL14}) holds, but with a non-zero right-hand side; our
conservation law is a 
special case of this.  The general dilation identity is of
considerable use in the analysis of non-linear wave equations.  It
is analogous to the ``almost-conservation 
law''\index{conservation law!almost@``almost''} derived
from scaling symmetry used to obtain lower
bounds on the area growth of minimal surfaces, as discussed in
\S\ref{Subsection:MinimalCLs}. 

\subsubsection{An Inversion Conservation Law}
\index{inversion}

We now consider the inversion symmetry corresponding to the conjugate
of time-translation by inversion in a unit (Minkowski) ``sphere''.  We
will follow the same procedure as for dilation symmetry, first
determining a vector field on $J^0(\L^{n+1},\R^+)$ generating this
inversion symmetry, then lifting it to a
contact-preserving\index{contact!line bundle} vector
field on $J^1(\L^{n+1},\R^+)$, and then applying the Noether
prescription to obtain the conserved density.

The conjugate by sphere-inversion of a time-translation in
$SO^o(n+1,2)$ is the matrix
$$
  \left(\begin{array}{cccc}  1 & b & 0 & \sf12b^2 \\
    0 & 1 & 0 & -b \\ 0 & 0 & I_n & 0 \\ 0 & 0 & 0 & 1
    \end{array}\right),
$$
and differentiating its projective linear action on
$$
  \left[z^\frac{2}{n-1}\left(\begin{array}{c}1 \\ t \\ y^i \\
    \sf12||(t,y)||^2\end{array}\right)\right]
$$
yields the vector field
$$
  \bar V_{\mathit{inv}} = \sf{n-1}{2}tz\sf{\p}{\p z}
    -\sf12(t^2+|y|^2)\sf{\p}{\p t} - ty^i\sf{\p}{\p y^i}.
$$
Again, the coefficients of $\frac{\p}{\p t}$ and $\frac{\p}{\p y^i}$
describe an infinitesimal conformal motion\index{Lorentz!conformal group} 
of Minkowski space---representing an element of the Lie algebra of the
Lorentz conformal group---and the coefficient of
$\frac{\p}{\p z}$ gives the induced action on a density line
bundle\index{density line bundle}.
We now look for coefficients $v_0\frac{\p}{\p p_0}+v_i\frac{\p}{\p
  p_i}$ to add to $\bar V_{\mathit{inv}}$ so that the new vector field
will preserve $\theta$ up to scaling, and the unique solution is
$$
  V_{\mathit{inv}} = \bar V_{\mathit{inv}} + 
   \left(\sf{n-1}{2}z + p_iy^i + \sf{n+1}{2}tp_0\right)\sf{\p}{\p p_0}
   + \left(p_0y^i + \sf{n+1}{2}tp_i\right)\sf{\p}{\p p_i}.
$$
In applying the Noether prescription to $\Lambda$ and
$V_{\mathit{inv}}$, it will turn out that we need a compensating term,
because $\mathcal L_{V_{\mathit{inv}}}\Lambda\not\equiv 0$ modulo
${\mathcal E}_\Lambda$.  However, instead of performing this tedious
calculation, we can simply test
\begin{eqnarray*}
  \tilde\varphi_{\mathit{inv}} & \stackrel{\mathit{def}}{=} &
    -V_{\mathit{inv}}\innerprod\Lambda \\
   & \equiv & L(t,y,z,p)(\sf12(t^2+|y|^2)dy -
   ty^idt\wedge dy_{(i)}) \\ 
  & & \quad +(\sf{n-1}{2}tz + \sf12p_0(t^2+|y|^2) + tp_iy^i)\wedge
    (p_0dy + p_jdt\wedge dy_{(j)}),
\end{eqnarray*}
We find that on solutions to (\ref{InvarWave14}),
$$
  d\tilde\varphi_{\mathit{inv}} = \sf{n-1}{2}p_0z\,dt\wedge dy
    = d(\sf{n-1}{4}z^2dy),
$$
and we therefore set
$$
  \varphi_{\mathit{inv}} = \tilde\varphi_{\mathit{inv}} -
    \sf{n-1}{4}z^2dy.
$$
As usual, we consider the restriction of this form to a hyperplane
$\R^n_t = \{t\}\times\R^n$, which gives
$$
  \varphi_{\mathit{inv}} \equiv
     \left(\sf12(L+p_0^2)(t^2+|y|^2) + \sf{n-1}{2}p_0tz +
     tp_0p_iy^i - \sf{n-1}{4}z^2\right)dy,
$$
modulo $\{dt, \theta\}$.

Again, for more general wave equations (\ref{NLWave14}),
this quantity gives not a conservation law, but a useful integral
identity.  The usefulness of the integrand follows largely from the fact
that after adding an exact $n$-form, the coefficient of $dy$ is
positive.  One notices this by expanding in terms of radial and
tangential derivatives
\begin{equation}
\begin{split}
  \varphi_{\mathit{inv}} &= 
    \left\{\frac14\left((p_0^2+p_\tau^2+p_r^2)(t^2+r^2) + 2(n-1)p_0tz
        + 4trp_0p_r - (n-1)z^2\right)\right. \\
   &\qquad \qquad  +\left.\frac12(t^2+r^2)F(z)\right\}dy, \\
\end{split}
\label{InvCL14}
\end{equation}
which suggests completing squares:
\begin{equation}
\begin{split}
  \varphi_{\mathit{inv}} &=
   \left\{ \frac14\left(|p_0y+tp|^2 + (tp_0+rp_r+(n-1)z)^2 +
    (rp_\tau)^2\right) \right. \\ & \qquad 
   + \frac12(t^2+r^2)F(z) - \left.
    \frac{n-1}{4}(nz^2+2rzp_r)\right\} dy. \\
\end{split}
\end{equation}
The last term is the divergence $d(\frac{n-1}{4}y^iz^2dy_{(i)})$ modulo
$\{dt,\theta\}$, and the remaining terms are positive.  The positive
expression 
$$
  e_c\stackrel{\mathit{def}}{=} \sf12\left(|p_0y+tp|^2 +
     (tp_0+rp_r+(n-1)z)^2 + (rp_\tau)^2\right) + (t^2+r^2)F(z)
$$
is sometimes called the {\em conformal energy}\index{conformal energy}
of the solution $z(t,y)$.

It is also sometimes convenient to set
$$
  e_d = rp_0p_r + \sf{n-1}{2}zp_0
$$
so that
$$
  \varphi_{\mathit{dil}} \equiv (e_d+te)dy \pmod{\{dt\}},
$$
and also
$$
  \varphi_{\mathit{inv}} \equiv
  (te_d+\sf12(t^2+|y|^2)e-\sf{n-1}{4}z^2)dy \pmod{\{dt\}}.
$$
The fact that the integrals of these quantities are constant in $t$
yields results about growth of solutions.  

In fact, in the analysis of wave equations
that are perturbations of the conformally invariant equation
(\ref{InvarWave14}), 
the most effective estimates pertain to the quantities appearing in
these conservation laws.  One can think of the conservation laws as
holding for our ``flat'' non-linear wave equation (\ref{InvarWave14}),
and then the estimates are their analogs in the ``curved'' setting.
\index{conservation law!for wave equations!for $\square z = Cz^{\frac{n+3}{n-1}}$|)}
\index{conservation law!for wave equations|)}
\index{Minkowski space|)}
 
\subsection{Energy in Three Space Dimensions}
\index{energy|(}

We conclude this section by discussing a few more properties that
involve the energy in three space dimensions.

The fact that the energy $E(t)$ is constant implies in particular that
$\int_{\R^n_t}z_tdy$ is bounded with respect to $t$.  This allows
us to consider the evolution of the spatial $L^2$-norm of a solution
to $\square z = z^3$ as follows:
\begin{eqnarray*}
  \frac{d^2}{dt^2}\left(\int_{\R^3_t}z^2dy\right)
   & = & 2\int_{\R^3_t}(z\,z_{tt}dy + z_t^2)dy \\
   & = & 2\int_{\R^3_t}(z\Delta z - z^4 + z_t^2)dy \\
   & = & -2\int_{\R^3_t}(|\nabla_yz|^2+z^4)dy 
      +2\int_{\R^n_t}z_t^2dy \\
   & \leq & 4E,
\end{eqnarray*}
where the second equality follows from the differential equation and
the third from Green's theorem\index{Green's theorem|(} 
(integration by parts).  The conclusion is that
$||z||^2_{L^2(\R^n_t)}$ grows at most quadratically,
\begin{equation}
  \int_{\R^3_t}z^2dy\leq 2Et^2 + C_1t + C_2,
\label{L2Decay14}
\end{equation}
and in particular $\int z^2dy = O(t^2)$.

\

The energy plays another interesting role in the equation
\begin{equation}
  \square z = -z^3
\label{BadEqn14}
\end{equation}
Here, it is possible for the energy to be negative:
\begin{equation}
  E = \int_{\R_t^3}\left(\sf12(z_t^2+||\nabla z||^2)-\sf{z^4}{4}\right)
    dy.
\label{BadEnergy14}
\end{equation}
We will prove that a solution to (\ref{BadEqn14}), with compactly
supported initial data $z(0,y)$,
$z_t(0,y)$ satisfying $E<0$, must blow up in
finite time (\index{Levine, H.|nn}\cite{Levine:Instability}). 
Notice that any non-trivial compactly supported initial
data may be scaled up to achieve $E<0$, and may be scaled down to
achieve $E>0$.

The idea is to show that the quantity
$$
  I(t) \stackrel{\mathit{def}}{=}\int_{\R^3_t}\sf12z^2dy
$$
becomes unbounded as $t\nearrow T$ for some finite time $T>0$.
We start by computing its derivatives
\begin{eqnarray*}
  I^\prime(t) & = & \int zz_tdy,  \\
  I^{\prime\prime}(t) & = & \int z_t^2dy +
         \int zz_{tt}dy \\
    & = & \int z_t^2dy - \int|\nabla z|^2dy + \int z^4dy.
\end{eqnarray*}
The last step uses Green's theorem, requiring the solution to have
compact $y$-support for each $t\geq 0$.  To dispose of the $\int z^4$
term, we add $4E$ to each side using (\ref{BadEnergy14}):
$$
  I^{\prime\prime}(t) + 4E =
    3\int z_t^2dy + \int|\nabla z|^2dy.
$$
We can discard from the right-hand side the positive gradient term,
and from the left-hand side the negative energy term, to obtain
$$
  I^{\prime\prime}(t) > 3\int z_t^2dy.
$$
To obtain a second-order differential inequality for $I$,
we multiply the last inequality by $I(t)$ to obtain
\begin{eqnarray*}
  I(t)I^{\prime\prime}(t) & > &
    \frac32\left(\int z^2dy\right)\left(\int z_t^2dy\right) \\
  & \geq & \sf32 I^\prime(t)^2.
\end{eqnarray*}
The last step follows from the Cauchy-Schwarz
inequality\index{Cauchy-Schwarz inequality|(}, and says
that $I(t)^{-1/2}$ has negative
second derivative.  We would like to use this to conclude that
$I^{-1/2}$ vanishes for some $T>0$ (which would imply that $I$
blows up), but for this we would need to know that
$(I^{-1/2})^\prime(0)<0$, or equivalently $I^\prime(0)>0$, which may
not hold. 

To rectify this, we shift $I$ to
$$
  J(t) = I(t)-\sf12 E(t+\tau)^2,
$$
with $\tau>0$ chosen so that $J^\prime(0)>0$.  We now mimic the
previous reasoning to show that $(J^{-1/2})^{\prime\prime}(t)<0$.
We have
\begin{eqnarray*}
  J^\prime(t) & = & \int zz_tdy - E(t+\tau), \\
  J^{\prime\prime}(t) & = & \int z_t^2dy + \int zz_{tt}dy - E \\
    & = & 3\int z_t^2dy + \int ||\nabla z||^2dy - 5E \\
    & > & 3\left(\int z_t^2dy - E\right).
\end{eqnarray*}
From this we obtain
$$
  J(t)J^{\prime\prime}(t) - \sf32 J^\prime(t)^2  > \textstyle
   \frac32\left[\left(\int z^2 - E(t+\tau)^2\right)
              \left(\int z_t^2 - E\right) - 
              \left(\int zz_t - E(t+\tau)\right)^2\right],
$$
which is positive, again by the Cauchy-Schwarz
inequality\index{Cauchy-Schwarz inequality|)}.  This means that
$(J^{-1/2})^{\prime\prime}(t)<0$.  Along with $J^{-1/2}(0)>0$ and
$(J^{-1/2})^\prime(0) <0$, this implies that for some $T>0$,
$J^{-1/2}(T)=0$, so $J(t)$ blows up.
\index{energy|)}

\

We conclude by noting that the qualitative behavior of solutions of
$\square z = f(z)$ depends quite sensitively on the choice of
non-linear term $f(z)$.  In contrast to the results for $\square z =
\pm z^3$ described above, we have for the equation
$$
  \square z = -z^2 \qquad (n=3),
$$
that every solution must blow up in finite time 
(\index{John, F.}\cite{John:Blowup}).
We will outline the proof in case the initial data are
compactly supported and satisfy
$$
  \int u(0,t)dy > 0,\quad
  \int u_t(0,t)dy > 0.
$$
Note that replacing $z$ by $-z$ gives the equation
$\square z = z^2$, which therefore has the same behavior.  

This proof is fairly similar to the previous one; we will derive
differential inequalities for
$$
  J(t) \stackrel{\mathit{def}}{=} \int_{\R^n_t}z\,dy
$$
which imply that this quantity blows up.  We use integration by
parts\index{Green's theorem|)} to obtain
$$
  J^{\prime\prime}(t) = \int z_{tt}dy =
    \int (\Delta z + z^2)dy = \int z^2,
$$
and using H\"older's inequality\index{H\"older's inequality} 
on $\mbox{Supp }z\subset\{|y|\leq
R_0+t\}$ in the form $||z||_{L^1}\leq||z||_{L^2}||1||_{L^2}$, this gives
\begin{equation}
  J^{\prime\prime}(t) \geq C\left(\int z\,dy\right)^2(R_0+t)^{-3}
    \geq C(1+t)^{-3}J(t)^2.
\label{BadSingI14}
\end{equation}
This is the first ingredient.  

Next, we use the fact that if
$z_0(y,t)$ is the {\em free solution} to the homogeneous wave equation
$\square z_0 = 0$, with the same initial data as our $z$, then
$$
  z(y,t) \geq z_0(y,t)
$$
for $t\geq 0$; this follows from a certain explicit integral
expression for the solution.  Note that if we set $J_0(t)  =
\int_{\R^n_t}z_0$, then it follows from the equation alone that 
$J_0^{\prime\prime}(t) = 0$,
and by the hypotheses on our initial data we have
$J_0(t) = C_0 + C_1t$ with $C_0,C_1 > 0$.  Another property of the
free solution is that its support at time $t$ lies in the annulus
$A_t=\{t-R_0\leq |y|\leq t+R_0\}$.
Now 
\begin{eqnarray*}
  C_0 + C_1t & \leq & \int_{A_t} z\,dy \\
    & \leq & ||z||_{L^1(A)} \\
    & \leq & ||z||_{L^2(A)}||1||_{L^2(A)} \\
    & \leq & C(1+t)\left(\int z^2dy\right)^2.
\end{eqnarray*}
This gives
$$
  J^{\prime\prime}(t) = \int z^2dy \geq
   \left(\frac{C_0+C_1t}{C(1+t)}\right)^{1/2},
$$
and in particular, $J^{\prime\prime} > 0$.  With the assumptions on
the initial data, this gives
\begin{eqnarray*}
  J^{\prime}(t) & > & 0 \label{BadSingII14} \\
  J(t) & \geq & C(1+t)^2 \label{BadSingIII14}.
\end{eqnarray*}

We can use (\ref{BadSingI14}, \ref{BadSingII14}, \ref{BadSingIII14})
to conclude that $J$ must blow up at some finite time.
This follows by writing
\begin{eqnarray*}
  J^{\prime\prime}(t) & \geq & C(1+t)^{-3}J^{3/2}J^{1/2} \\
    & \geq & C(1+t)^{-2}J^{3/2}.
\end{eqnarray*}
Multiply by $J^\prime$ and integrate to obtain
$$
  J^{\prime}(t) \geq C(1+t)^{-1}J(t)^{5/4}.
$$
Integrating again, we have
$$
  J(t) \geq \left(J(0)^{-1/4} - \sf14C\ln(1+t)\right)^{-4},
$$
and because $J(0) > 0$ and $C > 0$, $J(t)$ must blow up in finite time.
\index{wave equation!non-linear|)}
\index{wave equation!conformally invariant|)}

\chapter{Additional Topics}
\label{Chapter:Additional}

\def\sfddt{\sf{\partial}{\partial t}}
\def\fddt{\frac{\partial}{\partial t}}

\section{The Second Variation}
\label{Section:SecondVarn}
\index{second variation|(}

In this section, we will discuss the second variation of the Lagrangian
functionals considered in the preceding chapters.  We begin by giving
an invariant, coordinate-free calculation of the formula
(\ref{GoodSecondVarn}) for the
second derivative of a functional under fixed-boundary variations.
This formula has an interpretation in terms of conformal
structures\index{conformal!structure} 
induced on integral manifolds of the Euler-Lagrange
system\index{Euler-Lagrange!system}, 
which we will describe.  The role played by conformal geometry
here is not to be confused with that in the previous chapter, although
both situations seem to reflect the increasing importance of variational 
equations in conformal geometry.  

The usual integration by parts that
one uses to establish local minimality of a solution to the Euler-Lagrange
equations cannot generally be done in an invariant manner, and we
discuss a condition under which this difficulty can be overcome.  
We considered in \S\ref{Section:PrescribedH} the example of
prescribed mean curvature\index{mean curvature!prescribed} 
hypersurfaces in Euclidean space; we will
give an invariant calculation of the second variation formula and the
integration by parts for this example.  We conclude by discussing
various classical conditions under which an integral manifold of an
Euler-Lagrange system is locally minimizing, using the
Poincar\'e-Cartan form
\index{Poincar\'e-Cartan form|(} 
to express and prove some of these results in a
coordinate-free manner. 

\subsection{A Formula for the Second Variation}
\label{Subsection:Formula}

We start by reconsidering the situation of
\S\ref{Section:EulerLagrange}, in which we calculated the first
variation\index{first variation|(}
of a Lagrangian $\Lambda$ on a contact manifold.  This
amounted to taking the first derivative, at some fixed time, of the
values of the functional
$\mathcal F_\Lambda$ on a $1$-parameter family of Legendre
submanifolds\index{Legendre submanifold|(}. 
Our goal is to extend the calculation to give the second variation of
$\Lambda$, or equivalently, the second derivative of $\mathcal
F_\Lambda$ on a $1$-parameter family at a Legendre submanifold for
which the first variation vanishes.  The result appears in formula
(\ref{SecondVar}) below, and in a more geometric form in
(\ref{GoodSecondVarn}).  This process is formally analogous to
computing the Hessian matrix of a smooth function $f:\R^n\to\R$ at a
critical point, which is typically done with the goal of identifying
local extrema.

Let $(M,I)$ be a contact manifold, with contact form
$\theta\in\Gamma(I)$, and Lagrangian $\Lambda\in\Omega^n(M)$
normalized so that the Poincar\'e-Cartan form is given by
$\Pi=d\Lambda=\theta\wedge\Psi$.  Fix a compact
manifold $N^n$ with boundary $\partial N$, and a smooth map
$$
  F:N\times[0,1]\to M
$$
which is a Legendre submanifold $F_t$ for each fixed $t\in[0,1]$ and is
independent of $t$ on $\partial N\times[0,1]$.   Two
observations will be important:
\begin{itemize}
\item
  $F^*\theta = G\,dt$ for some function $G$ on $N\times[0,1]$,
  depending on the choice of generator $\theta\in\Gamma(I)$.
  This holds because each $F_t$ is a Legendre submanifold, meaning that
  $F_t^*\theta=0$. 
\item
  For every form $\alpha\in\Omega^*(M)$, and every boundary point
  $p\in\partial N$, we have
$$
  \left(\sfddt\innerprod F^*\alpha\right)
    (p,t)=0.
$$
  This is equivalent to the fixed-boundary condition; at each
  $p\in\partial N$, we have $F_*(\fddt)=0$.
\end{itemize}
We previously calculated the first variation\index{first variation|)} as 
(see \S\ref{Subsection:EulerLagrange}) 
$$
  \frac{d}{dt}\left(\int_{N_t}F_t^*\Lambda\right) =
     \int_{N_t}G\cdot F_t^*\Psi,
$$
where $\Pi=d\Lambda=\theta\wedge\Psi$ is the Poincar\'e-Cartan form for
$\Lambda$.  This holds for each $t\in[0,1]$.

We now assume that $F_0$ is stationary\index{stationary|(} 
for $\Lambda$; that is, $F$ is a
Legendre variation\index{Legendre variation|(} 
of an integral manifold $F_0:N\hookrightarrow M$ of the
Euler-Lagrange system\index{Euler-Lagrange!system|(} 
$\mathcal E_\Lambda = \{\theta,d\theta,\Psi\}$.
This is the situation in which we want to calculate the
second derivative:
\begin{eqnarray*}
  \delta^2(\mathcal F_\Lambda)_{N_0}(g) & = &
  \left.\frac{d^2}{dt^2}\right|_{t=0}\left(\int_{N_t}F_t^*\Lambda\right) 
   \\ & = &
    \left.\frac{d}{dt}\right|_{t=0}\int_{N_t}G\,F_t^*\Psi \\
  & = & \int_{N_0}\mathcal L_{\fddt}
          (G\,F^*\Psi) \\
  & = & \int_{N_0}g\,\mathcal L_{\fddt}
            (F^*\Psi),
\end{eqnarray*}
where $g=G|_{t=0}$, and
the last step uses the fact that $F^*_0\Psi = 0$.

To better understand the Lie derivative $\mathcal
L_{\fddt}(F^*\Psi)$, we use the results obtained via the
equivalence method\index{equivalence method} 
in \S\ref{Section:Bigequiv}.  This means that we are
restricting our attention to the case $n\geq 3$, with a Poincar\'e-Cartan
form that is neo-classical and 
definite.\index{Poincar\'e-Cartan form!neo-classical!definite}  
In this situation, we have a $G$-structure $B\to M$, where $G\subset
GL(2n+1,\R)$ has Lie algebra
\begin{equation}
  \lie{g} = \left\{\left(\begin{array}{ccc} (n-2)r & 0 & 0 \\
    0 & -2r\delta^i_j+a^i_j & 0 \\ d_i & s_{ij} & nr\delta^j_i-a^j_i
    \end{array}\right):a^i_j+a^j_i= s_{ij}-s_{ji}=s_{ii}=0\right\}.
\label{LieAlgebra15}
\end{equation} 
The sections of $B\to M$ are local
coframings $(\theta,\omega^i,\pi_i)$ of $M$ for which:
\begin{itemize}
\item
  $\theta$ generates the contact line bundle $I$,
\item
  the Poincar\'e-Cartan form is $\Pi = -\theta\wedge\pi_i\wedge
     \omega_{(i)}$,
\item
  there exists a $\lie{g}$-valued $1$-form
\begin{equation}
  \varphi = \left(\begin{array}{ccc}
    (n-2)\rho & 0 & 0 \\ 0 & -2\rho\delta^i_j+\alpha^i_j & 0 \\
      \delta_i & \sigma_{ij} & n\rho\delta^j_i-\alpha^j_i
  \end{array}\right),
\label{StreqnA15}
\end{equation}
satisfying a structure equation
\begin{equation}
  d\left(\begin{array}{c}\theta \\ \omega^i \\ \pi_i\end{array}\right)
  =-\varphi\wedge\left(\begin{array}{c}\theta \\ \omega^j \\
       \pi_j\end{array}\right)+\left(\begin{array}{c}
     -\pi_i\wedge\omega^i \\ \Omega^i \\ 0\end{array}\right)
\label{StreqnB15}
\end{equation}
where
\begin{equation}
  \Omega^i = T^{ijk}\pi_j\wedge\omega^k -
    (S^i_j\omega^j + U^{ij}\pi_j)\wedge\theta,
\label{StreqnC15}
\end{equation}
with $T^{ijk}=T^{jik}=T^{kji}$, $T^{iik}=0$; $U^{ij}=U^{ji}$;
$S^i_j=S^j_i$, $S^i_i=0$.  The pseudo-connection form $\varphi$ may be
chosen so that also (cf.~(\ref{MasterDRho}))
$$
  (n-2)d\rho = -\delta_i\wedge\omega^i-S^i_j\pi_i\wedge\omega^j
    +(\sf{n-2}{2n}U^{ij}\sigma_{ij}-t^i\pi_i)\wedge\theta
$$
for some functions $t^i$.
\end{itemize}
For any point $p\in N$, we consider a neighborhood $U\subset M$ of
$F_0(p)$ on which we can fix onesuch coframing $(\theta,\omega^i,\pi_i)$ with
pseudo-connection $\varphi$.  All of the forms and functions may be
pulled back to $W=F^{-1}(U)\subset N\times[0,1]$, which is the setting
for the calculation of $\mathcal L_{\frac{\partial}{\partial
    t}}(F^*\Psi)$.  From now on, we drop all $F^*$s.

We have $\Psi = -\pi_i\wedge\omega_{(i)}$, and we now have the
structure equations needed to differentiate $\Psi$, but it will
simplify matters if we further adapt the forms $(\theta,\omega^i,\pi_i)$ on
$W$ in a way that does not alter the structure equations.  Note first
that restricted to each $W_t=W\cap N_t$, we have $\bigwedge\omega^i\neq 0$,
so $(\omega^1,\ldots,\omega^n,dt)$ forms a coframing on $W$.  We can therefore write
$\pi_i=s_{ij}\omega^j+g_idt$ for some $s_{ij}$, $g_i$, and because
each $W_t$ is Legendre, we must have $s_{ij}=s_{ji}$.  The structure
group of $B\to M$ admits addition of a traceless, symmetric
combination of the $\omega^j$ to the $\pi_i$, so we replace
$$
  \pi_i\leadsto\pi_i-s_{ij}^o\omega^j,
$$
where $s_{ij}^0=s_{ij}-\frac1n\delta^i_js_{kk}$ is the traceless
part.  Now we have
$$
  \pi_i = s\omega^i + g_idt
$$
for some functions $s$, $g_i$ on $W$, so that
$$
  \Psi = -\pi_i\wedge\omega_{(i)} = -ns\,\omega-g_idt\wedge
    \omega_{(i)}.
$$
Note that because $F_0$ is assumed integral for the Euler-Lagrange
system, we have $s=0$ everywhere on $W_0\subset W$; in particular,
$\pi_i|_{W_0}=0$.

With our choice of $\pi_i$, recalling that along $W_0$, $\theta = g\,dt$ for
some function $g$ on $W_0$, and keeping in mind that restricted to
$W_0$ we have $\pi_i=0$ and $s=0$, we can calculate
\begin{eqnarray*}
  \left.\mathcal L_{\fddt}\Psi\right|_{t=0} 
  & = &
    -\sfddt\innerprod d(\pi_i\wedge\omega_{(i)})
    -d(\sfddt\innerprod(\pi_i\wedge\omega_{(i)})) \\
  & = & \sfddt\innerprod\left((\delta_i\wedge\theta+
           (n\rho\delta^j_i-\alpha^j_i)\wedge\pi_j)
           \wedge\omega_{(i)}\right. \\
   & & \left.\qquad\quad + (s\omega^i+g_idt)\wedge d\omega_{(i)}\right)
    \\ & &-d\left((g_i+s\sfddt\innerprod\omega^i)\omega_{(i)}
           -\pi_i\wedge(\sfddt\innerprod\omega_{(i)})\right) \\
  & = & -(g\delta_i+dg_i+n\rho g_i-g_j\alpha^j_i)\wedge\omega_{(i)}.
\end{eqnarray*}
This gives our desired formula:
\begin{equation}
  \left.\frac{d^2}{dt^2}\right|_{t=0}\left(\int_{N_t}\Lambda\right)
    = -\int_{N_0}g(dg_i+n\rho g_i-g_j\alpha^j_i+g\delta_i)\wedge
     \omega_{(i)}.
\label{SecondVar}
\end{equation}
Unfortunately, in its present form this is not very enlightening, and
our next task is to give a geometric interpretation of the formula.

\subsection{Relative Conformal Geometry}
\label{Subsection:Relative}

It is natural to ask what kind of geometric structure is induced on an
integral manifold $f:N\hookrightarrow M$ of an Euler-Lagrange system
$\mathcal E_\Lambda$.  What we find is:
\begin{quote}{\em
If $\Pi=d\Lambda$ is a
definite, neo-classical Poincar\'e-Cartan form, then an integral
manifold $N$ of its Euler-Lagrange system ${\mathcal E}_\Lambda$
has a natural 
conformal structure\index{conformal!structure|(}, invariant under
symmetries\index{symmetry} of
$(M,\Pi,N)$, even though there may be no invariant conformal structure
on the ambient $M$.}
\end{quote}
This is a simple pointwise phenomenon, in the sense that {\em any}
$n$-plane $V^n\subset T_pM^{2n+1}$ on which $\bigwedge\omega^i\neq 0$
has a canonical conformal inner-product\index{conformal!inner-product}
defined as follows.  Taking any 
section $(\theta,\omega^i,\pi_i)$ of $B\to M$, we can restrict the
induced quadratic form $\sum(\omega^i)^2$ on $TM$ to $V\subset T_pM$,
where it is positive definite,
and then the action of the structure group (\ref{LieAlgebra15}) on
$(\omega^i)$ shows that up to scaling, this quadratic form is
independent of our choice of section.  Alternatively, one can show
this infinitesimally by using the structure equations to compute on
$B$ the Lie derivative of
$\sum(\omega^i)^2$ along a vector field that is vertical for $B\to M$;
this Lie derivative is itself multiple of
$\sum(\omega^i)^2$.  Note that we have {\em not} restricted to the
conformal branch of the equivalence problem, characterized by
$T^{ijk}=U^{ij}=S^i_j=0$ and discussed in \S\ref{Section:ConfEquiv}.

In particular, any integral manifold $f:N\hookrightarrow M$ for the
Euler-Lagrange system $\mathcal E_\Lambda$ inherits a canonical
conformal structure $[ds^2]_f$.
We now want to develop the conformal structure equations for 
$(N,[ds^2]_f)$, in terms of the structure equations
on $B\to M$, and our procedure will work only for integral manifolds of
$\mathcal E_\Lambda$.  We first note that along our
integral manifold $N$ we can choose local sections
$(\theta,\omega^i,\pi_i)$ of $B_N$ which are adapted to $N$ in the
sense that 
$$
  T_pN = \{\theta,\pi_1,\cdots,\pi_n\}^\perp\subset T_pM,
$$
for each $p\in N$.  In fact, such sections define a reduction $B_f\to
N$ of the principal bundle $B_N\to N$, having Lie algebra defined as
in (\ref{LieAlgebra15}) by $s_{ij}=0$.

Restricted to $B_f$, we have the same structure equations as on $B$,
but with
$\theta=\pi_i=0$, and $d\theta=d\pi_i=0$.  Now observe that two of our
structure equations restrict to give
\begin{eqnarray*}
d\omega^i & = & (2\rho\delta^i_j-\alpha^i_j)
    \wedge\omega^j, \\ 
  d\rho & = & -\sf{1}{n-2}\delta_i\wedge\omega^i.
\end{eqnarray*}
We therefore have a situation similar to that in \S\ref{Section:ConfEquiv}
(cf.~(\ref{domega12}, \ref{AlmostConf12})), with equations formally
like those in the conformal geometry equivalence problem of
\S\ref{Subsection:ConfEquiv}.
We can mimic the derivation of conformal structure equations in the
present case by first setting $\beta_i = \sf2{n-2}\delta_i$,
and then we know that this results in an equation
\begin{equation}
  d\alpha^i_j+\alpha^i_k\wedge\alpha^k_j+\beta_i\wedge\omega^j
   -\beta_j\wedge\omega^i =
   \sf12A^i_{jkl}\omega^k\wedge\omega^l,
\label{TempConf15}
\end{equation}
where the quantity $(A^i_{jkl})$ has the symmetries of the Riemann
curvature tensor\index{Riemann curvature tensor}.  
Furthermore, we know that there are unique
functions $t_{ij}=t_{ji}$ such that replacing
$\beta_i\leadsto\beta_i+t_{ij}\omega^j$ will yield the
preceding equation with $A^l_{jkl}=0$.  However, it will simplify
matters later if we go back and replace instead
$$
  \delta_i\leadsto\delta_i+\sf{n-2}{2}t_{ij}^o\omega^j,
$$
where $t_{ij}^o=t_{ij}-\sf1n\delta_{ij}t_{kk}$ is the traceless part
(note that only a {\em traceless} addition to $\delta_i$ will preserve
the structure equations on $B\to M$).  In terms of the new
$\delta_i$, we define
\begin{equation}
  \beta_i \stackrel{\mathit{def}}{=} \sf2{n-2}(\delta_i-\sf1nK\omega^i),
\label{DefBeta15}
\end{equation}
where
$$
  K \stackrel{\mathit{def}}{=} -\sf{n-2}{2}t_{jj}.
$$
This $K$ was chosen so that defining $\beta_i$ by
(\ref{DefBeta15}) gives the correct conformal structure equation
(\ref{TempConf15}) with $A^l_{jkl}=0$; it reflects the difference
between the pseudo-connection forms for the Poincar\'e-Cartan form
and the Cartan connection\index{connection!Cartan} 
forms for the induced conformal geometry on
the submanifold.  One interpretation is the following.
\begin{quote}{\em
The function $K$ on $B_f$ is a fundamental invariant of a stationary
submanifold $(N,[ds^2]_f)\hookrightarrow(M,\Pi)$ of an Euler-Lagrange
system, and may be thought of as an extrinsic curvature
depending on up to third derivatives of the immersion $f$.}
\end{quote}
In the classical setting, $N$ is already the $1$-jet graph of a
solution $z(x)$ of an Euler-Lagrange equation, so an expression for
$K(x)$ depends on fourth derivatives of $z(x)$.

Now suppose that our integral manifold $f:N\hookrightarrow M$, with
$\delta_i$, $\beta_i$, and $K$ as above, is the initial manifold $F_0=f$
in a Legendre variation $F:N\times[0,1]\to M$.  Then we can rewrite
our formula (\ref{SecondVar}) for the second variation as
$$
  \delta^2(\mathcal F_\Lambda)_{N_0}(g) =  -\int_{N_0}
    \left(g(dg_i+n\rho g_i-g_j\alpha^j_i+\sf{n-2}{2}g\beta_i)
       \wedge\omega_{(i)}
    +g^2K\omega\right).
$$
Part of this integrand closely resembles the expression
(\ref{SecondDeriv10}) for the second covariant
derivative\index{covariant derivative|(} of a section
of a density line bundle\index{density line bundle}, 
discussed in constructing the conformal
Laplacian\index{conformal!Laplacian} 
in \S\ref{Subsection:ConfLaplacian}.  This suggests the
following computations.

First, consider the structure equation
$d\theta=-(n-2)\rho\wedge\theta-\pi_i\wedge\omega^i$.  Along $W_0$
(but not yet restricted to $W_0$), where
$\theta = g\,dt$ and $\pi_i=g_idt$, this reads
$$
  dg\wedge dt = -(n-2)\rho\wedge g\,dt-g_idt\wedge\omega^i,
$$
so that restricted to $W_0$, we have
$$
  dg + (n-2)\rho g= g_i\omega^i.
$$
This equation should be interpreted on $B_f$, which is identified with
the principal bundle associated to the conformal structure $[ds^2]_f$.
It says that $g$ is a section of the density bundle
$D^{\frac{n-2}{2n}}$, and that
$g_i$ are the components of its covariant derivative
(see (\ref{GenDensityInfEquivar}) ff.).  We can now write
$$
  dg_i+n\rho g_i -g_j\alpha^j_i+\sf{n-2}{2}g\beta_i = g_{ij}\omega^j,
$$
and by definition
$$
  \Delta_fg = g_{ii}\in\Gamma(D^{\frac{n+2}{2n}}),
$$
where $\Delta_f$ is the conformal Laplacian on $N$ induced
by $f=F_0:N\hookrightarrow M$.  We now have a more promising version of
the second variation formula:
\begin{equation}
   \boxed {\delta^2(\mathcal F_\Lambda)_{N_0}(g) =
   -\int_{N_0}(g\Delta_fg+g^2K)\omega. }
\label{GoodSecondVarn}
\end{equation}
It is worth noting that the sign of this integrand does not depend on
the sign of the 
variation's generating function
$g$, and if we fix an orientation of $N$,
then the integrand $K\omega$ on $B_f$ has a well-defined sign at
each point of $N$.

\subsection{Intrinsic Integration by Parts}

In order to detect local minima using the second variation formula
(\ref{GoodSecondVarn}), it is often helpful to
convert an integral like $\int g\Delta g\,dx$ into one like $-\int ||\nabla
g||^2dx$.  In the Euclidean setting, with either compact supports or with
boundary terms, this is done with integration by parts and is
straightforward; the two integrands differ by the divergence of
$g\nabla g$, whose integral depends only on boundary data.

We would like to perform a calculation like this on $N$, for an
arbitrary Legendre variation $g$, using only intrinsic data.  In other words,
we would like to associate to any $N$ and $g$ some
$\xi\in\Omega^{n-1}(N)$ such
that $d\xi$ is the difference between $(g\Delta_fg)\omega$ and some
quadratic expression $Q(\nabla g, \nabla g)\omega$, possibly with some
additional zero-order terms.  A natural expression to consider,
motivated by the flat case, is
\begin{equation}
  \xi = gg_i\omega_{(i)}.
\label{BoundaryForm}
\end{equation}
Here, in order for $g_i$ to make sense, we are assuming that we have
a local coframing on $M$ adapted to the integral manifold $N$, so that
at points of $N$, $\theta = g\,dt$, $\pi_i=g_idt$, and so that restricted to
$N$, $\theta|_N=\pi_i|_N=0$.  We can then compute:
\begin{eqnarray}
  d\xi & = & dg\wedge g_i\omega_{(i)} +g\,dg_i\wedge\omega_{(i)}
               +gg_id\omega^j\wedge\omega_{(ij)} \\
   & = & (dg+(n-2)\rho g)\wedge g_i\omega_{(i)} +
       g(dg_i+n\rho g_i-g_j\alpha^j_i)\wedge\omega_{(i)} \\
   & = & \left(\sum (g_i)^2 +g\Delta_fg\right)\omega-\sf{n-2}{2}g^2\beta_i
           \wedge\omega_{(i)}.
\label{Parts15}
\end{eqnarray}
Now, the first term is exactly what we are looking for, and second
is fairly harmless because it is of order zero in the variation $g$,
and in practice contributes only terms similar to $g^2K\omega$.

The problem is that $\xi$ (\ref{BoundaryForm}) 
is defined on the total space $B_f\to N$,
and although semibasic for this bundle, it is not basic; that is,
there is no form on $N$ that pulls back to $B_f$ to equal $\xi$, even
locally.  
The criterion for $\xi$ to be basic is that $d\xi$ be semibasic, and
this fails because of the appearance of $\beta_i$ in (\ref{Parts15}).

But suppose that we can find some canonical reduction of the ambient
$B\to M$ to a subbundle on which $\delta_i$ becomes semibasic over
$M$; in terms of the Lie algebra (\ref{LieAlgebra15}), this means that
we can reduce to
the subgroup having Lie algebra given by $\{d_i=0\}$.  In this case, each
$\beta_i=\sf{2}{n-2}(\delta_i-\sf{K}{n}\omega^i)$ is a linear
combination of $\theta$, $\omega^i$, $\pi_i$, and is in particular
semibasic over $N$.  Consequently, $\xi$ is basic
over $N$, and we can perform the integration by parts in an invariant
manner. 

Unfortunately, there are cases in which no such canonical reduction of
$B$ is possible.  An example is the homogeneous Laplace
equation\index{Laplace equation}
on $\R^n$,
$$
  \Delta z=0,
$$
which is preserved under an action of the conformal
group\index{conformal!group}
$SO^o(n+1,1)$.  The associated conformal geometry on the trivial
solution $z=0$ is flat, and our second variation formula reads
$$
  \delta^2(\mathcal F_\Lambda)_0(g) = \int_\Omega g\Delta g\,dx
$$
for a variation $g\in C^\infty_0(\Omega)$, $\Omega\subset\R^n$.
It follows from our construction that this integrand is invariant
under a suitable action of the conformal group.  However,
the tempting integration by parts
$$
  \int_\Omega g\Delta g\,dx = -\int_\Omega||\nabla g||^2dx
$$
leaves us with an integrand which is {\em not} conformally invariant.
It is this phenomenon that we would like to avoid.

To get a sense of when one might be able to find a canonical subbundle
of $B\to M$ on which the $\delta^i$ are semibasic, recall that
$$
  (n-2)d\rho\equiv-\delta_i\wedge\omega^i\pmod{\{\theta,\pi_i\}}.
$$
Working modulo $\{\theta,\pi_i\}$ essentially amounts to restricting
to integral manifolds of the Euler-Lagrange system.  The preceding
then says that a choice of subbundle $B^\prime\subset B$ on which
$\delta_i$ are semibasic gives a subbundle of each conformal bundle
$B_f\to N$ on which $d\rho$ is semibasic over $N$.  Now, typically
the role of $\rho$ in the Cartan connection for a conformal structure
is as a psuedo-connection in the density line bundle $D$; a
special reduction is required for $\rho$ to be a genuine connection,
and the latter requirement is equivalent to $d\rho$ being semibasic.
In other words, being able to integrate by parts in an invariant
manner as described above is equivalent to having a connection in $D$
represented by the pseudo-connection $\rho$.
One way to find a connection in $D$ is to suppose that $D$
has somehow been trivialized, and this is equivalent to choosing a
Riemannian metric\index{Riemannian!metric} 
representing the conformal class.  This suggests
that Euler-Lagrange systems whose integral manifolds have canonical
Riemannian metrics will have canonical reductions of this type.

In fact, we have seen an example of a Poincar\'e-Cartan form whose
geometry $B\to M$ displays this behavior.  This is
the system for Riemannian hypersurfaces having prescribed mean
curvature\index{mean curvature!prescribed|(}, 
characterized in \S\ref{Section:PrescribedH} in terms of
differential invariants of its neo-classical, definite
Poincar\'e-Cartan form.  To illustrate the preceding discussion, we
calculate the second variation formula for this system.  The reader
should note especially how use of the geometry of the
Poincar\'e-Cartan form gives a somewhat simpler derivation of the formula
than one finds in standard sources.\footnote{See for example
  pp.~513--539 of \cite{Spivak:ComprehensiveIV},
where the calculation is prefaced by a colorful warning about
  its difficulty.} 
\index{second variation|)}
\index{conformal!structure|)}

\subsection{Prescribed Mean Curvature, Revisited}
\index{second variation!for prescribed mean curvature|(}

In \S\ref{Section:PrescribedH}, we considered a definite,
neo-classical Poincar\'e-Cartan form $(M,\Pi)$ whose associated
geometry $(B\to M,\varphi)$ had invariants satisfying
$$
  T^{ijk}=0,\quad U^{ij}=\lambda\delta^{ij}.
$$
We further assumed the open condition
$$
  \lambda < 0,
$$
and this led to a series of reductions of $B\to M$, resulting in a
principal subbundle $B_3\to M$, having structure group with Lie algebra
$$
  \lie{g}_3 = \left\{\left(\begin{array}{ccc}
    0 & 0 & 0 \\ 0 & a^i_j & 0 \\ 0 & 0 & -a^j_i \end{array}
  \right):a^i_j+a^j_i=0\right\},
$$
on which the original structure equations
(\ref{StreqnA15}, \ref{StreqnB15}, \ref{StreqnC15}) hold, with
$U^{ij}=-\delta^{ij}$, $S^i_j=T^{ijk}=0$, $\rho = -\sf{H}{2n}\theta$
for a function $H$ on $B_3$, and
$d\rho=-\sf{1}{n-2}\delta_i\wedge\omega^j$, where
$\delta_i\equiv0\mbox{ (mod $\{\theta,\omega^j,\pi_j\}$)}$.  In this case we
computed that $dH\equiv0\mbox{ (mod $\{\theta,\omega^i\}$)}$, so that $H$
defines a function on the local quotient space $Q^{n+1}$, which also
inherits a Riemannian metric $\sum(\omega^i)^2$.  The contact
manifold $M$ can be locally identified with the bundle of tangent
hyperplanes of $Q$, and the integral manifolds of the Euler-Lagrange
differential system $\mathcal E_\Lambda$ are
the tangent loci of hypersurfaces
in $Q$ whose mean curvature coincides with the background function $H$.
In this case $B_3\to M\to Q$ is locally identified with the
orthonormal frame bundle\index{Riemannian!frame bundle|(} 
of the Riemannian manifold $Q$.

An important point here is that the Riemannian geometry associated to
$\Pi$ only appears after reducing to $B_3\to M$.  However, if our goal
is to see the formula for the second variation, then we face the
following difficulty.  That formula required the use of coframes of
$M$ adapted to a stationary submanifold $N\hookrightarrow M$ in a
certain way, but while adapted coframes can always be found in
$B\to M$, there is no guarantee that they can be found in the
subbundle $B_3\to M$, where the Riemannian geometry is visible.

We will overcome these difficulties and illustrate the invariant
calculation of the second variation by starting only with the Riemannian
geometry of $(Q,ds^2)$.  This is expressed in the Levi-Civita
connection\index{connection!Levi-Civita} 
in the orthonormal frame bundle, where we can also give the
Poincar\'e-Cartan form and Euler-Lagrange system for prescribed mean
curvature.  We then introduce higher-order data on a larger bundle,
which allows us to study the second fundamental 
form\index{second fundamental form}.  In fact, this
larger bundle corresponds to the partial reduction $B_2\to M$ on which
$\rho$ and $\delta_i$ are semibasic, but $\sigma_{ij}$ is not.  
The end result of our calculation is formula (\ref{SecondVarnArea15}).
In the following discussion, index ranges are $0\leq a,b,c\leq n$ and
$1\leq i,j,k\leq n$.

We begin with a generalization of 
the discussion in \S\ref{Section:Hypersurface}
of constant mean curvature\index{mean curvature!constant} 
hypersurfaces in Euclidean space.
Let $(Q,ds^2)$ be an oriented Riemannian manifold of dimension $n+1$.
A {\em frame} for $Q$ is a pair
$f = (q,e)$
consisting of a point $q\in Q$ and a positively-oriented orthonormal basis
$e=(e_0,\ldots,e_n)$ for $T_qQ$.  The set $\mathcal F$ of all such
frames is a manifold, and the right
$SO(n+1,\R)$-action
$$
  (q,(e_0,\ldots,e_n))\cdot(g^a_b) = (q,(\textstyle\sum e_ag^a_0,\ldots,
     \sum e_ag^a_n))
$$
gives the basepoint map
$$
  q:{\mathcal F}\to Q
$$
the structure of a principal bundle.
The unit sphere bundle
$$
  M^{2n+1}=\{(q,e_0):q\in Q,\ e_0\in T_qQ,\ ||e_0||=1\}
$$
is identified with the Grassmannian\index{Grassmannian} 
bundle of oriented tangent $n$-planes in $TQ$, and it
has a contact structure generated by the $1$-form
\begin{equation}
  \theta_{(q,e_0)}(v) = ds^2(e_0,q_*(v)),\qquad v\in 
    T_{(q,e_0)}M,
\label{DefContact}
\end{equation}
where $q:M\to Q$ is the projection.  An immersed oriented hypersurface
in $Q$ has a unit normal vector field, which may be thought of as a
$1$-jet lift of the submanifold to $M$.  The submanifold of $M$ thus
obtained is easily seen to be a Legendre submanifold for this contact
structure, and the transverse Legendre 
submanifold\index{Legendre submanifold!transverse} is locally of this
form.

To carry out calculations on $M$, and even to verify the
non-degeneracy of $\theta$, we will use the
projection $\mathcal F\to M$ defined by $(q,(e_0,\ldots,e_n))\mapsto(q,e_0)$.
Calculations can then be carried out using structure equations for the
canonical parallelization of $\mathcal F$, which we now introduce.
First, there are
the $n+1$ tautological $1$-forms
$$
  \varphi^a_{(q,e)} = ds^2(e_a, q_*(\cdot))\in\Omega^1(\mathcal F),
$$
which
form a basis for the semibasic $1$-forms over $Q$. 
Next, there are
globally defined, uniquely determined Levi-Civita
connection\index{connection!Levi-Civita} forms
$\varphi^a_b=-\varphi^b_a\in\Omega^1(\mathcal F)$ satisfying
\begin{equation}
  \left\{\begin{array}{l}
    d\varphi^a = -\varphi^a_b\wedge\varphi^b, \\
    d\varphi^a_b=-\varphi^a_c\wedge\varphi^c_b+\sf12R^a_{bcd}
      \varphi^c\wedge\varphi^d.
  \end{array}\right.
\label{RiemStreqns15}
\end{equation}
The functions $R^a_{bcd}$ on $\mathcal F$ are the
components of the Riemann curvature 
tensor\index{Riemann curvature tensor} with respect to different
orthonormal frames, and satisfy
$$
  R^a_{bcd}+R^a_{bdc}=R^a_{bcd}+R^b_{acd}=R^a_{bcd}+R^a_{cdb}
    +R^a_{dbc}=0.
$$
We now distinguish the $1$-form
$$
  \theta \stackrel{\mathit{def}}{=} \varphi^0\in\Omega^1(\mathcal F),
$$
which is the pullback via $\mathcal F\to M$ of the
contact $1$-form given the same name in (\ref{DefContact}).
One of our structure equations now reads
\begin{equation}
  d\theta = -\varphi^0_i\wedge\varphi^i.
\label{dContact15}
\end{equation}
This implies that the original $\theta\in\Omega^1(M)$ is actually a
contact form, and also that $(\theta,\varphi^i,\varphi^0_i)$ is a
basis for the semibasic $1$-forms for $\mathcal F\to M$.

At this point, we can give the Poincar\'e-Cartan form for the
prescribed mean curvature system.  Namely, let $H\in C^\infty(Q)$ be a
smooth function, and define on $\mathcal F$ the $(n+1)$-form
$$
  \Pi \stackrel{\mathit{def}}{=} -\theta\wedge(\varphi^0_i\wedge\varphi_{(i)}-H\varphi).
$$
Because $H$ is the pullback of a function on $Q$, its derivative is
of the form
$$
  dH = H_\nu\theta + H_i\varphi^i,
$$
and using this and the structure equations (\ref{RiemStreqns15}), one
can verify
that $\Pi$ is closed.  Because $\Pi$ is semibasic over $M$ and closed,
it is the pullback of a closed form on $M$, which is then a definite,
neo-classical Poincar\'e-Cartan form.  The associated Euler-Lagrange
differential system then pulls back to $\mathcal F$ as
$$
  \mathcal E_H = \{\theta,d\theta,
  \varphi^0_i\wedge\varphi_{(i)}-H\varphi\}.
$$
While (\ref{dContact15}) shows that generic Legendre $n$-planes in $M$ are
defined by equations
$$
  \theta = 0,\quad \varphi^0_i-h_{ij}\varphi^j = 0,
$$
with $h_{ij}=h_{ji}$, integral $n$-planes in $M$ for $\mathcal E_H$
are defined by the same equations, plus
$$
  h_{ii}=H.
$$
The functions $h_{ij}$ describing the tangent locus of a transverse Legendre
submanifold of $M$ are of course the coefficients of the second
fundamental form of the corresponding submanifold of $Q$.  Therefore, 
the transverse integral manifolds of $\mathcal E_H$ correspond locally
to hypersurfaces in $Q$ whose mean curvature $h_{ii}$ equals the background
function $H$.  This will appear quite explicitly in what follows.

To investigate these integral manifolds, we employ the following
apparatus.  First consider the product
$$
  \mathcal F\times \R^{n(n+1)/2},
$$
where $\R^{n(n+1)/2}$ has coordinates $h_{ij}=h_{ji}$, and inside this
product define the locus
$$
  {\mathcal F}^{(1)} \stackrel{\mathit{def}}{=}
    \left\{(f,h)\in{\mathcal F}\times\R^{n(n+1)/2}: h_{ii} = H\right\}.
$$
To perform calculations, we want to extend our parallelization of
$\mathcal F$ to $\mathcal F^{(1)}$.  With a view toward
reconstructing some of the bundle $B\to M$ associated to the
Poincar\'e-Cartan form $\Pi$, we do this in a way that is as
well-adapted to $\Pi$ as possible.

On $\mathcal F^{(1)}$, we continue to work with $\theta$, and define
$$
  \begin{array}{l}
  \omega^i \stackrel{\mathit{def}}{=} \varphi^i, \\
  \pi_i \stackrel{\mathit{def}}{=} \varphi^0_i-h_{ij}\varphi^j.
  \end{array}
$$
With these definitions, we have
$$
  \begin{array}{l}
  d\theta = -\pi_i\wedge\omega^i, \\
  \Pi = -\theta\wedge\pi_i\wedge\omega_{(i)}.
  \end{array}
$$
Motivated by Riemannian geometry, we set
$$
  Dh_{ij} \stackrel{\mathit{def}}{=} dh_{ij} - h_{kj}\alpha^k_i -
  h_{ik}\alpha^k_j, 
$$
so that in particular $Dh_{ij}=Dh_{ji}$ and $Dh_{ii}=dH$.  We also define
the traceless part
$$
  Dh^0_{ij} \stackrel{\mathit{def}}{=}Dh_{ij}-\sf1n\delta_{ij}dH.
$$
Direct computations show that we will have exactly the
structure equations (\ref{StreqnA15}, \ref{StreqnB15}, \ref{StreqnC15})
if we define
\begin{eqnarray*}
    \alpha^i_j & = & \varphi^i_j, \\
    \rho & = & -\sf1{2n}H\theta, \\
    \delta_i & = & -\sf12H\pi_i+h_{ij}\pi_j +
      (h_{ik}h_{kj}-R^0_{ij0}-\sf1n\delta_{ij}H_\nu)\omega^j, \\
    \sigma_{ij}\wedge\omega^j & = & (Dh^0_{ij}+\sf1n\delta_{ij}H_k\omega^k
      +\sf12R^0_{ijk}\omega^k)\wedge\omega^j, 
\end{eqnarray*}
with $\sigma_{ij}=\sigma_{ji}$, $\sigma_{ii}=0$.
The last item requires some comment.  Some linear algebra
involving a Koszul complex\index{Koszul complex} 
shows that for any tensor $V_{ijk}$ with
$V_{ijk}=-V_{ikj}$ (this will be applied to
$\sf1{2n}(\delta_{ij}H_k-\delta_{ik}H_j)+\sf12R^0_{ijk}$), there is
another tensor $W_{ijk}$, not unique, satisfying $W_{ijk}=W_{jik}$,
$W_{iik}=0$, and $\sf12(W_{ijk}-W_{ikj})=V_{ijk}$.  This justifies the
existence of $\sigma_{ij}$ satisfying our requirements.
The structure equations
(\ref{StreqnA15}, \ref{StreqnB15}, \ref{StreqnC15}) resulting from these
assignments have torsion coefficients
$$
  T^{ijk}=0,\quad U^{ij}=-\delta^{ij},\quad
    S^i_j = -h_{ij}+\sf1n\delta^i_jH.
$$
\index{Riemannian!frame bundle|)}

The general calculations of \S\ref{Subsection:Formula} for the second
variation can now be applied; note that we have the freedom to adapt
coframes to a single integral submanifold in $M$ of $\mathcal E_H$.
Repeating those calculations verbatim leads to
$$
  \left.\frac{d^2}{dt^2}\right|_{t=0}\left(\int_{N_t}\Lambda\right) =
    -\int_{N_0}g(dg_i+n\rho g_i-g_j\alpha^j_i+g\delta_i)\wedge
    \omega_{(i)},
$$
where $F:N\times[0,1]\to M$ is a Legendre variation
in $M$, $F_0$ is an integral manifold of $\mathcal E_H$, and the forms are
all pullbacks of forms on $F^*(\mathcal F^{(1)})$ by a section of
$F^*(\mathcal F^{(1)})\to N\times[0,1]$, adapted along $N_0$ in the
sense that
$$
  \theta_{N_0}=g\,dt,\ (\pi_i)_{N_0}=g_idt.
$$
These imply that {\em restricted} to $N_0$, we have
$\theta|_{N_0}=\pi_i|_{N_0}=0$, and the preceding formula becomes
$$
  \delta^2(\mathcal F_\Lambda)_{N_0}(g)
    = 
    -\int_{N_0}g(dg_i-g_j\alpha^j_i+g(h_{ik}h_{kj}-\sf1n
    \delta_{ij}H_\nu - R^0_{ij0})\omega^j)\wedge\omega_{(i)}.
$$
Recognizing that $g$ can be thought of as a section of the normal
bundle of the hypersurface $N\stackrel{F_0}{\hookrightarrow}M
\stackrel{q}{\to}Q$,
and that in this case $g_i$ are the coefficients of its covariant
derivative\index{covariant derivative|)}, this can be rewritten as
\begin{equation}
  \delta^2(\mathcal F_\Lambda)_{N_0}(g)
  = -\int_{N_0}(g\Delta g + g^2(||h||^2-H_\nu-R^0_{ii0}))\omega,
\label{SecondVarnArea15}
\end{equation}
where $\Delta$ is the Riemannian
Laplacian\index{Riemannian!Laplacian}, and $||h||^2=Tr(h^*h)$.
Here, the extrinsic curvature function $K$ appears as the quantity
$||h||^2-H_\nu-R^0_{ii0}$.
Notice that we have actually calculated this second variation without ever
determining the functional $\Lambda$.
In case the ambient manifold $Q$ is flat Euclidean space,
if the background function $H$ is a constant and the variation $g$ is
compactly supported in the interior of $N$, this simplifies to
\begin{eqnarray*}
  \left.\frac{d^2}{dt^2}\right|_{t=0}\left(\int_{N_t}\Lambda\right)
    & = & -\int_{N_0}(g\Delta g+g^2||h||^2)\omega \\
    & = & \int_{N_0}(||\nabla g||^2-g^2||h||^2)\omega.
\end{eqnarray*}
Even for the minimal surface\index{minimal surface} 
equation $H=0$, we cannot conclude from
this formula alone that a solution locally minimizes area.  
\index{second variation!for prescribed mean curvature|)}
\index{mean curvature!prescribed|)}

\subsection{Conditions for a Local Minimum}

We now discuss some conditions under which an integral manifold
$N\hookrightarrow M$ of
an Euler-Lagrange system $\mathcal E_\Lambda\subset\Omega^*(M)$ is a
local minimum for
the functional $\mathcal F_\Lambda$, in the sense that $\mathcal
F_\Lambda(N)<\mathcal F_\Lambda(N^\prime)$ for all Legendre
submanifolds $N^\prime$ near $N$.  
\index{Legendre submanifold|)}
However, there are two natural
meanings for ``near'' in this context,
and this will yield two notions of local
minimum.  Namely, we will say that $\mathcal F_\Lambda$ has a {\em
  strong}\index{strong local minimum|(} 
local minimum at $N$ if the preceding inequality holds
whenever $N^\prime$ is $C^0$-close to $N$, while $\mathcal F_\Lambda$
has a {\em weak} local minimum\index{weak local minimum|(} 
at $N$ if the preceding inequality holds only
among the narrower class of $N^\prime$ which are $C^1$-close to
$N$.\footnote{A thorough, coordinate-based discussion of the relevant
  analysis
  can be found in \cite{Giaquinta:Calculus}.}
\label{GHPage}

Our goal is to illustrate how the Poincar\'e-Cartan form may be
used to understand in a simple geometric manner some classical
conditions on extrema.  Specifically, we will introduce the notion of
a {\em calibration}\index{calibration|(} 
for an integral manifold of the Euler-Lagrange
system; its existence (under mild topological hypotheses) implies that
the integral manifold is a strong local minimum.  Under certain classical
conditions for a local minimum, we will use the
Poincar\'e-Cartan form to construct an analogous {\em weak calibration}.
Finally, our geometric description of the second variation formula
highlights the Jacobi operator\index{Jacobi operator} 
$Jg=-\Delta_cg+Kg\omega$, and some
linear analysis shows that the positivity of the first eigenvalue of
$J$ implies the classical conditions. 

Let $\Pi$ be
a neo-classical Poincar\'e-Cartan form $\Pi$ on a
contact manifold $(M,I)$.  There is a local foliation $M\to Q$,
and we can choose coordinates to have $(x^i,z)\in Q\subset
\R^n\times\R$, $(x^i,z,p_i)\in M\subset J^1(\R^n,\R)$,
$\theta=dz-p_idx^i\in\Gamma(I)$, and a Lagrangian
potential\index{Lagrangian potential}
$$
  \Lambda=L(x,z,p)dx+\theta\wedge L_{p_i}dx_{(i)}\in
    \Omega^n(M)
$$
whose Poincar\'e-Cartan form is
$$
  \Pi = d\Lambda = \theta\wedge(-dL_{p_i}\wedge dx_{(i)}+L_zdx).
$$
We will confine our discussion to a domain where this classical
description holds.  We may regard an integral manifold of the
Euler-Lagrange system $\mathcal E_\Lambda$ as a submanifold
$N_0\hookrightarrow Q$ given by the graph $\{(x,z_0(x)):x\in U\}$ of a
solution to the Euler-Lagrange equation over some open
$U\subset\R^n$.  It has a natural $1$-jet extension
$N_0^{(1)}\hookrightarrow M$, equal to $\{(x,z_0(x),\nabla
z_0(x)):x\in U\}$, which is an integral manifold of $\mathcal
E_\Lambda$ in the sense discussed previously.  We define a {\em strong
  neighborhood} of $N_0$ to be the collection of hypersurfaces in $Q$
lying in some open neighborhood of $N_0$ in $Q$, and a {\em weak
  neighborhood} of $N_0$ to be the collection of hypersurfaces
$N\hookrightarrow Q$ whose $1$-jet prolongations $N^{(1)}$ lie in some
open neighborhood of $N_0^{(1)}$ in $M$.  Whether or not a given
stationary submanifold $N_0$ is minimal depends on which of
these two classes of competing submanifolds one studies.

Starting with strong neighborhoods, we fix a neighborhood
$W\subset Q$ of a stationary submanifold $N_0\subset Q$ for $\Lambda$,
and introduce the following useful notion.
\begin{Definition}
A {\em calibration} for $(\Lambda, N_0)$ is an $n$-form
$\tilde\Lambda\in\Omega^n(W)$ satisfying
\begin{itemize}
\item $d\tilde\Lambda = 0$;
\item $\tilde\Lambda|_{N_0} = \Lambda|_{N^{(1)}_0}$;
\item $\tilde\Lambda|_E\leq\Lambda_E$ for each $n$-plane $E^n\subset
  T_qW$.
\end{itemize}
\end{Definition}
In the right-hand side of the last inequality, we are regarding the
$n$-plane $E$ as specifying a point of $M\subset G_n(TQ)$ over $q\in
Q$, and evaluating $\Lambda_{(q,E)}$ on any
tangent $n$-plane $E^\prime\subset T_{(q,E)}M$ projecting one-to-one
into $T_qQ$; the value is independent of the choice of $E^\prime$,
because $\Lambda$ is semibasic over $Q$.
In particular, the third condition says that the integral of
$\tilde\Lambda$ on any $N\subset W\subset Q$ will not exceed the integral of
$\Lambda$ on $N^{(1)}\subset M$.

In both the strong and weak settings, we will only have $N_0$ compete
against submanifolds having the same boundary.  For this reason, we
assume that
$W\supset N_0$ is chosen so that $\partial N_0=N_0\cap\partial W$, and
that $(N_0,\partial N_0)$ generates the relative homology 
$H_n(W,\partial W;\mathbf Z)$.
\begin{Proposition}
If there exists a calibration $\tilde\Lambda$ for $(\Lambda,N_0)$, 
then $\mathcal
F_\Lambda(N_0)\leq\mathcal F_\Lambda(N)$ for every hypersurface
$N\hookrightarrow W\subset Q$ satisfying $\partial N=\partial N_0$.
\label{CalProp15}
\end{Proposition}
We then say that $N_0$ is a {\em strong (but not {\em strict}!) local
  minimum} for $\mathcal F_\Lambda$.

\

\begin{Proof}
We simply calculate
\begin{eqnarray*}
  \mathcal F_\Lambda(N_0) & = & \int_{N^{(1)}_0}\Lambda \\
    & = & \int_{N_0}\tilde\Lambda \\
    & = & \int_N\tilde\Lambda \\
    & \leq & \int_{N^{(1)}}\Lambda \\
    & = & \mathcal F_\Lambda(N).
\end{eqnarray*}
The third equality uses Stokes' theorem\index{Stokes' theorem|(}, 
which applies because our
topological hypothesis on $W$ implies that the cycle $N-N_0$ in $W$ is
a boundary.
\end{Proof}

\

The question of when one can find a calibration naturally arises.  For
this, we use the following classical concept.
\begin{Definition}
A {\em field}\index{field|(} 
for $(\Lambda, N_0)$ is a neighborhood $W\subset Q$ of
$N_0$ with a smooth foliation by a $1$-parameter family
$F:N\times(-\varepsilon,\varepsilon)\to W$ of integral manifolds of
$\mathcal E_\Lambda$.
\end{Definition}
This family does {\em not} have a fixed boundary.
We retain the topological hypotheses on $W$ used in
Proposition~\ref{CalProp15}, and have the following.
\begin{Proposition}
If there exists a field for $(\Lambda, N_0)$, then there exists a
closed form $\tilde\Lambda\in\Omega^n(W)$ such that
$\tilde\Lambda|_{N_0} = \Lambda|_{N^{(1)}_0}$.
\label{FieldProp15}
\end{Proposition}
$\tilde\Lambda$ is then a calibration if it additionally
satisfies the third condition, $\tilde\Lambda|_E\leq\Lambda_E$.
In the proof, we will explicitly construct $\tilde\Lambda$ using the
Poincar\'e-Cartan form.
\index{Poincar\'e-Cartan form|)}

\

\begin{Proof}
The field $F:N\times(-\varepsilon,\varepsilon)\to Q$ may be thought of
as a family of graphs
$$
  N_t = \{(x,z(x,t))\},
$$
where each $z(\cdot,t)$ is a solution of the Euler-Lagrange equations,
and the domain of $z(\cdot,t)$ may depend on
$t\in(-\varepsilon,\varepsilon)$.
Because each point of $W$ lies on exactly one of these graphs, we
can define a $1$-jet lift $F^\prime:W\to M$, given by
$$
  (x,z)=(x,z(x,t))\mapsto (x,z(x,t),\nabla_xz(x,t)).
$$
Let $\tilde\Lambda=(F^\prime)^*\Lambda\in\Omega^n(W)$.  Then it is
clear that $\tilde\Lambda|_{N_0}=\Lambda|_{N^{(1)}_0}$, and to show
that $\tilde\Lambda$ is closed, we need to see that $(F^\prime)^*\Pi =
d\tilde\Lambda = 0$.  This holds because $\Pi=\theta\wedge\Psi$ is
quadratic in an ideal of forms vanishing on each leaf $F^{(1)}_t:N\to
M$; more concretely, each of $(F^\prime)^*\theta$ and
$(F^\prime)^*\Psi$ must be a multiple of $dt$, so their product vanishes.
\end{Proof}

\

General conditions for $\tilde\Lambda =(F^\prime)^*\Lambda$ to be a
calibration, and for $N_0$ to therefore 
be a strong local minimum, are not clear.  However, we can still use
the preceding to detect weak local minima.
\begin{Proposition}
Under the hypotheses of Propositions~\ref{CalProp15}
and~\ref{FieldProp15}, if
\begin{equation}
  L_{p_ip_j}(x,z_0(x),\nabla z_0(x))\xi_i\xi_j\geq c||\xi||^2,
\label{Hessian15}
\end{equation}
for some constant $c>0$ and all $(\xi_i)$, 
then $\tilde\Lambda|_E\leq\Lambda_E$ for all
$E\subset T_qQ$ sufficiently near $T_qN_0$, with equality if and only
if $E=T_qN_0$.  Furthermore, $\mathcal
F_\Lambda(N_0)<\mathcal F_\Lambda(N)$ for all $N\neq N_0$ in a weak
neighborhood of $N_0$.
\end{Proposition}
The first statement allows us to
think of $\tilde\Lambda$ as a {\em weak calibration} for
$(\Lambda,N_0)$.  The proof of the second statement from the first
will use Stokes' theorem\index{Stokes' theorem|)} in exactly the manner of
Proposition~\ref{CalProp15}.

\

\begin{Proof}
The positivity of the $\nabla z$-Hessian of $L$ suggests that we
define the {\em Weierstrass excess function}
\index{Weierstrass excess function|(}
$$
  E(x,z,p,q)\stackrel{\mathit{def}}{=}
   L(x,z,p)-L(x,z,q)-\ss(p_i-q_i)L_{q_i}(x,z,q),
$$
which is the second-order remainder in a Taylor series expansion
for $L$.  This function will appear in a more detailed expression
for $\pi^*\tilde\Lambda=(F^\prime\circ\pi)^*\Lambda\in\Omega^n(M)$, 
computed modulo $\{I\}$.  We write
$$
  F^\prime(x^i,z)=(x^i,z,q_i(x,z))\in M,
$$
where $(x^i,z,p_i)$ are the usual coordinates on $M$, and the
functions $q_i(x,z)$ are the partial derivatives of the field
elements $z(x,t)$.
We have
\begin{eqnarray*}
  \pi^*\tilde\Lambda & = & \pi^*\circ F^{\prime*}(L\,dx+\theta\wedge 
     L_{p_i}dx_{(i)}) \\
  & = & L(x^i,z,q_i(x,z))dx+(dz-q_i(x,z)dx^i)\wedge
     L_{p_i}(x^i,z,q_i(x,z))dx_{(i)} \\
  & \equiv & \left(L(x^i,z,q_i(x,z))+(p_i-q_i(x,z))
   L_{p_i}(x^i,z,q_i(x,z))\right)dx \pmod{\{I\}} \\
  & \equiv & \Lambda-E(x^i,z,p_i,q_i(x,z))dx \pmod{\{I\}}.
\end{eqnarray*}
The hypothesis (\ref{Hessian15}) on the Hessian $F_{p_ip_j}$ implies
that for each
$(x^i,z)$, and $p_i$ sufficiently close to $q_i(x,z)$, the
second-order remainder satisfies
$$ E(x^i,z,p_i,q_i(x,z))\geq0, $$
with equality if and only if $p_i=q_i(x,z)$.  
The congruence of $\pi^*\tilde\Lambda$ and $\Lambda-E\,dx$ modulo $\{I\}$
then implies our first statement.

For the second statement, we use the Stokes' theorem argument:
\begin{eqnarray*}
  \int_{N^{(1)}_0}\Lambda & = & \int_{N_0}\tilde\Lambda \\ & = &
    \int_N\tilde\Lambda \\
  & = & \int_{N^{(1)}}(\Lambda-E(x^i,z,p_i,q_i(x,z))dx) \\ & \leq &
  \int_{N^{(1)}}\Lambda,
\end{eqnarray*}
with equality in the last step if and only if $N_0=N$.
\end{Proof}
\index{calibration|)}

\

This proof shows additionally that if the Weiestrass excess function satisfies
$E(x^i,z,p_i,q_i)>0$ for all $p\neq q$, then $N_0$ is a strong (and
strict) local minimum for $\mathcal F_\Lambda$.
\index{Weierstrass excess function|)}
\index{strong local minimum|)}
\index{weak local minimum|)}

So far, we have shown that if we can cover
some neighborhood of a stationary submanifold $N_0$ with a field, then we can
construct an $n$-form $\tilde\Lambda$, whose calibration properties
imply extremal properties of $N_0$.
It is therefore natural to ask when there exists such a field, and the
answer to this involves some analysis of the {\em Jacobi
  operator}.\index{Jacobi operator|(}
We will describe the operator, and hint at the analysis.

The Jacobi operator acts on sections of a density line
bundle\index{density line bundle} on a
given integral manifold $N_0$ of the Euler-Lagrange system, with its
induced conformal structure\index{conformal!structure}.  Specifically, 
$$ J:D^{\frac{n-2}{2n}}\to D^{\frac{n+2}{2n}} $$ 
is the differential operator given by
$$
  Jg = -\Delta_cg-gK,
$$
where $\Delta_c$ is the conformal Laplacian, and $K$ is the
curvature invariant introduced in \S\ref{Subsection:Relative}. 
The second variation formula (\ref{GoodSecondVarn}) then reads
$$
  \delta^2(\mathcal F_\Lambda)_{N_0}(g) = \int_{N_0}g\,J\!g\,\omega.
$$
The main geometric fact is:
\begin{quote}
{\em The Jacobi operator gives the (linear) variational equations for
  integral manifolds of the Euler-Lagrange system $\mathcal
  E_\Lambda$.}
\end{quote}
\index{Euler-Lagrange!system|)}
This means the following.  Let $F:N\times[0,1]\to M$ be a Legendre
variation\index{Legendre variation|)} 
of the $\Lambda$-stationary submanifold\index{stationary|)} $F_0$---not
necessarily having fixed boundary---and let $g=\left(\sf{\partial
  F}{\partial t}|_{t=0}\right)\innerprod\theta$, as usual.  
Then our previous calculations imply that $Jg=0$ if and only if
$$
  \mathcal L_{\frac{\partial}{\partial t}}(F^*\Psi)|_{t=0}=0.
$$
We might express condition by saying that
$F_t$ is an integral manifold
for $\mathcal E_\Lambda=\mathcal I+\{\Psi\}$ modulo $O(t^2)$.

We now indicate how a condition on the Jacobi operator of $N_0$ can
imply the existence of a field near $N_0$.  Consider the eigenvalue
problem
$$
  Jg=-\Delta_cg-gK=\lambda g,\qquad g\in C^\infty_0(N),
$$
for smooth, fixed boundary variations.  It is well-known $J$ has a
discrete spectrum bounded from below, $\lambda_1<\lambda_2<\cdots$,
with $\lambda_k\to\infty$ and with finite-dimensional eigenspaces.
We consider the consequences of the assumption
$$
  \lambda_1>0.
$$
Because $\lambda_1=\inf\{\int_Ng\,J\!g\,\omega:||g||_{L^2}=1\}$, the assumption
$\lambda_1>0$ is equivalent to
$$
  \delta^2(\mathcal F_\Lambda)_N(g)>0,\qquad\mbox{for }g\neq 0.
$$
The main analytic result is the following.
\begin{Proposition}  If $\lambda_1>0$, then given $g_0\in
  C^\infty(\partial N)$, there is a unique solution $g\in C^\infty(N)$ to the boundary
  value problem
$$
  Jg=0,\qquad g|_{\partial N}=g_0.
$$
Furthermore, if $g_0>0$ on $\partial N$, then this solution satisfies
$g>0$ on $N$.
\end{Proposition}
The existence and uniqueness statements follow from standard elliptic
theory.  The point is that we can compare
the second variation
$\delta^2(\mathcal F_{\Lambda})_N(g)$ to the Sobolev norm
$||g||^2_1=\int_N(||\nabla g||^2+|g|^2)\omega$, and if $\lambda_1>0$,
then there are constants $c_1,c_2>0$ such that
$$
  c_1||g||_1^2\leq \int_Ng\,Jg\,\omega\leq c_2||g||_1^2.
$$
The Schauder theory gives existence and uniqueness in this situation.

Less standard is the positivity of the solution $g$ under the
assumption that $g|_{\partial N}>0$, and this is crucial for the
existence of a field.  Namely, a further implicit function argument
using elliptic theory guarantees that the variation $g$ is tangent to
an arc of integral manifolds of $\mathcal E$, and the fact that $g\neq
0$ implies that near the initial $N$, this arc defines a field.  For
the proof of the positivity of $g$, and details of all of the
analysis, see Giaquinta \& Hildebrandt (cit.~p.~\pageref{GHPage}n).
\index{field|)}
\index{Jacobi operator|)}

\section{Euler-Lagrange PDE Systems}
\label{Section:MultiPC}

Up to this point, we have studied geometric aspects of
first-order Lagrangian functionals
\begin{equation}
  {\mathcal F}_L(z) = \int_\Omega 
     L\left(x^i,z,\sf{\partial z}{\partial x^i}\right)dx,
  \qquad \Omega\subset\R^n,
\label{OneFunctional}
\end{equation}
where $x=(x^1,\ldots,x^n)$ and $z=z(x)$ is a scalar function.  In this
section, we consider the more general situation of functionals
\begin{equation}
  {\mathcal F}_L(z) =
  \int_\Omega L\left(x^i,z^\alpha(x),\sf{\partial z^\alpha}{\partial
    x^i}(x)\right)dx,\qquad \Omega\subset\R^n,
\label{MultiFunctional}
\end{equation}
where now $z(x)=(z^1(x),\ldots,z^s(x))$ is an $\R^s$-valued function
of $x=(x^i)$, and $L=L(x^i,z^\alpha,p^\alpha_i)$ is a smooth function
on $\R^{n+s+ns}$.  The Euler-Lagrange
equations\index{Euler-Lagrange!equation} describing maps 
$z:\Omega\to\R^s$ which are stationary for $\mathcal{F}_L$ under all
fixed-boundary variations form a PDE system
\begin{equation}
  \frac{\partial L}{\partial z^\alpha}-\sum_i\frac{d}{dx^i}
    \left(\frac{\partial L}{\partial p^\alpha_i}\right) = 0,\qquad
  \alpha=1,\ldots,s.
\label{ELEqn16}
\end{equation}
In the scalar case $s=1$, we have examined the geometry of the
equivalence class of $\mathcal F_L$ under contact
transformations\index{transformation!contact|(} and
found the canonically defined Poincar\'e-Cartan
form\index{Poincar\'e-Cartan form|(} to be of considerable
use.  In this section, we describe a generalization of the
Poincar\'e-Cartan form for $s\geq 1$.  Geometrically, we
study functionals on the space of compact submanifolds of codimension $s$, in
an $(n+s)$-dimensional manifold with local coordinates $(x^i,z^\alpha)$.

An immediate difference between the cases $s=1$ and $s\geq 2$ is that
in the latter case, there are no proper contact transformations of
$\R^{n+s+ns}$; that is, the only smooth maps
$x^\prime=x^\prime(x,z,p)$, $z^\prime=z^\prime(x,z,p)$,
$p^\prime=p^\prime(x,z,p)$ for which
$$
  \{dz^\alpha-\ss p^\alpha_idx^i\} =
    \{dz^{\alpha\prime}-\ss p^{\alpha\prime}_idx^{i\prime}\}
$$
are point transformations\index{transformation!point|(} 
$x^\prime=x^\prime(x,z)$,
$z^\prime=z^\prime(x,z)$, with $p^\prime=p^\prime(x,z,p)$ determined
by the chain rule.  We will explain why this is so, and later, we will
see that in case $s=1$ our original contact-invariant
Poincar\'e-Cartan form still appears naturally in the more limited
context of point transformations.  Our first task, however, is to
introduce the geometric setting for studying functionals
(\ref{MultiFunctional}) subject to point transformations, analogous to
our use of contact manifolds for (\ref{OneFunctional}).

Throughout this section, we have as always $n\geq 2$ and we use the index
ranges $1\leq i,j\leq n$, $1\leq\alpha,\beta\leq s$.

\subsection{Multi-contact Geometry}
\index{multi-contact geometry|(}

Having decided to apply point transformations to the functional
(\ref{MultiFunctional}), we interpret $z(x)=(z^\alpha(x^i))$
as corresponding to an $n$-dimensional submanifold of $\R^{n+s}$.  The first
derivatives $p^\alpha_i=\sf{\partial z^\alpha}{\partial x^i}$ specify
the tangent $n$-planes of this submanifold.  This suggests our first
level of geometric generalization.

Let $X$ be a manifold of dimension $n+s$, and let
$G_n(TX)\stackrel{\pi}{\to}X$ be
the Grassmannian\index{Grassmannian} 
bundle of $n$-dimensional subspaces of tangent spaces
of $X$; that is, a point of $G_n(TX)$ is of the form
$$
  m=(p,E),\qquad p\in X,\ E^n\subset T_pX.
$$
Any diffeomorphism of $X$ induces a diffeomorphism of $G_n(TX)$,
and either of these diffeomorphisms will be called a {\em point
  transformation}.

We can define on $G_n(TX)$ two Pfaffian systems $I\subset J\subset
T^*(G_n(TX))$, of ranks
$s$ and $n+s$, respectively, which are canonical in the sense that they are
preserved by any point transformation.  
First, $J=\pi^*(T^*X)$ consists of all forms that are semibasic over $X$; $J$
is integrable, and its maximal integral submanifolds are the fibers of
$G_n(TX)\to X$.
Second, we define $I$ at a point $(p,E)\in G_n(TX)$ to be 
$$
  I_{(p,E)} = \pi_p^*(E^\perp),
$$
where $E^\perp\subset T_p^*X$ is the $s$-dimensional annihilator of the
subspace $E\subset T_pX$.  
$I$ is not integrable, and to understand its integral submanifolds,
note that any $n$-dimensional immersion
$\iota:N\hookrightarrow X$ has a $1$-jet lift $\iota^{(1)}:N\hookrightarrow
G_n(TX)$.  In fact, such lifts are the transverse
integral submanifolds of the Pfaffian system $I\subset T^*(G_n(TX))$.

To see this explicitly, choose local coordinates $(x^i,z^\alpha)$
on $U\subset X^{n+s}$.  These induce local coordinates
$(x^i,z^\alpha,p^\alpha_i)$ corresponding to the $n$-plane $E\subset
T_{(x^i,z^\alpha)}U$ defined as
$$
  E=\{dz^1-p^1_idx^i,\ldots,dz^s-p^s_idx^i\}^\perp.
$$
These coordinates are defined on a dense open subset of
$\pi^{-1}(U)\subset G_n(TX)$, consisting of $n$-planes $E\subset TX$
for which $dx^1\wedge\cdots\wedge dx^n|_E\neq 0$.  In terms of these
local coordinates on $G_n(TX)$, our Pfaffian systems are
\begin{eqnarray*}
  J & = & \{dx^i,dz^\alpha\}, \\
  I & = & \{dz^\alpha-p^\alpha_idx^i\}.
\end{eqnarray*}
An immersed submanifold $N^n\hookrightarrow U$ for which
$dx^1\wedge\cdots\wedge dx^n|_N\neq 0$ may be regarded as a graph
$$
  N = \{(x^i,z^\alpha):z^\alpha = f^\alpha(x^1,\ldots,x^n)\}.
$$
Its lift to $N\hookrightarrow G_n(TX)$ lies in the
domain of the coordinates $(x^i,z^\alpha,p^\alpha_i)$, and equals the
$1$-jet graph
\begin{equation}
  N^{(1)} = \{(x^i,z^\alpha,p^\alpha_i):
    z^\alpha = f^\alpha(x^1,\ldots,x^n),\
    p^\alpha_i = \sf{\partial f^\alpha}{\partial x^i}
        (x^1,\ldots,x^n)\}.
\label{MultiGraph}
\end{equation}
Clearly this lift is an integral submanifold of $I$.  Conversely, a
submanifold $N^{(1)}\hookrightarrow\pi^{-1}(U)\subset G_n(TX)$ on
which $dz^\alpha-\sum p^\alpha_idx^i=0$ and $dx^1\wedge\cdots\wedge
dx^n\neq 0$ is necessarily given locally by a graph of the form
(\ref{MultiGraph}).  The manifold $M=G_n(TX)$ with its Pfaffian
systems $I\subset J$ is our standard example of a {\em multi-contact
  manifold}.  This notion will be defined shortly, in terms of the
following structural properties of the Pfaffian systems.

Consider on $G_n(TX)$ the differential ideal
$\mathcal{I}=\{I,dI\}\subset\Omega^*(G_n(TX))$ generated by $I\subset
T^*(G_n(TX))$. 
If we set 
$$
  \bar\theta^\alpha=dz^\alpha-p^\alpha_idx^i,\quad 
  \bar\omega^i=dx^i,\quad \bar\pi^\alpha_i=dp^\alpha_i, 
$$
then we have
the structure equations
\begin{equation}
  d\bar\theta^\alpha \equiv -\bar\pi^\alpha_i\wedge\bar\omega^i
     \pmod{\{I\}},\quad 1\leq\alpha\leq s.
\label{BasicMultiStructure}
\end{equation}
It is not difficult to verify that the set of all coframings
$(\theta^\alpha,\omega^i,\pi^\alpha_i)$ on $G_n(TX)$ for which
\begin{itemize}
\item $\theta^1,\ldots,\theta^s$ generate $I$,
\item $\theta^1,\ldots,\theta^s,\omega^1,\ldots,\omega^n$ generate
  $J$, and
\item $d\theta^\alpha\equiv-\pi^\alpha_i\wedge\omega^i \pmod{\{I\}}$
\end{itemize}
are the local sections of a
$G$-structure on $G_n(TX)$.   Here $G\subset GL(n+s+ns,\R)$ may be
represented as acting on $(\theta^\alpha,\omega^i,\pi^\alpha_i)$ by
\begin{equation}
  \left\{\begin{array}{l}
    \bar\theta^\alpha = a^\alpha_\beta\theta^\beta, \\
    \bar\omega^i = c^i_\beta\theta^\beta+b^i_j\omega^j, \\
    \bar\pi^\alpha_i = d^\alpha_{i\beta}\theta^\beta +
       e^\alpha_{kj}(b^{-1})^k_i\omega^j + a^\alpha_\beta
           \pi^\beta_j(b^{-1})^j_i,
  \end{array}\right.
\label{DefMultiG}
\end{equation}
where $(a^\alpha_\beta)\in GL(s,\R)$, $(b^i_j)\in GL(n,\R)$, and
$e^\alpha_{ij}=e^\alpha_{ji}$.
From these properties we make our definition.
\begin{Definition}  A {\em multi-contact manifold} is a
  manifold $M^{n+s+ns}$, with a $G$-structure as in
  (\ref{DefMultiG}), whose sections
  $(\theta^\alpha,\omega^i,\pi^\alpha_i)$ satisfy 
\begin{eqnarray}
  d\theta^\alpha  & \equiv &  -\pi^\alpha_i\wedge\omega^i
     \pmod{\{\theta^1,\ldots,\theta^s\}},
  \label{GenMultiStructure} \\
  d\omega^i  & \equiv &  0 \pmod{\{\theta^1,\ldots,\theta^s,\omega^1,\ldots,
    \omega^n\}}.
  \label{GenMultiInteg}
\end{eqnarray}
\label{DefMulti16}
\end{Definition}
Note that the $G$-structure determines Pfaffian systems
$I=\{\theta^1,\ldots,\theta^s\}$ and
$J=\{\theta^1,\ldots,\theta^s,\omega^1,\ldots,\omega^n\}$, and we may
often refer to $(M,I,J)$ as a multi-contact manifold, implicitly
assuming that $J$ is integrable and that there are coframings for
which the structure equations (\ref{BasicMultiStructure}) hold.
The integrability of $J$ implies that locally in $M$ one can define a
smooth leaf space $X^{n+s}$ and a surjective submersion $M\to X$ whose
fibers are integral manifolds of $J$.  When working locally in a
multi-contact manifold, we will often make reference to this
quotient $X$.

It is not difficult to show that any multi-contact structure $(M,I,J)$ is
locally equivalent to that of $G_n(TX)$ for a manifold $X^{n+s}$.  The
integrability of $J$ implies that there are local coordinates
$(x^i,z^\alpha,q^\alpha_i)$ for which $dx^i,dz^\alpha$ generate $J$.
We can relabel the $x^i,z^\alpha$ to assume that
$dz^\alpha-p^\alpha_idx^i$ generate $I$ for
some functions $p^\alpha_i(x,z,q)$.  The structure equations then
imply that $dx^i,dz^\alpha,dp^\alpha_i$ are linearly independent, so
on a possibly smaller neighborhood in $M$, we can replace the
coordinates $q^\alpha_i$ by $p^\alpha_i$, and this exhibits our
structure as equivalent to that of $G_n(TX)$.

We will see below that if $s\geq 2$, then the Pfaffian system
$I$ of a multi-contact manifold uniquely determines the larger system
$J$.
Also, if $s\geq 3$, the hypothesis (\ref{GenMultiInteg}) that
$J=\{\theta^\alpha,\omega^i\}$ is integrable is not necessary; it is
easily seen to be a consequence of the structure equation
(\ref{GenMultiStructure}).  However, in the case $s=2$, $J$ is
determined by $I$ but is not necessarily integrable; our study of
Euler-Lagrange systems will not involve this exceptional situation, so we
have ruled it out in our definition.

It is not at all obvious how one can
determine, given a Pfaffian system $I$ of rank $s$ on a manifold $M$
of dimension $n+s+ns$, whether $I$ comes from a multi-contact structure;
deciding whether structure equations (\ref{GenMultiStructure}) can be
satisfied for some generators of $I$ is a difficult problem.  Bryant has
given easily evaluated intrinsic criteria characterizing such $I$,
generalizing the Pfaff theorem's\index{Pfaff theorem} 
normal form for contact manifolds,
but we shall not need this here (see Ch.~II, \S4 of
  \cite{Bryant:Exterior}).

Aside from those of the form $G_n(TX)$, there are two other kinds of
multi-contact manifolds in common use.  One is $J^1(Y^n,Z^s)$, the
space of $1$-jets of maps from an $n$-manifold $Y$ to an $s$-manifold
$Z$.  The other is $J^1_\Gamma(E^{n+s},Y^n)$, the space of $1$-jets of
sections of a fiber bundle $E\to Y$ with base of dimension $n$ and fiber
$Z$ of dimension $s$.  These are distinguished by the kinds of coordinate
changes considered admissible in each case; to the space
$J^1(Y^n,Z^s)$, one would apply prolonged classical
transformations\index{transformation!classical}
$x^\prime(x)$, $z^\prime(z)$, while to 
$J^1_\Gamma(E^{n+s},Y^n)$, one would apply prolonged gauge
transformations\index{transformation!gauge}
$x^\prime(x)$, $z^\prime(x,z)$.  These are both smaller classes than
the point transformations $x^\prime(x,z)$, $z^\prime(x,z)$ that we
apply to $G_n(TX)$, the space of $1$-jets of $n$-submanifolds in
$X^{n+s}$.

Recall our claim that in the multi-contact case $s\geq 2$, every
contact transformation is a prolonged point transformation.  
This is the same as saying that any local diffeomorphism
of $M$ which preserves the Pfaffian system $I$ also preserves $J$; for
a local diffeomorphism of $M$ preserving $J$ must induce a
diffeomorphism of the local quotient space $X$, which in turn
uniquely determines the original local diffeomorphism of $M$.  
To see why a local diffeomorphism preserving $I$ must preserve $J$, we
will give an intrinsic construction of $J$ in terms of $I$ alone, for
the local model $G_n(TX)$.
First, define for any $2$-form $\Psi\in\bw{2}(T_m^*M)$ the space of
$1$-forms
$$
  \mathcal{C}(\Psi)=\{V\innerprod\Psi:V\in T_mM\}.
$$
This is a pointwise construction.  We apply it to each element of
the vector space
$$
  \{\lambda_\alpha d\theta^\alpha:(\lambda^\alpha)\in\R^s\},
$$
intrinsically given as the quotient of $\mathcal I_2$, the degree-$2$
part of the multi-contact differential ideal, by the subspace
$\{I\}_2$, the degree-$2$ part of the {\em algebraic}
ideal\index{ideal!algebraic} $\{I\}$.  For example,
$$
  \mathcal C(\overline{d\theta^\alpha}) \equiv
    \mbox{Span}\{\pi^\alpha_i,\omega^i\}\pmod{\{I\}}.
$$
The intersection
$$
  \bigcap_{\Theta\in\mathcal I_2/\{I\}_2}{\mathcal C}(\Theta)
$$
is a well-defined subbundle of $T^*M/I$.  If $s\geq 2$, then its
preimage in $T^*M$ is $J=\{\bar\theta^\alpha,\bar\omega^i\}$, as is
easily seen using the
structure equations (\ref{BasicMultiStructure}).
Any local diffeomorphism of $M$ preserving $I$ therefore preserves
$\mathcal I_2$, $\mathcal I_2/\{I\}_2$, $\bigcap\mathcal C(\Theta)$, and
finally $J$, which is what we wanted to prove.  Note that in the
contact case
$s=1$, $\bigcap\mathcal C(\Theta)\equiv\mathcal
C(d\theta)\equiv\{\pi_i,\omega^i\}$ modulo $\{I\}$, so instead of this
construction giving $J$, it gives all of $T^*M$.  
In this case, introducing $J$ in the definition of a multi-contact
manifold reduces our
pseudogroup from contact transformations to point transformations.
\index{transformation!contact|)}
\index{transformation!point|)}

We have given a generalization of the notion of a contact manifold to
accomodate the study of submanifolds of codimension greater than one.
There is a further generalization to higher-order contact geometry
which is the correct setting for studying higher-order Lagrangian
functionals, and we will consider it briefly in the next section.

In what follows, we will carry out the discussion of functionals
modelled on (\ref{MultiFunctional}) on a general multi-contact
manifold $(M,I,J)$, but the reader can concentrate on the case $M=G_n(TX)$.

\subsection{Functionals on Submanifolds of Higher Codimension}

Returning to our functional (\ref{MultiFunctional}), we
think of
the integrand $L(x^i,z^\alpha,p^i_\alpha)dx$
as an $n$-form on a dense open subset of the multi-contact
manifold $G_n(T\R^{n+s})$.  Note that this $n$-form is semibasic for
the projection $G_n(T\R^{n+s})\to\R^{n+s}$, and that any $n$-form
congruent to $L(x^i,z^\alpha,p^i_\alpha)dx$ modulo
$\{dz^\alpha-p^\alpha_idx^i\}$ gives the same classical
functional.  This suggests the following.

\begin{Definition}  A {\em Lagrangian}\index{Lagrangian|(}
 on a multi-contact manifold
  $(M,I,J)$ is a smooth section
  $\Lambda\in\Gamma(M,\bigwedge^nJ)\subset\Omega^n(M)$.  Two
  Lagrangians are {\em equivalent} if they are congruent modulo
  $\{I\}$.
\end{Definition}
An equivalence class $[\Lambda]$ of Lagrangians corresponds to a
section of the vector bundle $\bigwedge^n(J/I)$.  It also
defines a functional on the space of compact integral manifolds (possibly with
boundary) of the Pfaffian system $I$ by
$$
  \mathcal F_\Lambda(N)=\int_N\Lambda,
$$
where $\Lambda$ is any representative of the class.  The notion of
{\em divergence equivalence}\index{divergence equivalence} 
of Lagrangians will appear later.  In the
discussion in Chapter 1 of the scalar case $s=1$, we combined these
two types of equivalence by emphasizing a characteristic cohomology
class in $H^n(\Omega^*(M)/\mathcal I)$, and used facts about
symplectic linear algebra to investigate these classes.  However, the
analogous ``multi-symplectic'' linear algebra that is appropriate for
the study of multi-contact geometry is still
poorly understood.\footnote{The \label{Grassi} 
 recent work \cite{Grassi:Local} of M.~Grassi\index{Grassi, M.|nn}
    may illuminate this issue, along with some others that will come
    up in the following discussion.}

Our goal is to associate to any functional
$[\Lambda]\in\Gamma(M,\bigwedge^n(J/I))$ a Lagrangian
$\Lambda\in\Gamma(M,\bigwedge^nJ)\subset\Omega^n(M)$, not necessarily
uniquely determined, whose exterior
derivative $\Pi=d\Lambda$ has certain
favorable properties and {\em is} uniquely determined by $[\Lambda]$.  Among 
these properties are:
\begin{itemize}
\item $\Pi\equiv 0\mbox{ (mod ${\{I\}}$)}$;
\item $\Pi$ is preserved under any diffeomorphism of $M$ preserving
  $I$, $J$, and $[\Lambda]$;
\item $\Pi$ depends only on the
  divergence-equivalence\index{divergence equivalence} class of
  $[\Lambda]$;
\item $\Pi=0$ if and only if the Euler-Lagrange
  equations\index{Euler-Lagrange!equation} for
  $[\Lambda]$ are trivial.
\end{itemize}
Triviality of the Euler-Lagrange equations means that every compact
integral 
manifold of $\mathcal I\subset\Omega^*(M)$ is stationary for $\mathcal
F_\Lambda$ under fixed-boundary variations.
Some less obvious ways in which such $\Pi$ could be useful are the
following, based on our experience in the scalar case $s=1$:
\begin{itemize}
\item in Noether's theorem, where one would hope for $v\mapsto
  v\innerprod\Pi$ to give an isomorphism from a Lie algebra of
  symmetries to a space of conservation laws;
\item in the inverse problem, where one can try to detect equations
  that are locally of Euler-Lagrange type not by finding a Lagrangian,
  but by finding a Poincar\'e-Cartan form inducing the equations;
\item in the study of local minimization, where it could help one obtain a
    calibration\index{calibration} 
    in terms of a field of stationary submanifolds.
\end{itemize}

Recall that in the case of a contact manifold, we replaced any
Lagrangian $\Lambda\in\Omega^n(M)$ by
$$
  \Lambda -\theta\wedge\beta,
$$
the unique form congruent to $\Lambda\mbox{ (mod $\{I\}$)}$ with the property
that
$$
  d\Lambda\equiv 0\pmod{\{I\}}.
$$
What happens in the multi-contact case?  Any Lagrangian
$\Lambda_0\in\Gamma(M,\bigwedge^nJ)$ is congruent modulo $\{I\}$ to a form
(in local coordinates)
$$
  L(x^i,z^\alpha,p^\alpha_i)dx,
$$
and motivated by the scalar case, we consider the equivalent form
\begin{equation}
  \Lambda = L\,dx+\theta^\alpha\wedge\sf{\partial L}{\partial
    p^\alpha_i}dx_{(i)},
\label{FirstLift16}
\end{equation}
which has exterior derivative
\begin{equation}
  d\Lambda = \theta^\alpha\wedge\left(\sf{\partial L}{\partial
    z^\alpha}dx - d\left(\sf{\partial L}{\partial
    p^\alpha_i}\right)\wedge dx_{(i)}\right).
\label{FirstPC16}
\end{equation}
This suggests the following definition.
\begin{Definition}  An {\em admissible 
   lifting}\index{admissible lifting|(} of a functional
  $[\Lambda]\in\Gamma(M,\bigwedge^n(J/I))$ is a Lagrangian
  $\Lambda\in\Gamma(M,\bigwedge^nJ)$ representing the class
  $[\Lambda]$ and satisfying $d\Lambda\in\{I\}$.
\end{Definition}
The preceding calculation shows that locally, every functional
$[\Lambda]$ has an admissible lifting.
Unfortunately, the admissible lifting is generally not
unique.  This will be addressed below, but first we show that any
admissible lifting is adequate for calculating the first variation
and the Euler-Lagrange system of the functional $\mathcal F_\Lambda$.

We mimic the derivation in Chapter 1 of the Euler-Lagrange
differential system in the scalar case $s=1$.
Suppose that we have a $1$-parameter family $\{N_t\}$ of integral
manifolds of the multi-contact Pfaffian system $I$, given as a smooth
map
$$
  F:N\times[0,1]\to M,
$$
for which each $F_t=F|_{N\times\{t\}}:N\hookrightarrow M$ is an
integral manifold of $I$ and such that $F|_{\partial N\times[0,1]}$
is independent of $t$.  Then choosing generators
$\theta^\alpha\in\Gamma(I)$, $1\leq\alpha\leq s$, we have
$$
  F^*\theta^\alpha = G^\alpha dt
$$
for some functions $G^\alpha$ on $N\times[0,1]$.  As in the contact case, it
is not difficult to show that any collection of functions $g^\alpha$
supported in the interior of $N$ can be realized as $G^\alpha|_{t=0}$
for some $1$-parameter family $N_t$.

The hypothesis that $\Lambda$ is an admissible lifting means that we
can write
$$
  d\Lambda = \ss\theta^\alpha\wedge\Psi_\alpha
$$
for some $\Psi_\alpha\in\Omega^n(M)$.  Then we can proceed as in
\S\ref{Section:EulerLagrange} to calculate
\begin{eqnarray*}
  \left.\frac{d}{dt}\right|_{t=0}\left(\int_{N_t}F_t^*\Lambda\right)
    & = & \int_{N_0}\mathcal L_{\frac{\partial}{\partial t}}(F^*\Lambda) \\
    & = & \int_{N_0}\sf{\partial}{\partial t}\innerprod
            (\theta^\alpha\wedge\Psi_\alpha) +
          \int_{N_0}d(\sf{\partial}{\partial t}\innerprod\Lambda) \\
    & = & \int_{N_0}g^\alpha\Psi_\alpha,
\end{eqnarray*}
where in the last step we used the fixed-boundary condition, the vanishing
of $F_0^*\theta^\alpha$, and the definition
$g^\alpha=G^\alpha|_{N_0}$.  Now the same reasoning as in
\S\ref{Section:EulerLagrange} shows that $F_0:N\hookrightarrow M$ is
stationary for $\mathcal F_\Lambda$ under all fixed-boundary variations
if and only if $\Psi_\alpha|_{N_0}=0$ for all $\alpha=1,\ldots,s$.

We now have a differential system
$\{\theta^\alpha,d\theta^\alpha,\Psi_\alpha\}$  whose integral manifolds are
exactly the integral manifolds of $\I$ that are stationary for $[\Lambda]$, but
it is not clear that this system is uniquely determined by $[\Lambda]$
alone; we might get different systems for different admissible liftings.
To rule out this possibility, observe first that if $\Lambda$,
$\Lambda^\prime$ are any two admissible liftings of $[\Lambda]$, then the
condition 
$\Lambda-\Lambda^\prime\in\{I\}$ allows us to write
$\Lambda-\Lambda^\prime=\theta^\alpha\wedge\gamma_\alpha$, and then
the fact that $d\Lambda\equiv
d\Lambda^\prime\equiv 0\mbox{ (mod $\{I\}$)}$ along with the structure
equations (\ref{GenMultiStructure}) allows us to write
$$
  0\equiv d\theta^\alpha\wedge\gamma_\alpha\equiv
    -\pi^\alpha_i\wedge\omega^i\wedge\gamma_\alpha
    \pmod{\{I\}}.
$$
When $n\geq 2$, this implies that $\gamma_\alpha\equiv 0\mbox{ (mod
  $\{I\}$)}$, so while two general representatives of $[\Lambda]$ need
be congruent only modulo $\{I\}$, for admissible liftings we have the
following.
\begin{Proposition}
Two admissible liftings of the same
  $[\Lambda]\in\Gamma(\bigwedge^n(J/I))$ are congruent modulo
  $\{\bigwedge^2I\}$.
\label{LiftProp16}
\end{Proposition}
Of course, when $s=1$, $\bigwedge^2I=0$ and we have a unique
lifting, whose derivative is the familiar Poincar\'e-Cartan form.
This explains how the Poincar\'e-Cartan form occurs in the context of
point transformation as well as contact transformations.

We use the proposition as follows.
If we take two admissible liftings $\Lambda$, $\Lambda^\prime$ of the
same functional $[\Lambda]$, and write
$$
  \Lambda-\Lambda^\prime =
  \sf12\theta^\alpha\wedge\theta^\beta\wedge\gamma_{\alpha\beta},
$$
with $\gamma_{\alpha\beta}+\gamma_{\beta\alpha}=0$, then
$$
  \theta^\alpha\wedge(\Psi_\alpha-\Psi_\alpha^\prime)  =
  d(\Lambda-\Lambda^\prime) \equiv
    -\theta^\alpha\wedge d\theta^\beta\wedge\gamma_{\alpha\beta}
      \pmod{\{\textstyle\bigwedge^2I\}}.
$$
A consequence of this is that
$\Psi_\alpha-\Psi_\alpha^\prime\in\mathcal I$ for each $\alpha$, and
we can therefore give the following.

\begin{Definition} The {\em Euler-Lagrange system} $\mathcal E_{\Lambda}$ of
  $[\Lambda]\in\Gamma(M,\bigwedge^n(J/I))$ is the differential
  ideal\index{ideal!differential}
  on $M$ generated by $I$ and the $n$-forms
  $\{\Psi_1,\ldots,\Psi_s\}\subset\Omega^n(M)$, where $\Lambda$ is any
  admissible lifting of $[\Lambda]$ and
  $d\Lambda=\sum\theta^\alpha\wedge\Psi_\alpha$.
A {\em stationary Legendre submanifold}\index{Legendre submanifold} 
of $[\Lambda]$ is an integral
manifold of $\mathcal{E}_{\Lambda}$.
\end{Definition}
\index{multi-contact geometry|)}

\subsection{The Betounes and Poincar\'e-Cartan Forms}
\index{Betounes form|(}

For scalar variational problems, the Poincar\'e-Cartan
form $\Pi\in\Omega^{n+1}(M)$ on the contact manifold $(M,I)$ is an
object of central importance.  Some of its key features were outlined above.
Underlying its usefulness is the fact that we are associating
to a Lagrangian functional---a certain equivalence class of
differential forms---an object that is not merely an equivalence
class, but an actual differential form with which we can carry out
certain explicit computations.  We would like to construct an analogous object
in the multi-contact case.

We will do this by imposing pointwise algebraic conditions on
$\Pi\stackrel{\mathit{def}}{=}d\Lambda$.  Fix an admissible coframing
$(\theta^\alpha,\omega^i,\pi^\alpha_i)$ on a multi-contact manifold
as in Definition \ref{DefMulti16}.
Then any admissible lifting $\Lambda\in\Gamma(\bigwedge^nJ)$ of a
functional $[\Lambda]$ has the form
$$
  \Lambda = \sum_{k=0}^{\mathit{min}(n,s)}\left((k!)^{-2}
    \sum_{|A|=|I|=k}F^I_A\theta^A\wedge\omega_{(I)}\right),
$$
for some functions $F^A_I$, which are skew-symmetric with respect to
each set of indices.  Because $J$ is integrable,
$\Pi = d\Lambda$ lies in $\Omega^{n+1}(M)\cap\{\bigwedge^nJ\}$; the ``highest
weight'' part can be written as
\begin{equation}
\boxed{
  \Pi = \sum_{k=1}^{\mathit{min}(n,s)}\left((k!)^{-2}
   \sum_{\begin{smallmatrix} \alpha, i\\ |A|=|I|=k
     \end{smallmatrix} }
     H^{iI}_{\alpha A} \pi^\alpha_i
     \wedge\theta^A\wedge\omega_{(I)}\right)\pmod{\{\textstyle\bigwedge^{n+1}J\}}.
}
\label{Expansion16}
\end{equation}
The functions $H^{iI}_{\alpha A}$ are skew-symmetric in the
multi-indices $I$ and $A$.  Notice that the equation 
$d\Pi\equiv 0\ (\mbox{mod }\{\bigwedge^{n+1}J\})$ gives for the $k=1$ term
$$
  H^{ij}_{\alpha\beta}=H^{ji}_{\beta\alpha}.
$$

To understand the relevant linear algebra,
suppose that $V^n$ is a vector space with basis $\{v_i\}$, and
that $W^s$ is a vector space with basis $\{w_\alpha\}$ and dual basis
$\{w^\alpha\}$.  Then we have for $k\geq 2$ the
$GL(W)\times GL(V)$-equivariant exact sequence
$$
  0 \to U_k\to W^*\otimes V\otimes(\textstyle\bigwedge^kW^*\otimes\bigwedge^kV)
     \stackrel{\tau_k}{\to} \bigwedge^{k+1}W^*\otimes\bigwedge^{k+1}V\to 0.
$$
Here the surjection is the obvious skew-symmetrization map, and $U_k$
is by definition its kernel.  The term $k=1$ will be exceptional, and
we instead define
$$
  0 \to U_1\to \mbox{Sym}^2(W^*\otimes V)\stackrel{\tau}{\to}
    \textstyle\bigwedge^2W^*\otimes\bigwedge^2V\to 0,
$$
so that $U_1 = \mbox{Sym}^2W^*\otimes\mbox{Sym}^2V$.

Now we can regard our coefficients $H^{iI}_{\alpha A}$, with
$|I|=|A|=k$, at each point 
of $M$ as coefficients of an element 
$$ 
  H_k = H^{iI}_{\alpha A}w^\alpha\otimes v_i\otimes w^A\otimes v_I
    \in W^*\otimes V\otimes(\textstyle\bigwedge^kW^*\otimes\bigwedge^kV).
$$
\begin{Definition}
The form $\Pi\in\Omega^{n+1}(M)\cap\{\bigwedge^nJ\}$ is {\em
  symmetric}\index{symmetric form} 
  if its expansion (\ref{Expansion16}) has $H_k\in U_k$ for
  all $k\geq 1$. 
\end{Definition}
For $k=1$ the condition is
$$
  H^{ij}_{\alpha\beta}=H^{ji}_{\alpha\beta}=H^{ij}_{\beta\alpha}.
$$

We first need to show that the condition that a given $\Pi$ be
symmetric is independent of the choice of admissible coframe.
Equivalently, we can show that the symmetry condition is preserved
under the group of coframe changes of the form (\ref{DefMultiG}), and
we will show this under three subgroups generating the group.
First, it is obvious that a change
$$
    \bar\theta^\alpha = a^\alpha_\beta\theta^\beta, \quad
    \bar\omega^i = b^i_j\omega^j, \quad
    \bar\pi^\alpha_i = a^\alpha_\beta\pi^\beta_j(b^{-1})^j_i,
$$
preserves the symmetry condition, because of the equivariance of
the preceding exact sequences under $(a^\alpha_\beta)\times(b^i_j)\in
GL(W)\times GL(V)$.  Second, symmetry is preserved under
$$
    \bar\theta^\alpha = \theta^\alpha, \quad
    \bar\omega^i = \omega^i, \quad
    \bar\pi^\alpha_i = d^\alpha_{i\beta}\theta^\beta +
      e^\alpha_{ij}\omega^j +\pi^\alpha_i, 
$$
with $e^\alpha_{ij}=e^\alpha_{ji}$,
because such a change has no effect on the expression for $\Pi$ modulo
$\{\bigwedge^{n+1}J\}$.  Finally, consider a change of the form
$$
  \bar\theta^\alpha = \theta^\alpha,\quad
  \bar\omega^i = c^i_\beta\theta^\beta+\omega^i, \quad
  \bar\pi^\alpha_i = \pi^\alpha_i.
$$
We will prove the invariance of the symmetry condition
infinitesimally, writing instead of $\bar\omega^i$ the family
\begin{equation}
  \omega^i(\varepsilon) = \varepsilon c^i_\beta\theta^\beta + \omega^i.
\label{TransFamily16}
\end{equation}
This associates to each $H_k$ a tensor $H_{k+1}(\varepsilon)$ for each $k\geq 1$, and we
will show that
$\sf{d}{d\varepsilon}|_{\varepsilon=0}H_{k+1}(\varepsilon)\in U_{k+1}$.
This just amounts to looking at the terms linear in $\varepsilon$ when
(\ref{TransFamily16}) is substituted into (\ref{Expansion16}).
Noting that $C=(c^i_\beta)\in W^*\otimes V$, we consider the commutative
diagram
$$
 \begin{array}{ccc}
   W^*\otimes V\otimes\textstyle\bigwedge^kW^*\otimes\bigwedge^kV
     \otimes W^*\otimes V & \stackrel{\sigma}{\to} &
     W^*\otimes V\otimes\bigwedge^{k+1}W^*\otimes\bigwedge^{k+1}V \\
   \downarrow\tau_k\otimes 1 & & \downarrow\tau_{k+1} \\
  \bigwedge^{k+1}W^*\otimes\bigwedge^{k+1}V\otimes W^*\otimes V
    & \to & \bigwedge^{k+2}W^*\otimes\bigwedge^{k+2}V,
 \end{array}
$$
where $\sigma$ is skew-symmetrization with the latter $W^*\otimes V$,
and $\tau_k\otimes 1$ is an extension of the earlier
skew-symmetrization.  The point is that given $H_k\otimes C$ in the
upper-left space of this diagram, 
$$
  \sf{d}{d\varepsilon}|_{\varepsilon=0}H_{k+1}(\varepsilon) =
  \sigma(H_k\otimes C).
$$
So if we assume that $(H_k)\in U_k$, then $\tau_k(H_k)=0$, so
$\sigma(H_k\times C)\in U_{k+1}$, which is what we wanted to show.

This proves that the condition that the symmetry condition on $\Pi =
d\Lambda$ is independent of the choice of adapted coframe.  We can now
state the following. 
\begin{Theorem}
Given a functional $[\Lambda]\in\Gamma(M,\bigwedge^n(J/I))$, there is a
unique admissible lifting $\Lambda\in\Gamma(M,\bigwedge^n J)$ such that
$\Pi=d\Lambda\in\Omega^{n+1}(M)$ is symmetric.
\label{PCTheorem16}
\end{Theorem}
\begin{Proof}
We inductively construct $\Lambda =
\Lambda_0+\Lambda_1+\cdots+\Lambda_{\mathit{min}(n,s)}$, with each
$\Lambda_i\in\{\bigwedge^iI\}$ chosen to eliminate the fully
skew-symmetric part of 
$$
  \Pi_{i-1}\stackrel{\mathit{def}}{=}d(\Lambda_0+\cdots+\Lambda_{i-1}).
$$  
Initially,
$\Lambda_0=F\omega$ is the prescribed $X$-semibasic $n$-form modulo
$\{I\}$.  
We know from the existence of admissible liftings that there is some
$\Lambda_1\in\{I\}$
such that
$\Pi_1\stackrel{\mathit{def}}{=}d(\Lambda_0+\Lambda_1)\in\{I\}$; and
we know from
Proposition~\ref{LiftProp16} that $\Lambda_1$ is 
 uniquely determined modulo $\{\bigwedge^2I\}$.
Now let
$$
  \Pi_1 \equiv
  H^{ij}_{\alpha\beta}\pi^\alpha_i\wedge\theta^\beta\wedge\omega_{(j)}
  \pmod{\{\textstyle\bigwedge^2I\}+\{\bigwedge^{n+1}J\}}.
$$
If we add to $\Lambda_0+\Lambda_1$ the $I$-quadratic term
$$
  \Lambda_2 = \sf1{2!^2}F^{ij}_{\alpha\beta}
    \theta^\alpha\wedge\theta^\beta\wedge\omega_{(ij)},
$$
then the structure equation (\ref{GenMultiStructure}) shows that
this alters the $I$-linear term $\Pi_1$ only by
$$
  H^{ij}_{\alpha\beta}\leadsto
    H^{ij}_{\alpha\beta} + F^{ij}_{\alpha\beta}.
$$
Because
$F^{ij}_{\alpha\beta}=-F^{ji}_{\alpha\beta}=-F^{ij}_{\beta\alpha}$, we
see that $F^{ij}_{\alpha\beta}$ may be uniquely chosen so that the new
$H^{ij}_{\alpha\beta}$ lies in $U_1$.
\index{admissible lifting|)}

The inductive step is similar.  Suppose we have
$\Lambda_0+\cdots+\Lambda_l\in\Gamma(\bigwedge^nJ)$ such that
$\Pi_l =  d(\Lambda_0+\cdots+\Lambda_l)$
is symmetric modulo $\{\bigwedge^{l-1}I\}$.  Then the term of
$I$-degree $l$ is of the form
$$
    \Pi_l\equiv \sf{1}{l!^2}
     \sum H^{iI}_{\alpha A}\pi^\alpha_i
       \wedge\theta^A\wedge\omega_{(I)} \pmod{
       \{U_1\}+\cdots+\{U_{l-1}\} +\{\textstyle\bigwedge^{l+1}I\} }.
$$
for some $H^{iI}_{\alpha A}$.  There is a unique skew-symmetric term
$$
  \Lambda_{l+1}=\sf{1}{(l+1)!^2}\sum_{|I|=|A|=l+1}
    F^I_A\theta^A\wedge\omega_{(I)}
$$
which may be added so that
$$
    \Pi_{l+1}\in
    \{U_1\}+\cdots+\{U_l\}+\{\textstyle\bigwedge^{l+1}I\}.
$$
We can continue in this manner, up to $l=\mbox{min}(n,s)$.
\end{Proof}
\begin{Definition}
The unique $\Lambda$ in the preceding theorem is called the {\em
  Betounes} form for the functional $[\Lambda]$.\footnote{It was
  introduced in coordinates in \cite{Betounes:Extension},
and further discussed in \cite{Betounes:Differential}.}
Its derivative $\Pi 
= d\Lambda$ is the {\em Poincar\'e-Cartan} form for $[\Lambda]$.
\end{Definition}

The unique determination of $\Pi$, along with the invariance of the
symmetry condition under admissible coframe changes of $M$, implies
that $\Pi$ is globally defined and invariant under symmetries of
the functional $[\Lambda]$ and the multi-contact structure $(M,I,J)$.

It is instructive to see the first step of the preceding construction
in coordinates.  If our initial Lagrangian is
$$
  \Lambda_0 = L(x,z,p)dx,
$$
then we have already seen in (\ref{FirstLift16}) that 
$$
  \Lambda_0+\Lambda_1 = L\,dx+\theta^\alpha\wedge
    \sf{\partial L}{\partial p^\alpha_i}dx_{(i)}.
$$
The $H_1$-term of $d(\Lambda_0+\Lambda_1)$ (see (\ref{FirstPC16})) is
\begin{equation}
  \frac{\partial^2L}{\partial p^\alpha_i\partial p^\beta_j}dp^\beta_j
    \wedge\theta^\alpha\wedge dx_{(i)}.
\label{NonSymmetrized16}
\end{equation}
Of course $L_{p^\alpha_ip^\beta_j}=L_{p^\beta_jp^\alpha_i}$,
corresponding to the fact that $H_1\in\mbox{Sym}^2(W^*\otimes V)$
automatically.  The proof shows that we can add
$\Lambda_2\in\{\bigwedge^2I\}$ so that $\Pi_2$ instead includes
$$
   \sf12(L_{p^\alpha_ip^\beta_j}+L_{p^\alpha_jp^\beta_i})
     \pi^\alpha_i\wedge\theta^\beta\wedge\omega_{(j)} \in
     U_1=\mbox{Sym}^2(W^*)\otimes\mbox{Sym}^2V.
$$
In fact, this corresponds to 
the {\em principal symbol}\index{symbol} 
of the Euler-Lagrange PDE system
(\ref{ELEqn16}), given by the symmetric $s\times s$ matrix
\begin{eqnarray*}
  H_{\alpha\beta}(\xi) & = &
  \frac{\partial^2L}{\partial p^\alpha_i\partial p^\beta_j}\xi_i\xi_j
\\ & = & \sf12(L_{p^\alpha_ip^\beta_j}+L_{p^\alpha_jp^\beta_i})
      \xi_i\xi_j, \qquad \xi\in V^*.
\end{eqnarray*}
In light of this, it is not surprising to find that only the symmetric
part of (\ref{NonSymmetrized16}) has invariant meaning.

Note also that if Lagrangians $\Lambda$, $\Lambda^\prime$ differ by a
divergence,
$$
  \Lambda-\Lambda^\prime = d\lambda,\qquad
  \lambda\in\Gamma(M,\textstyle\bigwedge^{n-1}J),
$$
then the construction in the proof of Theorem~\ref{PCTheorem16} shows that the
Poincar\'e-Cartan forms are equal, though the Betounes forms may not be.
A related but more subtle property is the following.
\begin{Theorem}
For a functional $[\Lambda]\in\Gamma(M,\bigwedge^n(J/I))$, the
Poincar\'e-Cartan form $\Pi=0$ if and only if the Euler-Lagrange
system is trivial, $\mathcal E_\Lambda=\mathcal I$.
\end{Theorem}
\begin{Proof}
One direction is clear:  if $\Pi=0$, then the $n$-form generators
$\Psi_\alpha$ for $\mathcal E_\Lambda$ can be taken to be $0$, so that
$\mathcal E_\Lambda=\mathcal I$.  For the converse, we first consider
the $I$-linear term
$$
  \Pi_1 \equiv H^{ij}_{\alpha\beta}\pi^\alpha_i\wedge\theta^\beta\wedge
    \omega_{(j)}\pmod{\{\textstyle\bigwedge^2I\} +
       \{\textstyle\bigwedge^{n+1}J\}}.
$$
$\mathcal E_\Lambda$ is generated by $\mathcal I$
and $\Psi_\beta=H^{ij}_{\alpha\beta}\pi^\alpha_i\wedge\omega_{(j)}$,
$1\leq\beta\leq s$, and
our assumption $\mathcal E_\Lambda=\mathcal I$ then implies that
these $\Psi_\beta=0$; that is,
$$
  H_1 = H^{ij}_{\alpha\beta}\pi^\alpha_i\wedge\theta^\beta\wedge
     \omega_{(j)}=0.
$$

We will first show that this implies
$$
  H_2=H_3=\cdots=0
$$
as well,
which will imply $\Pi\in\{\bigwedge^{n+1}J\}$.  To see
this, suppose $H_l$ is the first non-zero term, having $I$-degree $l$.  Then we can consider
$$
  0 \equiv d\Pi
      \pmod{\{\textstyle\bigwedge^lI\}+\{\bigwedge^{n+1}J\}},
$$
and using the structure equations (\ref{GenMultiStructure}),
$$
  0 = \sum_{|I|=|A|=l-1}
    H^{ijI}_{\alpha\beta A}\pi^\alpha_i\wedge\pi^\beta_j.
$$
Written out fully, this says that
$$
  H^{i_1i_2\cdots i_{l+1}}_{\alpha_1\alpha_2\cdots\alpha_{l+1}}
    = H^{i_2i_1\cdots i_{l+1}}_{\alpha_2\alpha_1\cdots\alpha_{l+1}}.
$$
Also, $H^{iI}_{\alpha A}$ is fully skew-symmetric in $I$ and $A$.  But
together, these imply in that $H_l$ is fully skew-symmetric in all
upper and all lower indices, for
\begin{eqnarray*}
  H^{i_1i_2i_3\cdots}_{\alpha_1\alpha_2\alpha_3\cdots}
   & = & H^{i_2i_1i_3\cdots}_{\alpha_2\alpha_1\alpha_3\cdots} \\
   & = & -H^{i_2i_1i_3\cdots}_{\alpha_2\alpha_3\alpha_1\cdots} \\
   & = & -H^{i_1i_2i_3\cdots}_{\alpha_3\alpha_2\alpha_1\cdots} \\
   & = & H^{i_1i_2i_3\cdots}_{\alpha_3\alpha_1\alpha_2\cdots} \\ 
   & = & H^{i_2i_1i_3\cdots}_{\alpha_1\alpha_3\alpha_2\cdots} \\
   & = & -H^{i_2i_1i_3\cdots}_{\alpha_1\alpha_2\alpha_3\cdots}.
\end{eqnarray*}
This proves full skew-symmetry in the upper indices, and the proof for
lower indices is similar.  However, we constructed $\Pi$ so that each $H_k$
lies in the invariant complement of the fully skew-symmetric tensors,
so we must have $H_k=0$.

Now we have shown that if the Euler-Lagrange equations of $[\Lambda]$
are trivial, then $\Pi\in\{\bigwedge^{n+1}J\}$.  But that means that
the Betounes form $\Lambda$ is not merely semibasic over the quotient
space $X$, but actually basic.  We can then compute the (assumed
trivial) first variation down in $X$ instead of $M$, and find that for
any submanifold $N\hookrightarrow X$, and any vector field $v$ along
$N$ vanishing at $\partial N$, 
$$
  0 = \int_N v\innerprod d\Lambda.
$$
But this implies that $d\Lambda = 0$, which is what we wanted to
prove.
\end{Proof}

\

The preceding results indicate that $\Pi$ is a good generalization of
the classical Poincar\'e-Cartan form for second-order, scalar
Euler-Lagrange equations.
We note that for higher-order Lagrangian functionals on vector-valued
functions of one variable (i.e., functionals on curves), such a
generalization is known, and not difficult; but for functionals of
order $k\geq 2$ on vector-valued functions of several variables,
little is known.\footnote{But see Grassi\index{Grassi, M.|nn}, 
cit.~p.~\pageref{Grassi}n.}

We want to briefly mention a possible generalization to the multi-contact
case of Noether's theorem\index{Noether's theorem}, 
which gives an isomorphism from a Lie
algebra of symmetries to a space of conservation laws.
To avoid distracting global considerations, we will assume that
$H^q_{dR}(M)=0$ in all degrees $q>0$.
First, we have the space $\lie{g}_\Pi$, consisting of vector fields on
$M$ which preserve $I$ and $\Pi$,
$$
  \lie{g}_\Pi=\{v\in\mathcal V(M):\mathcal L_vI\subseteq I,\
    \mathcal L_v\Pi=0\}.
$$
Second, we have the space of conservation laws
$$
  \mathcal C=H^{n-1}(\Omega^*(M)/\mathcal E_\Lambda);
$$
under our topological assumption, this is identified with $H^n(\mathcal
E_\Lambda)$, and we need not introduce a notion of ``proper''
conservation law as in \S\ref{Section:Noether}.  In this situation,
Noether's theorem says the following.
\begin{quote}
{\em There is a map $\eta:\lie{g}_\Pi\to H^n(\mathcal E_\Lambda)$, defined by
  $v\mapsto [v\innerprod\Pi]$, which is an isomorphism if $\Pi$ is
  non-degenerate in a suitable sense.}
\end{quote}
The map is certainly well-defined; that is,
for any
$v\in\lie{g}_\Pi$, the form $v\innerprod\Pi$ is a closed section of
$\mathcal E_\Lambda$.  First,
$$
  v\innerprod\Pi = (v\innerprod\theta^\alpha)\Psi_\alpha-\theta^\alpha
    \wedge(v\innerprod\Psi_\alpha),
$$
so that $v\innerprod\Pi$ is a section of $\mathcal E_\Lambda$; and
second,
$$
  d(v\innerprod\Pi)=\mathcal L_v\Pi-v\innerprod d\Pi = 0,
$$
so that $v\innerprod\Pi$ is closed.
However, the proof that under the right conditions this map is an
isomorphism involves some rather sophisticated commutative algebra,
generalizing the symplectic linear 
algebra\index{symplectic!linear algebra} used in
Chapter~\ref{Chapter:Lagrangians}.
This will not be presented here.  

As in the scalar case, a simple prescription for the conserved density in
$H^{n-1}(\Omega^*(M)/\mathcal E_\Lambda)$ corresponding to
$v\in\lie{g}_\Pi$ is available when also
$$
  \mathcal L_v\Lambda = 0.
$$
One virtue of the Betounes form is that this holds for
infinitesimal multi-contact symmetries of $[\Lambda]$.
Assuming only that $d\Lambda = \Pi$ and $\mathcal L_v\Lambda = 0$, we
can calculate that
\begin{equation}
  d(-v\innerprod\Lambda) = -\mathcal L_v\Lambda + v\innerprod d\Lambda
    = v\innerprod \Pi.
\label{SimpleNoether17}
\end{equation}
Therefore, $-v\innerprod \Lambda\in\Omega^{n-1}(M)$ represents a class
in $\mathcal C=H^{n-1}(\Omega^*(M)/\mathcal E_\Lambda)$ corresponding
to $\eta(v) \in H^n(\mathcal E_\Lambda)$.  We will use this
prescription in the following.
\index{Betounes form|)}
\index{Poincar\'e-Cartan form|)}
\index{Lagrangian|)}

\subsection{Harmonic Maps of Riemannian Manifolds}
\index{harmonic!map|(}

The most familiar variational PDE systems in differential geometry are
those describing harmonic maps between Riemannian manifolds.

Let $P,Q$ be Riemannian manifolds of dimensions $n,s$.  We will
define a Lagrangian density on $P$, depending on a map $P\to Q$ and
its first derivatives, whose integral over $P$ may be thought of as
the {\em energy}\index{energy!of a map} 
of the map.  The appropriate multi-contact manifold
for this is the space of $1$-jets of maps $P\to Q$,
$$
  M = J^1(P,Q),
$$
whose multi-contact system will be described shortly.  We may also
think of $M$ as $\mbox{Hom}(TP,TQ)$, the total space of a rank-$ns$
vector bundle over $P\times Q$.  To carry out computations, it will be
most convenient to work on
$$
  \mathcal F \stackrel{\mathit{def}}{=} \mathcal F(P)\times\mathcal F(Q)
     \times\R^{ns},
$$
where $\mathcal F(P)$, $\mathcal F(Q)$ are the orthonormal frame
bundles\index{Riemannian!frame bundle}.  
These are parallelized in the usual manner by
$(\omega^i,\omega^i_j)$, $(\varphi^\alpha,\varphi^\alpha_\beta)$,
respectively, with structure equations
$$
  \left\{\begin{array}{ll}
    d\omega^i = -\omega^i_j\wedge\omega^j,&
      d\omega^i_j=-\omega^i_k\wedge\omega^k_j + \Omega^i_j, \\
    d\varphi^\alpha = -\varphi^\alpha_\beta\wedge\varphi^\beta,&
      d\varphi^\alpha_\beta = -\varphi^\alpha_\gamma\wedge
      \varphi^\gamma_\beta + \Phi^\alpha_\beta.
   \end{array}\right.
$$
These forms and structure equations will be considered pulled back
to $\mathcal F$.  To complete a coframing of $\mathcal F$, we take
linear fiber coordinates $p^\alpha_i$ on $\R^{ns}$, and define
$$
  \pi^\alpha_i = dp^\alpha_i+\varphi^\alpha_\beta p^\beta_i -
     p^\alpha_j\omega^j_i.
$$
The motivation here is that $\mbox{Hom}(TP,TQ)\to P\times Q$ is a vector
bundle associated to the principal $(O(n)\times O(s))$-bundle
$\mathcal F(P)\times\mathcal F(Q)\to P\times Q$, with the data
$((e^P_i),(e^Q_\alpha),(p^\alpha_i))\in\mathcal F$ defining the
homomorphism $e^P_i\mapsto e^Q_\alpha p^\alpha_i$.  Furthermore, if a
section $\sigma\in\Gamma(\mbox{Hom}(TP,TQ))$ is represented by an equivariant map
$(p^\alpha_i):\mathcal F(P)\times\mathcal F(Q)\to\R^{ns}$, then the
$\R^{ns}$-valued $1$-form $(\pi^\alpha_i)$
represents the covariant derivative\index{covariant derivative} 
of $\sigma$.

For our purposes, note that $M=\mbox{Hom}(TP,TQ)$ is the quotient of $\mathcal
F$ under a certain action of $O(n)\times O(s)$, and that the forms
semibasic for the projection
$\mathcal F\to M$ are generated by
$\omega^i,\varphi^\alpha,\pi^\alpha_i$.  A natural multi-contact
system on $M$ pulls back to $\mathcal F$ as the Pfaffian system $I$
generated by
$$
  \theta^\alpha \stackrel{\mathit{def}}{=} \varphi^\alpha-p^\alpha_i\omega^i,
$$
and the associated integrable Pfaffian system on $M$ pulls back to
$J=\{\varphi^\alpha,\omega^i\} = \{\theta^\alpha,\omega^i\}$.  The
structure equations on $\mathcal F$ adapted to these Pfaffian systems
are
\begin{equation}
  \left\{\begin{array}{l}
    d\theta^\alpha = -\pi^\alpha_i\wedge\omega^i
    -\varphi^\alpha_\beta\wedge\theta^\beta, \\
    d\omega^i = -\omega^i_j\wedge\omega^j, \\
    d\pi^\alpha_i = \Phi^\alpha_\beta p^\beta_i-p^\alpha_j\Omega^j_i
    -\varphi^\alpha_\beta\wedge\pi^\beta_i -
    \pi^\alpha_j\wedge\omega^j_i.
  \end{array}\right.
\label{HarmonicStreqns17}
\end{equation}

We now define the {\em energy Lagrangian}\index{energy!of a map}
$$
  \tilde\Lambda = \sf12||p||^2\omega\in\Gamma(\textstyle\bigwedge^nJ)\subset
     \Omega^n(\mathcal F),
$$
where the norm is
$$
  ||p||^2 = \mbox{Tr}(p^*p) = \textstyle\sum(p^\alpha_i)^2.
$$
Although this $\tilde\Lambda$ is not an admissible
lifting\index{admissible lifting} of its
induced functional $[\Lambda]$, a computation using the structure equations
(\ref{HarmonicStreqns17}) shows that
$$
  \Lambda\stackrel{\mathit{def}}{=} \sf12||p||^2\omega+p^\alpha_i\theta^\alpha
    \wedge\omega_{(i)}
$$
is admissible:
\begin{eqnarray*}
  d\Lambda & = & -\theta^\alpha\wedge\pi^\alpha_i\wedge\omega_{(i)}
   -p^\alpha_ip^\beta_i\varphi^\alpha_\beta\wedge\omega
   +p^\alpha_ip^\alpha_j\omega^j_i\wedge\omega \\
  & = & -\theta^\alpha\wedge\pi^\alpha_i\wedge\omega_{(i)},
\end{eqnarray*}
where the last step uses $\varphi^\alpha_\beta+\varphi^\beta_\alpha =
\omega^i_j+\omega^j_i = 0$.  Now we define
$$
  \Pi = -\theta^\alpha\wedge\pi^\alpha_i\wedge\omega_{(i)},
$$
and note that $\Pi$ is in fact the lift to ${\mathcal F}$ of a {\em
  symmetric} \index{symmetric form}form on
$M$, as defined earlier.  Therefore, we have found the Betounes form
and the Poincar\'e-Cartan form for the energy functional.  

The Euler-Lagrange system for $[\Lambda]$, pulled back to $\mathcal
F$, is
$$
  \mathcal E_\Lambda = \{\theta^\alpha,\ \pi^\alpha_i\wedge\omega^i,\
    \pi^\alpha_i\wedge\omega_{(i)}\}.
$$
A Legendre submanifold $N\hookrightarrow M=J^1(P,Q)$ on
which $\bigwedge\omega^i\neq 0$ is the $1$-jet
graph of a map $f:P\to Q$.  On the inverse image $\pi^{-1}(N)\subset\mathcal
F$, in addition to $\theta^\alpha = 0$, there are relations
$$
  \pi^\alpha_i= h^\alpha_{ij}\omega^j,\quad h^\alpha_{ij}=h^\alpha_{ji}.
$$
Differentiating this equation shows that the expression
$$
  h = h^\alpha_{ij}\omega^i\omega^j\otimes e^Q_\alpha 
$$
is invariant along fibers of $\pi^{-1}(N)\to N$, so it gives a
well-defined section of $\mbox{Sym}^2(T^*P)\otimes TQ$; this is called
the {\em second fundamental 
form}\index{second fundamental form!of a map} 
of the map $f:P\to Q$.  The
condition for $N$ to be an integral manifold of the Euler-Lagrange
system is then
$$
  \mbox{Tr}(h) = h^\alpha_{ii}=0\in\Gamma(M,f^*TQ),
$$
\begin{Definition}  A map $f:P\to Q$ between Riemannian manifolds is
  {\em harmonic} if the trace of its second fundemental form vanishes.
\end{Definition}
Expressed in coordinates on $P$ and $Q$, this
is a second-order PDE system for $f:P\to Q$.

We now consider conservation 
laws\index{conservation law!for harmonic maps|(} 
for the harmonic map system
$\mathcal E_\Lambda\subset\Omega^*(\mathcal F)$ corresponding to
infinitesimal isometries (Killing vector 
fields\index{Killing vector fields} of either $P$ or $Q$.  These are
symmetries
not only of $\Pi$ but of the Lagrangian $\Lambda$, so we can use the
simplified prescription (\ref{SimpleNoether17}) for a conserved
$(n-1)$-form.

First, an infinitesimal isometry of $P$ induces a unique vector field on
$M=J^1(P,Q)$ preserving $I$ and fixing $Q$.  This vector field
preserves $\Lambda$, $\Pi$, and $\mathcal E_\Lambda$, and has a
natural lift to $\mathcal F$ which does the same.  This vector field
$v$ on $\mathcal F$ satisfies
$$
  v\innerprod\omega^i = v^i,\quad
  v\innerprod\varphi^\alpha = 0\quad\Rightarrow\quad
  v\innerprod\theta^\alpha = -p^\alpha_iv^i,
$$
for some functions $v^i$.  We can then calculate
$$
  \varphi_v \stackrel{\mathit{def}}{=}
  v\innerprod\Lambda \equiv (\sf12||p||^2v^i-p^\alpha_ip^\alpha_jv^j)
    \omega_{(i)}\pmod{\{I\}}.
$$
As in Chapter~\ref{Chapter:Conformal}, it is useful to write this expression restricted to
the $1$-jet graph of a map $f:P\to Q$, which is
\begin{eqnarray*}
  v\innerprod\Lambda|_N & = &
    *_P(\sf12||p||^2v^i-p^\alpha_ip^\alpha_jv^j)\omega^i \\
  & = & *_P\sf12\left(v\innerprod(\sf12||df||^2\textstyle\sum(\omega^i)^2-
     f^*\sum(\varphi^\alpha)^2)\right),
\end{eqnarray*}
where we use $f^*\varphi^\alpha =p^\alpha_i\omega^i$,
and $*_P\omega_{(i)}=\omega^i$.
One might recognize the {\em stress-energy 
tensor}\index{stress-energy tensor}
$$
  S=\sf12||df||^2ds^2_P - f^*ds^2_Q,
$$
and write our conserved density as
\begin{equation}
  2(v\innerprod\Lambda) \equiv *_P(v\innerprod S)\pmod{\{I\}}.
\label{DefStress16}
\end{equation}
In fact, $S$ is traditionally defined as the unique symmetric
$2$-tensor on $P$ for which the preceding equation holds for arbitrary
$v\in\mathcal V(P)$ and $f:P\to Q$; 
then (\ref{DefStress16}) gives a conserved
density when $v$ is an infinitesimal isometry and $f$ is a harmonic map.
In this case
In fact, for any infinitesimal isometry $v$, a calculation gives
\begin{equation}
  d(*_P(v\innerprod S)) = (v\innerprod\mbox{div }S)\omega
\label{RiemID16}
\end{equation}
on the $1$-jet graph of any map.\footnote{The {\em divergence} of a
  symmetric $2$-form $S$ is the $1$-form $\mbox{div }S =
  \nabla_{e_i}S(e_i,\cdot)$, where $\nabla$ is the Levi-Civita
  covariant derivative and $(e_i)$ is any orthonormal frame.
  Equation (\ref{RiemID16}) is true of any symmetric $2$-form $S$ and
  infinitesimal isometry $v$.}

Now consider an infinitesimal isometry of $Q$, whose lift
$w\in\mathcal V(\mathcal F)$ satisfies
$$
  w\innerprod\omega^i = 0,\ w\innerprod\varphi^\alpha =
  w^i\quad\Rightarrow\quad w\innerprod\theta^\alpha = w^\alpha.
$$
Then
$$
  \varphi_w \stackrel{\mathit{def}}{=} w\innerprod\Lambda = w^\alpha
  p^\alpha_i\omega_{(i)}. 
$$
Given a map $f:P\to Q$, we can use $df\in\mbox{Hom}(TP,TQ)$ and
$ds^2_Q\in\mbox{Sym}^2(T^*Q)$ to regard $ds^2_Q(df(\cdot),w)$ as a
$1$-form on $P$, and then
$$
  \varphi_w = *_P(ds^2_Q(df(\cdot),w)).
$$
Because this expression depends linearly on $w$,
we can simplify further by letting $\lie{a}$ denote the Lie
algebra of infinitesimal symmetries of $Q$, and then the map $w\mapsto
w\innerprod \Lambda$
is an element of $\lie{a}^*\otimes\Omega^{n-1}(P)$.  If we define an
$\lie{a}^*$-valued $1$-form on $P$ by
$$
  \alpha(v) = ds^2_Q(df(v),\cdot),\quad v\in T_pP,
$$
then our conservation laws read
\begin{equation}
  d(*_P\alpha) = 0 \in \lie{a}^*\otimes\Omega^n(P).
\label{BadCLs17}
\end{equation}
The $\lie{a}^*$-valued $(n-1)$-form $*_P\alpha$ may be formed for any
map $f:P\to Q$, and it is closed if $f$ is harmonic.
In fact, if $Q$ is locally homogeneous, meaning that infinitesimal
isometries span each tangent space $T_qQ$, then (\ref{BadCLs17}) is
{\em equivalent} to the harmonicity of $f$.

An important special case of this last phenomenon is when $Q$ itself
is a Lie group $G$ with bi-invariant metric $ds^2_G$.  Examples are
compact semisimple Lie groups, such as $O(N)$ or $SU(N)$, with metric
induced by the Killing form on the Lie algebra $\lie{g}$.
Now a map $f:P\to G$ is uniquely determined up to
left-translation by the pullback $f^*\varphi$ of the left-invariant
$\lie{g}$-valued Maurer-Cartan\index{Maurer-Cartan!form} 
$1$-form $\varphi$.  Using the metric
to identify $\lie{g}\cong\lie{g}^*$, the conservation laws state that
if $f$ is harmonic, then $d(*_P(f^*\varphi)) = 0$.  Conversely, if $P$
is simply connected, then given a $\lie{g}$-valued $1$-form $\alpha$
on $P$ satisfying
$$
  \left\{\begin{array}{l}
    d\alpha + \sf12[\alpha,\alpha]=0, \\
    d(*_P\alpha) = 0,
  \end{array}\right.
$$
there is a harmonic map $f:P\to G$ with $f^*\varphi=\alpha$,
uniquely determined up to left-translation.  
This is the idea behind
the gauge-theoretic reformulation of certain harmonic map systems, for
which remarkable results have been obtained in the past
decade.\footnote{The literature on this subject is vast, but a good
  starting point is \cite{Wood:Harmonic}.}
Quite
generally, PDE systems that can be written as systems of conservation
laws have special properties; one typically exploits such expressions to
define weak solutions, derive integral identities, and
prove regularity theorems.
\index{conservation law!for harmonic maps|)}
\index{harmonic!map|)}

\section{Higher-Order Conservation Laws}

One sometimes encounters a conservation law for a PDE that involves
higher-order derivatives of the unknown function, but that cannot be
expressed in terms of derivatives of first-order conservation laws
considered up to this point.  An example is the $(1+1)$-dimensional
wave equation $-z_{tt}+z_{xx}=0$, for which $(z_{tt}^2+z_{tx}^2)dt +
2z_{tt}z_{tx}dx$ is closed on solutions, but cannot be obtained by
differentiating any conservation law on $J^1(\R^2,\R)$.  
In this section, we introduce the geometric framework in which such
conservation laws may be found, and we propose a version of Noether's
theorem appropriate to this setting.  While other general forms of
Noether's theorem have been stated and proved (e.g., see
\cite{Vinogradov:CSpectral} or \cite{Olver:Applications}), it is not
clear how they relate to that conjectured here.

We also discuss (independently from the preceding)
the higher-order relationship
between surfaces in Euclidean space with Gauss curvature $K=-1$ and
the sine-Gordon equation $z_{tx}=\sf12\sin(2z)$, in terms of exterior
differential systems.

\subsection{The Infinite Prolongation}

We begin by defining the {\em prolongation} of an exterior
differential system (EDS).  When this is applied to the EDS associated to a
PDE system, it gives the EDS associated to the PDE system
augmented by the first derivatives of the original equations.  This
construction then extends to that of the {\em infinite prolongation},
an EDS on an infinite-dimensional manifold which includes
information about derivatives of all orders.

The general definition of prolongation uses a construction introduced
in \S\ref{Section:MultiPC}, in the discussion of multi-contact manifolds.
Let $X^{n+s}$ be a manifold, and
$G_n(TX)\stackrel{\pi}{\to} X$ the bundle of tangent $n$-planes of
$X$; points of $G_n(TX)$
are of the form $(p,E)$, where $p\in X$ and $E\subset T_pX$ is a
vector subspace of dimension $n$.  
As discussed previously, there is a canonical Pfaffian
system $I\subset T^*G_n(TX)$ of 
rank $s$, defined at $(p,E)$ by
$$
  I_{(p,E)} \stackrel{\mathit{def}}{=} \pi^*(E^\perp).
$$
Given local coordinates $(x^i,z^\alpha)$ on $X$, there are induced
coordinates $(x^i,z^\alpha,p^\alpha_i)$ on $G_n(TX)$, in terms
of which $I$ is generated by the $1$-forms
\begin{equation}
  \theta^\alpha = dz^\alpha - p^\alpha_i dx^i.
\label{CoordMC17}
\end{equation}
We let $\mathcal I\subset\Omega^*(G_n(TX))$ be the differential
ideal generated by $I$.

Now let $(M,\mathcal E)$ be an exterior differential system; that is, $M$
is a manifold of dimension $m+s$ and $\mathcal E\subset\Omega^*(M)$
is a differential ideal for which we are interested in
$m$-dimensional integral manifolds.
We then define the locus $M^{(1)}\subset G_m(TM)$ to consist of the
{\em integral elements} of $\mathcal E\subset\Omega^*(M)$; that is,
$(p,E)\in M^{(1)}$ if and only if
$$
  \varphi_E = 0 \in \textstyle\bigwedge^*(E^*) \quad
    \mbox{for all }\varphi\in\mathcal E.
$$
We will assume from now on that
$M^{(1)}\stackrel{\iota}{\hookrightarrow} G_m(TM)$ is a smooth
submanifold.  Then we define
$$
  \mathcal E^{(1)} \stackrel{\mathit{def}}{=}
     \iota^*\mathcal I\subset\Omega^*(M^{(1)})
$$
as the
restriction to $M^{(1)}$ of the multi-contact differential ideal.
This is the same as the differential ideal generated by the Pfaffian system
$\iota^*I\subset T^*M^{(1)}$, and the {\em first prolongation} of
$(M,\mathcal E)$ is defined to be the exterior differential system
$(M^{(1)},\mathcal E^{(1)})$.  
Note that the first prolongation is always a Pfaffian system.
Furthermore, if $\pi:M^{(1)}\to M$ is the obvious
projection map, and assuming that $\mathcal E$ is a Pfaffian system,
then one can show that $\pi^*\mathcal
E\subseteq\mathcal E^{(1)}$.  However, the projection $\pi$ could be
quite complicated, and need not even be surjective.  Finally, note
that any integral manifold $f:N\hookrightarrow M$ of $\mathcal E$ lifts
to an integral manifold $f^{(1)}:N\hookrightarrow M^{(1)}$ of $\mathcal
E^{(1)}$, and that the transverse integral manifold of
$\mathcal E^{(1)}$ is locally of this form.

Inductively, the {\em $k$th prolongation}
$(M^{(k)},\mathcal E^{(k)})$ of $(M,\mathcal E)$ is the first
prolongation of the $(k-1)$st prolongation of $(M,\mathcal E)$.
This gives rise to the {\em prolongation tower}
$$
  \cdots\to M^{(k)}\to M^{(k-1)}\to\cdots\to M^{(1)}\to M.
$$
An integral manifold of $(M,\mathcal E)$ lifts to an integral manifold
of each $(M^{(k)},\mathcal E^{(k)})$ in this tower.

Two examples will help to clarify the construction.  The first 
is the prolongation tower of the multi-contact system
$(G_n(TX),\mathcal I)$ itself,
and this will give us more detailed information about the structure of
the ideals $\mathcal E^{(k)}$ for general $(M,\mathcal E)$.  The
second is the prolongation tower of the EDS associated to a first-order PDE
system, most of which we leave as an exercise.

\noindent
\textbf{Example 1.}
Consider the multi-contact ideal $\mathcal I$ on $G_n(TX)$, over a
manifold $X$ of dimension $n+s$ with local coordinates
$(x^i,z^\alpha)$.  We can see from
the coordinate expression (\ref{CoordMC17}) that its integral elements
over the dense open subset where
$\bigwedge_i dx^i\neq 0$ are exactly the $n$-planes of the form
$$
  E_{p^\alpha_{ij}} = \{dz^\alpha-p^\alpha_idx^i,\
    dp^\alpha_i - p^\alpha_{ij}dx^j\}^\perp\subset T(G_n(TX)),
$$
for some constants $p^\alpha_{ij}=p^\alpha_{ji}$.  These
$p^\alpha_{ij}$ are local fiber coordinates for the
prolongation $(G_n(TX)^{(1)},\mathcal I^{(1)})$.  Furthermore, with respect to
the full coordinates $(x^i,z^\alpha,p^\alpha_i,p^\alpha_{ij})$ for
$G_n(TX)^{(1)}\subset G_n(TG_n(TX))$, the $1$-jet graphs of integral manifolds
of $\mathcal I\subset\Omega^*(G_n(TX))$ satisfy
$$
  dz^\alpha-p^\alpha_idx^i = 0,\
  dp^\alpha_i-p^\alpha_{ij}dx^j = 0.
$$
It is these $s+ns$ $1$-forms that differentially generate the
prolonged Pfaffian system $\mathcal I^{(1)}$.  It is not difficult to
verify that we
have globally $G_n(TX)^{(1)} = G_{2,n}(X)$, the bundle of $2$-jets of
$n$-dimensional submanifolds of $X$, and that $\mathcal I^{(1)}\subset
\Omega^*(G_{2,n}(X))$ is the Pfaffian system whose transverse integral
manifolds are $2$-jet graphs of submanifolds of $X$.

More generally, let $G_k = G_{k,n}(X)\to X$ be the
bundle of $k$-jets of $n$-dimensional submanifolds of $X$.
Because a $1$-jet 
of a submanifold is the same as a tangent plane, $G_1=G_n(TX)$ is the
original space whose prolongation tower we are describing.
$G_k$ carries a canonical Pfaffian system
$\mathcal{I}_k\subset\Omega^*(G_k)$,
whose transverse integral manifolds are $k$-jet graphs
$f^{(k)}:N\hookrightarrow G_k$ of $n$-dimensional submanifolds
$f:N\hookrightarrow X$.  This is perhaps clearest in coordinates.
Letting $(x^i,z^\alpha)$ be coordinates on $X$, $G_k$ has induced
local coordinates $(x^i,z^\alpha,p^\alpha_i,\ldots,p^\alpha_I)$,
$|I|\leq k$, corresponding to the jet at $(x^i,z^\alpha)$ of the
submanifold 
$$
  \{(\bar x^i,\bar z^\alpha)\in X:\bar z^\alpha = z^\alpha +
    p^\alpha_i(\bar x^i-x^i)+\cdots+\sf1{I!}p^\alpha_I(\bar x-x)^I\}.
$$
In terms of these coordinates, the degree-$1$ part $I_k\subset
T^*(G_k)$ of the Pfaffian system $\mathcal{I}_k$ is generated by
\begin{equation}
  \begin{array}{l}
  \theta^\alpha  =  dz^\alpha -p^\alpha_idx^i, \\
  \theta^\alpha_i  =  dp^\alpha_i-p^\alpha_{ij}dx^j, \\
   \qquad \vdots  \\
  \theta^\alpha_I =  dp^\alpha_I - p^\alpha_{Ij}dx^j, \quad
    |I|=k-1.
  \end{array}
\label{ProlGen17}
\end{equation}
It is not hard to see that the transverse integral manifolds of this
$I_k$ are as described above.  
The point here is that $(G_k,\mathcal I_k)$ is the first prolongation of
$(G_{k-1},\mathcal I_{k-1})$ for each $k>1$, and is therefore the $(k-1)$st
prolongation of the original $(G_1,\mathcal I_1)=(G_n(TX),\mathcal I)$.

For future reference, we note the structure equations
\begin{equation}
  \begin{array}{ll}
  d\theta^\alpha  =  -\theta^\alpha_i\wedge dx^i, & \\
  d\theta^\alpha_i  =  -\theta^\alpha_{ij}\wedge dx^j, & \\
  \qquad \vdots & \\
  d\theta^\alpha_I =  -\theta^\alpha_{Ij}\wedge dx^j, &
    |I|=k-2, \\
  d\theta^\alpha_I  =  -dp^\alpha_{Ij}\wedge dx^j, &
    |I|=k-1.
  \end{array}
\label{ProlStreqn17}
\end{equation}
There is a tower
\begin{equation}
  \cdots\to G_k\to G_{k-1}\to\cdots\to G_1,
\label{MainTower16}
\end{equation}
and one can pull back to $G_k$ any functions or differential forms on
$G_{k^\prime}$, with $k^\prime<k$.  Under these maps, our
different uses of the coordinates $p^\alpha_I$ and forms
$\theta^\alpha_I$ are consistent, and we can also write
$$
  I_{k^\prime}\subset I_k\subset T^*G_k,\qquad
   \mbox{for }k^\prime<k.
$$
None of the $I_k$ is an integrable Pfaffian system.  In fact, the
filtration on $G_k$
$$
  I_k\supset I_{k-1}\supset\cdots\supset I_1\supset 0
$$
coincides with the {\em derived flag} of $I_k\subset T^*G_{k,n}$
(cf.~Ch.~II, \S4 of \cite{Bryant:Exterior}).

\noindent
\textbf{Example 2.}
Our second example of prolongation relates to a first-order PDE system
$F^a(x^i,z^\alpha(x),z^\alpha_{x^i}(x))=0$ for
some unknown functions $z^\alpha(x)$.  The equations
$F^a(x^i,z^\alpha,p^\alpha_i)$ define a locus $M_F$
in the space $J^1(\R^n,\R^s)$ of $1$-jets of maps $z:\R^n\to\R^s$, and
we will assume that this locus is a smooth submanifold which
submersively surjects 
onto $\R^n$.  The restriction to $M_F\subset J^1(\R^n,\R^s)$ of the
multi-contact Pfaffian system
$I_1=\{dz^\alpha-p^\alpha_idx^i\}$ generates an EDS $(M_F, \mathcal
I_F)$.  Now, the 
set of integral elements for $(M_F,\mathcal I_F)$ is a subset of the
set of integral 
elements for $\mathcal I_1$ in $J^1(\R^n,\R^s)$; it consists of those integral
elements of $\mathcal I_1$ which are tangent to $M_F\subset
J^1(\R^n,\R^s)$.  Just as in 
the preceding example, the integral elements of $\mathcal I_1$ may be
identified with elements of the space $J^2(\R^n,\R^s)$ of $2$-jets of maps.
The collection of $2$-jets which correspond to integral elements of
$(M_F, I_F)$ are exactly the $2$-jets satisfying the augmented PDE
system
\begin{eqnarray*}
 0 & = & F^a(x^i,z^\alpha(x), z^\alpha_{x^i}(x)), \\
 0 & = &  \frac{\partial F^a}{\partial x^i} +
   \frac{\partial F^a}{\partial z^\alpha}z^\alpha_{x^i} +
   \frac{\partial F^a}{\partial p^\alpha_j}z^\alpha_{x^ix^j}.
\end{eqnarray*}
Therefore, integral manifolds of the prolongation of the EDS
associated to a PDE system correspond to solutions of this augmented
system.  For this reason, 
prolongation may generally be thought of as adjoining
the derivatives of the original equations.

It is important to note that as the first prolongation of arbitrary
$(M,\mathcal E)$ is embedded in the canonical multi-contact system
$(G_n(TM),\mathcal I)$, so can all higher prolongations $(M^{(k)},\mathcal
E^{(k)})$ be embedded in the prolongations $(G_{k,n}(M), \mathcal I_k)$.  Among
other things, this implies that $\mathcal E^{(k)}$ is locally generated by
forms like (\ref{ProlGen17}), satisfying structure equations
(\ref{ProlStreqn17}), typically with additional linear-algebraic relations.

Of most interest to us is the {\em infinite prolongation}
$(M^{(\infty)},\mathcal
E^{(\infty)})$ of an EDS $(M,\mathcal E)$.  As a space, $M^{(\infty)}$ is
defined as the inverse limit of 
$$ 
  \cdots\stackrel{\pi_{k+1}}{\to}M^{(k)}\stackrel{\pi_k}{\to}
   \cdots\stackrel{\pi_1}{\to}
    M^{(0)}=M;
$$
that is,
$$
  M^{(\infty)} = \{(p_0,p_1,\ldots)\in M^{(0)}\times M^{(1)}\times\cdots:
   \pi_k(p_k)=p_{k-1}\mbox{ for each }k\geq 1\}.
$$
An element of $M^{(\infty)}$ may be thought of as a Taylor series
expansion for a possible integral manifold of $(M,\mathcal E)$.
$M^{(\infty)}$ is generally of infinite dimension, but its presentation as
an inverse limit will prevent us from facing analytic
difficulties.  In particular, smooth functions and differential forms
are by definition the corresponding objects on
some finite $M^{(k)}$, pulled up to $M^{(\infty)}$ by the projections.  
It therefore makes sense to define
$$
  \mathcal E^{(\infty)} = \bigcup_{k>0} \mathcal E^{(k)},
$$
which gives an EDS on $M^{(\infty)}$ whose transverse integral
manifolds are the infinite-jet graphs of integral manifolds of
$(M^{(0)},\mathcal E^{(0)})$.  
$\mathcal E^{(\infty)}$ is a Pfaffian system, differentially
generated by its degree-$1$ part $I^{(\infty)} = \bigcup I^{(k)}$, where
each $I^{(k)}$ is the degree-$1$ part of $\mathcal E^{(k)}$.  
In fact, we can see from (\ref{ProlStreqn17}) that $\mathcal
E^{(\infty)}$ is {\em algebraically} generated by 
$I^{(\infty)}$; that is, $I^{(\infty)}$ is a formally integrable
Pfaffian system, although this is not true of any finite $I^{(k)}$.
However, there is no analog of the Frobenius theorem for the
infinite-dimensional $M^{(\infty)}$, so we must be cautious about how we
use this fact.

Vector fields on $M^{(\infty)}$ are more subtle.  By definition,
$\mathcal V(M^{(\infty)})$ is the Lie algebra of derivations of the ring
$\mathcal R(M^{(\infty)})$ of smooth functions on $M^{(\infty)}$.
In case $M^{(\infty)} = J^\infty(\R^n,\R^s)$, a vector field is of the
form
$$
  v = v^i\sf{\partial}{\partial x^i} +
  v^\alpha_0\sf{\partial}{\partial z^\alpha}
   +\cdots+v^\alpha_I\sf{\partial}{\partial p^\alpha_I}+\cdots.
$$
Each coefficient $v^\alpha_I$ is a function on some $J^k(\R^n,\R^s)$,
possibly with $k>|I|$.  Although $v$ may have infinitely many terms,
only finitely many appear in its application to any particular 
$f\in\mathcal R(M^{(\infty)})$, so there are no issues of convergence.

\subsection{Noether's Theorem}

To give the desired generalization of Noether's theorem, we must first
discuss a generalization of the infinitesimal symmetries used in the
classical version.
For convenience, we change notation and let $(M,\mathcal E)$ denote the infinite prolongation of an
exterior differential system $(M^{(0)},\mathcal E^{(0)})$.
\begin{Definition}
A {\em generalized symmetry} of $(M^{(0)},\mathcal E^{(0)})$
is a vector field $v\in\mathcal V(M)$ such that $\mathcal L_v\mathcal
E\subseteq\mathcal E$.  A {\em trivial generalized symmetry} is a
vector field $v\in\mathcal V(M)$ such that $v\innerprod\mathcal
E\subseteq\mathcal E$.  The space $\lie{g}$ of {\em proper generalized
  symmetries} is the quotient of the space of generalized symmetries
by the subspace of trivial generalized symmetries.
\end{Definition}
Several remarks are in order.
\begin{itemize}
\item The Lie derivative in the definition of generalized symmetry
  is defined by the Cartan formula
$$
  \mathcal L_v\varphi = v\innerprod d\varphi + d(v\innerprod\varphi).
$$
The usual definition involves a flow along $v$, which may not exist in
this setting.
\item A trivial generalized symmetry is in fact a generalized
  symmetry; this is an immediate consequence of the fact that
  $\mathcal E$ is differentially closed.
\item The space of generalized symmetries has the obvious structure of
  a Lie algebra.
\item Using the fact that $\mathcal E$ is a formally integrable
  Pfaffian system, it is easy to show that the condition
  $v\innerprod\mathcal E\subseteq\mathcal E$ for $v$ to be a trivial
  generalized symmetry is equivalent to the condition $v\innerprod
  I=0$, where $I=\mathcal E\cap\Omega^1(M)$ is the degree-$1$ part
  of $\mathcal E$.
\item The vector subspace of trivial generalized symmetries is an
  ideal in the Lie algebra of generalized symmetries, so $\lie{g}$ is
  a Lie algebra as well.  The following proof of this fact uses the preceding
  characterization $v\innerprod I=0$ for trivial generalized symmetries:
if $\mathcal L_v\mathcal E\subset\mathcal E$, $w\innerprod I=0$, and
$\theta\in\Gamma(I)$, then
\begin{eqnarray*}
  [v,w]\innerprod\theta & = & -w\innerprod(v\innerprod d\theta)
    +v(w\innerprod\theta)-w(v\innerprod\theta) \\
  & = & -w\innerprod(\mathcal L_v\theta-d(v\innerprod\theta)) + 0
    -w\innerprod d(v\innerprod\theta) \\
  & = & -w\innerprod\mathcal L_v\theta \\
  & = & 0.
\end{eqnarray*}
\end{itemize}
The motivation for designating certain generalized symmetries as
trivial comes from a formal calculation which shows that a trivial
generalized symmetry is tangent to any integral manifold of the
formally integrable Pfaffian system $\mathcal E$.  Thus, the ``flow''
of a trivial generalized symmetry does not permute the integral
manifolds of $\mathcal E$, but instead acts by diffeomorphisms of each
``leaf''.

The following example is relevant
to what follows.  Let $M^{(0)}=J^1(\R^n,\R)$ be the standard contact
manifold of $1$-jets of functions, with global coordinates $(x^i,z,p_i)$ and
contact ideal $\mathcal E^{(0)}=\{dz-p_idx^i,\,dp_i\wedge dx^i\}$.  
The infinite prolongation of $(M^{(0)},\mathcal E^{(0)})$ is
$$
  M = J^\infty(\R^n,\R),\quad\mathcal E=\{\theta_I:|I|\geq 0\},
$$
where $M$ has coordinates $(x^i,z,p_i,p_{ij},\ldots)$, and
$\theta_I = dp_I-p_{Ij}dx^j$.  (For the empty index
$I=\emptyset$, we let $p=z$, so $\theta=dp-p_idx^i$ is the
original contact form.)  Then the trivial generalized symmetries of
$(M,\mathcal E)$ are
the total derivative vector fields
$$
  D_i = \sf\partial{\partial x^i} + p_i\sf{\partial}{\partial z} +
    \cdots + p_{Ii}\sf{\partial}{\partial p_I}+\cdots.
$$
We will determine the proper generalized symmetries of $(M,\mathcal
E)$ shortly.

There is another important feature of a vector field on the infinite
prolongation $(M,\mathcal E)$ of $(M^{(0)},\mathcal E^{(0)})$, which
is its {\em order}.  To introduce
this, first note that any vector field $v_0\in\mathcal V(M^{(0)})$ on
the original, finite-dimensional manifold induces a vector field and a
flow on each finite prolongation $M^{(k)}$, and therefore induces on
$M$ itself a vector field $v\in\mathcal V(M)$ having a flow.  A
further special property of $v\in\mathcal V(M)$ induced by
$v_0\in\mathcal V(M^{(0)})$ is that $\mathcal L_v(I_k)\subseteq I_k$
for each $k\geq 1$.  Though it is tempting to try to characterize
those $v\in\mathcal V(M)$ induced by such $v_0$ using this last
criterion, we ought not to do so, because this is not a criterion that
can be inherited by {\em proper} generalized symmetries of
$(M,\mathcal E)$.  Specifically, an arbitrary trivial generalized
symmetry $v\in\mathcal V(M)$ only satisfies
$$
  \mathcal L_v(I_k)\subseteq I_{k+1},
$$
so a generalized symmetry $v$ can be {\em equivalent} (modulo trivials) to one
induced by a $v_0\in\mathcal V(M^{(0)})$, without satisfying $\mathcal
L_v(I_k)\subseteq I_k$.  Instead, we have the following.
\begin{Definition} For a vector field $v\in\mathcal V(M)$, the {\em
    order} of $v$, written $o(v)$, is the minimal $k\geq 0$ such that
  $\mathcal L_v(I_0)\subseteq I_{k+1}$.
\end{Definition}
With the restriction $o(V)\geq 0$, the orders of
equivalent generalized symmetries of $\mathcal E$ are equal.
A vector field induced by $v_0\in\mathcal V(M^{(0)})$ has order $0$.
Further properties are:
\begin{itemize}
\item
  $o(v)=k$ if and only if for each $l\geq 0$, $\mathcal
  L_v(I_l)\subseteq I_{l+k+1}$;
\item letting $\lie{g}_k=\{v:o(v)\leq k\}$, we have
  $[\lie{g}_k,\lie{g}_l]\subseteq\lie{g}_{k+l}$.
\end{itemize}

We now investigate the generalized symmetries of the
prolonged contact system on $M=J^\infty(\R^n,\R)$.
The conclusion will be that the proper generalized symmetries
correspond to smooth functions on $M$; this is analogous to the
finite-dimensional contact case, in which we could locally associate
to each contact symmetry its generating function, and conversely.
Recall that we have a coframing $(dx^i,\theta_I)$ for $M$, satisfying
$d\theta_I=-\theta_{Ij}\wedge dx^j$ for all multi-indices $I$.  To describe
vector fields on $M$, we will work with the dual framing
$(D_i,\partial/\partial \theta_I)$, which in terms of the usual
framing $(\partial/\partial x^i,\partial/\partial p_I)$ is given by
\begin{eqnarray*}
  D_i & = & \frac{\partial}{\partial x^i} + \sum_{|I|\geq 0}
    p_{Ii}\frac{\partial}{\partial p_I}, \\
  \frac{\partial}{\partial\theta_I} & = & \frac{\partial}{\partial
    p_I}.
\end{eqnarray*}
The vector fields $D_i$ may be thought of as ``total derivative''
operators, and applied to a function $g(x^i,z,p_i,\ldots,p_I)$ on some
$J^k(\R^n,\R)$ give
$$
  (D_ig)(x^j,z,p_j,\ldots,p_I,p_{Ij}) = \sf{\partial g}{\partial
    x^i}+p_i\sf{\partial g}{\partial z}+p_{ij}\sf{\partial g}{\partial
    p_j}+\cdots+p_{Ii}\sf{\partial g}{\partial p_I},
$$
which will generally be defined only on $J^{k+1}(\R^n,\R)$ rather
than $J^k(\R^n,\R)$.  These operators can be composed, and we set
$$
  D_I = D_{i_1}\circ\cdots\circ D_{i_k},\qquad I=(i_1,\ldots,i_k).
$$
We do this because the proper generalized symmetries of $\mathcal
I=\{\theta_I:|I|\geq 0\}$ are uniquely represented by vector fields
\begin{equation}
  v = g\frac{\partial}{\partial\theta} + g_i\frac{\partial}{\partial
    \theta^i} +\cdots+ g_I\frac{\partial}{\partial \theta_I}+\cdots,
\label{FormVF17}
\end{equation}
where $g=v\innerprod\theta$ and
\begin{equation}
  g_I = D_Ig,\qquad |I|>0.
\label{Determine17}
\end{equation}
To see this, first note that any vector field is congruent modulo
trivial generalized symmetries to a unique one of the form
(\ref{FormVF17}).  It then follows from a straightforward calculation
that a vector field of the form (\ref{FormVF17}) is a generalized
symmetry of $\mathcal I$ if and only if it satisfies (\ref{Determine17}).
If one defines $\mathcal R_k\subset\mathcal R(M)$ to consist of
functions pulled back from $J^k(\R^n,\R)$, then one can verify that
for any proper generalized symmetry $v\in\lie{g}$,
$$
   o(v)\leq k\qquad\Longleftrightarrow \qquad
   g=v\innerprod\theta\in\mathcal R_{k+1}.
$$

The general version of Noether's theorem will involve proper
generalized symmetries.  However, recall that our first-order version
requires us to distinguish among the symmetries of an Euler-Lagrange
system the symmetries of the original variational problem; only the
latter give rise to conservation laws.  We therefore have to give the
appropriate corresponding notion for proper generalized symmetries.

For this purpose, we introduce the following algebraic apparatus.  We
filter the differential forms $\Omega^*(M)$ on the infinite
prolongation $(M,\mathcal{E})$ 
of an Euler-Lagrange system $(M^{(0)},\mathcal E^{(0)})$ by letting
\begin{equation}
   I^p\Omega^{p+q}(M) = \mbox{Image}
   (\underbrace{ \mathcal E\otimes \cdots\otimes
   \mathcal E}_p \otimes\Omega^*(M) \to
     \Omega^*(M))\cap \Omega^{p+q}(M).
\label{filtration}
\end{equation}
We define the associated graded objects
$$
  \Omega^{p,q}(M) = I^p\Omega^{p+q}(M)/I^{p+1}\Omega^{p+q}(M).
$$
Because $\mathcal E$ is formally integrable, the exterior derivative
$d$ preserves this filtration and its associated graded objects:
$$
  d:\Omega^{p,q}(M) \to \Omega^{p,q+1}(M).
$$
We define the cohomology
$$
  H_\Lambda^{p,q}(M) = \frac
    {\mbox{Ker}(d:\Omega^{p,q}(M)\to\Omega^{p,q+1}(M))}
    {\mbox{Im}(d:\Omega^{p,q-1}(M)\to\Omega^{p,q}(M))}.
$$
A simple diagram-chase shows that the exterior derivative operator
$d$ induces a map $d_1:H_\Lambda^{p,q}(M)\to
H_\Lambda^{p+1,q}(M)$.\footnote{Of course, 
  $H_\Lambda^{*,*}(M)$ is the $E_1$-term of a spectral sequence.
  Because we will not be using any of the higher terms, however, there
  is no reason to invoke this machinery.  Most of this theory was
  introduced in \cite{Vinogradov:CSpectral}.}
Now, the Poincar\'e-Cartan form $\Pi\in\Omega^{n+1}(M^{(0)})$
pulls back to an element $\Pi\in I^2\Omega^{n+1}(M)$
which is closed and therefore defines a class $[\Pi]\in
H_\Lambda^{2,n-1}(M)$.

It follows from the definition that a generalized symmetry of
$\mathcal E$ preserves the filtration $I^p\Omega^{p+q}(M)$, and
therefore acts on the cohomology group $H_\Lambda^{2,n-1}(M)$.  The
generalized symmetries appropriate for Noether's theorem are exactly
those generalized symmetries $v$ of $\mathcal{E}$ satisfying the
additional condition
$$
  \mathcal L_v[\Pi] = 0 \in H_\Lambda^{2,n-1}(M).
$$
In other words, $v$ is required to preserve $\Pi$ modulo (a) forms in
$I^3\Omega^{n+1}(M)$, and (b) derivatives of forms in
$I^2\Omega^n(M)$.  We also need to verify 
that a {\em trivial} generalized symmetry $v$ preserves the class $[\Pi]$;
this follows from the fact that $v\innerprod I = 0$, for then
$v\innerprod\Pi \in I^2\Omega^n(M)$, so that
$$
  \mathcal L_v[\Pi] = [d(I^2\Omega^n(M))] = 0.
$$
We now have the Lie subalgebra $\lie{g}_{[\Pi]}\subseteq\lie{g}$ of
proper generalized symmetries of the variational problem.  It is worth
noting that this requires only that we have $\Pi$ defined modulo
$I^3\Omega^{n+1}(M)$.  A consequence is that even in the most general
higher-order, multi-contact case where a canonical Poincar\'e-Cartan
form is not known to exist, there should be a version of
Noether's theorem that includes both the first-order
multi-contact version discussed in the previous section, and the
higher-order scalar version discussed below.  However, we will not
pursue this.

\index{conservation law!higher-order|(}
The other ingredient in Noether's theorem is a space of conservation
laws, defined by analogy with previous cases as
$$
  \mathcal C(\mathcal E) = H^{n-1}(\Omega^*/\mathcal E, \bar d)
     = H^{0,n-1}_\Lambda(M),
$$
where the last notation refers to the cohomology just introduced.
It is a substantial result (see~\cite{Bryant:CharacteristicI}) that
over contractible subsets of $M^{(0)}$ we can use the exterior
derivative to identify $\mathcal C(\mathcal E)$ with
$$
  \bar{\mathcal C}(\mathcal E) \stackrel{\mathit{def}}{=}
     \mbox{Ker} (d_1:
    H_\Lambda^{1,n-1}(M) \to H_\Lambda^{2,n-1}(M)).
$$
Now we can identify conservation laws as classes of $n$-forms, as in
the previous case of Noether's theorem, and we will do so without
comment in the following.

We define a Noether map 
$\lie{g}_{[\Pi]}\to\bar{\mathcal C}(\mathcal E)$ as 
$v\mapsto v\innerprod \Pi$.
To see that this is well-defined, first note that
$$
  v\innerprod\Pi \in I^1\Omega^n(M),
$$
so $v\innerprod\Pi$ represents an element of $\Omega^{1,n-1}(M)$,
which we shall also denote as $v\innerprod\Pi$.
Furthermore, its exterior derivative is
\begin{equation}
  d(v\innerprod\Pi) = \mathcal L_v\Pi,
\label{AnotherDerivative}
\end{equation}
and this lies in $I^2\Omega^{n+1}$, simply because $v$ preserves
$\mathcal E$ and therefore also the
filtration (\ref{filtration}).  Consequently, 
$$
  v\innerprod\Pi \in \mbox{Ker}(d:\Omega^{1,n-1}(M)\to
     \Omega^{1,n}(M)),
$$
and we therefore have an element 
$$
  [v\innerprod\Pi]\in H_\Lambda^{1,n-1}(M).
$$
Finally, we need to verify that
$$
  [v\innerprod\Pi]\in\mbox{Ker}(d_1:H_\Lambda^{1,n-1}(M)\to
       H_\Lambda^{2,n-1}(M)).
$$
This follows from the hypothesis that $v$ preserves not only
the Euler-Lagrange system $\mathcal{E}$ and associated
filtration (\ref{filtration}), but also the class
$[\Pi]$.  Specifically, the image
$$
  d_1([v\innerprod\Pi])\in H_\Lambda^{2,n-1}(M)
$$
is represented by the class (see (\ref{AnotherDerivative}))
$$
  [d(v\innerprod\Pi)] = [\mathcal L_v\Pi] = {\mathcal L}_v[\Pi] = 0.
$$
This proves that
$$
  v\mapsto [v\innerprod\Pi]
$$
defines a map between the appropriate spaces.

We can now make the following proposal for a general form of Noether's
theorem.
\begin{Conjecture}
Let $(M,\mathcal E)$ be the infinite prolongation of an
Euler-Lagrange system, and assume that the system is non-degenerate and
that $H^q_{dR}(M)=0$ for all $q>0$.  
Then the map $v\mapsto[v\innerprod\Pi]$ induces an isomorphism
$$
   \lie{g}_{[\Pi]}\stackrel{\sim}{\longrightarrow}
     \bar{\mathcal C}(\mathcal{E}).
$$
\end{Conjecture}
It is quite possible that this is already essentially proved
in~\cite{Vinogradov:CSpectral} or~\cite{Olver:Applications}, 
but we have not been able
to determine the relationship between their statements and ours.  In
any case, it would be illuminating to have a proof of the present
statement in a spirit similar to that of our 
Theorem~\ref{NoetherTheorem}. 

To clarify this, we will describe how it appears in coordinates.  First,
note that for the classical Lagrangian
$$
  L(x^i,z,p_i)dx,
$$
the Euler-Lagrange equation
$$
  E(x^i,z,p_i,p_{ij})=(\textstyle\sum D_jL_{p_j}-L_z)(x^i,z,p_i,p_{ij})
  = 0
$$
defines a locus $M^{(1)}\subset J^2(\R^n,\R)$, and the first
prolongation of the Euler-Lagrange system $(J^1(\R^n,\R),\mathcal
E_L)$ discussed previously is given by the restriction of the
second-order contact Pfaffian system on $J^2(\R^n,\R)$ to this locus.
Higher prolongations are defined by setting
$$
   M^{(k)}=\{E_I\stackrel{\mathit{def}}{=}D_IE=0,\
    |I|\leq k-1\}\subset J^{k+1}(\R^n,\R)
$$
and restricting the $(k+1)$st-order contact system $\mathcal
I^{(k+1)}$.  We will consider generalized symmetries of
$(M^{(\infty)},\mathcal E^{(\infty)})$ which arise as restrictions of
those generalized symmetries of $(J^\infty(\R^n,\R),\mathcal I)$
which are also tangent to $M^{(\infty)}\subset J^\infty(\R^n,\R)$.  This
simplifies matters insofar as we can understand generalized symmetries
of $\mathcal I$ by their generating functions.
The tangency condition is
\begin{equation}
  \mathcal L_v(E_I)|_{M^{(\infty)}}=0,\qquad |I|\geq 0.
\label{InfiniteTangency17}
\end{equation}
This Lie derivative is just the action of a vector field as a
derivation on functions.  Now, for a generalized symmetry $v$ of the
infinite-order contact system, all of
the conditions (\ref{InfiniteTangency17}) follow from just the first one,
$$
  \mathcal L_v(E)|_{M^{(\infty)}}=0.
$$
If we let $v$ have generating function $g=v\innerprod\theta\in\mathcal
R(J^\infty(\R^n,\R))$, then we can see from
(\ref{FormVF17}, \ref{Determine17}) that this 
condition on $g$ is
\begin{equation}
  \sum_{|I|\geq 0} \sf{\partial E}{\partial p_I}D_Ig =0
   \mbox{ on }M^{(\infty)}.
\label{SymCondn17}
\end{equation}
We are again using $p=z$ for convenience.
For instance,  $E=\sum p_{ii}-f(z)$ defines the Poisson equation
$\Delta z=f(z)$, and the preceding condition is
\begin{equation}
  \sum_{i}D_i^2g-f^\prime(z)g=0 \mbox{ on }M^{(\infty)}.
\label{PoissonGF17}
\end{equation}

We now consider
the Noether map for $M^{(\infty)}\subset J^\infty(\R^n,\R)$.
We write the Poincar\'e-Cartan form pulled back to
$M^{(\infty)}$ using coframes adapted to this infinite prolongation,
starting with
$$
  dL_{p_i} = D_j(L_{p_i})dx^j + L_{p_iz}\theta+L_{p_ip_j}\theta_j,
$$
and then the Poincar\'e-Cartan form on $J^\infty(\R^n,\R)$ is
\begin{eqnarray*}
  \Pi & = & d(L\,dx+\theta\wedge L_{p_i}dx_{(i)}) \\
    & = & \theta\wedge\left((-D_i(L_{p_i})+L_z)dx-\theta_j\wedge
           L_{p_ip_j}dx_{(i)}\right).
\end{eqnarray*}
Restriction to $M^{(\infty)}\subset J^\infty(\R^n,\R)$ kills the first
term, and we have
$$
  \Pi = -L_{p_ip_j}\theta\wedge\theta_j\wedge dx_{(i)}.
$$
We then apply a vector field $v_g$ with generating function $g$, and
obtain
$$
  v_g\innerprod\Pi = -gL_{p_ip_j}\theta_j\wedge dx_{(i)} +
    (D_jg)L_{p_ip_j}\theta\wedge dx_{(i)}.
$$
This will be the ``differentiated form'' of a conservation law
precisely if ${\mathcal L}_{v_g}[\Pi]=0$, that is, if
$$
  d(v_g\innerprod\Pi)\equiv 0
    \pmod{I^3\Omega^{n+1}(M)+dI^2\Omega^n(M)}.
$$

\

Concerning generalized symmetries of a PDE, note that in the condition 
(\ref{SymCondn17}) for $g=g(x^i,p,p_i,\ldots,p_I)\in\mathcal
R_k=C^\infty(J^k(\R^n,\R))$, the variables $p_I$ with $|I|>k$ appear only
polynomially upon taking the total derivatives $D_Jg$.  In other
words, the condition on $g$ is polynomial in the variables $p_I$ for $|I|>k$.
Equating coefficients of these polynomials gives a PDE system to be
satisfied by a generalized symmetry of an Euler-Lagrange equation.
With some effort, one can analyze the situation for our Poisson
equation $\Delta z = f(z)$ and find the following.
\begin{Proposition}
 If $n\geq 3$, then a solution $g=g(x^i,p,p_i,\ldots,p_I)$ of order
 $k$ to (\ref{PoissonGF17}) is equal on $M^{(\infty)}$ to a
 function that is
linear in the variables $p_J$ with $|J|\geq k-2$.  
If in addition $f^{\prime\prime}(z)\neq 0$, so that the
Poisson equation is non-linear, then every solution's restriction to
$M^{(\infty)}$ is the pullback of a function on $M^{(0)}\subset
 J^2(\R^n,\R)$, which generates a classical symmetry of the equation.
\label{NoCLs17}
\end{Proposition}
In other words, a non-linear Poisson equation in $n\geq 3$ independent
variables has
no non-classical generalized symmetries, and consequently no
higher-order conservation laws.

\

\begin{Proof}
Because the notation involved here becomes rather tedious, we will
sketch the proof and leave it to the reader to verify the
calculations.  We previously hinted at the main idea:  the
condition (\ref{PoissonGF17}) on a generating function
$g=g(x^i,p,\ldots,p_I)$, $|I| = k$, is polynomial in the highest order
variables with coefficients depending on partial derivatives of $g$.
To isolate these terms, we filter the functions on $M^{(\infty)}$ by
letting $\overline{\mathcal R_l}$ be the image of $\mathcal R_l$ under
restriction to $M^{(\infty)}$; in other words, $\overline{\mathcal R_l}$
consists of functions which can be expressed as functions of
$x^i,p_J$, $|J|\leq l$, {\em after} substituting the defining
relations of $M^{(\infty)}$,
$$
  p_{Jii} = \frac{d^{|J|}f}{dx^J}.
$$
To calculate in $\overline{\mathcal R_l}$ we will need to define
variables $q_J$ to be the harmonic parts of $p_J$; that is,
\begin{eqnarray*}
  q_i  & = & p_i, \\
  q_{ij} & = & p_{ij} - \sf1n\delta_{ij}p_{ll} \\
    & = & p_{ij}-\sf1n\delta_{ij}f(p), \\
  q_{ijk} & = & p_{ijk} - \sf1{n+2}(\delta_{ij}p_{kll}
    +\delta_{jk}p_{ill} + \delta_{ki}p_{jll}) \\
   & = & p_{ijk}-\sf1{n+2}(\delta_{ij}p_k + \delta_{jk}p_i
    + \delta_{ki}p_j)f^\prime(p), \quad\mbox{\&c}.
\end{eqnarray*}
These, along with $x_i$ and $p$, give coordinates on $M^{(\infty)}$.
In addition to working modulo various $\overline{\mathcal R_l}$ to isolate
terms with higher-order derivatives, we will also at times work modulo
functions that are {\em linear} in the $q_I$.  In what follows, we use
the following index conventions:  $p^{(l)} = (p,p_j,\ldots,p_J)$
(with $|J|=l$) denotes the derivative variables up to order $l$, and
the multi-indices $I$, $K$, $A$, satisfy $|I|=k$, $|K|=k-1$,
$|A|=k-2$. 

Now, starting with $g=g(x^j,p^{(k)})\in\mathcal R_k$, we note that
$$
  0 = \sum_i D_i^2g - f^\prime g \in
    \overline{\mathcal R_{k+1}};
$$
that is, the possible order-$(k+2)$ term resulting from two
differentiations of $g$ already drops to order $k$ when restricted to the
equation manifold.
We consider this expression modulo $\overline{\mathcal R_k}$\footnote{In
  this context, ``modulo'' refers to quotients of vector spaces by
  subspaces, not of rings by ideals, as in exterior algebra.}, and
obtain a quadratic polynomial in $q_{Ij}$ with coefficients
in $\overline{\mathcal R_k}$.  We consider only the quadratic terms of this
polynomial, which are
\begin{equation}
  0 \equiv
   \sum_{\begin{smallmatrix} |I|,|I^\prime|=k \\ 1\leq i\leq n
     \end{smallmatrix}}
    \frac{\partial^2g}{\partial p_I\partial p_{I^\prime}}
    q_{Ii}q_{I^\prime i}.
\label{HarmonicOrtho17}
\end{equation}
To draw conclusions about $\frac{\partial^2g}{\partial p_I\partial
  p_{I^\prime}}$ from this, we need the following fundamental lemma,
  in which the difference between the cases $n=2$ and $n\geq 3$
  appears:
\begin{quote}
{\em If ${\mathbf H}_l\subset \mathit{Sym}^l(\R^n)^*$ denotes the space of
  degree-$l$ homogeneous harmonic polynomials on $\R^n$, $n\geq 3$, then the
  $O(n)$-equivariant contraction map
$ {\mathbf H}_{l+1}\otimes{\mathbf H}_{m+1} \to
    {\mathbf H}_l\otimes{\mathbf H}_m $ given by
$$
  \quad X_Cz^C\otimes Y_Qz^Q \mapsto \sum_i 
         X_{Bi}z^B\otimes Y_{Pi}z^P
$$
is surjective; here, $|C|-1=|B|=l$, $|Q|-1=|P|=m$.
}
\end{quote}
We will apply this in situations where a given $g^{BP}\in
\mathit{Sym}^l(\R^n)\otimes \mathit{Sym}^m(\R^n)$ is known to
annihilate all $X_{Bi}Y_{Pi}$
with $X_C$, $Y_Q$ harmonic, for then we have $g^{BP}$ orthogonal to
$\mathbf{H}_l\otimes\mathbf{H}_m\subset \mathit{Sym}^l(\R^n)^*\otimes
\mathit{Sym}^m(\R^n)^*$.  In particular, from (\ref{HarmonicOrtho17})
we have 
$$
  \frac{\partial^2g}{\partial p_I\partial p_J}q_Iq_J \equiv 0.
$$
This means that the restriction of the function $g$ to the hyperplanes
$p_{Aii} = \frac{d^{|A|}f}{dx^A}$ is linear in the highest $p_I$;
in other words, we can write
$$
  g(x^i,p^{(k)}) = h(x,p^{(k-1)}) + h^I(x,p^{(k-1)})p_I
$$
for some $h^I\in\mathcal{R}_{k-1}$.
We can further assume that all $h^{Aii}=0$, where $|A|=k-2$.
This completes the first step.

The second step is to simplify the functions $h^I(x,p^{(k-1)})$,
substituting our new form of $g$ into the condition (\ref{PoissonGF17}).
Again, the ``highest'' terms appear modulo $\overline{\mathcal R_k}$,
and are
$$
  0 \equiv 2\frac{\partial h^I}{\partial p_K}q_{Ki}q_{Ii}
     + 2\frac{\partial h^I}{\partial p_A}q_{Ai}q_{Ii};
$$
here we recall our index convention $|A| = k-2$, $|K|=k-1$, $|I|=k$.
Both terms must vanish separately, and for the first,
our lemma on harmonic polynomials gives that
$[h]^{IK}\stackrel{\mathit{def}}{=}\frac{\partial h^I}{\partial p_K}$ is
orthogonal to harmonics; but then our normalization hypothesis
$h^{Aii}=0$ gives that $[h]^{IK}=0$.  We conclude that
$$
  g(x^i,p^{(k)}) = h(x,p^{(k-1)}) + h^I(x,p^{(k-2)})p_I, \qquad
    |I|=k.
$$
For the second term, our lemma gives similarly that $\frac{\partial
  h^I}{\partial p_A} = 0$, so have
$$
  g(x^i,p^{(k)}) = h(x,p^{(k-1)}) + h^I(x,p^{(k-3)})p_I.
$$
This completes the second step.

For the third step, we again substitute the latest form of $g$ into
the condition (\ref{PoissonGF17}), and now work modulo
$\overline{\mathcal  R_{k-1}}$.  The only term non-linear in the $p_I$ is
$$
  0 \equiv \frac{\partial^2h}{\partial p_K\partial
      p_{K^\prime}}p_{Ki}p_{K^\prime i}, 
$$
and as before, the lemma implies that $h$ is linear in $p_K$.  Now we
have
$$
  g(x^i,p^{(k)}) = h(x,p^{(k-2)}) + h^K(x,p^{(k-2)})p_K +
    h^I(x,p^{(k-3)})p_I.
$$

Again working modulo $\mathcal R_{k-1}$, the only term that is
non-linear in 
$(p_K,p_I)$ is $\frac{\partial h^K}{\partial p_A}p_{Ai}p_{Ki}$, so as
before, we can assume $\frac{\partial h^K}{\partial p_A}=0$ and write
$$
  g(x^i,p^{(k)}) = h(x,p^{(k-2)}) + h^K(x,p^{(k-3)})p_K +
    h^I(x,p^{(k-3)})p_I.
$$

The final step is similar, and gives that $h$ is linear in $p_A$.
This yields
\begin{equation}
  g(x^i,p^{(k)}) = h + h^Ap_A + h^Kp_K + h^Ip_I,
\label{FinalGen17}
\end{equation}
where each of $h, h^A, h^K, h^I$ is a function of $(x^i,p^{(k-3)})$.
This is the first statement of the proposition.

To derive the second statement, we use the form (\ref{FinalGen17}) in
the condition (\ref{PoissonGF17}) modulo $\overline{\mathcal R_{k-2}}$, in
which the only term {\em non-linear} in $(p_K,p_I)$ is
\begin{eqnarray*}
  0 & \equiv & h^I(x^i,p^{(k-3)})\frac{d^{|I|}f}{dx^I} \\
    & \equiv & h^I(x^i,p^{(k-3)})f^{\prime\prime}(p)
     \sum_{i\in I}p_{I\backslash i}p_i.
\end{eqnarray*}
In particular, if $f^{\prime\prime}(p)\neq 0$, then we must have $h^I
= 0$.  Therefore $g\in\mathcal R_k$ actually lies in
$\mathcal R_{k-1}$, and we can induct downward on $k$, eventually
proving that $g\in\mathcal R_1$, as desired.
\end{Proof}

\

We mention two situations which contrast sharply with that of the
non-linear Poisson equation in dimension $n\geq 3$.  First, in the
case of a {\em linear} Poisson equation $\Delta z = f(z)$,
$f^{\prime\prime}(z)=0$ (still in $n\geq 3$ dimensions), one can
extend the preceding argument to show 
that a generating function for a conservation law is linear in all of
the derivative variables $p_I$.  In particular, the infinite
collection of conservation laws for the Laplace equation $\Delta z =
0$ can be determined in this manner; it is interesting to see how all
of these disappear upon the addition of a non-linear term to the
equation.

Second, in dimension $n=2$, there are well-known non-linear Poisson
equations $\Delta u = \sinh u$ and $\Delta u = e^u$ having infinitely
many higher-order conservation laws, but we will not discuss these.


\subsection{The $K=-1$ Surface System}
\label{Subsection:Gauss}
\index{Gauss curvature|(}

In \S\ref{Section:Hypersurface}, we constructed Monge-Ampere systems
on the contact manifold
$M^5$ of oriented tangent planes to Euclidean space $\mathbf E^3$, whose
integral manifolds corresponded to linear Weingarten
surfaces\index{Weingarten equation}.  We
briefly recall this setup for the case of surfaces with Gauss
curvature $K=-1$.  Our index ranges are now $1\leq a,b,c\leq 3$,
$1\leq i,j,k\leq 2$.

Let $\mathcal F\to\mathbf E^3$ be the Euclidean frame
bundle\index{Euclidean!frame bundle}, with its
global coframing $\omega^a$, $\omega^a_b=-\omega^b_a$
satisfying the structure equations
$$
  d\omega^a = -\omega^a_b\wedge\omega^b,\quad
  d\omega^a_b = -\omega^a_c\wedge\omega^c_b.
$$
We set $\theta = \omega^3$, $\pi_i=\omega^3_i$, and then
$\theta\in\Omega^1(\mathcal F)$ is the pullback of a global contact
form on $M=G_2(T\mathbf E^3)$.  The forms that are semibasic over $M$
are generated by $\theta,\omega^i,\pi_i\in\Omega^1(\mathcal F)$.

We define the $2$-forms on $\mathcal F$
$$
  \begin{array}{l}
  \Theta  =  d\theta  =  -\pi_1\wedge\omega^1-\pi_2\wedge\omega^2, \\
  \Psi  =  \pi_1\wedge\pi_2+\omega^1\wedge\omega^2,
  \end{array}
$$
which are pullbacks of uniquely determined forms on $M$.
On a transverse integral element $E^2\subset T_mM$ of the contact system
$\mathcal I=\{\theta,\Theta\}$,
on which $\omega^1\wedge\omega^2\neq 0$, there are relations
$$
  \pi_i = h_{ij}\omega^j,\quad h_{ij}=h_{ji}.
$$
In this case,
$$
  \pi_1\wedge\pi_2 = K\omega^1\wedge\omega^2,
$$
where $K=h_{11}h_{22}-h_{12}h_{21}$ is the Gauss curvature of any
surface $N^2\hookrightarrow\mathbf E^3$ whose $1$-jet graph in $M^5$
is tangent to $E\subset T_mM$.
Therefore, transverse integral manifolds of the EDS
\begin{equation}
  \mathcal E = \{\theta,\Theta,\Psi\}
\label{KSystem17}
\end{equation}
correspond locally to surfaces in $\mathbf E^3$ with constant Gauss
curvature $K=-1$.

The EDS $(M,\mathcal E)$ is an example of a hyperbolic Monge-Ampere
system\index{Monge-Ampere system!hyperbolic}; 
this notion appeared in \S\ref{Section:SmallEquiv},
where we used it to specify a branch of the equivalence problem for
Poincar\'e-Cartan forms on contact $5$-manifolds.
The defining property of a hyperbolic Monge-Ampere system $\mathcal
E=\{\theta,\Theta,\Psi\}$ is that modulo the algebraic ideal
$\{\theta\}$, $\mathcal E$ contains two distinct (modulo scaling) {\em
  decomposable}\index{decomposable form} 
$2$-forms; that is, one can find two non-trivial linear
combinations of the form
\begin{eqnarray*}
  \lambda_1\Theta+\mu_1\Psi & = & \alpha_1\wedge\beta_1, \\
  \lambda_2\Theta+\mu_2\Psi & = & \alpha_2\wedge\beta_2.
\end{eqnarray*}
This exhibits two rank-$2$ Pfaffian systems
$I_i=\{\alpha_i,\beta_i\}$, called the {\em characteristic
  systems}\index{characteristic system|(} of
$\mathcal E$, which are easily seen to be independent of choices
(except for which one is $I_1$ and which one is $I_2$).
The relationship between the geometry of the characteristic systems
and that of the original hyperbolic Monge-Ampere system is very
rich (see \cite{Bryant:Hyperbolic}).
Of particular
interest are those hyperbolic systems whose characteristic systems
each contain a non-trivial conservation 
law\index{conservation law!for K@for $K=-1$ surfaces|(}.  We will show that this
holds for the $K=-1$ system introduced above, but only after one
prolongation.  In other words, for the prolonged system $\mathcal E^{(1)}$,
there is also a notion of characteristic systems $I^{(1)}_i$ which
restrict to any integral surface as the original
$I_i$, and each of these $I^{(1)}_i$
contains a non-trivial conservation law for $\mathcal E^{(1)}$.

Returning to the discussion of integral elements of $\mathcal
E=\{\theta,\Theta,\Psi\}$, note that for any integral element
$E\subset T_{(p,H)}M$, given by equations
$$
  \pi_i - h_{ij}\omega^j = 0,
$$ 
there is a unique frame
$(p,(e_1,e_2,e_1\times e_2))\in\mathcal F$ over $(p,e_1\wedge e_2)\in
M$ for which the second fundamental 
form\index{second fundamental form} is normalized as
$$
  h_{11}=a>0,\quad h_{22}=-\sf1a,\quad h_{12}=h_{21}=0.
$$
The tangent lines in $\mathbf E^3$ spanned by these $e_1,e_2$ are the
{\em principal directions}\index{principal direction} 
at $p$ of any surface whose $1$-jet graph is
tangent to $E$; they define an orthonormal frame in which the second
fundamental form is diagonal.  The number $a>0$ is determined by the
plane $E\subset TM$, so to study integral elements of
$(M,\mathcal E)$, and in particular to calculate on its first
prolongation\index{prolongation!of an EDS}, we introduce
$$
  \mathcal F^{(1)} = \mathcal F\times\R^*,
$$
where $\R^*$ has the coordinate $a>0$.  There is a projection
$\mathcal F^{(1)}\to M^{(1)}$, mapping
$$
  (p,e,a)\mapsto
  (p,e_1\wedge e_2,\{\pi_1-a\omega^1,\pi_2+\sf1a\omega^2\}^\perp).
$$
We define on $\mathcal F^{(1)}$ the forms
$$
  \begin{array}{l} \theta_1 = \pi_1-a\omega^1, \\
    \theta_2 = \pi_2+\sf1a\omega^2, \end{array}
$$
which are semibasic for $\mathcal F^{(1)}\to M^{(1)}$. 
The first
prolongation of the system $\mathcal E$ on $M$ is a Pfaffian system on
$M^{(1)}$, which pulls back to $\mathcal F^{(1)}$ as
$$
  \mathcal E^{(1)}  = 
  \{\theta,\theta_1,\theta_2,d\theta_1,d\theta_2\}.
$$
We have structure equations
$$
  \left.\begin{array}{l}
    d\theta_{\ } \equiv -\theta_1\wedge\omega^1-\theta_2\wedge\omega^2
       \equiv 0 \\
    d\theta_1 \equiv
    -da\wedge\omega^1+\sf{1+a^2}{a}\omega^1_2\wedge\omega^2 \\
    d\theta_2 \equiv \sf{1}{a^2}da\wedge\omega^2 +
    \sf{1+a^2}{a}\omega^1_2\wedge\omega^1 
  \end{array} \right\}\pmod{\{\theta,\theta_1,\theta_2\}},
$$
and in particular, we have the decomposable linear combinations
\begin{equation}
  \begin{array}{l}
    -d\theta_1-a\,d\theta_2 = 
      (da-(1+a^2)\omega^1_2)\wedge(\omega^1+\sf1a\omega^2), \\
    -d\theta_1+a\,d\theta_2 =
      (da+(1+a^2)\omega^1_2)\wedge(\omega^1-\sf1a\omega^2).
  \end{array}
\label{FirstProl17}
\end{equation}
The EDS $\mathcal E^{(1)}$ is algebraically generated by
$\theta,\theta_1,\theta_2$, and these two decomposable $2$-forms.  The
characteristic systems are by definition differentially generated by
\begin{equation}
  \begin{array}{l}
    I_1^{(1)} = \{\theta,\ \theta_1,\ \theta_2,\
      da-(1+a^2)\omega^1_2,\ \omega^1+\sf1a\omega^2\}, \\
    I_2^{(1)} = \{\theta,\ \theta_1,\ \theta_2,\
      da+(1+a^2)\omega^1_2,\ \omega^1-\sf1a\omega^2\}.
  \end{array}
\label{DefChar17}
\end{equation}
Now, the ``universal'' second fundamental form can be factored as
\begin{equation}
  II = a(\omega^1)^2-\sf1a(\omega^2)^2 =
     a(\omega^1+\sf1a\omega^2)(\omega^1-\sf1a\omega^2).
\label{IIFactors17}
\end{equation}
These linear factors, restricted any $K=-1$ surface,
define the {\em asymptotic curves}\index{asymptotic curve} 
of that surface, so by comparing
(\ref{DefChar17}) and (\ref{IIFactors17}) we find that:
\begin{quote}
{\em On a $K=-1$ surface, the integral curves of the characteristic
  systems are the asymptotic curves.}
\end{quote}

Now we look for Euclidean-invariant conservation laws in each
$I^{(1)}_i$.  Instead of using Noether's 
theorem\index{Noether's theorem}, we work directly.  We start by setting
$$
  \varphi_1 = f(a)(\omega^1+\sf1a\omega^2)\in I^{(1)}_1,
$$
and seek conditions on $f(a)$ to have $d\varphi_1\in\mathcal
E^{(1)}$.  A short computation using the structure equations gives
$$
  d\varphi_1\equiv (f^\prime(a)(1+a^2)-f(a)a)\omega^1_2\wedge
    (\omega^1+\sf1a\omega^2)\pmod{\mathcal E^{(1)}},
$$
so the condition for $\varphi_1$ to be a conserved $1$-form is
$$
  \frac{f^\prime(a)}{f(a)} = \frac{a}{1+a^2}.
$$
A solution is 
$$
  \varphi_1 = \sf12\sqrt{1+a^2}(\omega^1+\sf1a\omega^2);
$$
the choice of multiplicative constant $\sf12$ will simplify later
computations. 
A similar computation, seeking an appropriate multiple of
$\omega^1-\sf1a\omega^2$, yields the conserved $1$-form
$$
  \varphi_2 = \sf12\sqrt{1+a^2}(\omega^1-\sf1a\omega^2).
$$
\index{characteristic system|)}

On any simply connected integral surface of the $K=-1$
system, there are coordinate functions $s,t$ such that
$$
  \varphi_1 = ds,\qquad \varphi_2 = dt.
$$
If we work in these coordinates, and in particular use the
non-orthonormal coframing $(\varphi_1,\varphi_2)$
then we can write
\begin{eqnarray}
  \omega^1 & = & \left(\frac1{\sqrt{1+a^2}}\right)
                   (\varphi_1+\varphi_2), \\
  \omega^2 & = & \left(\frac a{\sqrt{1+a^2}}\right)
                   (\varphi_1-\varphi_2), \\
\label{FirstFF17}
  I & = & (\omega^1)^2+(\omega^2)^2 =
     (\varphi_1)^2+2\left(\frac{1-a^2}{1+a^2}\right)
        \varphi_1\varphi_2+(\varphi_2)^2, \\
\label{SecondFF17}
  II & = & \left(\frac{4a}{1+a^2}\right)\varphi_1\varphi_2.
\end{eqnarray}
These expressions suggest that we define
$$
  a = \tan z,
$$
where $a>0$ means that we can smoothly choose
$z=\tan^{-1}a\in(0,\pi/2)$.  Note that
$2z$ is the angle measure between the asymptotic directions
$\varphi_1^\perp,\varphi_2^\perp$, and that
\begin{equation}
  \omega^1 = (\cos z)(\varphi_1+\varphi_2),\quad
  \omega^2 = (\sin z)(\varphi_1-\varphi_2).
\label{FrameChange17}
\end{equation}

The following is fundamental in the study of $K=-1$ surfaces.
\begin{Proposition}  On an immersed surface in $\mathbf E^3$ with
  constant Gauss  curvature $K=-1$, the associated function $z$,
  expressed in terms of asymptotic coordinates $s,t$, satisfies the
  {\em sine-Gordon equation}\index{sine-Gordon equation|(} 
\begin{equation}
  z_{st}=\sf12\sin(2z).
\label{sG17}
\end{equation}
\end{Proposition}
One can prove this by a direct computation, but we will instead
highlight certain general EDS constructions which relate
the $K=-1$ differential system to a hyperbolic Monge-Ampere system
associated to the sine-Gordon equation.  One of these is the notion of
an {\em integrable extension}\index{integrable extension|(} 
of an exterior differential system,
which we have not yet encountered.  This is a
device that handles a forseeable difficulty; namely, the sine-Gordon
equation is expressed in terms of the variables $s$ and $t$, but for
the $K=-1$ system, these are primitives of a conservation law, defined
on integral manifolds of the system
only up to addition of integration constants.  One can think of an
integrable extension as a device for appending the primitives of
conserved $1$-forms.  More precisely, an integrable
extention of an EDS $(M,\mathcal E)$ is given by a submersion
$M^\prime\stackrel{\pi}{\to}M$, 
with a differential ideal $\mathcal E^\prime$ on
$M^\prime$ generated {\em algebraically} by $\pi^*\mathcal E$ and some
$1$-forms on $M^\prime$.  In this case, the preimage in $M^\prime$ of
an integral manifold of $\mathcal E$ is foliated by integral manifolds
of ${\mathcal E}^\prime$.
For example, if a $1$-form $\varphi\in\Omega^1(M)$ is a conservation
law for $\mathcal E$, then one can take $M^\prime=M\times\R$, and let
$\mathcal E^\prime\subset\Omega^*(M^\prime)$ be generated by $\mathcal
E$ and $\tilde\varphi = \varphi-ds$, where $s$ is a fiber coordinate
on $\R$.  Then the preimage in $M^\prime$ of any integral manifold of
$(M,\mathcal E)$ is foliated by a $1$-parameter family of integral
manifolds of $(M^\prime,\mathcal E^\prime)$, where the parameter
corresponds to a choice of integration constant for
$\varphi$.\footnote{For more information about integrable extensions,
  see \S6 of \cite{Bryant:CharacteristicII}.}

\

\begin{Proof}
Because $z$ is defined only on $M^{(1)}$, it is clear that we will
have to prolong once more to study $z_{st}$.  (Twice is unnecessary,
because of the Monge-Ampere form of (\ref{sG17}).)  
From (\ref{FirstProl17}), we see that integral elements with
$\varphi_1\wedge\varphi_2\neq 0$ satisfy
$$
  dz-\omega^1_2 = 2p\varphi_1,\quad
  dz+\omega^1_2 = 2q\varphi_2,
$$
for some $p,q$.  These $p,q$ can be taken as fiber coordinates on the
second prolongation
$$
  M^{(2)}=M^{(1)}\times\R^2.
$$
Let 
$$
  \begin{array}{l} \theta_3 = dz-p\varphi_1-q\varphi_2,\\
    \theta_4=\omega^1_2+p\varphi_1-q\varphi_2,
  \end{array}
$$
and then the prolonged differential system is
$$
  \mathcal E^{(2)} = \{\theta,\theta_1,\ldots,\theta_4,
    d\theta_3, d\theta_4\}.
$$
Integral manifolds for the original $K=-1$ system correspond to
integral manifolds of $\mathcal E^{(2)}$; in particular, on such an
integral manifold $f^{(2)}:N \hookrightarrow M^{(2)}$ we have
$$
  0  =  d\theta_3|_{N^{(2)}}
     =  -dp\wedge\varphi_1-dq\wedge\varphi_2
$$
and
\begin{eqnarray*}
  0 & = & d\theta_4|_{N^{(2)}} \\
    & = & K\omega^1\wedge\omega^2+dp\wedge\varphi_1-dq\wedge\varphi_2
       \\
    & = & (-1)(\cos z)(\varphi_1+\varphi_2)\wedge
            (\sin z)(\varphi_1-\varphi_2) + dp\wedge\varphi_1-
            dq\wedge\varphi_2 \\
    & = & \sin(2z)\varphi_1\wedge\varphi_2+dp\wedge\varphi_1
        -dq\wedge\varphi_2.
\end{eqnarray*}

Now we define the integrable extension
$$
  M^{(2)\prime} = M^{(2)}\times\R^2,
$$
where $\R^2$ has coordinates $s,t$, and on $M^{(2)\prime}$ we define
the EDS $\mathcal E^{(2)\prime}$ to be generated by $\mathcal
E^{(2)}$, along with the $1$-forms $\varphi_1-ds$, $\varphi_2-dt$.
An integral manifold $f^{(2)}:N\hookrightarrow M^{(2)}$ of $\mathcal
E^{(2)}$ gives a $2$-parameter family of integral manifolds
$f^{(2)}_{s_0,t_0}:N \hookrightarrow M^{(2)\prime}$ of $\mathcal
E^{(2)\prime}$.  On any of these, the functions $s,t$ will be local
coordinates, and we will have
$$
  \begin{array}{l}
    0 = dz-p\,ds-q\,dt,\\
    0 = -dp\wedge ds-dq\wedge dt,\\
    0 = \sin(2z)ds\wedge dt-ds\wedge dp-dq\wedge dt.
  \end{array}
$$
These three clearly imply that $z(s,t)$ satisfies (\ref{sG17}).
\end{Proof}

\

Note that we can start from the other side, defining the
differential system for the sine-Gordon equation as
$$
  \mathcal E_{sG} = \{dz-p\,ds-q\,dt,
      -dp\wedge ds-dq\wedge dt,
      ds\wedge dp+dq\wedge dt - \sin(2z)ds\wedge dt\},
$$
which is a Monge-Ampere system on the contact manifold
$J^1(\R^2,\R)$.  One can form a ``non-abelian'' integrable extension
$$
  \mathcal P = J^1(\R^2,\R)\times\mathcal F
$$
of $(J^1(\R^2,\R),\mathcal E_{sG})$ by taking
$$
  \mathcal E_{\mathcal P}=\mathcal E_{sG}+
    \left\{\begin{array}{l}
      \omega^1 - (\cos z)(ds+dt), \\
      \omega^2 - (\sin z)(ds+dt), \\
      \omega^3, \\
      \omega^1_2+p\,ds-q\,dt, \\
      \omega^3_1 - (\sin z)(ds+dt), \\
      \omega^3_2 + (\cos z)(ds-dt). \end{array}\right\}.
$$
This system is differentially closed, as one can see by using the
structure equations for $\mathcal F$ and assuming that $z$ satisfies
the sine-Gordon equation.  Though the following diagram is
complicated, it sums up the whole story:
$$
  
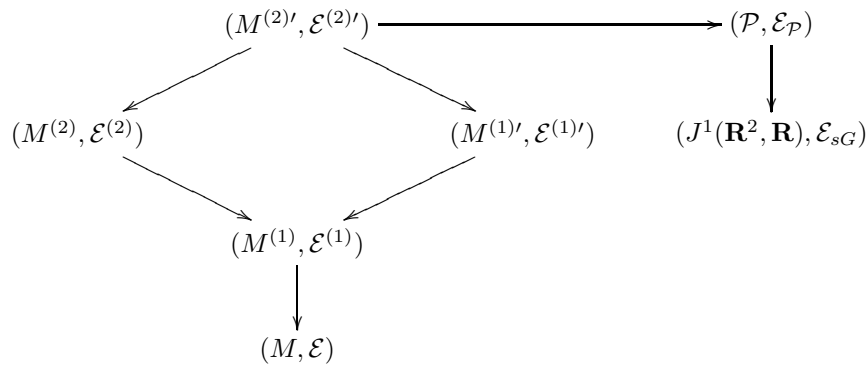
\xymatrix{ 
     & (M^{(2)\prime},\mathcal E^{(2)\prime})\ar[rr]\ar[dl]\ar[dr] & 
          & (\mathcal P,\mathcal E_{\mathcal P})\ar[d] \\
   (M^{(2)},\mathcal E^{(2)})\ar[dr] & & 
         (M^{(1)\prime},\mathcal E^{(1)\prime})\ar[dl] & 
         (J^1(\R^2,\R),\mathcal E_{sG}) \\
    & (M^{(1)},\mathcal E^{(1)})\ar[d] & & \\
    & (M,\mathcal E) & & }
$$
The system $(M^{(1)\prime},\mathcal E^{(1)\prime})$ was not introduced
in the proof; it is the integrable extension of $(M^{(1)},\mathcal
E^{(1)})$ formed by adjoining primitives $s,t$ for the conserved
$1$-forms $\varphi_1,\varphi_2$, and its prolongation turns out to be
$(M^{(2)\prime},\mathcal E^{(2)\prime})$.  In other words, starting on
$M^{(1)}$, one can first prolong and then adjoin primitives, or vice
versa.

The main point of this diagram is:
\begin{quote}
{\em The identification $(M^{(2)\prime},\mathcal E^{(2)\prime})\longrightarrow
  (\mathcal P,\mathcal E_{\mathcal P})$ is an isomorphism of exterior
  differential systems.  In other words, modulo prolongations and
  integrable extensions, the $K=-1$ system and the sine-Gordon system
  are equivalent.}
\end{quote}

Note that while conservation laws are preserved under prolongation,
there is an additional subtlety for integrable extensions.  In this
case, only those conservation laws for the sine-Gordon system that are
invariant under $s,t$-translation give conservation laws for the
$K=-1$ system.  Conversely, only those conservation laws for the
$K=-1$ system that are invariant under Euclidean motions (i.e.,
translations in $\mathcal F$) give conservation laws for the
sine-Gordon system.  There is even a difficulty involving
trivial\index{conservation law!trivial}
conservation laws; namely, the non-trivial, Euclidean-invariant
conservation laws $\varphi_1,\varphi_2$ for the $K=-1$ system induce
the trivial conservation laws $ds,dt$ for the sine-Gordon system.
\index{integrable extension|)}

However, we do have two conservation laws for sine-Gordon, obtained
via Noether's theorem applied to $s,t$-translations and the Lagrangian
$$
  \Lambda_{sG} = (pq-\sf12(\cos(2z)-1))ds\wedge dt;
$$
they are
$$
  \begin{array}{l}
    \psi_1 = \sf12p^2ds-\sf14(\cos(2z)-1)dt, \\
    \psi_2 = \sf12q^2dt+\sf14(\cos(2z)-1)ds.
  \end{array}
$$
The corresponding conserved $1$-forms on $(M^{(2)},\mathcal E^{(2)})$
are
$$
  \begin{array}{l}
    \varphi_3 = \sf12p^2\varphi_1 - \sf14(\cos(2z)-1)\varphi_2, \\
    \varphi_4 = \sf12q^2\varphi_2 + \sf14(\cos(2z)-1)\varphi_1.
  \end{array}
$$
In the next section, we will introduce a B\"acklund
transformation\index{B\"acklund transformation} for
the sine-Gordon equation, which can be used to generate an
infinite double-sequence $\psi_{2k-1},\psi_{2k}$ of conservation
laws (see \cite{Anderson:Lie}).
These in turn give an infinite double-sequence
$\varphi_{2k+1},\varphi_{2k+1}$ of independent conservation laws for $K=-1$,
extending the two pairs that we already have.
Although we will not discuss these, it is worth pointing out that the
generalized symmetries\index{symmetry!generalized} on
$(M^{(\infty)},\mathcal E^{(\infty)})$ to which they correspond under
Noether's
theorem are {\em not} induced by symmetries at any finite
prolongation.  For this reason, they are called {\em hidden
  symmetries}\index{symmetry!hidden}.  
\index{sine-Gordon equation|)}

Finally, we mention the following non-existence result, which is 
similar to the
result in Proposition~\ref{NoCLs17} for higher-dimensional non-linear
Poisson equations\index{Poisson equation!non-linear}.
\begin{Proposition}
In dimension $n\geq 3$, there are {\em no} second-order Euclidean-invariant
conservation laws for the
linear Weingarten system for hypersurfaces in $\mathbf{E}^{n+1}$ with
Gauss curvature $K=-1$.  
\end{Proposition}
\begin{Proof}
In contrast to our proof of the analogous statement for non-linear
Poisson equations, we give here a direct argument not
appealing to generating functions.  We work on the product
$$
  \mathcal{F}^{(1)} = \mathcal F\times\{(a_1,\ldots,a_n)\in\R^n:
    \textstyle\prod a_i = -1\displaystyle\},
$$
where $\mathcal{F}$ is the Euclidean frame bundle for
$\mathbf{E}^{n+1}$, and the other factor parameterizes eigenvalues of
admissible second fundamental forms.  We use the usual structure
equations, but without the sum convention:
\begin{eqnarray*}
  d\omega_i & = & -\sum_j\omega_{ij}\wedge\omega_j +
    \pi_i\wedge\theta, \\ 
  d\theta & = & -\sum_j\pi_j\wedge\omega_j \\
  d\omega_{ij} & = & -\sum_k\omega_{ik}\wedge\omega_{kj} +
    \pi_i\wedge\pi_j, \\
  d\pi_i & = & -\sum_j\pi_j\wedge\omega_{ji}.
\end{eqnarray*}
There is a Pfaffian system $\mathcal{I}$ on $\mathcal{F}^{(1)}$ whose transverse
$n$-dimensional integral
manifolds correspond to $K=-1$ hypersurfaces; it is differentially
generated by $\theta$ and the $n$ $1$-forms
$$
  \theta_i = \pi_i - a_i\omega_i.
$$
A conservation law for this system is an
$(n-1)$-form on $\mathcal{F}^{(1)}$ whose exterior derivative vanishes
on any integral manifold of $\mathcal{I}$.  A conservation law is
invariant under Euclidean motions if its restriction has the form
$$
  \varphi = \sum f_i(a_1,\ldots,a_n)\omega_{(i)}.
$$
What we will show is that for $n\geq 3$, such a form cannot be closed
modulo $\mathcal I$ unless it equals $0$. 

First, we calculate using the structure equations that
\begin{equation}
  d\varphi \equiv \sum_i(df_i-\sum_jf_j\omega_{ij})\wedge
    \omega_{(i)} \pmod{\mathcal{I}}.
\label{HalfThere17}
\end{equation}
To proceed further, we want an expression for
$\omega_{ij}\wedge\omega_{(i)}$, which we obtain first by computing
\begin{eqnarray*}
  0 & \equiv & d\theta_i \pmod{\mathcal{I}} \\
    & \equiv & -\sum_j\pi_k\wedge\omega_{ki} - da_i\wedge\omega_i
      + a_i\sum_k\omega_{ik}\wedge\omega_k \\
    & \equiv & -da_i\wedge\omega_i + 
      \sum_k(a_i-a_k)\omega_{ik}\wedge\omega_k.
\end{eqnarray*}
and then by multiplying the last result by $\omega_{(ij)}$:
\begin{eqnarray*}
  0 \equiv da_i\wedge\omega_{(j)} +
    (a_i-a_j)\omega_{ij}\wedge\omega_{(i)} \pmod{\mathcal{I}}.
\end{eqnarray*}
Using this in (\ref{HalfThere17}), we obtain
$$
  d\varphi \equiv  \sum_i\left(df_i + 
    \sum_j\frac{f_i}{a_j-a_i}da_j\right)\wedge\omega_{(i)}
    \pmod{\mathcal{I}}.
$$
So $\varphi$ is a conservation law only if for each $i$,
$$
  df_i = f_i\sum_j\frac{da_j}{a_i-a_j}.
$$
Keep in mind that we are requiring this equation to hold on the locus
$\mathcal{F}^{(1)}$ where $\Pi a_i = -1$.  Wherever at least one
$f_i(a_1,\ldots,a_n)$ is non-zero, we must have
$$
  0 = d\left(\sum_j\frac{da_j}{a_i-a_j}\right)
    = -\sum_j\frac{da_i\wedge da_j}{(a_i-a_j)^2}.
$$
However, when $n\geq 3$, the summands in the expression are linearly
independent $2$-forms on $\mathcal{F}^{(1)}$.
\end{Proof}

\

It is worth noting that when $n=2$, this $2$-form {\em does}
vanish, and we can solve for $f_1(a_1,a_2)$, $f_2(a_1,a_2)$ to obtain
the conservation laws $\varphi_1$, $\varphi_2$ discussed earlier.
The same elementary method can be used to analyze second-order
conservation laws for more general Weingarten equations; in this way,
one can obtain a full classification of those few Wiengarten equations
possessing higher-order conservation laws.
\index{conservation law!higher-order|)}
\index{conservation law!for K@for $K=-1$ surfaces|)}

\subsection{Two B\"acklund Transformations}
\label{Subsection:Backlund}
\index{B\"acklund transformation|(}

We have seen a relationship between the $K=-1$ surface system, and
the sine-Gordon equation\index{sine-Gordon equation|(}
\begin{equation}
  z_{xy}=\sf12\sin(2z).
\label{sineGordon}
\end{equation}
Namely, the half-angle measure between the asymptotic
directions\index{asymptotic curve} on a $K=-1$
surface, when expressed in asymptotic coordinates, satisfies the
sine-Gordon equation.  We have also interpreted this relationship in
terms of important EDS constructions.  In this section, we will explain
how this relationship connects the {\em B\"acklund transformations}
associated to each of these systems.

There are many definitions of B\"acklund transformation in the
literature, and instead of trying to give an all-encompassing
definition, we will restrict attention to Monge-Ampere
systems
$$
  \mathcal E = \{\theta,\Theta,\Psi\},
$$
where $\theta$ is a contact form on a manifold $(M^5,I)$, and
$\Theta,\Psi\in\Omega^2(M)$ are linearly independent modulo $\{I\}$.
Suppose that $(M,\mathcal E)$ and $(\bar M,\bar{\mathcal E})$ are
two Monge-Ampere systems, with
$$
  \mathcal E=\{\theta,\Theta,\Psi\},\quad
  \bar{\mathcal E} = \{\bar\theta,\bar\Theta,\bar\Psi\}.
$$
A {\em B\"acklund transformation} between $(M,\mathcal E)$ and $(\bar
M, \bar{\mathcal E})$, is a $6$-dimensional submanifold $B\subset
M\times \bar M$ such that in the diagram
\begin{equation}
  \xymatrix{ & B \ar[dl]_\pi \ar[dr]^{\bar\pi} & \\ M & & \bar M,}
\label{BacklundDiagram17}
\end{equation}
\begin{itemize}
\item each projection $B\to M$, $B\to\bar M$ is a submersion; and
\item pulled back to $B$, we have
$$
  \{\Psi,\bar\Psi\} \equiv
  \{\Theta,\bar\Theta\}\pmod{\{\theta,\bar\theta\}}. 
$$
\end{itemize}
The second condition implies that the dimension of the space of
$2$-forms spanned by
$\{\Theta,\Psi,\bar\Theta,\bar\Psi\}$ modulo 
$\{\theta,\bar\theta\}$ is at
most $2$.  Therefore,
$$
  \{\Theta,\Psi\}\equiv\{\bar\Theta,\bar\Psi\}\pmod
    {\{\theta,\bar\theta\}}.
$$
This consequence is what we really want, but the original formulation
has the extra benefit of ruling out linear dependence between
$\Theta$ and $\bar\Theta$, which would lead to a triviality in what
follows.

A B\"acklund transformation allows one to find a family of integral
manifolds of
$(\bar M,\bar{\mathcal E})$ from one integral manifold
$N^2\hookrightarrow M$ of $(M,\mathcal E)$, as follows.  
On the $3$-dimensional preimage $\pi^{-1}(N)\subset B$, the
restriction $\bar\pi^*\bar{\mathcal E}$ is algebraically generated by
$\bar\theta$ alone, and is therefore an integrable Pfaffian system.  Its
integral manifolds can therefore be found by ODE methods, and they foliate
$\pi^{-1}(N)$ into a $1$-parameter family of surfaces which project
by $\bar\pi$ to integral manifolds of $(\bar M,\bar{\mathcal E})$.  In
each of the following two examples, $(M,\mathcal E)$ and $(\bar
M,\bar{\mathcal E})$ are
equal, so one can generate from one known solution many others.

\subsubsection{Example 1: B\"acklund transformation for the
  sine-Gordon equation.} 
The primary example concerns the sine-Gordon equation (\ref{sineGordon}).  
The well-known coordinate phenomenon is that if two functions $u(x,y)$,
$\bar u(x,y)$ satisfy the first-order PDE system
\begin{equation}
  \left\{\begin{array}{l}
    u_x-\bar u_x = \lambda\sin(u+\bar u), \\
    u_y+\bar u_y = \sf1\lambda\sin(u-\bar u),
  \end{array}\right.
\label{sGsystem}
\end{equation}
where $\lambda\neq 0$ is any constant,\footnote{This $\lambda$ will
  {\em not} correspond to the integration parameter in the B\"acklund
  transformation.  It plays a role only in the relation to the $K=-1$
  system, to be discussed shortly.} then each of $u(x,y)$ and $\bar
u(x,y)$ satisfies (\ref{sineGordon}).  Conversely, given a function
$\bar u(x,y)$, the overdetermined system (\ref{sGsystem}) for unknown
$u(x,y)$ is compatible, and can therefore be reduced to an ODE system,
if and only if $\bar u(x,y)$ satisfies
(\ref{sineGordon}).  This indicates that given one solution of the
sine-Gordon equation, ODE methods give a family of additional
solutions. 

We fit this example into our definition of a B\"acklund transformation
as follows.
Start with two copies of the sine-Gordon
Monge-Ampere system, one on $M=\{(x,y,u,p,q)\}$ generated by
$$
   \mathcal E = \left\{\begin{array}{l}
     \theta = du-p\,dx-q\,dy, \\
     \Theta = d\theta = -dp\wedge dx-dq\wedge dy, \\
     \Psi = dx\wedge dp + dq\wedge dy-\sin(2u)dx\wedge dy
   \end{array}\right\},
$$
the other on $\bar M=\{(\bar x,\bar y, \bar u,\bar p,\bar q)\}$
generated by
$$
  \bar{\mathcal E} = \left\{\begin{array}{l}
     \bar\theta = d\bar u-\bar p\,d\bar x-\bar q\,d\bar y, \\
     \bar\Theta = d\bar\theta = -d\bar p\wedge d\bar x-d\bar q
                    \wedge d\bar y, \\
     \bar\Psi = d\bar x\wedge d\bar p + d\bar q\wedge d\bar y
                 -\sin(2\bar u)d\bar x\wedge d\bar y
   \end{array}\right\}.
$$
One can verify that the submanifold $B\subset M\times\bar M$ defined by
$$
  \left\{\begin{array}{l}
  \bar x = x, \quad \bar y = y, \\
  p-\bar p  =  \lambda\sin(u+\bar u), \\
  q+\bar q  =  \sf1\lambda\sin(u-\bar u),
  \end{array}\right.
$$
satisfies the criteria for a B\"acklund transformation, and that the
process of solving the overdetermined system (\ref{sGsystem}) for
$u(x,y)$ corresponds to integrating the Frobenius system as described
previously.  

For example, the solution $\bar u(x,y)=0$ of sine-Gordon corresponds
to the integral manifold $N=\{(x,y,0,0,0)\}\subset M$, whose preimage
in $B\subset M\times\bar M$ has coordinates $(x,y,u)$ and satisfies
$$
  \bar u = \bar p = \bar q = 0,\
  p = \lambda\sin(u),\ q = \sf1\lambda\sin(u).
$$
The system $\mathcal E$ is algebraically generated by the form
$$
  \theta|_{\bar\pi^{-1}(N)} = du - \lambda\sin(u)dx -
        \sf1\lambda\sin(u)dy.
$$
The problem of finding $u(x,y)$ on which this $\theta$ vanishes is the
same as solving the overdetermined system (\ref{sGsystem}) with $\bar
u=0$.  It is obtained by integrating
$$
  \frac{du}{\sin u}-\lambda\,dx-\sf1\lambda\,dy = 0,
$$
which has the implicit solution
$$
  -\ln(\csc u +\cot u)-\lambda x-\sf1\lambda y = c,
$$
where $c$ is the integration constant.  This can be solved for $u$ to
obtain 
$$
  u(x,y) = 2\tan^{-1}(e^{\lambda x+\frac{1}{\lambda} y+c}).
$$
One can verify that this is indeed a solution to the sine-Gordon
equation.  In principle, we could rename this as $\bar u$, and repeat
the process to obtain more solutions.

\subsubsection{Example 2: B\"acklund transformation for the $K=-1$ system.}
Suppose that $f,\bar f:N\hookrightarrow\mathbf E^3$ are two immersions
of a 
surface into Euclidean space.
We say that there is a {\em pseudospherical line 
congruence}\index{pseudospherical line congruence|(} between
$f,\bar f$ if for each $p\in N$:
\begin{enumerate}
\item the line through $f(p)$ and $\bar f(p)$ in $\mathbf E^3$ is tangent
  to each surface at these points (we assume $f(p)\neq \bar f(p)$);
\item the distance $r=||f(p)-\bar f(p)||$ is constant;
\item the angle $\tau$ between the normals $\nu(p)$ and $\bar\nu(p)$ is
  constant.
\end{enumerate}
This relationship between $f,\bar f$ will play a role analogous to
that of
the system (\ref{sGsystem}).  We prove the following theorem of
Bianchi.\index{Bianchi, L.}
\begin{Theorem}  If there is a pseudospherical line congruence between
  $f,\bar f:N\hookrightarrow \mathbf E^3$, then each of $f$ and
  $\bar f$ has constant negative Gauss curvature
$$
  K = -\frac{\sin^2(\tau)}{r^2}.
$$
\label{Theorem:Congruence}
\end{Theorem}
It is also true that given one surface $\bar f$, there locally exists a
surface $f$ sharing a pseudospherical line congruence with $\bar f$ if
and only if $\bar f$ has constant negative Gauss curvature.  
We will partly verify this claim, after proving Bianchi's theorem.

\

\begin{Proof}
Choose Euclidean frame fields $F,\bar
F:N\to \mathcal F$ which are adapted to the pair of surfaces in the
sense that 
\begin{equation}
  \bar e_1=e_1,
\label{FirstFrameReln17}
\end{equation}
made possible by condition 1 above.
Also, as usual, we let $e_3,\bar e_3$ be unit normals to $f,\bar
f$, respectively, which must then satisfy
\begin{equation}
   \left\{\begin{array}{l}
            \bar e_2 = (\cos\tau)e_2+(\sin\tau)e_3,
      \\ \bar e_3 = (-\sin\tau)e_2+(\cos\tau)e_3, \end{array}
   \right.
\label{FrameRelns17}
\end{equation}
with $\tau$ constant by condition 3.  Now condition 2 says that
\begin{equation}
  \bar f(p) = f(p) + re_1(p)
\label{PointReln17}
\end{equation}
for fixed $r$.  We can use the structure equations
$$
    df= e_i\cdot\omega^i,\qquad
    de_i = e_j\cdot\omega^j_i,
$$
and similar for $d\bar f$, $d\bar e_i$, to obtain relations among the
pullbacks by $F$ and $\bar F$ of the canonical forms on $\mathcal F$.
Namely,
$$
  d\bar f = e_1\bar\omega^1 + ((\cos\tau)e_2+(\sin\tau)e_3)\bar\omega^2,
$$
and also
$$
  d\bar f = d(f+re_1) =
    e_1\omega^1+e_2(\omega^2+r\omega^2_1)+e_3(r\omega^3_1),
$$
so that
\begin{equation}
  \left\{\begin{array}{l}
    \bar\omega^1 = \omega^1, \\
    (\cos\tau)\bar\omega^2 = \omega^2+r\omega^2_1, \\
    (\sin\tau)\bar\omega^2 = r\omega^3_1.
  \end{array}\right.
\label{FirstRelns17}
\end{equation}
Note that
$\omega^2+r\omega^2_1=(r\cot\tau)\omega^3_1$, and this gives a necessary
condition on $f$ alone to share a pseudospherical line congruence.
Similar calculations using $e_1=\bar e_1$ yield
$$
  \left\{\begin{array}{l}
   \bar\omega^2_1 = (\cos\tau)\omega^2_1+(\sin\tau)\omega^3_1, \\
   \bar\omega^3_1 = -(\sin\tau)\omega^2_1+(\cos\tau)\omega^3_1,
  \end{array}\right.
$$
and differentiating the remaining relations (\ref{FrameRelns17}) gives
$$
  \bar\omega^3_2 = \omega^3_2,
$$
giving complete expressions for $\bar F^*\omega$ in terms of
$F^*\omega$.
Note in particular that
$$
  \bar\omega^3_1 = \frac{\sin\tau}{r}\omega^2.
$$

Now we can consider the curvature, expanding both sides of the
definition
\begin{equation}
  d\bar\omega^2_1 = -\bar K\bar\omega^1\wedge\bar\omega^2.
\label{GaussCurv17}
\end{equation}
First,
\begin{eqnarray*}
  d\bar\omega^2_1 & = & -\bar\omega^2_3\wedge\bar\omega^3_1 \\
    & = & -\omega^2_3\wedge\frac{\sin\tau}{r}\omega^2 \\
    & = & h_{21}\frac{\sin\tau}{r}\omega^1\wedge\omega^2,
\end{eqnarray*}
where $h_{21}$ is part of the second fundamental form of
$F:N\hookrightarrow M$, defined by $\omega^3_i = h_{ij}\omega^j$.
Note in particular that if $h_{12}=0$, then $\omega^3_1$ is a multiple
of $\omega^1$, so by (\ref{FirstRelns17}),
$\bar\omega^1\wedge\bar\omega^2=0$, a 
contradiction; we can now assume that $h_{12}=h_{21}\neq 0$.
On the right-hand side of (\ref{GaussCurv17}),
\begin{eqnarray*}
  -\bar K\bar\omega^1\wedge\bar\omega^2 & = &
    -\bar K\omega^1\wedge\frac{r}{\sin\tau}\omega^3_1 \\
    & = & -\bar Kh_{12}\frac{r}{\sin\tau}\omega^1\wedge\omega^2.
\end{eqnarray*}
Equating these expressions, we have
$$
  \bar K = -\frac{\sin^2\tau}{r^2}
$$
as claimed.
\end{Proof}

\

Now suppose given a surface $f:N\hookrightarrow \mathbf E^3$ with
constant negative Gauss curvature $K=-1$.  We are interested in
finding $\bar f$ which shares with $f$ a pseudospherical line
congruence.

We start with local coordinates $(s,t)$ on $N$ whose coordinate lines
$ds=0$, $dt=0$ define the asymptotic curves of $f$.  It will be
convenient to instead have orthogonal coordinate lines, so we define
the coordinates $x=s+t$, $y=s-t$, for which $dx=0$ and $dy=0$ define
the principal curves of $f$.  We have seen in
(\ref{FirstFF17}, \ref{SecondFF17}) that the first and second
fundamental forms are given by
\begin{eqnarray}
  I & = & ds^2 + 2\cos(2z)ds\,dt+dt^2 \\
     & = & \cos^2z\,dx^2+\sin^2z\,dy^2, \\
\label{SecondAsymp17}
  II & = & 2\sin(2z)ds\,dt \\
     & = & \sin(z)\cos(z)(dx^2-dy^2),
\end{eqnarray}
where $2z$ is half of the angle measure between the asymptotic directions and
satisfies the sine-Gordon equation. 
One orthonormal coframing is given by $(\cos(z)dx,\sin(z)dy)$;
we consider an orthonormal coframing differing from this one
by rotation by some $\alpha$:
\begin{eqnarray}
  \left(\begin{array}{c}\omega^1 \\ \omega^2\end{array}\right) & = &
    \left(\begin{array}{cc}\cos\alpha & \sin\alpha \\ -\sin\alpha &
            \cos\alpha \end{array}\right)\left(\begin{array}{c}
       \cos(z)dx \\ \sin(z)dy\end{array}\right) \\ & = &
  \left(\begin{array}{cc}\cos(\alpha-z) & \cos(\alpha+z) \\
    -\sin(\alpha-z) & -\sin(\alpha+z)\end{array}\right)\left(\begin{array}{c}
      ds \\ dt\end{array}\right).
\label{NewFrame17}
\end{eqnarray}
The idea here is that we are looking for a function $\alpha$ on $N$
for which this coframing could be part of that induced by a
pseudospherical line congruence.  
The main compatibility condition, derived from (\ref{FirstRelns17}), is
\begin{equation}
  \omega^2+r\omega^2_1 = r(\cot\tau)\omega^3_1.
\label{PseudoCompat17}
\end{equation}
We can compute the Levi-Civita
connection form $\omega^2_1$ using the structure equations
$d\omega^i=-\omega^i_j\wedge\omega^j$, and find
$$
  \omega^2_1 = (\alpha_s+z_s)ds+(\alpha_t-z_t)dt.
$$
Similarly, we can compute from (\ref{SecondAsymp17}) and (\ref{NewFrame17}) the coefficients
of the second fundamental form with respect to the coframe
$(\omega^1,\omega^2)$, and find
$$
  \begin{array}{l}
    \omega^3_1 = \sin(\alpha-z)ds-\sin(\alpha+z)dt, \\
    \omega^3_2 = \cos(\alpha-z)ds-\cos(\alpha+z)dt.
  \end{array}
$$
Substituting these into the compatibility condition
(\ref{PseudoCompat17}), we obtain an equation of $1$-forms whose $ds$,
$dt$ coefficients are
\begin{eqnarray*}
    \alpha_s+z_s & = & (\csc\tau+\cot\tau)\sin(\alpha-z), \\
    \alpha_t-z_t & = & \left(\frac{1}{\csc\tau+\cot\tau}\right)\sin(\alpha+z).
\end{eqnarray*}
We compare this to (\ref{sGsystem}), and conclude that the local existence of a
solution $\alpha$ is equivalent to having $z$ satisfy the sine-Gordon
equation.  Note that the role played by $\lambda$ in (\ref{sGsystem})
is similar to that played by the angle $\tau$ in the pseudospherical
line congruence.
\index{pseudospherical line congruence|)}
\index{sine-Gordon equation|)}
We conclude by exhibiting a surprising use of the B\"acklund
transformation for the $K=-1$ system.  This starts with an
integral manifold in $M^5$ of $\mathcal E$ (see (\ref{KSystem17}))
that is {\em not} transverse as a Legendre 
submanifold,\index{Legendre submanifold!non-transverse} in the sense of
being a $1$-jet lift of an immersed surface in $\mathbf E^3$.
Instead, regarding the contact manifold $M$ as the unit sphere bundle
over $\mathbf E^3$, $N\hookrightarrow M$ consists of the unit normal
bundle of the line $\{(0,0,w):w\in\R\}\subset\mathbf E^3$.
This Legendre surface is topologically a
cylinder.  To study its geometry, we will work in the circle bundle
$\mathcal F\to M$.  Its preimage there is parameterized by
\begin{equation}
  (u,v,w)\mapsto \left((0,0,w),\left(\begin{array}{l}
    e_1 = (\sin u\cos v, -\sin u\sin v, \cos u), \\
    e_2 = (\cos u\cos v, -\cos u\sin v, -\sin u), \\
    e_3 = (\sin v, \cos v, 0)\end{array}\right)\right),
\label{Degenerate17}
\end{equation}
where $(u,v,w)\in S^1\times S^1\times\R$.  It is easily verified that
this is an integral manifold for the pullback of $\mathcal E$ by
$\mathcal F\to M$.  We will apply the B\"acklund transformation to
this degenerate integral manifold, and obtain a non-trivial surface in
${\mathbf E}^3$ with Gauss curvature $K=-1$.

We will take the B\"acklund transformation to be the submanifold
$B\subset\mathcal F\times\bar{\mathcal F}$ ($\bar{\mathcal F}$ is another
copy of $\mathcal F$) defined by (\ref{FirstFrameReln17},
\ref{FrameRelns17}, \ref{PointReln17}); this is a lift of the original
picture (\ref{BacklundDiagram17}) from $M$ to $\mathcal F$.  We
fix the constants of the line congruence to be
$$
  \tau = \sf{\pi}{2},\quad r=1.
$$
As a consequence, if our B\"acklund transformation gives a
transverse Legendre submanifold, then Theorem~\ref{Theorem:Congruence}
states that the corresponding surface in $\mathbf E^3$ will have
Gauss curvature $K=-1$.

Now, the definition of $B\subset\mathcal F\times\bar{\mathcal F}$ provides
a unique lift $\pi^{-1}(N)$ of our degenerate integral manifold
(\ref{Degenerate17}) to $B$.  The ``other'' $K=-1$ system
$\bar{\mathcal E}$ should restrict to $\pi^{-1}(N)\subset
B\subset\mathcal F\times\bar{\mathcal F}$ to be algebraically generated
by the $1$-form $\bar\theta$, and then $\pi^{-1}(N)$ will be
foliated into surfaces which project into integral manifolds of
$\bar{\mathcal E}$.  So we compute $\bar\theta=\bar\omega^3$:
\begin{eqnarray*}
  \bar\omega^3 & = & \langle d\bar x,\bar e_3\rangle \\
    & = & \langle dx+de_1, -e_2\rangle \\
    & = & \sin u\,dw-du.
\end{eqnarray*}
Indeed, this $1$-form is integrable, and its integral manifolds are of
the form
\begin{equation}
  u = 2\tan^{-1}(\exp(w+c)),
\label{BacklundLeaf17}
\end{equation}
where $c$ is an integration constant.  We will consider the integral
manifold corresponding to $c=0$.  The Euclidean surface that we are
trying to construct is now parameterized by $\bar x(u,v,w)$,
constrained by (\ref{BacklundLeaf17}).  We obtain
\begin{eqnarray*}
  \bar x(u,v,w) & = & x(u,v,w)+e_1(u,v,w) \\
    & = & (0,0,w) + (\sin u\cos v, -\sin u\sin v, \cos u) \\
    & = & \left(\sf{2e^w}{1+e^{2w}}\cos v,
          -\sf{2e^w}{1+e^{2w}}\sin v,
          w+\sf{1-e^{2w}}{1+e^{2w}}\right) \\
    & = & \left(\mbox{sech}\,w\cos v,-\mbox{sech}\,w\sin v,
             w-\mbox{tanh}\,w\right).
\end{eqnarray*}
This surface in $\mathbf E^3$ is the pseudosphere, the most familiar
surface of constant negative Gauss curvature; we introduced both it
and the ``framed line'' in \S\ref{Section:Hypersurface},
as examples of smooth but non-transverse Legendre submanifolds of
the unit sphere bundle $M\to\mathbf E^3$.  In principle, we could
iterate this B\"acklund
transformation, obtaining arbitrarily many examples of $K=-1$ surfaces.
\index{Gauss curvature|)}
\index{B\"acklund transformation|)}

\backmatter






\bibliography{pcfull}
\bibliographystyle{amsalpha}

\printindex

\end{document}